%% file: CFU_2023_DyadicKDE.tex
\documentclass{article}
\pdfoutput=1

\input{preamble_CFU.tex}

\begin{document}


\title{Uniform Inference for Kernel Density Estimators with Dyadic Data}

\author{
  Matias D.\ Cattaneo\textsuperscript{1}
  \and
  Yingjie Feng\textsuperscript{2}
  \and
  William G.\ Underwood\textsuperscript{1*}
}

\makeatletter
\def\nmfootnote{\gdef\@thefnmark{}\@footnotetext}
\makeatother

\maketitle

\footnotetext[1]{
  Department of Operations Research
  and Financial Engineering,
  Princeton University
}
\footnotetext[2]{
  School of Economics and Management,
  Tsinghua University
}
\nmfootnote{
  \textsuperscript{*}Corresponding author:
  \href{mailto:wgu2@princeton.edu}{\texttt{wgu2@princeton.edu}}
}
\setcounter{footnote}{3}

\setcounter{page}{0}\thispagestyle{empty}

\begin{abstract}

  Dyadic data is often encountered when quantities of interest
  are associated with the edges of a network.
  As such it plays an important role in statistics,
  econometrics and many other data science disciplines.
  %
  We consider the problem of uniformly estimating
  a dyadic Lebesgue density function,
  focusing on nonparametric kernel-based estimators taking the form of
  dyadic empirical processes.
  %
  Our main contributions include the minimax-optimal uniform convergence rate
  of the dyadic kernel density estimator,
  along with strong approximation results
  for the associated standardized and Studentized $t$-processes.
  A consistent variance estimator
  enables the construction of valid and feasible uniform confidence bands for
  the unknown density function.
  We showcase the broad applicability of our results by
  developing novel counterfactual density estimation
  and inference methodology for dyadic data,
  which can be used for causal inference and program evaluation.
  %
  A crucial feature of dyadic distributions
  is that they may be ``degenerate'' at certain
  points in the support of the data,
  a property making our analysis somewhat delicate.
  Nonetheless our methods for uniform inference
  remain robust to the potential presence of such points.
  %
  For implementation purposes,
  we discuss inference procedures based on
  positive semi-definite covariance estimators,
  mean squared error optimal bandwidth selectors
  and robust bias correction techniques.
  We illustrate the empirical finite-sample performance of our
  methods both in simulations and with real-world trade data,
  for which we make comparisons between observed and counterfactual
  trade distributions in different years.
  Our technical results concerning strong approximations
  and maximal inequalities are of potential independent interest.

\end{abstract}


\clearpage

\tableofcontents
\pagebreak


\section{Introduction}
\label{sec:introduction}

Dyadic data, also known as graphon data,
plays an important role in the statistical, social,
behavioral and biomedical sciences.
In network settings,
this type of dependent data captures
interactions between the units of study,
and its analysis is
of interest in statistics \citep{kolaczyk2009statistical}, economics
\citep{graham2020network}, psychology \citep{kenny2020dyadic}, public health
\citep{luke2007network} and many other data science disciplines.
For $n \geq 2$, a dyadic data set contains
$\frac{1}{2}n(n-1)$
observed real-valued random variables
\begin{align*}
  \bW_n = (W_{ij}:1\leq i<j \leq n),
  \quad\qquad W_{ij}
  &= W(A_i,A_j,V_{ij}),
\end{align*}
where $W$ is an unknown function, $\bA_n=(A_{i}:1\leq i \leq n)$
are independent and identically distributed (i.i.d.)\ latent random variables,
and $\bV_n=(V_{ij}:1\leq i<j \leq n)$ are i.i.d.\
latent random variables independent of $\bA_n$.
A natural interpretation of
this data is as a complete
undirected network on $n$ vertices,
with the latent variable $A_i$ associated with node $i$
and the observed variable $W_{ij}$ associated
with the edge between nodes $i$ and $j$.
The data generating process above is justified
without loss of generality by the celebrated
Aldous--Hoover representation theorem
for exchangeable arrays
\citep{aldous1981representations, hoover1979relations}.

Various distributional features of dyadic data are of interest in applications.
Most of the statistical literature focuses on parametric analysis, almost
exclusively considering moments of (transformations of)
the identically distributed $W_{ij}$.
See \citet{davezies2021exchangeable}, \citet{gao2021minimax},
\citet{MatsushitaOtsu2021}, and references therein, for contemporary
contributions and overviews. More recently, however, a few nonparametric
procedures for dyadic data have been proposed in the literature
\citep{graham2021dyadicregression,graham2022kernel}.

With the aim of estimating density-like functions
associated with $W_{ij}$ using
nonparametric kernel-based methods,
we investigate the statistical properties of a class
of local stochastic processes given by
\begin{align}\label{eq:estimator}
  w \mapsto \widehat{f}_W(w)
  = \frac{2}{n(n-1)} \sum_{i=1}^{n-1} \sum_{j=i+1}^n k_h(W_{ij},w),
\end{align}
where $k_h(\cdot,w)$ is a kernel function
that can change with the $n$-varying bandwidth parameter
$h=h(n)$ and the evaluation point
$w \in \cW\subseteq \R$. For each $w\in\cW$ and with an appropriate choice of
the kernel function
(e.g.\ $k_h(\cdot,w)=K((\cdot-w)/h)/h$ for an interior point $w$
of $\cW$ and a fixed symmetric integrable kernel function $K$),
the statistic $\widehat{f}_W(w)$ becomes a kernel density estimator
for the Lebesgue density function
$f_W(w) = \E\big[f_{W \mid AA}(w \mid A_i,A_j)\big]$,
where $f_{W \mid AA}(w \mid A_i,A_j)$ denotes the conditional Lebesgue density
of $W_{ij}$ given $A_i$ and $A_j$.
Setting
$k_h(\cdot,w)=K((\cdot-w)/h)/h$, \citet{graham2022kernel}
recently introduced the dyadic point estimator $\widehat{f}_W(w)$ and studied
its large sample properties pointwise in $w\in\cW=\mathbb{R}$, while
\citet{chiang2020empirical} established its rate of convergence uniformly in
$w\in\cW$ for a compact interval $\cW$ strictly contained in the support of
the dyadic data $W_{ij}$.
\citet{chiang2022inference} obtained a distributional approximation for the
supremum statistic $\sup_{w\in\cW}\big|\widehat{f}_W(w)\big|$
over a finite collection $\cW$ of design points.
More generally, as we discuss below,
the estimand $f_W(w)$ is useful in
different applications because it forms the basis for
counterfactual distributional analysis
(Section~\ref{sec:counterfactual})
and other nonparametric and semiparametric methods
(Section~\ref{sec:future}).
We remark that while we assume throughout that the network is complete,
our approach generalizes in a straightforward way to networks
with missing edges, as in Section~\ref{sec:trade_data}.
This can be seen by setting $W_{ij} = -\infty$
whenever the edge $\{i, j\}$ is not present, so that the law of
$W_{ij}$ is a mixture between a continuous distribution and a point mass
at $-\infty$.
We then apply our methodology to recover the
continuous component of this distribution,
following \citet{chiang2022inference}.

We contribute to the emerging
literature on nonparametric smoothing methods
for dyadic data with two main technical
results. Firstly, we derive the minimax rate of uniform convergence for density
estimation with dyadic data and show that the estimator
$\widehat{f}_W$ in \eqref{eq:estimator}
is minimax-optimal under appropriate conditions.
Secondly, we present a set of
uniform distributional approximation results
for the \emph{entire} stochastic process $\big(\widehat{f}_W(w):w\in\cW\big)$.
Furthermore, we illustrate the usefulness of our main results
with two distinct substantive statistical applications:
\begin{inlineroman}
  \item
  confidence bands for $f_W$
  (Section~\ref{sec:implementation}), and
  \item
  estimation and inference for counterfactual
  dyadic distributions (Section~\ref{sec:counterfactual}).
\end{inlineroman}
Our main results also lay the
foundation for studying the uniform distributional properties of other
nonparametric and semiparametric tests and estimators based on
dyadic data (Section~\ref{sec:future}).
Importantly, our inference
results cannot be deduced from the existing U-statistic, empirical process and
U-process theory available in the literature
\citep{van1996weak,gine2021mathematical} because,
as explained in detail below,
$\widehat{f}_W(w)$ is not a standard U-statistic,
nor is the stochastic process $\widehat{f}_W$ Donsker in general,
and the underlying dyadic data $\bW_n$ exhibits
statistical dependence due to its network structure.

Section~\ref{sec:setup} outlines the setup and presents the main assumptions
imposed throughout the paper. We first discuss a Hoeffding-type decomposition of
the U-statistic-like $\widehat{f}_W$ which is more general than the standard
Hoeffding decomposition for second-order U-statistics due
to its dyadic data structure.
In particular,
\eqref{eq:h-decomposition} shows that $\widehat{f}_W(w)$ decomposes into a sum
of the four terms $B_n(w)$, $L_n(w)$, $E_n(w)$ and $Q_n(w)$,
where $E_n(w)$ is not present in the classical second-order U-statistic theory.
The first term $B_n(w)$ captures the usual smoothing bias, the second term
$L_n(w)$ is akin to the H\'{a}jek projection for second-order U-statistics, the
third term $E_n(w)$ is a mean-zero double
average of conditionally independent terms,
and the fourth term $Q_n(w)$ is a negligible totally degenerate
second-order U-process.
The leading stochastic
fluctuations of the process $\widehat{f}_W$
are captured by $L_n$ and $E_n$,
both of which are known to be asymptotically distributed
as Gaussian random variables pointwise in $w\in\cW$ \citep{graham2022kernel}.
However, the H\'{a}jek projection term $L_n$ will often
be ``degenerate'' at some
or possibly all evaluation points $w\in\cW$.

Section~\ref{sec:point_estimation} studies minimax convergence rates for point
estimation of $f_W$ uniformly over $\cW$ and gives precise conditions under
which the estimator $\widehat{f}_W$ is minimax-optimal.
Firstly, in Theorem~\ref{thm:uniform_consistency} we establish the uniform rate
of convergence of $\widehat{f}_W$ for $f_W$. This result improves upon the
recent paper of \citet{chiang2020empirical} by allowing
for compactly supported dyadic data and
generic kernel-like functions $k_h$ (including boundary-adaptive kernels), while
also explicitly accounting for possible degeneracy of the
H\'{a}jek projection term $L_n$ at some or possibly all points $w\in\cW$.
Secondly, in Theorem~\ref{thm:minimax} we derive the
minimax uniform convergence rate for estimating
$f_W$, again allowing for possible degeneracy,
and verify that it is achieved by $\widehat f_W$.
This result appears to be new to the literature,
complementing recent work on parametric moment estimation
using graphon data \citep{gao2021minimax} and on nonparametric
kernel-based regression using dyadic data \citep{graham2021dyadicregression}.

Section~\ref{sec:inference} presents a distributional analysis of
the stochastic process $\widehat{f}_W$ uniformly in $w \in \cW$. Because
$\widehat{f}_W$ is not asymptotically tight in general, it
does not converge
weakly in the space of uniformly bounded real functions supported on $\cW$ and
equipped with the uniform norm \citep{van1996weak},
and hence is non-Donsker.
To circumvent this problem,
we employ strong approximation methods to characterize
its distributional properties.
Up to the smoothing bias term $B_n$ and the negligible term $Q_n$,
it is enough to consider
the stochastic process $w \mapsto L_n(w)+E_n(w)$.
Since $L_n$ can be degenerate at some or possibly all points $w\in\cW$, and also
because under some bandwidth choices both $L_n$ and $E_n$ can be of comparable
order, it is crucial to analyze the joint distributional properties of $L_n$ and
$E_n$. To do so, we employ a carefully crafted conditioning approach where we
first establish an unconditional strong approximation for $L_n$ and a
conditional-on-$\bA_n$ strong approximation for $E_n$.
We then combine these to obtain a strong approximation for $L_n+E_n$.

The stochastic process $L_n$ is an empirical process indexed by an $n$-varying
class of functions depending only on the i.i.d.\ random variables $\bA_n$. Thus
we use the celebrated Hungarian
construction \citep{komlos1975approximation}, building on ideas in
\citet{gine2004kernel} and \citet{gine2010confidence}. The resulting rate of
strong approximation is optimal, and follows from a generic strong approximation
result of potential independent interest given in
Section~SA3 of
the online supplemental appendix.
Our main result for $L_n$ is given as Lemma~\ref{lem:strong_approx_Ln},
and makes explicit the potential presence of degenerate points.

The stochastic process $E_n$ is an empirical process depending on the dyadic
variables $W_{ij}$ and indexed by an $n$-varying class of functions.
When conditioning on $\bA_n$, the variables $W_{ij}$ are independent but not
necessarily identically distributed (i.n.i.d.), and thus we establish a
conditional-on-$\bA_n$ strong
approximation for $E_n$ based on the Yurinskii coupling
\citep{yurinskii1978error}, leveraging a recent refinement obtained by
\citet*[Lemma~38]{belloni2019conditional}.
This result follows from a generic
strong approximation result which
gives a novel rate of strong approximation for
(local) empirical processes based on i.n.i.d. data,
given in
Section~SA3 of
the online supplemental appendix.
Lemma~\ref{lem:conditional_strong_approx_En}
gives our conditional strong approximation for $E_n$.

Once the unconditional strong approximation for $L_n$ and the
conditional-on-$\bA_n$ strong approximation for $E_n$ are established,
we show how to properly ``glue'' them together to deduce a final unconditional
strong approximation for $L_n+E_n$ and hence also for $\widehat{f}_W$ and its
associated $t$-process. This final step requires some additional technical work.
Firstly, building on our conditional strong approximation for $E_n$, we
establish an unconditional strong approximation for $E_n$ in
Lemma~\ref{lem:unconditional_strong_approx_En}. We then employ a
generalization of the celebrated Vorob'ev--Berkes--Philipp theorem
\citep{dudley1999uniform}, given in
Section~SA3 of
the online supplemental appendix,
to deduce a \emph{joint} strong
approximation for $(L_n,E_n)$ and, in particular, for $L_n+E_n$.
Thus we obtain our main result in Theorem~\ref{thm:strong_approx_fW},
which establishes a valid strong approximation for $\widehat{f}_W$
and its associated $t$-process.
This uniform inference result complements the recent contribution of
\citet{davezies2021exchangeable}, which is not applicable
here because $\widehat{f}_W$ is non-Donsker in general.

We illustrate the applicability of our strong approximation result for
$\widehat{f}_W$ and its associated $t$-process by constructing valid
standardized uniform confidence bands for
the unknown density function $f_W$.
Instead of relying on extreme value theory \citep*[e.g.][]{gine2004kernel}, we
employ anti-concentration methods, following \citet{chernozhukov2014anti}.
This illustration improves on the recent work of \citet{chiang2022inference},
which obtained simultaneous confidence intervals for the dyadic density $f_W$
based on a high-dimensional central limit theorem over rectangles,
following prior work by \cite{Chernozhukov-Chetverikov-Kato_2017_AoP}.
The distributional
approximation therein is applied to the
H\'{a}jek projection term $L_n$ only, whereas
our main construction leading to Theorem~\ref{thm:strong_approx_fW} gives a
strong approximation for the entire U-process-like
$\widehat{f}_W$ and its associated $t$-process, uniformly on $\cW$.
As a consequence, our uniform inference theory is robust to potential
unknown degeneracies in $L_n$ by virtue of our strong
approximation of $L_n+E_n$ and the use of proper standardization,
delivering a ``rate-adaptive'' inference procedure.
Our result appears to be the first to provide
confidence bands that are valid uniformly
over $w \in \cW$ rather than over some finite collection of design points.
Moreover, they provide distributional approximations for
the whole $t$-statistic process,
which can be useful in applications where functionals other
than the supremum are of interest.

Section~\ref{sec:implementation} addresses outstanding issues of
implementation. Firstly, we discuss estimation of the covariance function of the
Gaussian process underlying our strong approximation results. We present two
estimators, one based on the plug-in method, and the other based on a
positive semi-definite regularization thereof
\citep{laurent2005semidefinite}. We derive the uniform convergence rates for
both estimators in Lemma~\ref{lem:SDP}, which we then use to justify
Studentization of $\widehat{f}_W$ and a
feasible simulation-based approximation of the infeasible Gaussian process
underlying our strong approximation results.
Secondly, we discuss integrated mean
squared error (IMSE) bandwidth selection and provide a simple rule-of-thumb
implementation for applications \citep{wand1994kernel,simonoff2012smoothing}.
Thirdly, we provide
feasible, valid uniform inference methods for
$f_W$ by employing robust bias correction
\citep{Calonico-Cattaneo-Farrell_2018_JASA,
  Calonico-Cattaneo-Farrell_2022_Bernoulli}.
Algorithm~\ref{alg:method} summarizes our entire feasible methodology.

Section~\ref{sec:simulations} reports
empirical evidence for our proposed
feasible robust bias-corrected confidence bands for $f_W$.
We use simulations to show that these confidence bands are robust to potential
unknown degenerate points in the underlying dyadic distribution.

Section~\ref{sec:counterfactual} presents novel results for counterfactual
dyadic density estimation and inference,
offering an application of our general theory to a
substantive problem in statistics and other data science disciplines.
Counterfactual distributions are important for causal inference
and policy evaluation
\citep{dinardo1996distribution,chernozhukov2013inference},
and in the context of network data,
such analysis can be used to answer empirical questions
such as ``what would the international trade distribution in one year have been
if the gross domestic product (GDP) of the countries had
remained the same as in a previous year?''
We formally show how our theory for kernel-based dyadic estimators
can be used
to infer the counterfactual density function
of dyadic data had some monadic covariates followed a
different distribution.
We propose a two-step semiparametric reweighting approach in which we
first estimate the Radon--Nikodym derivative between the
observed and counterfactual
covariate distributions using a simple parametric estimator,
and then use this to construct a weighted dyadic kernel density estimator.
We present uniform consistency, strong approximation
and feasible inference results
for this dyadic counterfactual density estimator.
Finally, we also illustrate our methods
with a real dyadic data set recording bilateral trade between countries,
using GDP as a covariate for the counterfactual analysis.

Section~\ref{sec:future} discusses further statistical
applications of our main results,
including dyadic density hypothesis testing and nonparametric
and semiparametric dyadic regression.
Section~\ref{sec:conclusion} concludes the paper.
The online supplemental appendix includes other technical and methodological
results, proofs and additional details omitted here to conserve space.
Section~SA3 may be of independent interest, containing
two generic strong approximation theorems for empirical processes,
a generalized Vorob'ev--Berkes--Philipp theorem
and a maximal inequality for i.n.i.d.\ random variables.

\subsection{Notation}

The total variation norm of a
real-valued function $g$ of a single real variable is
defined to be
$\|g\|_\TV = \sup_{n \geq 1} \sup_{x_1 \leq \cdots \leq x_n}
\sum_{i=1}^{n-1} |g(x_{i+1}) - g(x_i)|$.
For an integer $m\geq 0$,
denote by $\mathcal{C}^m(\mathcal{X})$
the space of all $m$-times continuously differentiable functions
on $\mathcal{X}$.
For $\beta > 0$ and $C>0$,
define the H\"{o}lder class on $\mathcal{X}$ to be
$\cH^\beta_C(\cX)
=
\big\{
g \in \cC^{\flbeta}(\cX): \
\max_{1 \leq r \leq \flbeta}
\big| g^{(r)}(x) \big| \leq C
\text{ and }
\big| g^{(\flbeta)}(x) - g^{(\flbeta)}(x') \big|
\leq C |x-x'|^{\beta - \flbeta}, \
\forall x, x' \in \cX
\big\},
$
where $\flbeta$ denotes the largest integer
which is strictly less than $\beta$.
For $a \in \R$ and $b \geq 0$,
we write $[a \pm b]$ for the interval $[a-b, a+b]$.
For non-negative sequences
$a_n$ and $b_n$, write
$a_n \lesssim b_n$ or $a_n = O(b_n)$
to indicate that
$a_n / b_n$ is bounded for $n\geq 1$.
Write $a_n \ll b_n$ or $a_n = o(b_n)$ if $a_n / b_n \to 0$.
If $a_n \lesssim b_n \lesssim a_n$,
write $a_n \asymp b_n$.
For random non-negative sequences
$A_n$ and $B_n$, write
$A_n \lesssim_\P B_n$ or $A_n = O_\P(B_n)$ if
$A_n / B_n$ is bounded in probability.
Write $A_n = o_\P(B_n)$ if $A_n / B_n \to 0$ in probability.
For $a,b \in \R$, define $a\wedge b=\min\{a,b\}$ and
$a \vee b = \max\{a,b\}$.

\section{Setup}\label{sec:setup}

We impose the following two assumptions throughout this paper.
\begin{assumption}[Data generation]\label{ass:data}\onehalfspacing
  %
  Let $\bA_n = (A_i: 1 \leq i \leq n)$ be i.i.d.\
  random variables supported on $\cA \subseteq \R$
  and let $\bV_n = (V_{ij}: 1 \leq i < j \leq n)$
  be i.i.d.\ random variables with a Lebesgue density $f_V$ on $\R$,
  with $\bA_n$ independent of $\bV_n$.
  %
  Let $W_{ij} = W(A_i, A_j, V_{ij})$
  and $\bW_n = (W_{ij}: 1 \leq i < j \leq n)$,
  where $W$ is an unknown real-valued function
  which is symmetric in its first two arguments.
  Let $\cW \subseteq \R$ be a compact interval
  with positive Lebesgue measure $\Leb(\cW)$.
  The conditional distribution
  of $W_{ij}$ given $A_i$ and $A_j$
  admits a Lebesgue density $f_{W \mid AA}(w \mid A_i, A_j)$.
  For $\CH > 0$ and $\beta \geq 1$,
  take $f_W \in \cH^\beta_{\CH}(\cW)$ where
  $f_{W}(w) = \E\left[f_{W \mid AA}(w \mid A_i,A_j)\right]$
  and
  $f_{W \mid AA}(\cdot \mid a, a') \in \cH^1_{\CH}(\cW)$
  for all $a,a' \in \cA$.
  Suppose
  $\sup_{w \in \cW} \|f_{W \mid A}(w \mid \cdot\,)\|_\TV <\infty$ where
  $f_{W \mid A}(w \mid a) = \E\left[f_{W \mid AA}(w \mid A_i,a)\right]$.
\end{assumption}

In Assumption~\ref{ass:data} we require
the density $f_W$ be in a $\beta$-smooth H\"older class of
functions on the compact interval $\cW$.
H\"older classes are well-established in the minimax estimation literature
\citep{gine2021mathematical}, with the smoothness parameter $\beta$ appearing
in the minimax-optimal rate of convergence.
If the H\"older condition is satisfied only piecewise,
then our results remain valid provided that the boundaries
between the pieces are known and treated as boundary points.


\begin{assumption}[Kernels and bandwidth]\label{ass:kernel_bandwidth}%
  \onehalfspacing

  Let $h = h(n) > 0$
  be a sequence of bandwidths
  satisfying $h \log n \to 0$
  and $\frac{\log n}{n^2h} \to 0$.
  For each $w \in \cW$, let $k_h(\cdot, w)$
  be a real-valued function supported on
  $[w \pm h] \cap \cW$.
  For an integer $p \geq 1$,
  let $k_h$ belong to a family of
  boundary bias-corrected kernels of order $p$, i.e.,
  \begin{align*}
    \int_{\cW}
    (s-w)^r k_h(s,w) \diff{s}
    \quad
    \begin{cases}
      \begin{alignedat}{2}
        &= 1 &\qquad &\text{for all } w \in \cW \text{ if }\, r = 0, \\
        &= 0 & &\text{for all } w \in \cW \text{ if }\, 1 \leq r \leq p-1, \\
        &\neq 0 & &\text{for some } w \in \cW \text{ if }\, r = p.
      \end{alignedat}
    \end{cases}
  \end{align*}
  Also, for $\CL > 0$,
  suppose $k_h(s, \cdot) \in \cH^1_{\CL h^{-2}}(\cW)$
  for all $s \in \cW$.
\end{assumption}

This assumption allows for all standard compactly supported
and possibly boundary-corrected kernel functions
\citep{wand1994kernel,simonoff2012smoothing}.
Assumption~\ref{ass:kernel_bandwidth} implies that
if $h \leq 1$ then
$k_h$ is uniformly bounded by
$\Ck h^{-1}$ where $\Ck \vcentcolon = 2 \CL + 1 + 1/\Leb(\cW)$.

\subsection{Hoeffding-type decomposition and degeneracy}
\label{sec:degeneracy}

The estimator $\widehat{f}_W(w)$ is akin to a
U-statistic and thus admits a Hoeffding-type decomposition
which is the starting point for our analysis. We have
\begin{align}\label{eq:h-decomposition}
  \widehat{f}_W(w) - f_W(w) = B_n(w) + L_n(w) + E_n(w) + Q_n(w)
\end{align}
with
$B_n(w) = \E\big[\widehat{f}_W(w)\big] - f_W(w)$
and
\begin{align*}
  L_n(w)
  &= \frac{2}{n} \sum_{i=1}^n l_i(w),
  &
  E_n(w)
  &= \frac{2}{n(n-1)} \sum_{i=1}^{n-1} \sum_{j=i+1}^{n} e_{ij}(w),
  &
  Q_n(w)
  &= \frac{2}{n(n-1)} \sum_{i=1}^{n-1} \sum_{j=i+1}^{n} q_{ij}(w),
\end{align*}
where
$l_i(w) = \E[k_h(W_{ij},w) \mid A_i] - \E[k_h(W_{ij},w)]$,
$e_{ij}(w) = k_h(W_{ij},w) - \E[k_h(W_{ij},w) \mid A_i, A_j]$
and
$q_{ij}(w) = \E[k_h(W_{ij},w) \mid A_i, A_j] - \E[k_h(W_{ij},w) \mid A_i]
- \E[k_h(W_{ij},w) \mid A_j] + \E[k_h(W_{ij},w)]$.
The non-random term $B_n$ captures the smoothing (or misspecification) bias,
while the three stochastic processes $L_n$, $E_n$ and $Q_n$
capture the variance of the estimator.
These processes are
mean-zero with $\E[L_n(w)] = \E[Q_n(w)] = \E[E_n(w)] = 0$ for all $w \in \cW$,
and mutually orthogonal in $L^2(\P)$ since
$\E[L_n(w) Q_n(w')] = \E[L_n(w) E_n(w')] = \E[Q_n(w) E_n(w')] = 0$
for all $w, w' \in \cW$.

The stochastic process $L_n$ is akin to the H\'{a}jek projection of a U-process,
which can (and often will) exhibit degeneracy at some
or possibly all points $w \in \cW$.
To characterize different types of degeneracy,
we introduce the following non-negative lower and upper degeneracy constants:
\[
  \Dl^2 := \inf_{w \in \cW} \Var\left[f_{W \mid A}(w \mid A_i)\right]
  \qquad \text{ and } \qquad
  \Du^2 := \sup_{w \in \cW} \Var\left[f_{W \mid A}(w \mid A_i)\right].
\]
The following lemma describes the stochastic order of
different terms in the Hoeffding-type decomposition,
explicitly accounting for potential degeneracy.

\begin{lemma}[Bias and variance] \label{lem:variance}
  Suppose that Assumptions \ref{ass:data}~and~\ref{ass:kernel_bandwidth} hold.
  Then the bias term satisfies
  $\sup_{w \in \cW} \big| B_n(w) \big| \lesssim h^{p\wedge\beta}$
  and the variance terms satisfy
  \begin{align*}
    \E\left[ \sup_{w \in \cW} |L_n(w)| \right]
    &\lesssim
    \frac{\Du}{\sqrt n},
    &
    \E\left[ \sup_{w \in \cW} |E_n(w)| \right]
    &\lesssim \sqrt{\frac{\log n}{n^2h}},
    &
    \E\left[ \sup_{w \in \cW} |Q_n(w)| \right]
    \lesssim \frac{1}{n}.
  \end{align*}
\end{lemma}

Lemma~\ref{lem:variance} captures the potential total degeneracy of $L_n$ by
showing that if
$\Du=0$ then $L_n=0$ everywhere on $\cW$ almost surely.
The following lemma
captures the potential partial degeneracy of $L_n$, where
$\Du > \Dl = 0$.
For $w,w' \in \cW$,
define the covariance function of the dyadic kernel density estimator as
\begin{equation*}
  \Sigma_n(w,w')
  =
  \E\Big[
    \Big(
    \widehat f_W(w)
    - \E\big[\widehat f_W(w)\big]
    \Big)
    \Big(
    \widehat f_W(w')
    - \E\big[\widehat f_W(w')\big]
    \Big)
    \Big].
\end{equation*}
\begin{lemma}[Variance bounds]
  \label{lem:variance_bounds}
  Suppose that Assumptions~\ref{ass:data}
  and~\ref{ass:kernel_bandwidth} hold.
  Then for sufficiently large $n$,
  \begin{align*}
    \frac{\Dl^2}{n} + \frac{1}{n^2h}
    \inf_{w \in \cW} f_W(w)
    &\lesssim
    \inf_{w \in \cW} \Sigma_n(w,w)
    \leq
    \sup_{w \in \cW} \Sigma_n(w,w)
    \lesssim
    \frac{\Du^2}{n} + \frac{1}{n^2h}.
  \end{align*}
\end{lemma}

Combining Lemmas~\ref{lem:variance} and~\ref{lem:variance_bounds},
we have the following trichotomy
for degeneracy of dyadic distributions
based on $\Dl$ and $\Du$:
\begin{inlineroman}
  \item
  total degeneracy if $\Du = \Dl = 0$,
  \item
  partial degeneracy if $\Du > \Dl = 0$,
  \item
  no degeneracy if $\Dl > 0$.
\end{inlineroman}
In the case of no degeneracy,
it can be shown that
$\inf_{w \in \cW} \Var[L_n(w)] \gtrsim n^{-1}$,
while in the case of total degeneracy,
$L_n(w) = 0$ for all $w \in \cW$ almost surely.
When the dyadic distribution is partially degenerate,
there exists at least one point
$w\in\cW$ such that
$\Var\left[f_{W \mid A}(w \mid A_i)\right] = 0$
and $\Var[L_n(w)]\lesssim h n^{-1}$,
and there also exists at least one point
$w'\in \cW$ such that
$\Var\left[f_{W \mid A}(w' \mid A_i)\right] > 0$ and
$\Var[L_n(w')] \gtrsim \frac{1}{n}$.
We say $w$ is a \emph{degenerate point}
if $\Var\left[f_{W \mid A}(w \mid A_i)\right] = 0$,
and otherwise say it is a
\emph{non-degenerate point}.

As a simple example,
consider the family of dyadic distributions $\P_{\pi}$
indexed by $\pi = (\pi_1, \pi_2, \pi_3)$
with $\sum_{i=1}^3 \pi_i = 1$ and $\pi_i \geq 0$,
generated by $W_{ij} = A_i A_j + V_{ij}$,
where $A_i$ equals $-1$ with probability $\pi_1$,
equals $0$ with probability $\pi_2$
and equals $+1$ with probability $\pi_3$,
and $V_{ij}$ is standard Gaussian.
This model induces a latent
``community structure'' where community membership
is determined by the value of $A_i$ for each node $i$,
and the interaction outcome $W_{ij}$ is a function
only of the communities which $i$ and $j$ belong to and some
idiosyncratic noise.
Unlike the stochastic block model
\citep{kolaczyk2009statistical},
our setup assumes that
community membership has no impact on edge existence,
as we work with fully connected networks;
see Section~\ref{sec:trade_data} for a discussion of
how to handle missing edges in practice.
Also note that the parameter of interest in this paper is the Lebesgue density
of a continuous random variable $W_{ij}$
rather than the probability of network edge existence,
which is the focus of graphon estimation
literature \citep{gao2021minimax}.

\begin{figure}[H]
  \captionsetup[subfigure]{justification=centering}
  \centering
  \begin{subfigure}{0.32\textwidth}
    \centering
    \includegraphics[scale=0.52]{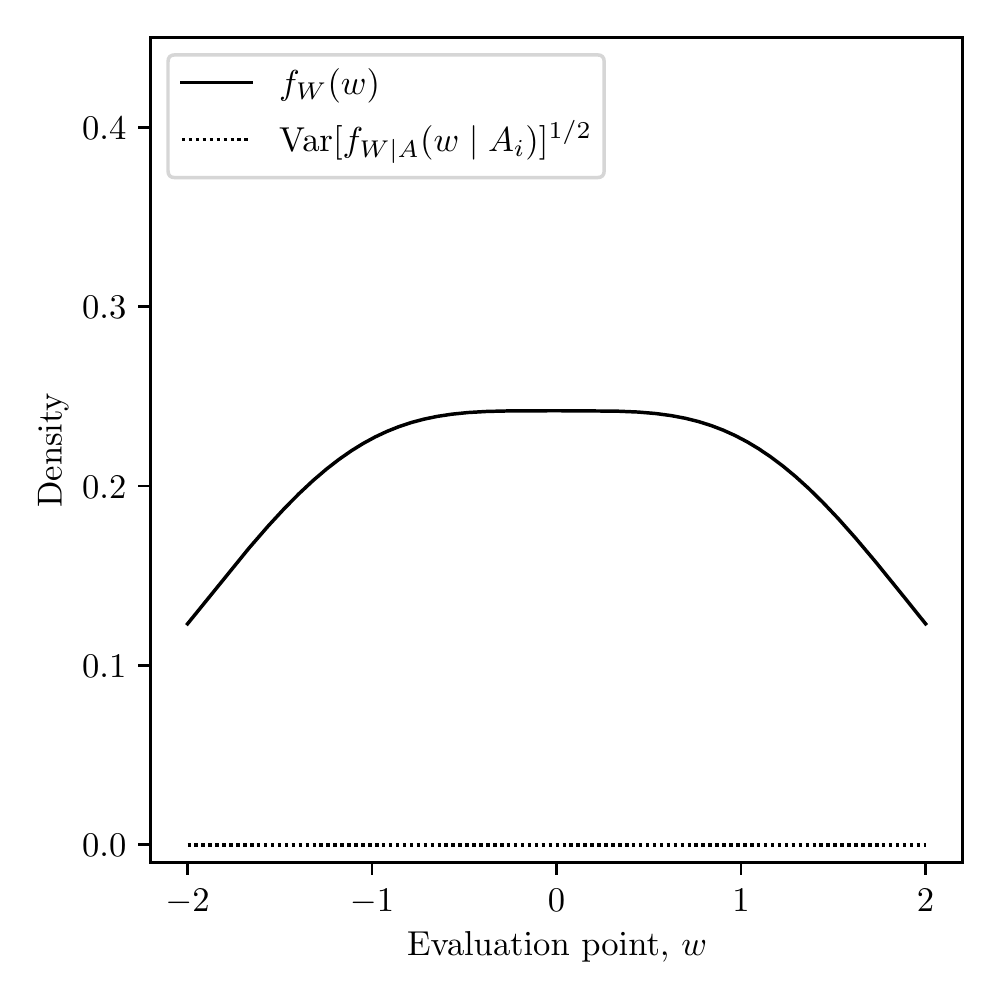}
    \caption{
      Total degeneracy, \\
      $\pi = \left( \frac{1}{2}, 0, \frac{1}{2} \right)$
    }
  \end{subfigure}
  \begin{subfigure}{0.32\textwidth}
    \centering
    \includegraphics[scale=0.52]{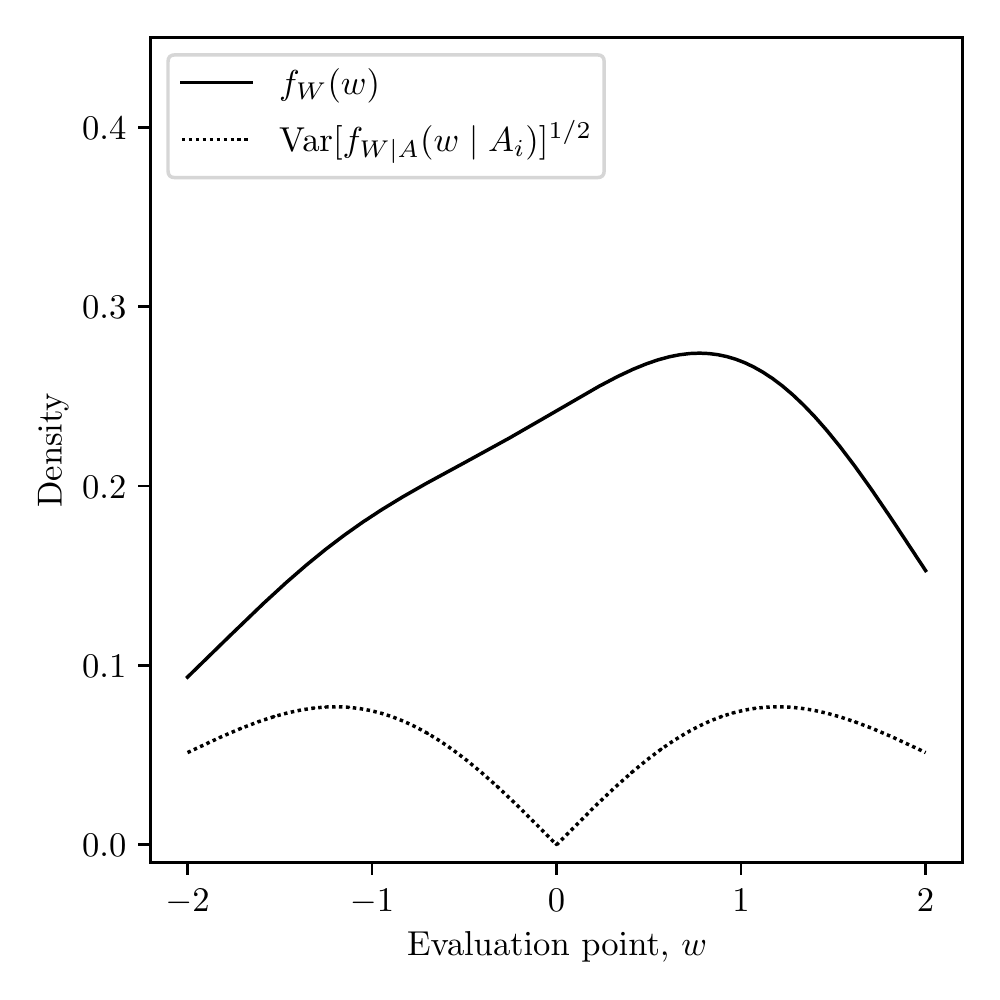}
    \caption{
      Partial degeneracy, \\
      $\pi = \left( \frac{1}{4}, 0, \frac{3}{4} \right)$
    }
  \end{subfigure}
  \begin{subfigure}{0.32\textwidth}
    \centering
    \includegraphics[scale=0.52]{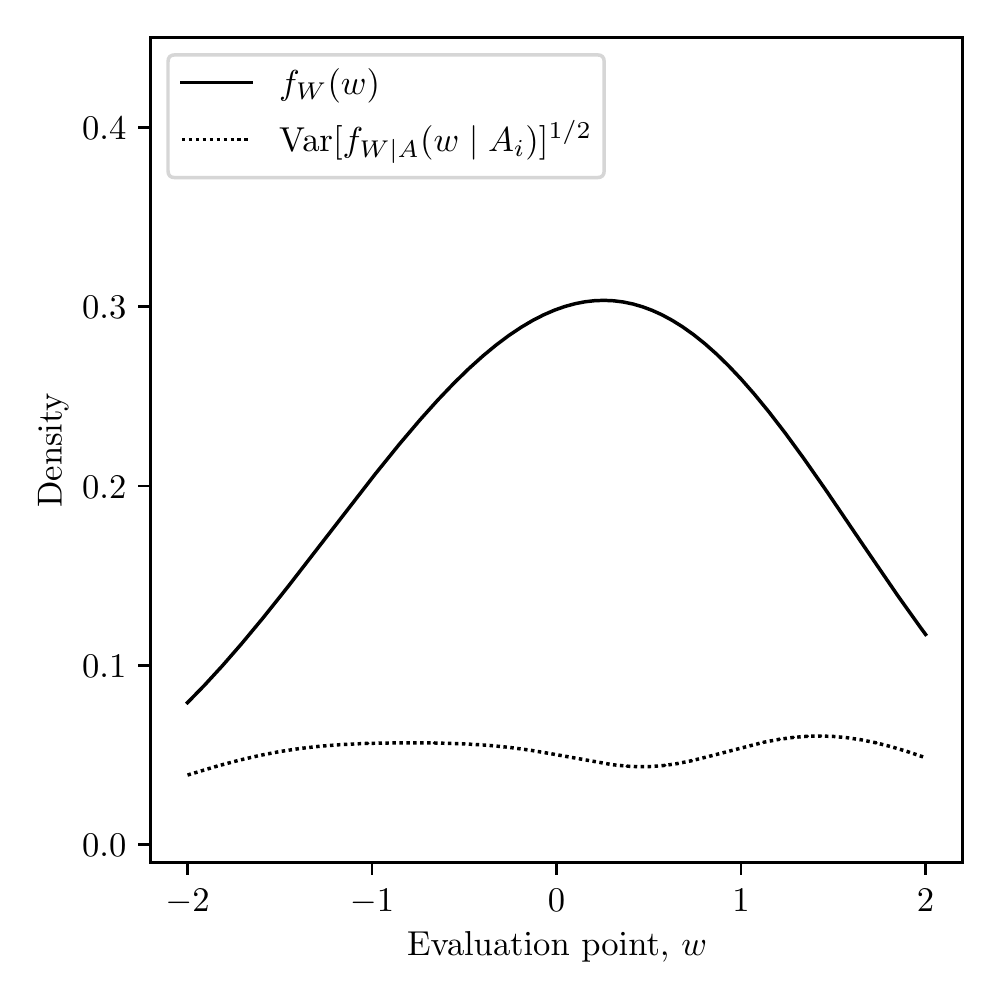}
    \caption{
      No degeneracy, \\
      $\pi = \left( \frac{1}{5}, \frac{1}{5}, \frac{3}{5} \right)$
    }
  \end{subfigure}
  \caption{
    Density $f_W$ and standard deviation
    of $f_{W|A}$ for the family of distributions $\P_\pi$.\\
  }
  \label{fig:distribution}
\end{figure}

In line with Assumption~\ref{ass:data}, $\bA_n$ and $\bV_n$ are i.i.d.\
sequences independent of each other.
Then
$f_{W \mid AA}(w \mid A_i, A_j) = \phi(w - A_i A_j)$,\,
$f_{W \mid A}(w \mid A_i) = \pi_1 \phi(w + A_i) + \pi_2 \phi(w)
+ \pi_3 \phi(w - A_i)$
and
$f_W(w) = (\pi_1^2 + \pi_3^2) \phi(w-1) + \pi_2 (2 - \pi_2) \phi(w) + 2
\pi_1 \pi_3 \phi(w+1),$
where $\phi$ denotes the probability density function of
the standard normal distribution.
Note that $f_W(w)$ is strictly positive for all $w \in \R$.
Consider the parameter choices:
\begin{enumerate}[label=(\roman*)]\onehalfspacing

  \item $\pi = \left( \frac{1}{2}, 0, \frac{1}{2} \right)$:\quad
    $\P_\pi$ is degenerate at all $w \in \R$,

  \item $\pi = \left( \frac{1}{4}, 0, \frac{3}{4} \right)$:\quad
    $\P_\pi$ is degenerate only at $w=0$,

  \item $\pi = \left( \frac{1}{5}, \frac{1}{5}, \frac{3}{5} \right)$:\quad
    $\P_\pi$ is non-degenerate for all $w \in \R$.

\end{enumerate}
Figure~\ref{fig:distribution} demonstrates these phenomena,
plotting the unconditional density $f_W$
and the standard deviation of the conditional density $f_{W|A}$
over $\cW = [-2,2]$ for each choice of the parameter $\pi$.

The trichotomy of total/partial/no  degeneracy is useful for understanding the
distributional properties of the dyadic kernel density estimator
$\widehat{f}_W(w)$. Crucially, our need for uniformity in $w$ complicates the
simpler degeneracy/no degeneracy dichotomy observed previously in the literature
\citep{graham2022kernel}.
More specifically,
from a pointwise-in-$w$ perspective,
partial degeneracy causes no issues,
while it is a fundamental problem
when conducting inference uniformly over $w \in \cW$.
We develop inference methods
that are valid uniformly over $w \in \cW$,
regardless of the presence of partial or total degeneracy.

\section{Point estimation results}\label{sec:point_estimation}

Using Lemma~\ref{lem:variance},
the next theorem establishes the uniform
convergence rate of $\widehat{f}_W$.
\begin{theorem}[Uniform convergence rate]%
  \label{thm:uniform_consistency}%
  Suppose that Assumptions
  \ref{ass:data}~and~\ref{ass:kernel_bandwidth}
  hold. Then
  \begin{align*}
    \E\left[
      \sup_{w \in \cW}
      \big|\widehat{f}_W(w) - f_W(w)\big|
    \right]
    \lesssim
    h^{p\wedge\beta}
    + \frac{\Du}{\sqrt n}
    + \sqrt{\frac{\log n}{n^2h}}.
  \end{align*}
\end{theorem}
The constant in Theorem~\ref{thm:uniform_consistency}
depends only on $\cW$, $\beta$, $\CH$
and the choice of kernel.
We interpret this result in light of the degeneracy trichotomy.
\begin{enumerate}[label=(\roman*)]\onehalfspacing
  \item Partial or no degeneracy:
    $\Du > 0$.
    Any bandwidths satisfying
    $n^{-1} \log n \lesssim h \lesssim n^{-\frac{1}{2(p\wedge\beta)}}$
    yield
    $\E\big[\sup_{w \in \cW}\big|\widehat f_W(w)
      - f_W(w)\big| \big] \lesssim \frac{1}{\sqrt n}$,
    the ``parametric'' bandwidth-independent rate noted by
    \citet{graham2022kernel}.

  \item Total degeneracy:
    $\Du = 0$.
    Minimizing the bound in
    Theorem~\ref{thm:uniform_consistency} with
    $h \asymp \left( \frac{\log n}{n^2} \right)^{\frac{1}{2(p\wedge\beta)+1}}$
    yields
    $\E\big[ \sup_{w \in \cW} \big|\widehat f_W(w) - f_W(w)\big| \big]
    \lesssim
    \big(\frac{\log n}{n^2} \big)^{\frac{p\wedge\beta}{2(p\wedge\beta)+1}}$.
\end{enumerate}

These results generalize \citet*[Theorem~1]{chiang2020empirical}
by allowing for compactly supported data and more general kernel-like
functions $k_h(\cdot,w)$, enabling boundary-adaptive density estimation.

\subsection{Minimax optimality}

We establish the minimax rate under the
supremum norm for density estimation
with dyadic data.
This implies minimax optimality of the kernel density estimator $\widehat
f_W$, regardless of the degeneracy type of the dyadic distribution.

\begin{theorem}[Uniform minimax rate] \label{thm:minimax}

  Fix $\beta \geq 1$ and $\CH > 0$,
  and let $\cW$ be a compact interval with positive Lebesgue measure.
  Define $\cP = \cP(\cW, \beta, \CH)$
  as the class of dyadic distributions
  satisfying Assumption~\ref{ass:data}.
  Define $\cPd$ as the subclass of $\cP$
  containing only those dyadic distributions
  which are totally degenerate on $\cW$ in the sense that
  $\sup_{w \in \cW} \Var\left[f_{W \mid A}(w \mid A_i)\right] = 0$.
  Then we have
  $\inf_{\widetilde f_W} \sup_{\P \in \cP}
  \E_\P\big[\sup_{w \in \cW} \big| \widetilde f_W(w) - f_W(w) \big| \big]
  \asymp \frac{1}{\sqrt n}$
  and
  $\inf_{\widetilde f_W} \sup_{\P \in \cPd}
  \E_\P\big[\sup_{w \in \cW} \big| \widetilde f_W(w) - f_W(w) \big| \big]
  \asymp \big( \frac{\log n}{n^2} \big)^{\frac{\beta}{2\beta+1}}$,
  where $\widetilde f_W$ is any estimator depending only on
  the data $\bW_n = (W_{ij}: 1 \leq i < j \leq n)$
  distributed according to the dyadic distribution $\P$.
  The constants underlying $\asymp$ depend only on
  $\cW$, $\beta$ and $\CH$.

\end{theorem}

Theorem~\ref{thm:minimax} shows that the uniform convergence rate of $n^{-1/2}$
obtained in Theorem~\ref{thm:uniform_consistency}
(coming from the $L_n$ term) is minimax-optimal in general.
When attention is restricted to totally degenerate dyadic distributions,
$\widehat f_W$ also achieves the minimax rate of uniform convergence
(assuming a kernel of sufficiently high order $p \geq \beta$),
which is on the order of
$\left(\frac{\log n}{n^2}\right)^{\frac{\beta}{2\beta+1}}$ and
is determined by the bias $B_n$ and the leading variance term $E_n$ in
\eqref{eq:h-decomposition}.

Combining Theorems \ref{thm:uniform_consistency}~and~\ref{thm:minimax}, we
conclude that the estimator $\widehat{f}_W(w)$ achieves the minimax-optimal rate
of uniform convergence for estimating $f_W(w)$ if
$h \asymp \left( \frac{\log n}{n^2} \right)^{\frac{1}{2\beta+1}}$
and $p \geq \beta$,
whether or not there are
any degenerate points in the underlying data generating process.
This result appears to be new to the literature on
nonparametric estimation with dyadic data.
See \citet{gao2021minimax} for a contemporaneous review.

\section{Distributional results}\label{sec:inference}

We investigate the distributional properties of the
standardized $t$-statistic process
\[ T_n(w) = \frac{\widehat{f}_W(w) - f_W(w)}{\sqrt{\Sigma_n(w,w)}}
  , \qquad w\in\cW,\]
which is not necessarily asymptotically tight.
Therefore, to approximate the distribution of the entire $t$-statistic process,
as well as specific functionals thereof, we rely on a novel strong approximation
approach outlined in this section.
Our results can be used to perform valid uniform inference irrespective of the
degeneracy type.

This section is largely concerned with distributional properties and thus
frequently requires copies of stochastic processes.
For succinctness of notation, we will not differentiate between a process
and its copy, but details are available in
Section~SA3 of the supplemental appendix.

\subsection{Strong approximation}

By the Hoeffding-type decomposition
\eqref{eq:h-decomposition} and Lemma~\ref{lem:variance},
it suffices to consider the distributional properties
of the stochastic process $(L_n(w)+E_n(w):w\in\cW)$.
Our approach combines the K\'omlos--Major--Tusn\'ady (KMT) approximation
\citep{komlos1975approximation} to obtain a strong approximation of
$L_n$ with a Yurinskii approximation
\citep{yurinskii1978error} to obtain a
\emph{conditional} (on $\bA_n$) strong approximation of
$E_n$.
The latter is necessary because
$E_n$ is akin to a local empirical process of i.n.i.d.\
random variables, conditional on $\bA_n$,
and therefore the KMT approximation is not applicable.
These approximations are then combined to give a final
(unconditional) strong approximation for $L_n+E_n$,
and thus for the $t$-statistic process
$(T_n(w):w\in\cW)$.

The following lemma is an application of our generic
KMT approximation result for empirical processes,
given in Section~SA3 of
the online supplemental appendix,
which builds on earlier work by \citet{gine2004kernel}
and \citet{gine2010confidence} and may be of independent interest.

\begin{lemma}[Strong approximation of $L_n$]
  \label{lem:strong_approx_Ln}
  Suppose that Assumptions
  \ref{ass:data}~and~\ref{ass:kernel_bandwidth} hold.
  For each $n$ there exists
  a mean-zero Gaussian process
  $Z^L_n$ indexed on $\cW$ satisfying
  $\E\big[ \sup_{w \in \cW} \big| \sqrt{n} L_n(w) -  Z_n^L(w) \big| \big]
  \lesssim \frac{\Du \log n}{\sqrt{n}}$,
  where $\E[Z_n^L(w)Z_n^L(w')] =  n\E[L_n(w)L_n(w')]$ for all $w, w' \in \cW$.
  The process $Z_n^L$
  is a function only of $\bA_n$ and some random noise
  independent of $(\bA_n, \bV_n)$.
\end{lemma}

The strong approximation result in Lemma~\ref{lem:strong_approx_Ln} would be
sufficient to develop valid and even optimal
uniform inference procedures whenever
(i) $\Dl > 0$ (no degeneracy in $L_n$) and
(ii) $n h \gg \log n$ ($L_n$ is leading).
In this special case, the recent Donsker-type results of
\citet{davezies2021exchangeable}
can be applied to analyze the limiting distribution
of the stochastic process $\widehat{f}_W$.
Alternatively,
again only when $L_n$ is non-degenerate and leading,
standard empirical process methods could also be used.
However, even in the special case when $\widehat{f}_W(w)$ is asymptotically
Donsker, our result in Lemma~\ref{lem:strong_approx_Ln}
improves upon the literature by providing a rate-optimal strong approximation
for $\widehat{f}_W$ as opposed to only a weak convergence result.
See Theorem \ref{thm:infeasible_ucb} and the subsequent
discussion below.

More importantly, as illustrated above,
it is common in the literature to find dyadic
distributions which exhibit partial or total degeneracy, making the process
$\widehat{f}_W$ non-Donsker.
Thus approximating only $L_n$ is in general insufficient
for valid uniform inference,
and it is necessary to capture the distributional properties of $E_n$ as well.
The following lemma is an application of
our strong approximation result for
empirical processes based on the Yurinskii approximation,
which builds on a refinement by \citet{belloni2019conditional}.

\begin{lemma}[Conditional strong approximation of $E_n$]
  \label{lem:conditional_strong_approx_En}

  Suppose Assumptions
  \ref{ass:data}~and~\ref{ass:kernel_bandwidth} hold
  and take any $R_n \to \infty$.
  For each $n$ there exists
  $\widetilde Z^E_n$
  which is a mean-zero Gaussian process
  conditional on $\bA_n$ satisfying
  $\sup_{w \in \cW}
  \big| \sqrt{n^2h} E_n(w) - \widetilde Z_n^E(w) \big|
  \lesssim_\P \frac{(\log n)^{3/8} R_n}{n^{1/4}h^{3/8}}$,
  where $\E[\widetilde Z_n^E(w)\widetilde Z_n^E(w')\bigm\vert \bA_n]
  =n^2h\E[E_n(w)E_n(w')\bigm\vert \bA_n]$
  for all $w, w' \in \cW$.
\end{lemma}

The process $\widetilde Z_n^E$
is a Gaussian process conditional on $\bA_n$ but is not in general a
Gaussian process unconditionally.
The following lemma further constructs an unconditional Gaussian process
$Z_n^E$ that approximates $\widetilde Z_n^E$.

\begin{lemma}[Unconditional strong approximation of $E_n$]
  \label{lem:unconditional_strong_approx_En}

  Suppose that Assumptions
  \ref{ass:data}~and~\ref{ass:kernel_bandwidth} hold.
  For each $n$ there exists
  a mean-zero Gaussian process $Z^E_n$ satisfying
  $\E\big[ \sup_{w \in \cW} \big|\widetilde Z_n^E(w) - Z_n^E(w)\big| \big]
  \lesssim \frac{(\log n)^{2/3}}{n^{1/6}}$,
  where
  $Z_n^E$ is independent of $\bA_n$ and
  $\E[Z_n^E(w)Z_n^E(w')]=\E[\widetilde Z_n^E(w)\widetilde Z_n^E(w')]
  = n^2h \, \E[E_n(w)E_n(w')]$
  for all $w, w' \in \cW$.
\end{lemma}


Combining Lemmas \ref{lem:conditional_strong_approx_En}
and~\ref{lem:unconditional_strong_approx_En},
we obtain an unconditional strong approximation for $E_n$.
The resulting rate of approximation may not be optimal,
due to the Yurinskii coupling, but to the best of our knowledge
it is the first in the literature for the process
$E_n$, and hence for $\widehat{f}_W$ and its associated $t$-process in the
context of dyadic data.
The approximation rate is
sufficiently fast to allow for optimal bandwidth choices; see Section
\ref{sec:implementation} for more details.
Strong approximation results for local empirical processes
(e.g.\ \citealp{gine2010confidence})
are not applicable here because the summands in the
non-negligible $E_n$
are not (conditionally) i.i.d.
Likewise, neither standard empirical process and U-process theory
\citep{van1996weak,gine2021mathematical} nor the recent results in
\citet{davezies2021exchangeable} are applicable
to the non-Donsker process $E_n$.

The previous lemmas showed that $L_n$ is
$\sqrt{n}$-consistent while $E_n$ is $\sqrt{n^2h}$-consistent
(pointwise in $w$),
showcasing the importance of careful standardization
(cf.\ Studentization in Section~\ref{sec:implementation})
for the purpose of rate
adaptivity to the unknown degeneracy type.
In other
words, a challenge in conducting uniform inference is that the
finite-dimensional distributions of the stochastic process $L_n+E_n$,
and hence those of $\widehat{f}_W$ and its associated $t$-process $T_n$,
may converge at different rates at different points $w\in\cW$.
The following theorem provides an
(infeasible)
inference procedure which is fully adaptive to
such potential unknown degeneracy.

\begin{theorem}[Strong approximation of $T_n$]
  \label{thm:strong_approx_fW}

  Suppose that Assumptions~\ref{ass:data}
  and~\ref{ass:kernel_bandwidth} hold and
  $f_W(w) > 0$ on $\cW$, and take any $R_n \to \infty$.
  Then for each $n$ there exists a
  centered Gaussian process $Z_n^{T}$ such that
  \begin{align*}
    &\sup_{w \in \cW}
    \left|
    T_n(w)
    -
    Z_n^{T}(w)
    \right|
    \lesssim_\P
    \frac{
      n^{-1} \log n
      + n^{-5/4} h^{-7/8} (\log n)^{3/8} R_n
      + n^{-7/6} h^{-1/2} (\log n)^{2/3}
      + h^{p\wedge\beta}}
    {\Dl/\sqrt{n} + 1/\sqrt{n^2h}},
  \end{align*}
  where $\E[Z_n^T(w)Z_n^T(w')] = \E[T_n(w)T_n(w')]$
  for all $w,w' \in \cW$.
\end{theorem}

The first term in the numerator corresponds to
the strong approximation error for $L_n$ in Lemma~\ref{lem:strong_approx_Ln}
and the error introduced by $Q_n$.
The second and third terms correspond to the conditional
and unconditional strong approximation errors for
$E_n$ in Lemmas \ref{lem:conditional_strong_approx_En}
and \ref{lem:unconditional_strong_approx_En}, respectively.
The fourth term is from the smoothing bias result
in Lemma~\ref{lem:variance}.
The denominator is the lower bound on the standard
deviation $\Sigma_n(w,w)^{1/2}$
formulated in Lemma~\ref{lem:variance_bounds}.

In the absence of degenerate points ($\Dl > 0$) and if
$n h^{7/2}\gtrsim 1$, Theorem~\ref{thm:strong_approx_fW}
offers a strong approximation of the $t$-process at the rate
$(\log n)/\sqrt{n}+\sqrt{n}h^{p\wedge\beta}$, which matches the celebrated KMT
approximation rate for i.i.d.\ data plus the
smoothing bias.
Therefore, our novel $t$-process strong
approximation can achieve the optimal KMT rate for
non-degenerate dyadic
distributions provided that $p\wedge\beta \geq 3.5$.
This is achievable
if a fourth-order (boundary-adaptive) kernel is used
and $f_W$ is sufficiently smooth.

In the presence of partial or total degeneracy ($\Dl =0$),
Theorem~\ref{thm:strong_approx_fW} provides a strong approximation for the
$t$-process at the rate
$\sqrt{h}\log n + n^{-1/4}h^{-3/8}(\log n)^{3/8} R_n +
n^{-1/6}(\log n)^{2/3} + n h^{1/2+p\wedge\beta}$.
If, for example, $n h^{p\wedge\beta}\lesssim 1$,
then our result can achieve a strong
approximation rate of $n^{-1/7}$
up to $\log n $ terms.
Theorem~\ref{thm:strong_approx_fW}
appears to be the first in the dyadic literature
which is also robust to the presence of (unknown)
degenerate points in the underlying dyadic distribution.

\subsection{Application: confidence bands}

Theorem~\ref{thm:infeasible_ucb} constructs standardized confidence
bands for $f_W$ which are infeasible as they depend on the unknown population
variance $\Sigma_n$.
In Section~\ref{sec:implementation}
we will make this inference procedure
feasible by proposing a valid estimator of the covariance function $\Sigma_n$
for Studentization, as well as developing
bandwidth selection and robust bias correction methods.
Before presenting our result on valid infeasible uniform confidence bands,
we first impose in Assumption~\ref{ass:rates}
some extra restrictions on the bandwidth sequence,
which depend on the degeneracy type of the dyadic distribution,
to ensure the coverage rate converges in large samples.

\begin{assumption}[Rate restriction for uniform confidence bands]
  \label{ass:rates}\onehalfspacing

  Assume that one of the following holds:

  \begin{enumerate}[label=(\roman*)]

    \item
      No degeneracy
      ($\Dl > 0$):
      $n^{-6/7} \log n \ll h
      \ll (n \log n)^{-\frac{1}{2(p \wedge \beta)}}$,

    \item
      Partial or total degeneracy
      ($\Dl = 0$):
      $n^{-2/3} (\log n)^{7/3} \ll h
      \ll (n^2 \log n)^{-\frac{1}{2(p \wedge \beta) + 1}}$.
  \end{enumerate}
\end{assumption}

We now construct the infeasible uniform confidence bands.
For $\alpha \in (0,1)$, let $q_{1-\alpha}$ be the quantile satisfying
$ \P\left(\sup_{w \in \cW}
  \left| Z_n^T(w) \right|
  \leq q_{1-\alpha} \right)
= 1 - \alpha$.
The following result employs the anti-concentration idea due to
\citet{chernozhukov2014anti} to deduce valid standardized confidence bands,
where we approximate the quantile of the unknown finite sample distribution of
$\sup_{w\in\cW} |T_n(w)|$ by the quantile $q_{1-\alpha}$ of
$\sup_{w\in\cW}|Z_n^T(w)|$. This approach offers a better rate of
convergence than relying on extreme value theory
for the distributional approximation,
hence improving the finite sample performance of the proposed confidence bands.

\begin{theorem}[Infeasible uniform confidence bands]
  \label{thm:infeasible_ucb}

  Suppose that Assumptions~\ref{ass:data},~\ref{ass:kernel_bandwidth}
  and~\ref{ass:rates} hold and
  $f_W(w) > 0$ on $\cW$.
  Then
  \begin{align*}
    \P\left(
      f_W(w)
      \in
      \left[
        \widehat f_W(w)
        \pm
        q_{1-\alpha}
        \sqrt{\Sigma_n(w,w)}
        \, \right]
      \, \textup{for all }
      w \in \cW
    \right)
    \to 1 - \alpha.
  \end{align*}
\end{theorem}

By Theorem~\ref{thm:uniform_consistency},
the asymptotically optimal choice of bandwidth
for uniform convergence is
$h \asymp ((\log n)/n^2)^{\frac{1}{2(p \wedge \beta)+1}}$.
As discussed in the next section, the approximate
IMSE-optimal bandwidth is
$h \asymp (1/n^2)^{\frac{1}{2(p \wedge \beta)+1}}$.
Both bandwidth choices satisfy Assumption~\ref{ass:rates}
only in the case of no degeneracy.
The degenerate cases in Assumption~\ref{ass:rates}(ii),
which require $p \wedge \beta > 1$,
exhibit behavior more similar to that of standard
nonparametric kernel-based estimation
and so the aforementioned optimal bandwidth
choices will lead to a non-negligible smoothing bias in the distributional
approximation for $T_n$.
Different approaches are available in the literature to address this issue,
including undersmoothing or ignoring the bias \citep{hall2001bootstrapping},
bias correction \citep{hall1992effect}, robust bias correction
\citep{Calonico-Cattaneo-Farrell_2018_JASA,
  Calonico-Cattaneo-Farrell_2022_Bernoulli}
and Lepskii's method \citep{lepskii1992asymptotically,birge2001alternative},
among others.
In the next section we develop a feasible uniform inference procedure,
based on robust bias correction methods,
which amounts to first selecting an optimal bandwidth for the point estimator
$\widehat{f}_W$ using a $p$th-order kernel,
and then correcting the bias of the point estimator while also adjusting
the standardization (Studentization) when forming the $t$-statistic $T_n$.

Importantly, regardless of the specific implementation details,
Theorem~\ref{thm:infeasible_ucb} shows that any bandwidth sequence $h$
satisfying both (i) and (ii)
in Assumption~\ref{ass:rates} leads to valid uniform inference which is robust
and adaptive to the (unknown) degeneracy type.

\section{Implementation}\label{sec:implementation}

We address outstanding implementation details to make our main
uniform inference results feasible.
In Section~\ref{sec:covariance_estimation} we propose a covariance estimator
along with a modified version which is guaranteed to be positive semi-definite.
This allows for the construction of fully feasible confidence bands
in Section~\ref{sec:feasible_confidence_bands}.
In Section~\ref{sec:bandwidth_selection} we discuss bandwidth selection
and formalize our procedure for
robust bias correction inference.

\subsection{Covariance function estimation}
\label{sec:covariance_estimation}

Define the following plug-in covariance function
estimator of $\Sigma_n$: for $w, w' \in \cW$,
\begin{align*}
  \widehat \Sigma_n(w,w')
  &=
  \frac{4}{n^2}
  \hspace*{-0.5mm}
  \sum_{i=1}^n
  \hspace*{-0.5mm}
  S_i(w) S_i(w')
  \hspace*{-0.5mm}
  -
  \hspace*{-0.5mm}
  \frac{4}{n^2(n-1)^2}
  \hspace*{-0.5mm}
  \sum_{i<j}
  \hspace*{-0.5mm}
  k_h(W_{ij},w)
  k_h(W_{ij},w')
  \hspace*{-0.5mm}
  -
  \hspace*{-0.5mm}
  \frac{4n-6}{n(n-1)}
  \widehat f_W(w)
  \widehat f_W(w'),
\end{align*}
where
$S_i(w) = \frac{1}{n-1} \big( \sum_{j = 1}^{i-1} k_h(W_{j i}, w)
+ \sum_{j = i+1}^n k_h(W_{ij}, w) \big)$
is an ``estimator'' of
$\E[k_h(W_{ij},w) \mid A_i]$.
Though $\widehat\Sigma_n(w,w')$ is consistent in an appropriate sense
as shown in Lemma~\ref{lem:SDP} below,
it is not necessarily positive semi-definite,
even in the limit.
We therefore propose a modified covariance estimator
which is guaranteed to be positive semi-definite.
Specifically, consider the following optimization problem
where $\Ck$ and $\CL$ are as in Section~\ref{sec:setup}.
\begin{equation}
  \label{eq:SDP}
  \begin{aligned}
    \minimize
    \qquad
     & \sup_{w,w' \in \cW}
    \left|
    \frac{M(w,w') - \widehat\Sigma_n(w,w')}
    {\sqrt{\widehat \Sigma_n(w,w) + \widehat \Sigma_n(w',w')}}
    \right|
    \quad \textup{ over } M: \cW \times \cW \to \R
    \\
    \subjectto
    \qquad
     & M \textup{ is symmetric and positive semi-definite}, \\
     & \big|M(w,w') - M(w, w'')\big|
    \leq \frac{4}{n h^3}
    \Ck \CL
    |w'-w''|
    \textup{ for all }
    w, w', w'' \in \cW.
  \end{aligned}
\end{equation}

Denote by $\widehat\Sigma_n^+$ any (approximately) optimal solution to
(\ref{eq:SDP}). The following lemma establishes uniform convergence rates
for both $\widehat \Sigma_n$ and $\widehat \Sigma_n^+$. It allows us to use
these estimators to construct feasible versions of $T_n$ and its associated
Gaussian approximation $Z_n^{T}$ defined in Theorem~\ref{thm:strong_approx_fW}.
\begin{lemma}[Consistency of $\widehat \Sigma_n$ and $\widehat \Sigma_n^+$]
  \label{lem:SDP}
  Suppose that Assumptions~\ref{ass:data}
  and~\ref{ass:kernel_bandwidth} hold
  and that
  $n h \gtrsim \log n$ and
  $f_W(w) > 0$ on $\cW$.
  Then
  \begin{align*}
    \sup_{w,w' \in \cW}
    \left|
    \frac{\widehat \Sigma_n(w,w') - \Sigma_n(w,w')}
    {\sqrt{\Sigma_n(w,w) + \Sigma_n(w',w')}}
    \right|
    &\lesssim_\P
    \frac{\sqrt{\log n}}{n}.
  \end{align*}
  Also, the optimization problem (\ref{eq:SDP})
  is a semi-definite program
  \citep[SDP,][]{laurent2005semidefinite}
  and has an approximately
  optimal solution
  $\widehat\Sigma_n^+$
  satisfying
  \begin{align*}
    \sup_{w,w' \in \cW}
    \left|
    \frac{\widehat \Sigma_n^+(w,w') - \Sigma_n(w,w')}
    {\sqrt{\Sigma_n(w,w) + \Sigma_n(w',w')}}
    \right|
    &\lesssim_\P
    \frac{\sqrt{\log n}}{n}.
  \end{align*}
\end{lemma}

In practice we take $w, w' \in \cW_d$ where $\cW_d$ is a finite
subset of $\cW$, typically taken to be an equally-spaced grid.
This yields finite-dimensional covariance matrices,
for which \eqref{eq:SDP}
can be solved in polynomial time in $|\cW_d|$ using
a general-purpose SDP solver
\citep[e.g.\ by interior point methods,][]{laurent2005semidefinite}.
The number of points in $\cW_d$ should be taken as large as is
computationally practical in order to generate
confidence bands rather than merely simultaneous confidence intervals.
It is worth noting that the complexity of solving \eqref{eq:SDP}
does not depend on the number of vertices $n$,
and so does not influence the ability
of our methodology to handle large and possibly sparse networks.

The bias-corrected variance estimator in
\citet[Section~3.2]{MatsushitaOtsu2021} takes a similar form
to our estimator $\widehat\Sigma_n$ but in the parametric setting,
and is therefore also not guaranteed
to be positive semi-definite in finite samples.
Our approach addresses this issue,
ensuring a positive semi-definite estimator
$\widehat\Sigma_n^+$ is always available.

\subsection{Feasible confidence bands}
\label{sec:feasible_confidence_bands}

Given a choice of the kernel order $p$
and a bandwidth $h$, we construct a valid confidence
band that is implementable in practice.
Define the Studentized $t$-statistic process
\[ \widehat T_n(w) = \frac{\widehat{f}_W(w) - f_W(w)}{\sqrt{\widehat
      \Sigma_n^+(w,w)}} , \qquad w\in\cW.\]
Let $\widehat Z_n^T(w)$ be a process which,
conditional on the data $\bW_n$,
is mean-zero and Gaussian, whose
conditional covariance structure is
$\E\big[ \widehat Z_n^T(w) \widehat Z_n^T(w') \bigm\vert \bW_n \big]
= \frac{\widehat \Sigma_n^+(w,w')}
{\sqrt{\widehat \Sigma_n^+(w,w) \widehat \Sigma_n^+(w',w')}}$.
For $\alpha \in (0,1)$, let $\widehat q_{1-\alpha}$ be the
conditional quantile satisfying
$ \P\big(\sup_{w \in \cW}
\big| \widehat Z_n^T(w) \big|
\leq \widehat q_{1-\alpha}
\bigm\vert \bW_n \big)
= 1 - \alpha$,
which is shown to be well-defined in
Section~SA5 of
the online supplemental appendix.

\begin{theorem}[Feasible uniform confidence bands]
  \label{thm:ucb}

  Suppose that Assumptions \ref{ass:data}, \ref{ass:kernel_bandwidth}
  and \ref{ass:rates} hold and $f_W(w) > 0$ on $\cW$. Then
  \begin{align*}
    \P\left(
      f_W(w)
      \in
      \left[
        \widehat f_W(w)
        \pm
        \widehat q_{1-\alpha}
        \sqrt{\widehat\Sigma_n^+(w,w)}
        \,\right]
      \,\textup{for all }
      w \in \cW
    \right)
    \to 1 - \alpha.
  \end{align*}
\end{theorem}

Recently, \cite{chiang2022inference} derived high-dimensional central limit
theorems over rectangles for exchangeable arrays and applied them to construct
simultaneous confidence intervals for a sequence of design points. Their
inference procedure relies on the multiplier bootstrap,
and their conditions for valid inference
depend on the number of design points considered.
In contrast, Theorem~\ref{thm:ucb} constructs a feasible uniform confidence band
over the entire domain of inference $\cW$ based on our strong
approximation results for the whole $t$-statistic process
and the covariance estimator $\widehat\Sigma_n^+$.
The required rate condition specified in
Assumption~\ref{ass:rates} does not depend on the number of design points.
Furthermore, our proposed inference methods are robust to potential
unknown degenerate points in the underlying dyadic data generating process.

In practice, suprema over $\cW$ can be
replaced by maxima over sufficiently many design points in $\cW$.
The conditional quantile $\widehat q_{1-\alpha}$
can be estimated by Monte Carlo simulation,
resampling from the Gaussian process defined by the law of
$\widehat Z_n^T \mid \bW_n$.

The bandwidth restrictions in Theorem~\ref{thm:ucb} are the same
as those required for the infeasible version
given in Theorem~\ref{thm:infeasible_ucb}, namely those imposed in Assumption
\ref{ass:rates}. This follows from the rates of convergence obtained in
Lemma~\ref{lem:SDP}, coupled with some careful technical work given in the
supplemental appendix to handle the potential presence
of degenerate points in $\Sigma_n$.

\subsection{Bandwidth selection and robust bias-corrected inference}
\label{sec:bandwidth_selection}

We give some practical suggestions for selecting
the bandwidth parameter $h$.
Let $\nu(w)$ be a non-negative real-valued function on $\cW$ and
suppose we use a kernel of order $p < \beta$
of the form $k_h(s,w) = K\big((s-w) / h\big)/h$.
Then the $\nu$-weighted asymptotic IMSE (AIMSE) is minimized by
\begin{align*}
  h^*_{\AIMSE}
  &=
  \left(
    \frac{
    p!(p-1)!
    \Big(\int_\cW f_W(w) \nu(w) \diff{w}\Big)
    \Big(\int_\R K(w)^2 \diff{w}\Big)
    }
    {
    2
    \Big(
    \int_{\cW}
    f_W^{(p)}(w)^2
    \nu(w)
    \diff{w}
    \Big)
    \Big(
    \int_\R
    w^p K(w)
    \diff{w}
    \Big)^2
    }
  \right)^{\frac{1}{2p+1}}
  \left(
    \frac{n(n-1)}{2}
  \right)^{-\frac{1}{2p+1}}.
\end{align*}
This is akin to the AIMSE-optimal bandwidth choice for traditional
monadic kernel density estimation
with a sample size of $\frac{1}{2}n(n-1)$.
The choice $h^*_{\AIMSE}$ is slightly undersmoothed
(up to a polynomial $\log n$ factor)
relative to the uniform minimax-optimal bandwidth choice discussed in
Section~\ref{sec:point_estimation}, but it is easier to implement in practice.

To implement the AIMSE-optimal bandwidth choice,
we propose a simple \emph{rule-of-thumb} (ROT)
approach based on Silverman's rule.
Suppose $p\wedge\beta=2$ and let $\widehat\sigma^2$ and $\widehat\IQR$
be the sample variance and sample interquartile range respectively
of the data $\bW_n$. Then
$\widehat{h}_{\ROT} = C(K) \big( \widehat\sigma \wedge
\frac{\widehat\IQR}{1.349} \big) \big(\frac{n(n-1)}{2} \big)^{-1/5}$,
where $C(K)=2.576$ for the triangular kernel $K(w) = (1 - |w|) \vee 0$,
and $C(K)=2.435$ for the Epanechnikov kernel
$K(w) = \frac{3}{4}(1 - w^2) \vee 0$.

The AIMSE-optimal bandwidth selector $h^*_{\AIMSE}\asymp n^{-\frac{2}{2p+1}}$
and any of its feasible estimators only satisfy
Assumption~\ref{ass:rates} in the case of no degeneracy ($\Dl>0$).
Under partial or total degeneracy,
such bandwidths are not valid due to the usual leading
smoothing
(or misspecification) bias of the distributional approximation. To
circumvent this problem and construct feasible uniform confidence bands
for $f_W$, we employ the following robust bias correction approach.

Firstly, estimate the bandwidth $h^*_{\AIMSE}\asymp n^{-\frac{2}{2p+1}}$ using a
kernel of order $p$, which leads to an AIMSE-optimal point estimator
$\widehat{f}_W$ in an $L^2(\nu)$ sense. Then use this bandwidth and a kernel of
order $p' > p$ to construct the statistic $\widehat T_n$ and the confidence band
as detailed in Section~\ref{sec:feasible_confidence_bands}.
Importantly, both $\widehat{f}_W$ and $\widehat{\Sigma}^+_n$ are
recomputed with the new higher-order kernel.
The change in centering is equivalent to a bias correction
of the original AIMSE-optimal point estimator,
while the change in scale captures the additional variability
introduced by the bias correction itself.
As shown formally in
\citet{Calonico-Cattaneo-Farrell_2018_JASA,
  Calonico-Cattaneo-Farrell_2022_Bernoulli}
for the case of kernel-based density estimation with i.i.d.\ data,
this approach leads to higher-order refinements in the distributional
approximation whenever additional smoothness is available ($p'\leq\beta$).
In the present dyadic setting, this procedure is valid so long as
$n^{-2/3} (\log n)^{7/3} \ll n^{-\frac{2}{2p+1}}
\ll (n^2 \log n)^{-\frac{1}{2p' + 1}}$,
which is equivalent to
$2 \leq p < p'$.
For concreteness, we recommend taking
$p = 2$ and $p' = 4$,
and using the rule-of-thumb bandwidth choice $\widehat{h}_{\ROT}$ defined above.
In particular, this approach automatically delivers a KMT-optimal
strong approximation whenever there are no degeneracies in the
underlying dyadic data generating process.

Our feasible robust bias correction method based on AIMSE-optimal dyadic
kernel density estimation for constructing uniform confidence bands
for $f_W$ is summarized in Algorithm~\ref{alg:method}.

\begin{algorithm}[ht]
  \caption{Feasible uniform confidence bands for
    dyadic kernel density estimation}
  \label{alg:method}

  Choose a kernel $k_h$ of order $p \geq 2$ satisfying
  Assumption~\ref{ass:kernel_bandwidth}. \\

  Select a bandwidth $h \approx h^*_{\AIMSE}$ for $k_h$
  as in Section~\ref{sec:bandwidth_selection},
  perhaps using $h = \widehat{h}_{\ROT}$. \\

  Choose another kernel $k_h'$ of order $p'>p$ satisfying
  Assumption~\ref{ass:kernel_bandwidth}.

  For $d \geq 1$, choose a set of $d$
  distinct evaluation points $\cW_d$. \\

  For each $w \in \cW_d$, construct the density estimate $\widehat f_W(w)$
  using $k'_{h}$ as in Section~\ref{sec:introduction}. \\

  For $w, w' \in \cW_d$, construct the
  covariance estimate $\widehat \Sigma_n(w,w')$
  using $k'_{h}$ as in Section~\ref{sec:covariance_estimation}. \\

  Construct the $d \times d$ positive semi-definite
  covariance estimate $\widehat \Sigma_n^+$
  as in Section~\ref{sec:covariance_estimation}. \\

  For $B \geq 1$,
  let $(\widehat Z_{n,r}^T: 1\leq r\leq B)$ be i.i.d.\ Gaussian vectors
  from $\widehat{Z}_n^T$ as in
  Section~\ref{sec:feasible_confidence_bands}. \\

  For $\alpha \in (0,1)$, set
  $\widehat q_{1-\alpha} = \inf_{q \in \R}
  \{ q :
  \# \{r: \max_{w\in\cW_d}|\widehat Z_{n,r}^T(w)|
  \leq q \}
  \geq B(1-\alpha) \}$. \\

  Construct
  $ \big[\widehat f_W(w) \pm
    \widehat q_{1-\alpha} \widehat\Sigma_n^+(w,w)^{1/2} \big]$
  for each $w \in \cW_d$.
\end{algorithm}

\section{Simulations}
\label{sec:simulations}

We investigate the
empirical finite-sample performance of
the kernel density estimator with dyadic data
using simulations.
The family of dyadic distributions
defined in Section~\ref{sec:degeneracy},
along with its three parametrizations,
is used to generate data sets with different degeneracy types.

We use two different boundary bias-corrected
Epanechnikov kernels of orders $p=2$ and $p=4$ respectively,
on the inference domain $\cW = [-2,2]$.
We select an optimal bandwidth for $p=2$ as recommended in
Section~\ref{sec:bandwidth_selection},
using the rule-of-thumb with $C(K) = 2.435$.
The semi-definite program in Section~\ref{sec:covariance_estimation}
is solved with the MOSEK interior point optimizer \citep{mosek},
ensuring positive semi-definite covariance estimates.
Gaussian vectors are resampled
$B = 10\,000$ times.

\begin{figure}[H]
  \centering
  \begin{subfigure}{0.32\textwidth}
    \centering
    \includegraphics[scale=0.52]{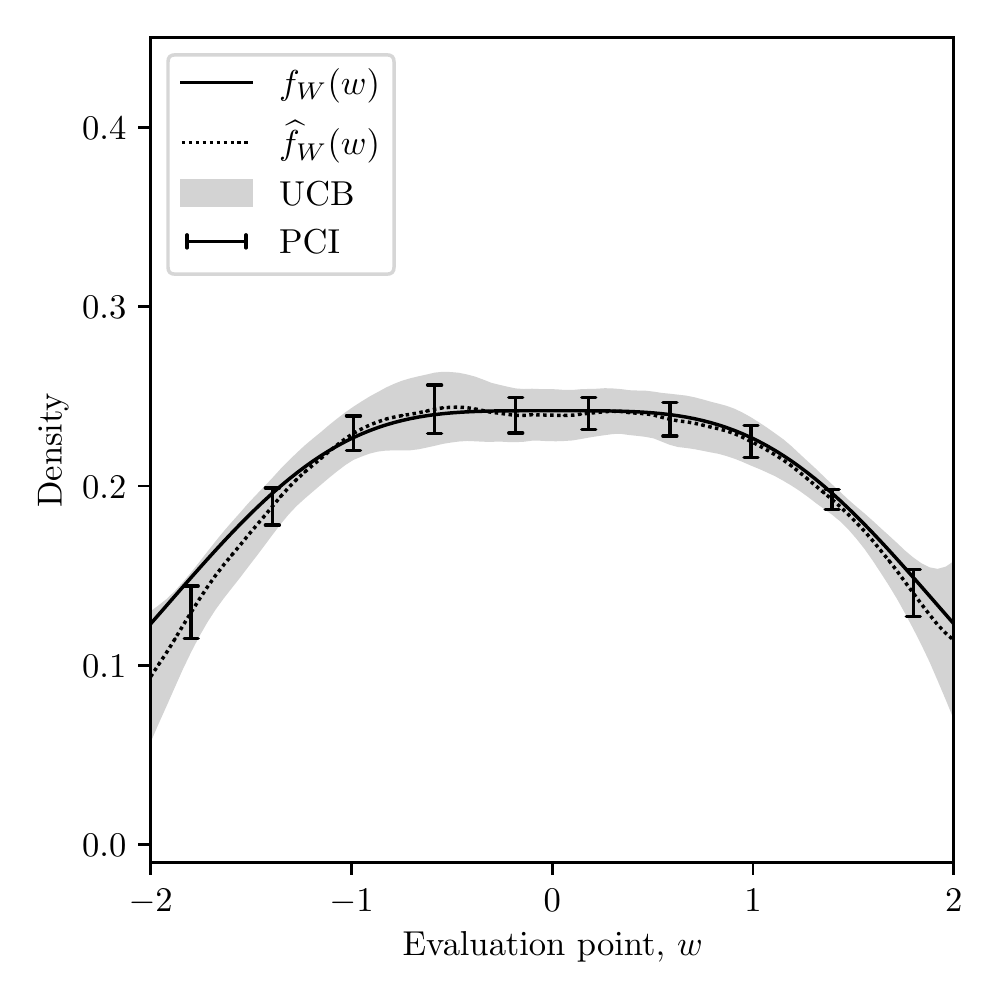}
    \caption{
      Total degeneracy, \\
      $\pi = \left( \frac{1}{2}, 0, \frac{1}{2} \right)$
    }
  \end{subfigure}
  \begin{subfigure}{0.32\textwidth}
    \centering
    \includegraphics[scale=0.52]{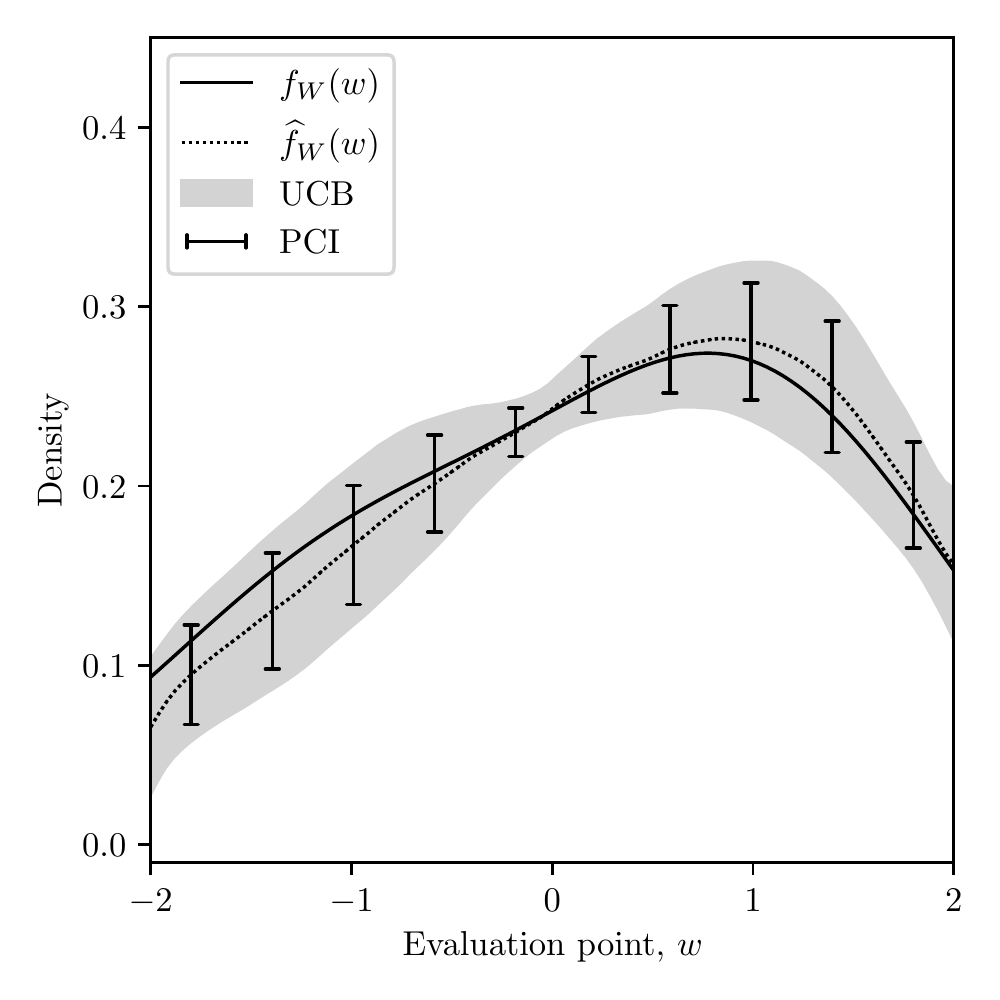}
    \caption{
      Partial degeneracy, \\
      $\pi = \left( \frac{1}{4}, 0, \frac{3}{4} \right)$
    }
  \end{subfigure}
  \begin{subfigure}{0.32\textwidth}
    \centering
    \includegraphics[scale=0.52]{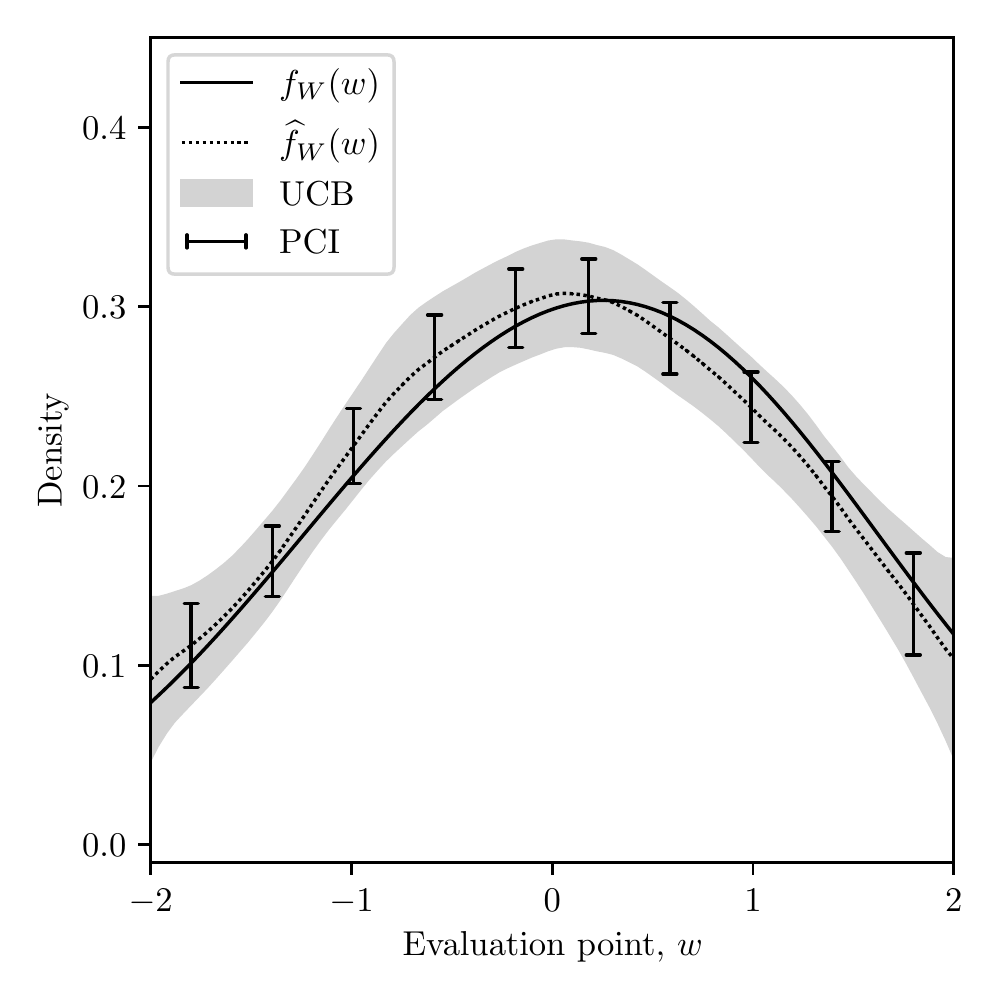}
    \caption{
      No degeneracy, \\
      $\pi = \left( \frac{1}{5}, \frac{1}{5}, \frac{3}{5} \right)$
    }
  \end{subfigure}
  \caption{
    Typical outcomes for three different values of the parameter $\pi$.
  }
  \label{fig:results}
  \vspace{-.8em}
\end{figure}

In Figure~\ref{fig:results} we plot a typical outcome
for each of the three degeneracy types (total, partial, none),
using the Epanechnikov kernel of order $p=2$,
with sample size $n=100$
(so $N=4950$ pairs of nodes)
and with $d=100$ equally-spaced evaluation points.
Each plot contains the true density function $f_W$,
the dyadic kernel density estimate $\widehat f_W$
and two different approximate $95\%$ confidence bands for $f_W$.
The first is the uniform confidence band (UCB)
constructed using one of our main results, Theorem~\ref{thm:ucb}.
The second is a sequence of pointwise confidence intervals (PCI)
constructed by finding a confidence interval for
each evaluation point separately.
We show only $10$ pointwise confidence intervals for clarity.
In general, the PCIs are too narrow as they fail to
provide simultaneous (uniform) coverage over the evaluation points.
Note that
under partial degeneracy the confidence band narrows near the
degenerate point $w = 0$.

Next, Table~\ref{tab:results} presents numerical results.
For each degeneracy type (total, partial, none)
and each kernel order ($p=2$, $p=4$),
we run $2000$ repeats
with sample size $n=3000$ (giving $N=4\,498\,500$ pairs of nodes)
and with $d=50$ equally-spaced evaluation points.
We record the average rule-of-thumb bandwidth $\widehat{h}_{\ROT}$
and the average root integrated mean squared error (RIMSE).
For both the uniform confidence bands (UCB)
and the pointwise confidence intervals (PCI),
we report the coverage rate (CR) and the average width (AW).
The lower-order kernel ($p=2$) ignores the bias, leading to
good RIMSE performance and acceptable UCB coverage
under partial or no degeneracy,
but gives invalid inference under total degeneracy.
In contrast, the higher-order kernel ($p=4$)
provides robust bias correction
and hence improves the coverage of the UCB in every regime,
particularly under total degeneracy,
at the cost of increasing both the RIMSE
and the average widths of the confidence bands.
As expected, the pointwise (in $w\in\cW$) confidence intervals (PCIs)
severely undercover in every regime.
Thus our simulation results show that the proposed feasible
inference methods based on robust bias correction and proper Studentization
deliver valid uniform inference
which is robust to unknown degenerate points
in the underlying dyadic distribution.

\begin{table}[H]
  \caption{Numerical results for three values of the parameter $\pi$.}
  \vspace{-1em}
  \input{table.tex}
  \label{tab:results}
\end{table}

\section{Application: counterfactual dyadic density estimation}
\label{sec:counterfactual}

To further showcase the applicability of our main results,
we develop a kernel density estimator for
dyadic counterfactual distributions.
The aim of such counterfactual analysis is to estimate the
distribution of an outcome variable had some covariates
followed a distribution different from the actual one,
and it is important in causal inference and program evaluation settings
\citep{dinardo1996distribution,chernozhukov2013inference}.

For each $r \in \{0,1\}$,
let $\bW_n^r$, $\bA_n^r$ and $\bV_n^r$ be random variables
as defined in Assumption~\ref{ass:data} and
$\bX_n^r = (X_1^r, \ldots, X_n^r)$ be some covariates.
We assume that $(A_i^r, X_i^r)$ are
independent over $1 \leq i \leq n$
and that $\bX_n^r$ is independent of $\bV_n^r$,
that $W_{ij}^r \mid X_i^r, X_j^r$
has a conditional Lebesgue density
$f_{W \mid XX}^r(\,\cdot \mid x_1, x_2) \in \cH^\beta_{\CH}(\cW)$,
that $X_i^r$ follows a distribution function $F_X^r$
on a common support $\cX$, and that
$(\bA_n^0, \bV_n^0, \bX_n^0)$
is independent of $(\bA_n^1, \bV_n^1, \bX_n^1)$.

We interpret $r$ as an index for two populations,
labeled $0$ and $1$.
The counterfactual density of the outcome of
population $1$ had it followed the same covariate distribution
as population $0$ is
\begin{align*}
  f_W^{1 \triangleright 0}(w)
  &=
  \E\left[
    f_{W \mid XX}^1\big(w \mid X_1^0, X_2^0\big)
  \right]
  = \int_{\cX}
  \int_{\cX}
  f_{W \mid XX}^{1}(w \mid x_1, x_2)
  \psi(x_1)
  \psi(x_2)
  \diff F_X^{1}(x_1)
  \diff F_X^{1}(x_2),
\end{align*}
where $\psi(x) = \mathrm{d} F_X^0(x) / \mathrm{d} F_X^1(x)$ for $x \in \cX$
is a Radon--Nikodym derivative.
If $X^0_i$ and $X^1_i$ have Lebesgue densities,
it is natural to consider a parametric model of the form
$\mathrm{d} F_X^{r}(x)=f_X^r(x;\theta)\diff x$
for some finite-dimensional parameter $\theta$.
Alternatively, if the covariates $X_n^r$ are discrete and have a positive
probability mass function $p_X^r(x)$ on a finite
support $\cX$, the object of interest becomes
$f_W^{1 \triangleright 0}(w)
= \sum_{x_1 \in \cX} \sum_{x_2 \in \cX}
f_{W \mid XX}^{1}(w \mid x_1, x_2) \psi(x_1) \psi(x_2)
p_X^{1}(x_1) p_X^{1}(x_2)$,
where $\psi(x) = p_X^0(x)/p_X^1(x)$ for $x \in \cX$.
We consider discrete covariates for simplicity,
and hence the counterfactual dyadic kernel density estimator is
\begin{align*}
  \widehat f_W^{\,1 \triangleright 0}(w)
  &=
  \frac{2}{n(n-1)}
  \sum_{i=1}^{n-1}
  \sum_{j=i+1}^n
  \widehat \psi(X_i^1)
  \widehat \psi(X_j^1)
  k_h(W_{ij}^1, w),
\end{align*}
where $\widehat\psi(x) = \widehat p_X^{\,0}(x) / \widehat p_X^{\,1}(x)$ and
$\widehat p_X^{\,r}(x) = \frac{1}{n}\sum_{i = 1}^n \I\{X_i^r = x\}$,
with $\I$ the indicator function.

Section~SA2.10 of the online supplemental appendix provides technical details:
we show how an asymptotic
linear representation for $\widehat\psi(x)$ leads to a
modified Hoeffding-type decomposition of
$\widehat f_W^{\,1 \triangleright 0}(w)$,
which is then used to establish that $\widehat f_W^{\,1 \triangleright 0}$
is uniformly consistent for $f_W^{\,1 \triangleright 0}(w)$
and also admits a Gaussian strong approximation,
with the same rates of convergence
as for the standard density estimator.
Furthermore, define the covariance function of
$\widehat f_W^{\,1 \triangleright 0}(w)$ as
$\Sigma_n^{1 \triangleright 0}(w,w') = \Cov\big[
\widehat f_W^{\,1 \triangleright 0}(w),
\widehat f_W^{\,1 \triangleright 0}(w') \big]$,
which can be estimated as follows.
First let
$\widehat\kappa(X_i^0, X_i^1, x)
= \frac{\I\{X_i^0 = x\} - \widehat p_X^0(x)}{\widehat p_X^1(x)}
- \frac{\widehat p_X^0(x)}{\widehat p_X^1(x)} \frac{\I\{X_i^1 = x\} - \widehat
  p_X^1(x)}{\widehat p_X^1(x)}$
be a plug-in estimate of the influence function for $\widehat\psi(x)$
and define the leave-one-out
conditional expectation estimators\newline
$S_i^{1 \triangleright 0}(w)
= \frac{1}{n-1} \big( \sum_{j=1}^{i-1} k_h(W_{j i}^1,w) \widehat\psi(X_j^1)
+ \sum_{j=i+1}^n k_h(W_{ij}^1,w) \widehat\psi(X_j^1) \big)$
and\newline
$\widetilde S_i^{1 \triangleright 0}(w)
= \frac{1}{n-1} \sum_{j=1}^n \I\{j \neq i\}
\widehat\kappa(X_i^0, X_i^1, X_j^1) S_j^{1 \triangleright 0}(w)$.
Then define the covariance estimator
\begin{align*}
  \widehat\Sigma_n^{1 \triangleright 0}(w,w')
  &=
  \frac{4}{n^2}
  \sum_{i=1}^n
  \big(
  \widehat\psi(X_i^1)
  S_i^{1 \triangleright 0}(w)
  + \widetilde S_i^{1 \triangleright 0}(w)
  \big)
  \big(
  \widehat\psi(X_i^1)
  S_i^{1 \triangleright 0}(w')
  + \widetilde S_i^{1 \triangleright 0}(w')
  \big) \\
  &\quad-
  \frac{4}{n^3(n-1)}
  \sum_{i<j}
  k_h(W_{ij}^1, w)
  k_h(W_{ij}^1, w')
  \widehat\psi(X_i^1)^2
  \widehat\psi(X_j^1)^2
  - \frac{4}{n}
  \widehat f_W^{\,1 \triangleright 0}(w)
  \widehat f_W^{\,1 \triangleright 0}(w').
\end{align*}
We use a positive semi-definite approximation to
$\widehat\Sigma_n^{1 \triangleright 0}$, denoted by
$\widehat\Sigma_n^{+, 1 \triangleright 0}$,
as in Section~\ref{sec:covariance_estimation}.
To construct feasible uniform confidence bands,
define a process $\widehat Z_n^{T, 1 \triangleright 0}(w)$ which
is conditionally mean-zero and conditionally Gaussian
given the data $\bW_n^1$, $\bX_n^0$ and $\bX_n^1$
and whose conditional covariance structure is
$\E\big[\widehat Z_n^{T, 1 \triangleright 0}(w)
\widehat Z_n^{T, 1 \triangleright 0}(w')
\bigm| \bW_n^1, \bX_n^0, \bX_n^1 \big]
= \frac{\widehat \Sigma_n^{+, 1 \triangleright 0}(w,w')}
{\sqrt{\widehat \Sigma_n^{+, 1 \triangleright 0}(w,w)
    \widehat \Sigma_n^{+, 1 \triangleright 0}(w',w')}}$.
For $\alpha \in (0,1)$, define
$\widehat q_{1-\alpha}^{\,1 \triangleright 0}$
as the conditional quantile satisfying
$\P\big(\sup_{w \in \cW}\big| \widehat Z_n^{T, 1 \triangleright 0}(w) \big|
\leq \widehat q_{1-\alpha}^{\,1 \triangleright 0}
\bigm\vert \bW_n^1, \bX_n^0, \bX_n^1 \big)
= 1 - \alpha$.
Then, assuming that the covariance estimator is appropriately consistent,
\begin{align*}
  \P\left(
    f_W^{1 \triangleright 0}(w)
    \in
    \left[
      \widehat f_W^{\,1 \triangleright 0}(w)
      \pm
      \widehat q^{\,1 \triangleright 0}_{1-\alpha}
      \sqrt{\widehat\Sigma_n^{+, 1 \triangleright 0}(w,w)}
      \,\right]
    \,\textup{for all }
    w \in \cW
  \right)
  \to 1 - \alpha,
\end{align*}
giving feasible uniform inference methods,
which are robust to unknown degeneracies,
for counterfactual distribution analysis in dyadic data settings.

\subsection{Application to trade data}
\label{sec:trade_data}

We illustrate the performance of our
estimation and inference methods with a real-world data set.
We use international bilateral trade data
from the International Monetary Fund's
Direction of Trade Statistics (DOTS),
previously analyzed by \citet{head2014gravity}
and \citet{chiang2022inference}.
This data set contains information about
the yearly trade flows among
$n = 207$ economies ($N = 21\,321$ pairs),
and we focus on the years
$1995$, $2000$ and $2005$.

We define the \emph{trade volume} between countries
$i$ and $j$ as the logarithm of the sum of
the trade flow (in billions of US dollars)
from $i$ to $j$
and the trade flow from $j$ to $i$.
In each year several pairs of countries
did not trade directly, yielding trade flows of zero
and hence a trade volume of $-\infty$.
We therefore assume that the distribution of trade volumes
is a mixture of a point mass at $-\infty$ and a
Lebesgue density on $\R$.
The local nature of our estimator means that
observations taking the value of $-\infty$
can simply be removed from the data set.
Table~\ref{tab:trade_network_stats}
gives summary statistics for these trade networks,
and shows how the networks tend to become more connected over time,
with edge density, average degree and clustering coefficient
all increasing.
\begin{table}[H]
  \caption{Summary statistics for the DOTS trade networks.}
  \label{tab:trade_network_stats}
  \vspace{-1em}
  \begin{center}
    \begin{tabular}{|c|c|c|c|c|c|}
      \hline
      Year
      & Nodes
      & Edges
      & Edge density
      & Average degree
      & Clustering coefficient \\
      \hline
      1995 & 207 & 11\,603 & 0.5442 & 112.1 & 0.7250 \\
      2000 & 207 & 12\,528 & 0.5876 & 121.0 & 0.7674 \\
      2005 & 207 & 12\,807 & 0.6007 & 123.7 & 0.7745 \\
      \hline
    \end{tabular}
  \end{center}
  \vspace{-2em}
\end{table}

For counterfactual analysis we use the
gross domestic product (GDP) of each country as a
covariate, using $10\%$-percentiles to group the values into
$10$ different levels for ease of estimation.
This allows for a comparison of the observed distribution of trade at each year
with, for example, the counterfactual distribution of trade
had the GDP distribution remained as it was in $1995$.
As such we can measure how much of the change in trade
distribution is attributable
to a shift in the GDP distribution.

\begin{figure}[H]
  \centering
  \begin{subfigure}{0.32\textwidth}
    \centering
    \includegraphics[scale=0.52]{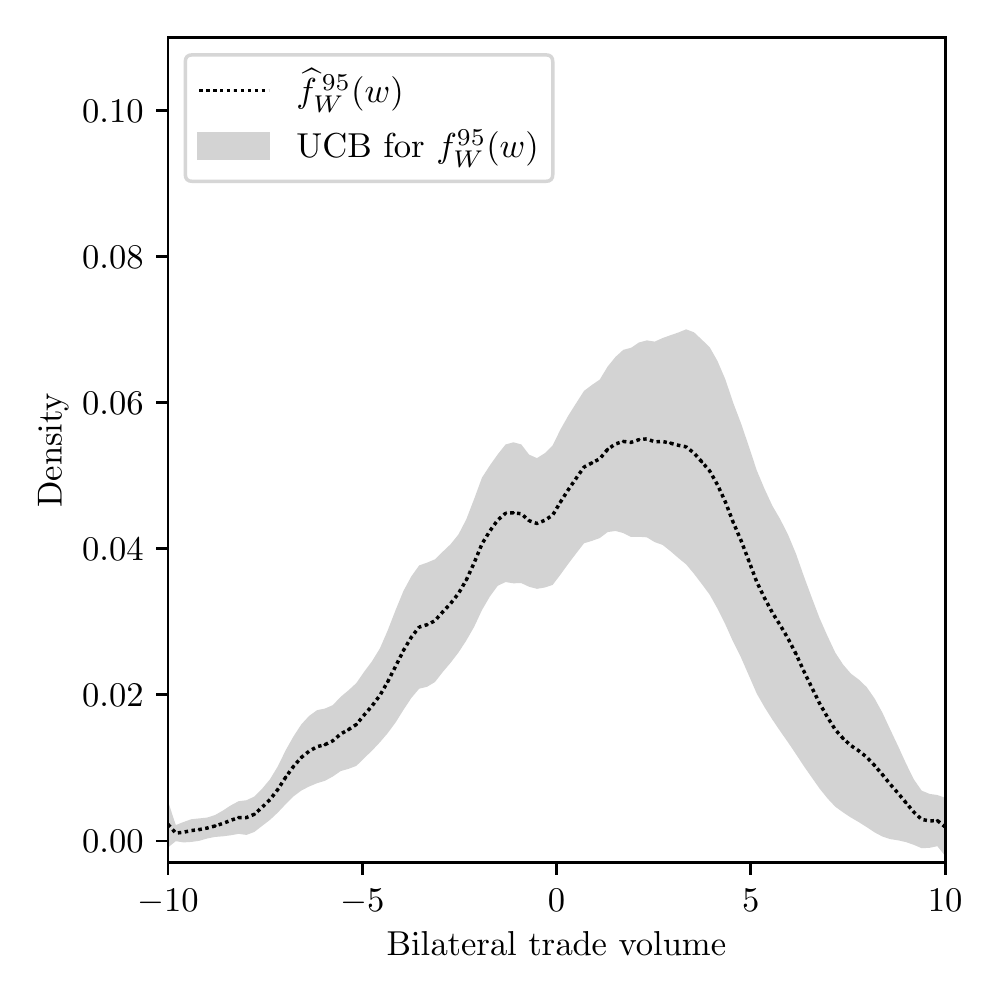}
    \caption{Year 1995, $\widehat h_{\ROT} = 1.27$}
  \end{subfigure}
  \begin{subfigure}{0.32\textwidth}
    \centering
    \includegraphics[scale=0.52]{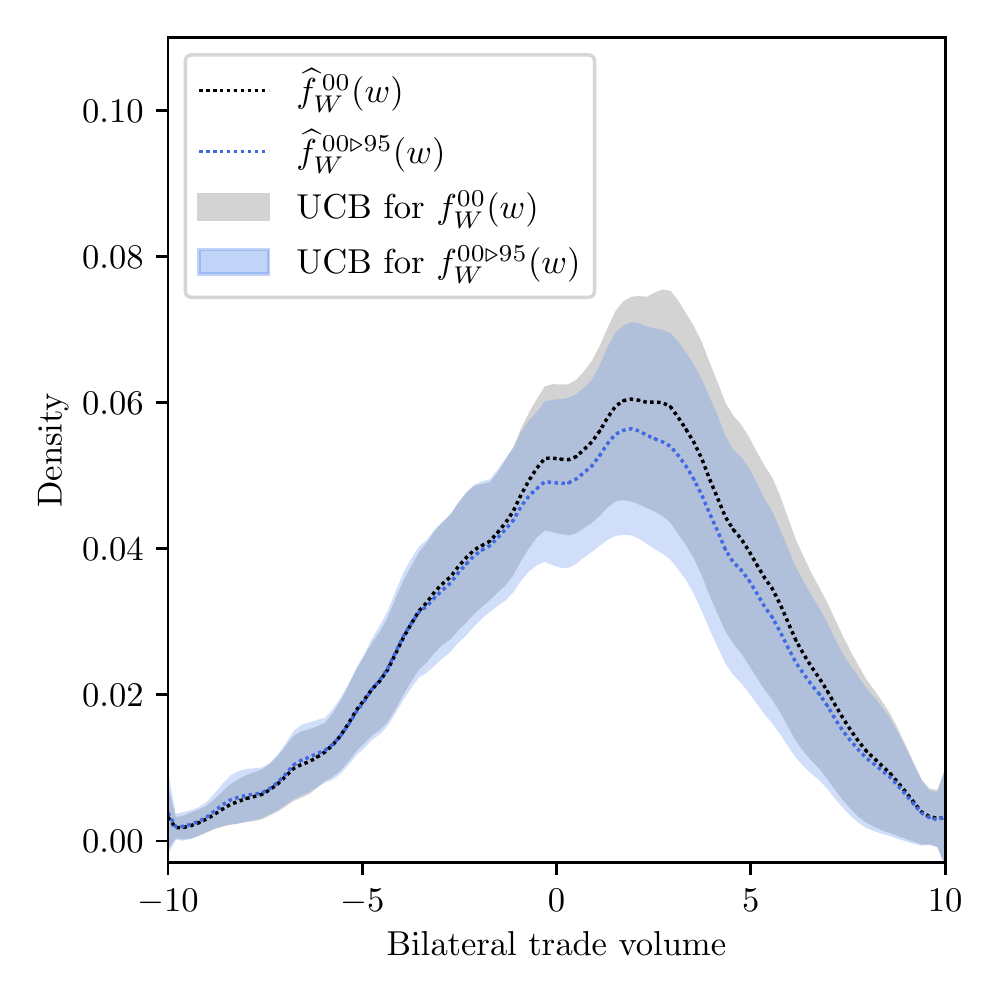}
    \caption{Year 2000, $\widehat h_{\ROT} = 1.31$}
  \end{subfigure}
  \begin{subfigure}{0.32\textwidth}
    \centering
    \includegraphics[scale=0.52]{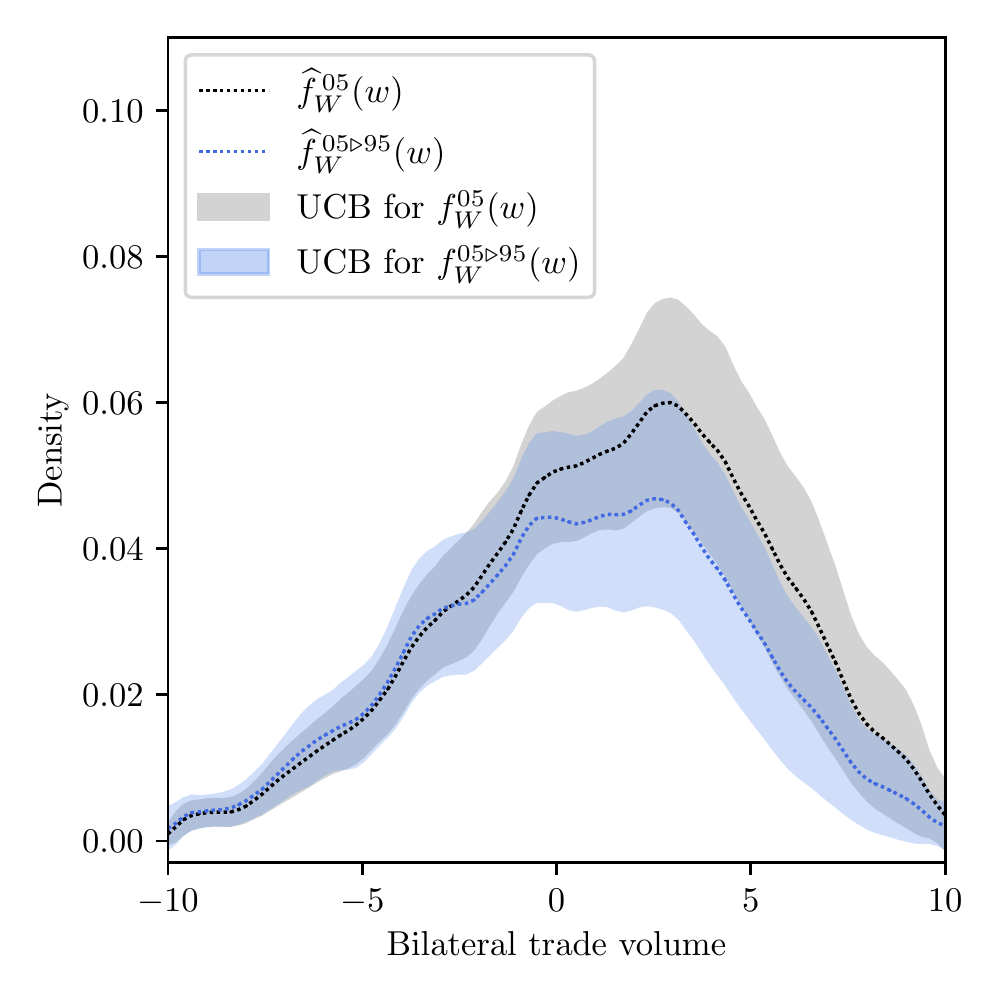}
    \caption{Year 2005, $\widehat h_{\ROT} = 1.37$}
  \end{subfigure}
  \caption{Real and counterfactual density estimates
    and confidence bands for the DOTS data.}
  \label{fig:trade}
  %
\end{figure}

To estimate the trade volume density function
we use Algorithm~\ref{alg:method}
with $d=100$ equally-spaced evaluation points in $[-10,10]$,
using the rule-of-thumb bandwidth selector $\widehat h_{\ROT}$
from Section~\ref{sec:bandwidth_selection} with
$p=2$ and $C(K) = 2.435$.
For inference we use an
Epanechnikov kernel of order $p=4$ and
resample the Gaussian process
$B = 10\,000$ times.
We also estimate the counterfactual trade distributions
in 2000 and 2005 respectively,
replacing the GDP distribution with that from 1995.
For each year, Figure~\ref{fig:trade}
plots the real and counterfactual density estimates
along with their respective uniform confidence bands (UCB)
at the nominal coverage rate of $95\%$.
Our empirical results show that
the counterfactual distribution drifts further from the truth
in 2005 compared with 2000, indicating a more significant shift
in the GDP distribution.
In the online supplemental appendix we repeat the procedure
using a parametric (log-normal maximum likelihood) preliminary estimate
of the GDP distribution, and observe that the results
are qualitatively similar.

\section{Other applications and future work}
\label{sec:future}

To emphasize the broad applicability of our methods
to network science problems,
we present three application scenarios.
The first concerns comparison of networks
\citep{kolaczyk2009statistical},
while the second and third involve
nonparametric and semiparametric dyadic regression
respectively.

Firstly, consider the setting where there are
two independent networks with
continuous dyadic covariates $\bW_n^0$ and $\bW_m^1$ respectively.
Practitioners may wish to test if these two dyadic distributions
are the same, that is, whether their density functions
$f_W^0$ and $f_W^1$ are equal on their
common support $\cW \subseteq \R$.
We present a family of hypothesis tests for this scenario
based on dyadic kernel density estimation.
Let $\widehat f_W^{\,0}(w)$ and $\widehat f_W^{\,1}(w)$
be the associated
(bias-corrected) dyadic kernel density estimators.
Consider the test statistics
$\tau_p$ for $1 \leq p \leq \infty$ where
\begin{align}
  \label{eq:hypothesis_test}
  \tau_p
  &=
  \left(
    \int_{-\infty}^{\infty}
    \left|
    \widehat f_W^{\,1}(w)
    - \widehat f_W^{\,0}(w)
    \right|^p
    \diff w
  \right)^{1/p}
  \text{ for } p < \infty, \text{ and }\,
  \tau_\infty
  = \sup_{w \in \cW}
  \left|
  \widehat f_W^{\,1}(w)
  - \widehat f_W^{\,0}(w)
  \right|.
\end{align}
Clearly we should reject the null hypothesis that $f_W^0 = f_W^1$ whenever
the test statistic $\tau_p$ is sufficiently large.
To estimate the critical value, let
$\widehat\Sigma_n^{+,0}(w, w')$ and $\widehat\Sigma_m^{+,1}(w, w')$
be the positive semi-definite estimators defined in
Section~\ref{sec:covariance_estimation} and
let $\widehat Z^0_n(w)$ and $\widehat Z^1_m(w)$ be zero-mean
Gaussian processes with covariance
structures $\widehat\Sigma_n^{+,0}(w, w')$
and $\widehat\Sigma_m^{+,1}(w, w')$
respectively, which are independent conditional on the data.
Define the approximate null test statistic
$\widehat \tau_p$ by replacing
$\widehat f_W^{\,0}(w)$ and $\widehat f_W^{\,1}(w)$
with $\widehat Z^0_n(w)$ and $\widehat Z^1_m(w)$ respectively
in \eqref{eq:hypothesis_test}.
For a significance level $\alpha \in (0,1)$,
the critical value is $\widehat C_\alpha$ where
$\P \big(
\widehat \tau_p \geq \widehat C_\alpha \bigm\vert \bW_n^0, \bW_n^1
\big) = \alpha$.
This is estimated by Monte Carlo simulation, resampling
from the conditional law of
$\widehat Z^0_n(w)$ and $\widehat Z^1_m(w)$
and replacing integrals and suprema by sums and maxima over a
finite partition of $\cW$.

While our focus has been on density estimation with dyadic data,
our uniform dyadic estimation and inference results are readily
applicable to the settings of nonparametric
and semiparametric dyadic regression.
For our second example, suppose that $Y_{ij} = Y(X_i, X_j, A_i, A_j, V_{ij})$,
where only $\bX_n$ and $\bY_n$ are observed and
$\bV_n$ is independent of $(\bX_n, \bA_n)$,
with $\bX_n = (X_i : 1 \leq i \leq n)$,
$\bA_n = (A_i : 1 \leq i \leq n)$,
$\bY_n = (Y_{ij}:1\leq i<j\leq n)$
and $\bV_n = (V_{ij} : 1 \leq i < j \leq n)$.
A parameter of interest is the regression function
$\mu(x_1, x_2) = \E[Y_{ij} \mid X_i=x_1, X_j=x_2]$,
which can be used to analyze average or partial effects
of changing the node attributes $X_i$ and $X_j$
on the edge variable $Y_{ij}$.
This conditional expectation could be estimated
using local polynomial methods:
suppose that $X_i$ takes values in $\R^m$ and
let $r(x_1, x_2)$ be a monomial basis up to degree
$\gamma \geq 0$ on $\R^m \times \R^m$. Then,
for some bandwidth $h > 0$ and
a kernel function $k_h$ on $\R^m \times \R^m$,
the local polynomial regression estimator of
$\mu(x_1, x_2)$ is
$\widehat\mu(x_1, x_2)
= e_1^\T \widehat\beta(x_1, x_2)$ where
$e_1$ is the first standard
unit vector in $\R^q$ for $q=\binom{2m+\gamma}{\gamma}$ and
\begin{align}
  \nonumber
  \widehat{\beta}(x_1, x_2)
  &=
  \argmin_{\beta \in \R^q}
  \sum_{i=1}^{n-1}
  \sum_{j=i+1}^n
  \left(
    Y_{ij} - r(X_i-x_1, X_j-x_2)^\T \beta
  \right)^2
  k_h(X_i-x_1, X_j-x_2) \\
  \label{eq:locpol}
  &=
  \left(
    \sum_{i=1}^{n-1}
    \sum_{j=i+1}^n
    k_{ij} r_{ij} r_{ij}^\T
  \right)^{-1}
  \left(
    \sum_{i=1}^{n-1}
    \sum_{j=i+1}^n
    k_{ij} r_{ij} Y_{ij}
  \right),
\end{align}
with $k_{ij} = k_h(X_i-x_1, X_j-x_2)$
and
$r_{ij} = r(X_i-x_1, X_j-x_2)$.
\citet{graham2021dyadicregression}
established pointwise distribution theory
for the special case of the dyadic Nadaraya--Watson
kernel regression estimator ($\gamma=0$), but no uniform analogues
have yet been given.
It can be shown that the
``denominator'' matrix
in \eqref{eq:locpol} converges uniformly
to its expectation,
while the U-process-like ``numerator'' matrix
can be handled the same way as we analyzed
$\widehat f_W(w)$ in this paper,
through a Hoeffding-type decomposition and
strong approximation methods,
along with standard bias calculations.
Such distributional approximation results can be used to construct valid
uniform confidence bands for the regression
function $\mu(x_1, x_2)$,
as well as to conduct hypothesis testing for
parametric specifications or shape constraints.

As a third example, we consider applying our results to
semiparametric semi-linear regression problems.
The dyadic semi-linear regression model is
$\E[Y_{ij} \mid W_{ij}, X_i, X_j]
= \theta^\T W_{ij} + g(X_i, X_j)$
where $\theta$ is the finite-dimensional parameter
of interest
and $g(X_i, X_j)$ is an unknown
function of the covariates
$(X_i, X_j)$.
Local polynomial (or other) methods can be used to
estimate $\theta$ and $g$,
where the estimator of the nonparametric component
$g$ takes a similar form to \eqref{eq:locpol},
that is, a ratio of two kernel-based
estimators as in \eqref{eq:estimator}.
Consequently, our strong approximation techniques
presented in this paper can be appropriately modified
to develop valid uniform inference procedures for
$g$ and $\E[Y_{ij} \mid W_{ij}=w, X_i=x_1, X_j=x_2]$,
as well as functionals thereof.

\section{Conclusion}
\label{sec:conclusion}

We studied the uniform estimation and inference properties of the
dyadic kernel density estimator
$\widehat{f}_W$ given in \eqref{eq:estimator},
which forms a class
of U-process-like estimators
indexed by the $n$-varying kernel function $k_h$ on $\cW$.
We established uniform minimax-optimal point estimation results
and uniform distributional approximations for this estimator
based on novel strong approximation strategies.
We then applied these results to derive valid
and feasible uniform confidence bands for
the dyadic density estimand $f_W$,
and also developed a substantive application of our theory to
counterfactual dyadic density analysis.
We gave some other statistical
applications of our methodology as well as
potential avenues for future research.
From a technical perspective,
the online supplemental appendix contains several generic results
concerning strong approximation methods and maximal inequalities
for empirical processes that may be of independent interest.

\section*{Acknowledgments}

We thank the Co-Editor, Associate Editor, and three reviewers,
along with
Harold Chiang,
Laurent Davezies,
Xavier D'Haultf{\oe}uille,
Jianqing Fan,
Yannick Guyonvarch,
Kengo Kato,
Jason Klusowski,
Ricardo Masini,
Yuya Sasaki
and Boris Shigida
for useful comments.
Cattaneo was supported through
National Science Foundation grants SES-1947805 and DMS-2210561.
Feng was supported by the
National Natural Science Foundation of China (NSFC)
under grants 72203122 and 72133002.

\section*{Supplemental material}

A supplemental appendix containing technical and methodological details
as well as proofs and additional empirical results is available at
\href{https://arxiv.org/abs/2201.05967}%
{\texttt{https://arxiv.org/abs/2201.05967}}.
Replication files for the empirical studies are provided at
\href{https://github.com/WGUNDERWOOD/DyadicKDE.jl}%
{\texttt{https://github.com/wgunderwood/DyadicKDE.jl}}.


\appendix

\bibliographystyle{hapalike}
\bibliography{CFU_2023_DyadicKDE--bib}

\clearpage

\end{document}


\def\arraystretch{1.2}

\maketitle

\footnotetext[1]{
  Department of Operations Research
  and Financial Engineering,
  Princeton University
}
\footnotetext[2]{
  School of Economics and Management,
  Tsinghua University
}
\begin{NoHyper}
  \nmfootnote{\textsuperscript{*}Corresponding author:
    \href{mailto:wgu2@princeton.edu}{\texttt{wgu2@princeton.edu}}}
\end{NoHyper}
\setcounter{footnote}{3}

\tableofcontents

\setcounter{page}{0}\thispagestyle{empty}

\clearpage
\pagebreak

\section{Introduction}

We describe the setup of
dyadic density estimation,
define the dyadic kernel density estimator,
give some notation
and state our assumptions.

\subsection{Setup and estimator}
\label{sec:setup_estimator}

Fix $n \geq 2$ and suppose there is
a probability space carrying the latent random variables
$\bA_n=(A_i: 1\leq i\leq n)$ and $\bV_n=(V_{ij}: 1\leq i<j\leq n)$.
Suppose that the $A_i$ are i.i.d., the $V_{ij}$ are i.i.d., and
that $\bA_n$ is independent of $\bV_n$.
Define the observable dyadic random variables
$W_{ij} = W(A_i, A_j, V_{ij})$
where $W$ is some unknown real-valued function.
There are $\frac{1}{2}n(n-1)$ such variables;
one for each unordered pair of distinct indices $i<j$.
Note that if $i<j$ and $i'<j'$ are all distinct,
then $W_{ij}$ is independent of $W_{i' j'}$.
However, $W_{ij}$ is
not in general independent of $W_{ij'}$,
as they may both depend on the latent variable $A_i$.
This data generating process is justified by the
Aldous--Hoover representation theorem for exchangeable arrays
\citep{aldous1981representations, hoover1979relations}.

Denote the kernel weight of the data point $W_{ij}$
at the evaluation point $w$
with bandwidth $h$
(see Section~\ref{sec:assumptions} for details)
by
$k_h(W_{ij}, w)$.
Then the dyadic
kernel density estimator
is defined as
\begin{align*}
  \widehat f_W(w) = \frac{2}{n(n-1)}
  \sum_{i=1}^{n-1} \sum_{j=i+1}^{n} k_h(W_{ij}, w).
\end{align*}

\subsection{Notation}
\label{sec:notation}

\subsubsection{Norms}

For real vectors,
$\|\cdot\|_p$ is the standard $L^p$
norm defined for $p \in [1, \infty]$.
For real square matrices,
$\|\cdot\|_p$ is the operator
norm induced by the corresponding vector norm.
In particular,
$\|\cdot\|_1$
is the maximum absolute column sum,
$\|\cdot\|_\infty$
is the maximum absolute row sum,
and
$\|\cdot\|_2$
is the maximum singular value.
For real symmetric matrices,
$\|\cdot\|_2$
coincides with the maximum absolute eigenvalue.
We use $\|\cdot\|_\maxi$
to denote the largest absolute entry of a real matrix.
For real-valued functions,
$\|\cdot\|_\infty$
denotes the (essential) supremum norm.
The total variation norm of a
real-valued function of a single real variable is
$\|g\|_\TV = \sup_{n \geq 1} \sup_{x_1 \leq \cdots \leq x_n}
\sum_{i=1}^{n-1} |g(x_{i+1}) - g(x_i)|$.

\subsubsection{Inequalities}

For deterministic non-negative sequences
$a_n$ and $b_n$,
write
$a_n \lesssim b_n$
or $a_n = O(b_n)$
to indicate that there exists a positive constant
$C$ which does not depend on $n$
(although might depend on other quantities,
depending on context)
satisfying
$a_n \leq C b_n$
for all sufficiently large $n$.
Write $a_n \ll b_n$
or $a_n = o(b_n)$
to indicate that
$a_n / b_n \to 0$.
If $a_n \lesssim b_n \lesssim a_n$,
write $a_n \asymp b_n$.
For random non-negative sequences
$A_n$ and $B_n$,
write
$A_n \lesssim_\P B_n$
or $A_n = O_\P(B_n)$
to indicate that for any $\varepsilon > 0$
there exists a deterministic positive constant
$C_\varepsilon$ satisfying
$\P(A_n \leq C_\varepsilon B_n) \geq 1 - \varepsilon$
for all sufficiently large $n$.
Write $A_n = o_\P(B_n)$ if
$A_n / B_n \to 0$ in probability.

\subsubsection{Sets}

For $x \in \R$ and $a \geq 0$,
we use $[x \pm a]$
to denote the compact interval
$[x-a, x+a]$.
For a bounded set $\cX \subseteq \R$ and $a \geq 0$
we use $[\cX \pm a]$ to denote the compact interval
$[\inf \cX - a, \ \sup \cX + a]$.
For measurable subsets of $\R^d$
we use $\Leb$ to denote the Lebesgue measure,
and for finite sets we use $|\cdot|$
for the cardinality.

\subsubsection{Sums}

We use $\sum_i$
to indicate $\sum_{i=1}^n$
when clear from context.
Similarly we use $\sum_{i<j}$
for $\sum_{i=1}^{n-1} \sum_{j=i+1}^n$
and $\sum_{i<j<r}$
for $\sum_{i=1}^{n-2} \sum_{j=i+1}^{n-1} \sum_{r=j+1}^n$.

\subsubsection{Function classes}

Let $\cX \subseteq \R$ be an interval and $\beta > 0$.
Define $\flbeta$ as the largest integer
which is \emph{strictly} smaller than $\beta$.
Let $\cC^{\flbeta}(\cX)$ be the class
of functions from $\R$ to $\R$ which are
$\flbeta$ times continuously differentiable on $\cX$.
Note that $\cC^0(\cX)$ is the class of functions
which are continuous on $\cX$.
For $C > 0$, define the H\"older class with smoothness
$\beta > 0$ by
\begin{align*}
  \cH^\beta_C(\cX)
  &=
  \Big\{
  g \in \cC^{\flbeta}(\cX): \
  \max_{1 \leq r \leq \flbeta}
  \big| g^{(r)}(x) \big| \leq C
  \text{ and }
  \big| g^{(\flbeta)}(x) - g^{(\flbeta)}(x') \big|
  \leq C |x-x'|^{\beta - \flbeta}, \
  \text{ for all }
  x, x' \in \cX
  \Big\}.
\end{align*}
Note that $\cH^1_C(\cX)$
is the class of functions which are $C$-Lipschitz on $\cX$,
and observe that the functions in
$\cH^\beta_C(\cX)$ are not uniformly bounded on $\cX$.

\subsection{Assumptions}
\label{sec:assumptions}

\begin{assumption}[Data generation]
  \label{ass:data}

  Fix $n \geq 2$ and
  let $\bA_n = (A_i: 1 \leq i \leq n)$
  be i.i.d.\
  real-valued random variables supported on $\cA \subseteq \R$.
  Let $\bV_n = (V_{ij}: 1 \leq i < j \leq n)$
  be i.i.d.\
  real-valued random variables with a Lebesgue density $f_V$ on $\R$.
  Suppose that $\bA_n$ is independent of $\bV_n$.
  Let $W_{ij} = W(A_i, A_j, V_{ij})$
  and $\bW_n = (W_{ij}: 1 \leq i < j \leq n)$,
  where $W$ is some unknown real-valued function
  which is symmetric in its first two arguments.
  Let $\cW \subseteq \R$ be a compact interval
  with positive Lebesgue measure $\Leb(\cW)$.
  Assume that the conditional distribution
  of $W_{ij}$ given $A_i$ and $A_j$
  admits a Lebesgue density
  denoted $f_{W \mid AA}(w \mid A_i, A_j)$,
  and define
  $f_{W \mid A}(w \mid a) = \E\left[f_{W \mid AA}(w \mid A_i,a)\right]$
  and
  $f_W(w) = \E\left[f_{W \mid AA}(w \mid A_i,A_j)\right]$.
  Take $\CH > 0$ and $\beta \geq 1$,
  and suppose that $f_W \in \cH^\beta_{\CH}(\cW)$ and that
  $f_{W \mid AA}(\cdot \mid a, a') \in \cH^1_{\CH}(\cW)$
  for all $a,a' \in \cA$.
  Assume that
  $\sup_{w \in \cW} \|f_{W \mid A}(w \mid \cdot\,)\|_\TV <\infty$.
\end{assumption}

\begin{remark}
  If $W(a_1, a_2, v)$
  is strictly monotonic
  and continuously differentiable
  in its third argument,
  we can give the conditional
  density of $W_{ij}$ explicitly
  using the usual change-of-variables formula:
  with $w=W(a_1,a_2,v)$, we have
  $f_{W \mid AA}(w \mid a_1,a_2)
  = f_V(v)
  \Big|
  \frac{\partial W(a_1,a_2,v)}{\partial v}
  \Big|^{-1}$.
\end{remark}

\begin{remark}
  By Lemma~\ref{lem:lipschitz_kernels_bounded},
  Assumption~\ref{ass:data} implies that the densities
  $f_W$, $f_{W \mid A}$ and $f_{W \mid AA}$ are all uniformly bounded by
  $\Cd \vcentcolon= 2 \sqrt{\CH} + 1/\Leb(\cW)$.
\end{remark}

\begin{assumption}[Kernels and bandwidth]
  \label{ass:kernel_bandwidth}

  Let $h = h(n) > 0$
  be a sequence of bandwidths
  satisfying $h \log n \to 0$
  and $\frac{\log n}{n^2h} \to 0$.
  For each $w \in \cW$ let $k_h(\cdot, w)$
  be a real-valued function supported on
  $[w \pm h] \cap \cW$.
  Let $p \geq 1$ be an integer and
  suppose that $k_h$ belongs to a family of
  boundary bias-corrected kernels of order $p$,
  which is to say that
  \begin{align*}
    \int_{\cW}
    (s-w)^r k_h(s,w) \diff{s}
    \quad
    \begin{cases}
      \begin{alignedat}{2}
        &= 1 &\qquad &\text{for all } w \in \cW \text{ if }\, r = 0, \\
        &= 0 & &\text{for all } w \in \cW \text{ if }\, 1 \leq r \leq p-1, \\
        &\neq 0 & &\text{for some } w \in \cW \text{ if }\, r = p.
      \end{alignedat}
    \end{cases}
  \end{align*}
  Suppose also that for some $\CL > 0$,
  the kernels satisfy
  $k_h(s, \cdot) \in \cH^1_{\CL/h^2}(\cW)$
  for all $s \in \cW$.

\end{assumption}

\begin{remark}
  The kernel functions required by
  Assumption~\ref{ass:kernel_bandwidth}
  can be constructed using
  polynomials on $[w \pm h] \cap \cW$,
  solving a family of linear systems to find the coefficients.
\end{remark}

\begin{remark}
  By Lemma~\ref{lem:lipschitz_kernels_bounded},
  Assumption~\ref{ass:kernel_bandwidth} implies that
  if $h \leq 1$ then
  $k_h$ is uniformly bounded by
  $\Ck / h$ where $\Ck \vcentcolon = 2 \CL + 1 + 1/\Leb(\cW)$.
\end{remark}
\section{Main results}

\subsection{Bias}

Lemma~\ref{lem:bias}
is a standard result in kernel density estimation
with boundary bias correction,
and does not rely on the dyadic structure of the data.

\begin{lemma}[Bias of $\widehat f_W$]
  \label{lem:bias}

  Suppose that Assumptions \ref{ass:data} and
  \ref{ass:kernel_bandwidth} hold.
  For $w \in \cW$
  define the leading bias term as
  \begin{align*}
    B_p(w)
    &=
    \frac{f_W^{(p)}(w)}{p!}
    \int_{\cW}
    k_h(s,w)
    \left(
      \frac{s-w}{h}
    \right)^p
    \diff{s}.
  \end{align*}
  for $1 \leq p \leq \flbeta$.
  Then we have the following bias bounds.
  \begin{enumerate}[label=(\roman*)]

    \item If $p \leq \flbeta - 1$,
      \begin{align*}
        \sup_{w \in \cW}
        \big|
        \E\big[\widehat f_W(w)\big]
        - f_W(w)
        - h^p B_p(w)
        \big|
        \leq
        \frac{2\Ck \CH}{(p+1)!}
        h^{p+1}.
      \end{align*}

    \item If $p = \flbeta$,
      \begin{align*}
        \sup_{w \in \cW}
        \big|
        \E\big[\widehat f_W(w)\big]
        - f_W(w)
        - h^p B_p(w)
        \big|
        \leq
        \frac{2 \Ck \CH}{\flbeta !}
        h^\beta.
      \end{align*}

    \item If $p \geq \flbeta+1$,
      \begin{align*}
        \sup_{w \in \cW}
        \big|
        \E\big[\widehat f_W(w)\big]
        - f_W(w)
        \big|
        \leq
        \frac{2 \Ck \CH}{\flbeta !}
        h^\beta.
      \end{align*}

  \end{enumerate}
  Noting that
  $\sup_{\cW} |B_p(w)| \leq 2 \Ck \CH / p!$,
  we deduce that for $h \leq 1$,
  \begin{align*}
    \sup_{w \in \cW}
    \big|
    \E\big[\widehat f_W(w)\big]
    - f_W(w)
    \big|
    \leq
    \frac{4 \Ck \CH}{(p \wedge \flbeta)!}
    h^{p \wedge \beta}
    \lesssim
    h^{p \wedge \beta}.
  \end{align*}

\end{lemma}

\subsection{Uniform consistency}

In this section we demonstrate uniform consistency of the
dyadic kernel density estimator.
Lemma~\ref{lem:hoeffding}
provides a U-statistic decomposition of the estimator
and Lemma~\ref{lem:uniform_concentration}
employs this decomposition to establish
uniform concentration.
Theorem~\ref{thm:uniform_consistency}
then combines this with the bias result from Lemma~\ref{lem:bias}
to show uniform consistency.
Lemma~\ref{lem:trichotomy} provides a useful trichotomy for
interpreting our results in various classes of data distributions.

\begin{lemma}[Hoeffding-type decomposition for $\widehat f_W$]
  \label{lem:hoeffding}

  Suppose that Assumptions~\ref{ass:data}
  and~\ref{ass:kernel_bandwidth}
  hold.
  Define the linear term
  (H\'{a}jek projection),
  quadratic term
  and error term
  of $\widehat f_W(w)$ as
  \begin{align*}
    L_n(w)
    &=
    \frac{2}{n} \sum_{i=1}^n
    l_i(w),
    &Q_n(w)
    &=
    \frac{2}{n(n-1)}
    \sum_{i=1}^{n-1}
    \sum_{j=i+1}^{n}
    q_{ij}(w),
    &E_n(w)
    &=
    \frac{2}{n(n-1)}
    \sum_{i=1}^{n-1}
    \sum_{j=i+1}^{n}
    e_{ij}(w)
  \end{align*}
  respectively, where
  \begin{align*}
    l_i(w)
    &=
    \E\left[k_h(W_{ij},w) \mid A_i\right] - \E\left[k_h(W_{ij},w)\right], \\
    q_{ij}(w)
    &=
    \E\left[k_h(W_{ij},w) \mid A_i, A_j\right]
    - \E\left[k_h(W_{ij},w) \mid A_i\right]
    - \E\left[k_h(W_{ij},w) \mid A_j\right]
    + \E\left[k_h(W_{ij},w)\right], \\
    e_{ij}(w)
    &=
    k_h(W_{ij},w)
    - \E\left[k_h(W_{ij},w) \mid A_i, A_j\right].
  \end{align*}
  Then the following Hoeffding-type decomposition holds:
  \begin{align*}
    \widehat f_W(w)
    = \E\big[\widehat f_W(w)\big]
    + L_n(w)
    + Q_n(w)
    + E_n(w).
  \end{align*}
  Further,
  the stochastic processes
  $L_n$, $Q_n$ and $E_n$
  are all mean-zero, since
  \begin{equation*}
    \E\big[L_n(w)\big]
    = \E\big[Q_n(w)\big]
    = \E\big[E_n(w)\big]
    = 0
  \end{equation*}
  for all $w \in \cW$.
  Also they are mutually orthogonal
  in $L^2(\P)$ as
  \begin{equation*}
    \E\big[
      L_n(w) Q_n(w')
      \big]
    = \E\big[
      L_n(w) E_n(w')
      \big]
    = \E\big[
      Q_n(w) E_n(w')
      \big]
    = 0
  \end{equation*}
  for all $w, w' \in \cW$.

\end{lemma}

\begin{lemma}[Trichotomy of degeneracy]
  \label{lem:trichotomy}

  Suppose that Assumptions~\ref{ass:data}
  and~\ref{ass:kernel_bandwidth} hold,
  and define the non-negative
  upper and lower degeneracy constants
  \begin{align*}
    \Du^2
    &=
    \sup_{w \in \cW}
    \Var\left[
      f_{W \mid A}(w \mid A_i)
    \right],
    &\Dl^2
    &=
    \inf_{w \in \cW}
    \Var\left[
      f_{W \mid A}(w \mid A_i)
    \right]
  \end{align*}
  respectively.
  Then precisely one of the following three statements must hold.
  \begin{enumerate}[label=(\roman*)]

    \item Total degeneracy:
      $\Du = \Dl = 0$.
      Then
      $L_n(w) = 0$ for all $w \in \cW$
      almost surely.

    \item No degeneracy:
      $\Dl > 0$. Then
      $\inf_{w \in \cW} \Var[L_n(w)] \geq \frac{2 \Dl}{n}$
      for all large enough n.

    \item Partial degeneracy:
      $\Du > \Dl = 0$.
      There exists $w \in \cW$ with
      $\Var\left[f_{W \mid A}(w \mid A_i)\right] = 0$;
      such a point is labelled \emph{degenerate} and satisfies
      $\Var[L_n(w)] \leq 64 \Ck \CH \Cd \frac{h}{n}$.
      There also exists a point $w' \in \cW$ with
      $\Var\left[f_{W \mid A}(w' \mid A_i)\right] > 0$;
      such a point is labelled \emph{non-degenerate} and satisfies
      $\Var[L_n(w')] \geq
      \frac{2}{n} \Var\left[f_{W \mid A}(w' \mid A_i)\right]$
      for all large enough $n$.

  \end{enumerate}

\end{lemma}

\begin{remark}

  The trichotomy of total/partial/no degeneracy
  given in Lemma~\ref{lem:trichotomy} is useful for understanding
  the asymptotic behavior of the dyadic kernel density estimator.
  Note that our need for uniformity in $w$
  complicates the simpler degeneracy/no degeneracy
  dichotomy observed for pointwise results by
  \citet{graham2019kernel}.

\end{remark}

\begin{lemma}[Uniform concentration of $\widehat f_W$]
  \label{lem:uniform_concentration}

  Suppose Assumptions
  \ref{ass:data}
  and
  \ref{ass:kernel_bandwidth}
  hold.
  Then with $L_n$, $Q_n$ and $E_n$ defined as in
  Lemma~\ref{lem:hoeffding},
  we have
  \begin{align*}
    \E\left[ \sup_{w \in \cW} |L_n(w)| \right]
    &\lesssim
    \frac{\Du}{\sqrt n},
    &\E\left[ \sup_{w \in \cW} |Q_n(w)| \right]
    &\lesssim
    \frac{1}{n},
    &\E\left[ \sup_{w \in \cW} |E_n(w)| \right]
    &\lesssim
    \sqrt{\frac{\log n}{n^2h}},
  \end{align*}
  where $\lesssim$ is up to constants which
  depend on the underlying data distribution
  and the choice of kernel.
  Note that the $Q_n$
  term is dominated by the $L_n$
  term uniformly in the bandwidth $h$.
  Therefore by Lemma~\ref{lem:hoeffding}
  \begin{align*}
    \E\left[
      \sup_{w \in \cW}
      \big|\widehat f_W(w) - \E[\widehat f_W(w)]\big|
    \right]
    &\lesssim
    \frac{\Du}{\sqrt n}
    + \sqrt{\frac{\log n}{n^2h}}.
  \end{align*}

\end{lemma}

\begin{theorem}[Uniform consistency of $\widehat f_W$]
  \label{thm:uniform_consistency}

  Suppose Assumptions \ref{ass:data}
  and \ref{ass:kernel_bandwidth} hold. Then
  \begin{align*}
    \E\left[
      \sup_{w \in \cW}
      \big|\widehat f_W(w) - f_W(w)\big|
    \right]
    &\lesssim
    h^{p \wedge \beta}
    + \frac{\Du}{\sqrt n}
    + \sqrt{\frac{\log n}{n^2h}},
  \end{align*}
  where $\lesssim$ is up to constants which
  depend on the underlying data distribution
  and the choice of kernel.

\end{theorem}

\begin{remark}

  In light of the degeneracy trichotomy
  in Lemma~\ref{lem:trichotomy},
  we interpret Theorem~\ref{thm:uniform_consistency}.
  \begin{enumerate}[label=(\roman*)]

    \item Partial or no degeneracy:
      when $\Du > 0$,
      any bandwidth sequence satisfying
      $\frac{\log n}{n} \lesssim h \lesssim n^{-\frac{1}{2(p \wedge \beta)}}$
      gives the bandwidth-independent ``parametric'' rate noted by
      \citet{graham2019kernel}:
      \begin{align*}
        \E\left[
          \sup_{w \in \cW}
          \big|\widehat f_W(w) - \E[\widehat f_W(w)]\big|
        \right]
        \lesssim
        \frac{1}{\sqrt n},
      \end{align*}

    \item Total degeneracy:
      when $\Du = 0$,
      minimizing the upper bound by setting
      $h \asymp
      \left( \frac{\log n}{n^2} \right)^{\frac{1}{2(p \wedge \beta)+1}}$
      yields
      \begin{align*}
        \E\left[
          \sup_{w \in \cW}
          \big|\widehat f_W(w) - f_W(w)\big|
        \right]
        &\lesssim
        h^{p \wedge \beta}
        + \sqrt{\frac{\log n}{n^2h}}
        \lesssim
        \left(
          \frac{\log n}{n^2}
        \right)^{\frac{p \wedge \beta}{2(p \wedge \beta)+1}},
      \end{align*}
  \end{enumerate}

\end{remark}

\subsection{Minimax optimality}
\label{sec:minimax}

In this section we demonstrate
minimax optimality of our estimator under uniform convergence
by providing upper and lower bounds in expectation uniformly over
some classes of dyadic distributions.
\begin{theorem}[Minimax optimality of $\widehat f_W$]
  \label{thm:minimax}

  Fix $\beta \geq 1$ and $\CH > 0$,
  and take $\cW$ a compact interval with
  positive Lebesgue measure.
  Define $\cP = \cP(\cW, \beta, \CH)$
  as the class of dyadic distributions
  satisfying Assumption~\ref{ass:data}.
  Define $\cPd$ as the subclass of $\cP$
  containing only those distributions
  which are totally degenerate on $\cW$ in the sense that
  $\sup_{w \in \cW} \Var\left[f_{W \mid A}(w \mid A_i)\right] = 0$.
  Then
  \begin{align*}
    \inf_{\widetilde f_W}
    \sup_{\P \in \cP}
    \E_\P\left[
      \sup_{w \in \cW}
      \big|
      \widetilde f_W(w) - f_W(w)
      \big|
    \right]
    &\asymp
    \frac{1}{\sqrt n}, \\
    \inf_{\widetilde f_W}
    \sup_{\P \in \cPd}
    \E_\P\left[
      \sup_{w \in \cW}
      \big|
      \widetilde f_W(w) - f_W(w)
      \big|
    \right]
    &\asymp
    \left(
      \frac{\log n}{n^2}
    \right)^{\frac{\beta}{2\beta+1}},
  \end{align*}
  where $\widetilde f_W$ is any estimator depending only on
  the data $\bW_n = (W_{ij}: 1 \leq i < j \leq n)$
  distributed according to the dyadic law $\P$.
  The constants in $\asymp$ depend only on
  $\cW$, $\beta$ and $\CH$.

\end{theorem}

\begin{remark}

  Theorem~\ref{thm:minimax}
  verifies that the rates of uniform consistency derived
  in Theorem~\ref{thm:uniform_consistency}
  are minimax-optimal when using a kernel of sufficiently high order
  $(p \geq \beta)$.
  It also shows that both the
  $L_n$ and $E_n$ terms are important
  in the consistency of $\widehat f_W$,
  and that their relative magnitude
  in the supremum norm is determined
  by the degeneracy type of the underlying distribution.

\end{remark}

\subsection{Covariance structure}

\begin{lemma}[Covariance structure]
  \label{lem:covariance_structure}

  Suppose Assumptions~\ref{ass:data}
  and~\ref{ass:kernel_bandwidth} hold.
  Define the covariance function of the
  dyadic kernel density estimator by
  \begin{align*}
    \Sigma_n(w,w')
    &=
    \Cov\Big[
      \widehat f_W(w),
      \widehat f_W(w')
      \Big]
  \end{align*}
  for $w,w' \in \cW$.
  Then $\Sigma_n$ admits the following
  representations.
  \begin{align*}
    \Sigma_n(w,w')
    &=
    \frac{2}{n(n-1)}
    \Cov\hspace*{-0.7mm}\big[
      k_h(W_{ij},w),
      k_h(W_{ij},w')
      \big]
    +
    \frac{4(n-2)}{n(n-1)}
    \Cov\hspace*{-0.7mm}\big[
      k_h(W_{ij},w),
      k_h(W_{ir},w')
      \big] \\
    &=
    \frac{2}{n(n-1)}
    \Cov\hspace*{-0.7mm}\big[
      k_h(W_{ij},w),
      k_h(W_{ij},w')
      \big]
    +
    \frac{4(n-2)}{n(n-1)}
    \Cov\hspace*{-0.7mm}\big[
      \E[k_h(W_{ij},w) \mid A_i],
      \E[k_h(W_{ij},w') \mid A_i]
      \big],
  \end{align*}
  where $1 \leq i < j < r \leq n$.

\end{lemma}

\begin{lemma}[Variance bounds]
  \label{lem:variance_bounds}

  Suppose that Assumptions~\ref{ass:data}
  and~\ref{ass:kernel_bandwidth} hold.
  Then for all large enough $n$,
  \begin{align*}
    \frac{\Dl^2}{n}
    + \frac{1}{n^2h}
    \inf_{w \in \cW} f_W(w)
    &\lesssim
    \inf_{w \in \cW} \Sigma_n(w,w)
    \leq
    \sup_{w \in \cW} \Sigma_n(w,w)
    \lesssim
    \frac{\Du^2}{n}
    + \frac{1}{n^2h}.
  \end{align*}

\end{lemma}

\subsection{Strong approximation}
\label{sec:strong_approx}

In this section we give a strong approximation
for the empirical process $\widehat f_W$.
We begin by using the
K\'omlos--Major--Tusn\'ady (KMT) approximation
to obtain a strong approximation for $L_n$
in Lemma~\ref{lem:strong_approx_Ln}.
Since $E_n$ is an empirical process of i.n.i.d.\ variables,
the KMT approximation is not valid.
Instead we apply a conditional version of
Yurinskii's coupling to obtain a
conditional strong approximation for $E_n$
in Lemma~\ref{lem:conditional_strong_approx_En},
and then construct an unconditional
strong approximation for $E_n$
in Lemma~\ref{lem:unconditional_strong_approx_En}.
These approximations are combined to give a
strong approximation for $\widehat f_W$
in Theorem~\ref{thm:strong_approx_fW}.
We do not need to construct a strong approximation
for the negligible $Q_n$.

This section is largely concerned with
distributional properties,
and as such will frequently involve
\emph{copies} of processes.
We say that $X'$ is a copy of a random variable $X$
if they have the same distribution,
though they may be defined on different probability spaces.
To ensure that all of the joint distributional properties of
such processes are preserved,
we also carry over a copy of the latent variables
$(\bA_n, \bV_n)$
to the new space.

Many of the technical details regarding
the copying and embedding of stochastic processes
are covered by the
Vorob'ev--Berkes--Philipp Theorem,
which is stated and discussed in
Lemma~\ref{lem:vbp}.
In particular, this theorem can be used
for random vectors or for stochastic processes
indexed by a compact rectangle in $\R^d$
with a.s.\ continuous sample paths.

\begin{lemma}[Strong approximation of $L_n$]
  \label{lem:strong_approx_Ln}

  Suppose that Assumptions
  \ref{ass:data}
  and
  \ref{ass:kernel_bandwidth} hold.
  For each $n \geq 2$
  there exists
  on some probability space
  a copy of $\big(\bA_n, \bV_n, L_n\big)$,
  denoted $\big(\bA_n', \bV_n', L_n'\big)$,
  and a mean-zero Gaussian process
  $Z^{L\prime}_n$
  indexed on $\cW$ satisfying
  \begin{align*}
    \P\left(
      \sup_{w \in \cW}
      \big| \sqrt{n} L'_n(w) -  Z_n^{L\prime}(w)\big|
      > \Du
      \frac{t + C_1 \log n}{\sqrt{n}}
    \right)
    &\leq C_2 e^{-C_3 t},
  \end{align*}
  for some positive constants
  $C_1$, $C_2$, $C_3$
  and for all $t > 0$.
  By integration of tail probabilities,
  \begin{align*}
    \E\left[
      \sup_{w \in \cW}
      \big| \sqrt{n} L_n'(w) -  Z_n^{L\prime}(w)\big|
    \right]
    &\lesssim
    \frac{\Du \log n}{\sqrt{n}}.
  \end{align*}
  Further,
  $Z_n^{L\prime}$ has the same covariance structure as
  $\sqrt{n} L_n'$ in the sense that for all $w, w' \in \cW$,
  \begin{align*}
    \E\left[
      Z_n^{L\prime}(w)
      Z_n^{L\prime}(w')
    \right]
    &=
    n
    \E\left[
      L_n'(w)
      L_n'(w')
    \right].
  \end{align*}
  It also satisfies the following
  trajectory regularity property
  for any $\delta_n \in (0, 1/2]$:
  \begin{align*}
    \E\left[
      \sup_{|w-w'| \leq \delta_n}
      \big|
      Z_n^{L\prime}(w)
      - Z_n^{L\prime}(w')
      \big|
    \right]
    &\lesssim
    \Du
    \delta_n \sqrt{\log 1/\delta_n},
  \end{align*}
  and has continuous trajectories.
  The process $Z_n^{L\prime}$
  is a function only of $\bA_n'$
  and some random noise
  which is independent of $(\bA_n', \bV_n')$.

\end{lemma}

\begin{lemma}[Conditional strong approximation of $E_n$]
  \label{lem:conditional_strong_approx_En}

  Suppose that Assumptions
  \ref{ass:data} and \ref{ass:kernel_bandwidth} hold.
  For each $n \geq 2$
  and $t_n > 0$ with $\left|\log t_n\right| \lesssim \log n$,
  there exists on some probability space
  a copy of
  $\big(\bA_n, \bV_n, E_n\big)$,
  denoted
  $\big(\bA_n', \bV_n', E_n'\big)$,
  and a process
  $\widetilde Z^{E\prime}_n$
  which is Gaussian conditional on $\bA_n'$
  and mean-zero conditional on $\bA_n'$,
  satisfying
  \begin{align*}
    \P\left(
      \sup_{w \in \cW}
      \big|
      \sqrt{n^2h} E_n'(w) - \widetilde Z_n^{E\prime}(w)
      \big|
      > t_n
      \Bigm\vert \bA_n'
    \right)
    &\leq
    C_1
    t_n^{-2}
    n^{-1/2}
    h^{-3/4}
    (\log n)^{3/4},
  \end{align*}
  $\bA_n'$-almost surely
  for some constant $C_1 > 0$.
  Setting $t_n = n^{-1/4} h^{-3/8} (\log n)^{3/8} R_n$
  for any sequence $R_n \to \infty$
  and taking an expectation gives
  \begin{align*}
    \sup_{w \in \cW}
    \big|
    \sqrt{n^2h} E_n'(w) - \widetilde Z_n^{E\prime}(w)
    \big|
    &\lesssim_\P
    n^{-1/4}
    h^{-3/8} (\log n)^{3/8} R_n.
  \end{align*}
  Further,
  $\widetilde Z_n^{E\prime}$ has the same
  conditional covariance structure as
  $\sqrt{n^2h} E_n'$ in the sense that for all $w, w' \in \cW$,
  \begin{align*}
    \E\left[
      \widetilde Z_n^{E\prime}(w)
      \widetilde Z_n^{E\prime}(w')
      \bigm\vert \bA_n'
    \right]
    &=
    n^2h
    \E\left[
      E_n'(w)
      E_n'(w')
      \bigm\vert \bA_n'
    \right].
  \end{align*}
  It also satisfies the following
  trajectory regularity property
  for any $\delta_n \in (0, 1/(2h)]$:
  \begin{align*}
    \E\left[
      \sup_{|w-w'| \leq \delta_n}
      \big|
      \widetilde Z_n^{E\prime}(w)
      - \widetilde Z_n^{E\prime}(w')
      \big|
    \right]
    &\lesssim
    \frac{\delta_n}{h}
    \sqrt{\log \frac{1}{h\delta_n}},
  \end{align*}
  and has continuous trajectories.

\end{lemma}

\begin{lemma}[Unconditional strong approximation of $E_n$]
  \label{lem:unconditional_strong_approx_En}

  Suppose that Assumptions
  \ref{ass:data} and \ref{ass:kernel_bandwidth} hold.
  Let
  $\big(\bA_n', \bV_n', \widetilde Z_n^{E\prime}\big)$
  be defined as in
  Lemma~\ref{lem:conditional_strong_approx_En}.
  For each $n \geq 2$
  there exists
  (on some probability space)
  a copy of
  $\big(\bA_n', \bV_n', \widetilde Z_n^{E\prime}\big)$,
  denoted
  $\big(\bA_n'', \bV_n'', \widetilde Z_n^{E\dprime}\big)$,
  and a centered
  Gaussian process
  $Z^{E\dprime}_n$
  satisfying
  \begin{align*}
    \E\left[
      \sup_{w \in \cW}
      \big|\widetilde Z_n^{E\dprime}(w) - Z_n^{E\dprime}(w)\big|
    \right]
    &\lesssim
    n^{-1/6} (\log n)^{2/3}.
  \end{align*}
  Further,
  $Z_n^{E\dprime}$ has the same
  (unconditional) covariance structure as
  $\widetilde Z_n^{E\dprime}$
  and
  $\sqrt{n^2h} E_n$
  in the sense that for all $w, w' \in \cW$,
  \begin{align*}
    \E\left[
      Z_n^{E\dprime}(w)
      Z_n^{E\dprime}(w')
    \right]
    &=
    \E\left[
      \widetilde Z_n^{E\dprime}(w)
      \widetilde Z_n^{E\dprime}(w')
    \right]
    =
    n^2h \,
    \E\left[
      E_n(w)
      E_n(w')
    \right].
  \end{align*}
  It also satisfies the following
  trajectory regularity property
  for any $\delta_n \in (0, 1/(2h)]$:
  \begin{align*}
    \E\left[
      \sup_{|w-w'| \leq \delta_n}
      \big|
      Z_n^{E\dprime}(w)
      - Z_n^{E\dprime}(w')
      \big|
    \right]
    &\lesssim
    \frac{\delta_n}{h}
    \sqrt{\log \frac{1}{h\delta_n}}.
  \end{align*}
  Finally, $Z_n^{E\dprime}$ is independent of $\bA_n''$
  and has continuous trajectories.

\end{lemma}

\begin{remark}

  Note that the process $\widetilde Z_n^{E\prime}$,
  constructed in Lemma~\ref{lem:conditional_strong_approx_En},
  is a conditionally Gaussian process but is
  not in general a Gaussian process.
  The process $Z_n^{E\dprime}$,
  constructed in Lemma~\ref{lem:unconditional_strong_approx_En},
  is a true Gaussian process.

\end{remark}

\begin{theorem}[Strong approximation of $\widehat f_W$]
  \label{thm:strong_approx_fW}

  Suppose that Assumptions
  \ref{ass:data} and \ref{ass:kernel_bandwidth} hold.
  For each $n \geq 2$
  and any sequence $R_n \to \infty$
  there exists on some probability space
  a centered Gaussian process $Z_n^{f\prime}$
  and a copy of $\widehat f_W$,
  denoted $\widehat f_W'$,
  satisfying
  \begin{align*}
    \sup_{w \in \cW}
    \Big|
    \widehat f_W'(w) - \E[\widehat f_W'(w)]
    - Z_n^{f\prime}(w)
    \Big|
    &\lesssim_\P
    n^{-1} \log n
    + n^{-5/4} h^{-7/8} (\log n)^{3/8} R_n
    + n^{-7/6} h^{-1/2} (\log n)^{2/3}.
  \end{align*}
  Further, $Z_n^{f\prime}$ has the same covariance
  structure as
  $\widehat f_W'(w)$
  in the sense that for all
  $w, w' \in \cW$,
  \begin{align*}
    \E\big[Z_n^{f\prime}(w) Z_n^{f\prime}(w')\big]
    &=
    \Cov\Big[
      \widehat f_W'(w),
      \widehat f_W'(w')
      \Big]
    = \Sigma_n(w,w').
  \end{align*}
  It also has continuous trajectories satisfying
  the following trajectory regularity property
  for any $\delta_n \in (0, 1/2]$:
  \begin{align*}
    \E\left[
      \sup_{|w-w'| \leq \delta_n}
      \Big|
      Z_n^{f\prime}(w)
      - Z_n^{f\prime}(w')
      \Big|
    \right]
    &\lesssim
    \frac{\Du}{\sqrt n} \delta_n
    \sqrt{\log \frac{1}{\delta_n}}
    + \frac{1}{\sqrt{n^2h}}
    \frac{\delta_n}{h}
    \sqrt{\log \frac{1}{h\delta_n}}.
  \end{align*}
\end{theorem}

\begin{remark}
  The interpretation of
  Theorem~\ref{thm:strong_approx_fW}
  is deferred to
  Section~\ref{sec:infeasible_confidence_bands},
  in which we scale the processes by their pointwise variance
  in order to better understand the role of degeneracy
  on strong approximation rates.
\end{remark}

\subsection{Infeasible uniform confidence bands}
\label{sec:infeasible_confidence_bands}

We use the strong approximation and bias results
from Theorem~\ref{thm:strong_approx_fW}
and Lemma~\ref{lem:bias}
respectively to construct uniform confidence bands for
the true density function $f_W$.
From now on we will drop the prime notation for copies of processes
in the interest of clarity.
In this section we will assume oracle knowledge of the true
covariance function $\Sigma_n$,
which is not typically available in practice.
For feasible versions of these results which use a covariance estimator,
see Section~\ref{sec:feasible_confidence_bands}.
We also assume that the true density $f_W$
is bounded away from zero on the domain of inference,
which is a standard assumption when constructing confidence bands.

\begin{lemma}[Infeasible Gaussian approximation
    of the standardized $t$-statistic]
  \label{lem:infeasible_t_statistic}

  Let Assumptions \ref{ass:data} and
  \ref{ass:kernel_bandwidth}
  hold and suppose that
  $f_W(w) > 0$ on $\cW$.
  Define for $w \in \cW$
  \begin{align*}
    T_n(w) = \frac{\widehat f_W(w) - f_W(w)}{\sqrt{\Sigma_n(w,w)}}
    \quad\text{and}\quad
    Z_n^T(w) = \frac{Z_n^f(w)}{\sqrt{\Sigma_n(w,w)}}.
  \end{align*}
  Then with $R_n \to \infty$ as in
  Theorem~\ref{thm:strong_approx_fW},
  \begin{align*}
    \sup_{w \in \cW} \left| T_n(w) - Z_n^T(w) \right|
    &\lesssim_\P
    \frac{
      n^{-1/2} \log n
      + n^{-3/4} h^{-7/8} (\log n)^{3/8} R_n
      + n^{-2/3} h^{-1/2} (\log n)^{2/3}
      + n^{1/2} h^{p \wedge \beta}}
    {\Dl + 1/\sqrt{nh}}.
  \end{align*}

\end{lemma}

For the coverage rate error
in the upcoming Theorem~\ref{thm:infeasible_ucb}
to converge to zero in large samples,
we require further restrictions on the bandwidth sequence.
These restrictions depend on the degeneracy type
of the dyadic distribution,
and are given as Assumption~\ref{ass:rates}.

\begin{assumption}[Rate restriction for uniform confidence bands]
  \label{ass:rates}

  Suppose that one of the following holds.

  \begin{enumerate}[label=(\roman*)]

    \item
      No degeneracy:
      $\Dl > 0$ and
      $n^{-6/7} \log n \ll h
      \ll (n \log n)^{-\frac{1}{2(p \wedge \beta)}}$.

    \item
      Partial or total degeneracy:
      $\Dl = 0$ and
      $n^{-2/3} (\log n)^{7/3} \ll h
      \ll (n^2 \log n)^{-\frac{1}{2(p \wedge \beta) + 1}}$.

  \end{enumerate}

\end{assumption}

\begin{theorem}[Infeasible uniform confidence bands]
  \label{thm:infeasible_ucb}

  Let Assumptions \ref{ass:data}, \ref{ass:kernel_bandwidth}
  and \ref{ass:rates}
  hold and suppose that $f_W(w) > 0$ on $\cW$.
  Let $\alpha \in (0,1)$ be a confidence level
  and define $q_{1-\alpha}$ as the quantile
  satisfying
  \begin{align*}
    \P\left(
      \sup_{w \in \cW}
      \left| Z_n^T(w) \right|
      \leq q_{1-\alpha}
    \right)
    &=
    1 - \alpha.
  \end{align*}
  Then
  \begin{align*}
    \P\left(
      f_W(w)
      \in
      \left[
        \widehat f_W(w)
        \pm
        q_{1-\alpha}
        \sqrt{\Sigma_n(w,w)}
        \, \right]
      \, \textup{for all }
      w \in \cW
    \right)
    \to 1 - \alpha.
  \end{align*}
\end{theorem}

\begin{remark}

  By Theorem~\ref{thm:uniform_consistency},
  the asymptotically optimal bandwidth choice
  under uniform convergence is given by
  $h \asymp (n^{-2} \log n)^{\frac{1}{2(p \wedge \beta)+1}}$.
  This bandwidth satisfies Assumption~\ref{ass:rates}
  only in the case of no degeneracy.
  Thus in degenerate cases, one must undersmooth by choosing
  a bandwidth which is smaller than the optimal bandwidth.
  This robust bias correction can be achieved in practice by selecting an
  approximately optimal bandwidth
  for a kernel of order $p$, but then using a kernel of higher order
  $p' > p$ for constructing the confidence bands.

\end{remark}

\subsection{Covariance estimation}
\label{sec:covariance_estimation}

In this section
we provide a consistent estimator for the
covariance function $\Sigma_n$.
In Lemma~\ref{lem:covariance_estimation}
we define the estimator and demonstrate that it
converges in probability in a suitable sense
In Lemma~\ref{lem:alternative_covariance_estimator}
we give an alternative representation
which is more amenable to computation.

\begin{lemma}[Covariance estimation]
  \label{lem:covariance_estimation}

  Let Assumptions \ref{ass:data} and \ref{ass:kernel_bandwidth}
  hold and suppose that $nh \gtrsim \log n$
  and $f_W(w) > 0$ on $\cW$.
  For $w,w' \in \cW$ define
  \begin{align*}
    \widehat \Sigma_n(w,w')
    &=
    \frac{4}{n^2(n-1)^2}
    \sum_{i<j}
    k_h(W_{ij},w)
    k_h(W_{ij},w')
    +
    \frac{24}{n^2(n-1)^2}
    \sum_{i<j<r}
    S_{ijr}(w,w')
    -
    \frac{4n-6}{n(n-1)}
    \widehat f_W(w)
    \widehat f_W(w'),
  \end{align*}
  where
  \begin{align*}
    S_{ijr}(w,w')
    &=
    \frac{1}{6}
    \Big(
    k_h(W_{ij},w)
    k_h(W_{ir},w')
    + k_h(W_{ij},w)
    k_h(W_{jr},w')
    + k_h(W_{ir},w)
    k_h(W_{ij},w') \\
    &\qquad+
    k_h(W_{ir},w)
    k_h(W_{jr},w')
    + k_h(W_{jr},w)
    k_h(W_{ij},w')
    + k_h(W_{jr},w)
    k_h(W_{ir},w')
    \Big).
  \end{align*}
  Then $\widehat \Sigma_n$
  is uniformly entrywise-consistent in the sense that
  \begin{align*}
    \sup_{w,w' \in \cW}
    \left|
    \frac{\widehat \Sigma_n(w,w') - \Sigma_n(w,w')}
    {\sqrt{\Sigma_n(w,w) + \Sigma_n(w',w')}}
    \right|
    &\lesssim_\P
    \frac{\sqrt{\log n}}{n}.
  \end{align*}

\end{lemma}

\begin{lemma}[Alternative covariance estimator representation]
  \label{lem:alternative_covariance_estimator}

  Suppose that Assumptions~\ref{ass:data}
  and~\ref{ass:kernel_bandwidth} hold
  and let $\widehat \Sigma_n$
  be the covariance estimator defined
  in Lemma~\ref{lem:covariance_estimation}.
  Then the following alternative representation
  for $\widehat \Sigma_n$ holds,
  which may be easier to compute
  as it does not involve any triple summations
  over the data.
  \begin{align*}
    \widehat \Sigma_n(w,w')
    &=
    \frac{4}{n^2}
    \sum_{i=1}^n
    S_i(w) S_i(w')
    - \frac{4}{n^2(n-1)^2}
    \sum_{i<j}
    k_h(W_{ij},w)
    k_h(W_{ij},w')
    -
    \frac{4n-6}{n(n-1)}
    \widehat f_W(w)
    \widehat f_W(w'),
  \end{align*}
  where
  \begin{align*}
    S_i(w)
    &=
    \frac{1}{n-1}
    \left(
      \sum_{j = 1}^{i-1}
      k_h(W_{ji}, w)
      + \sum_{j = i+1}^n
      k_h(W_{ij}, w)
    \right)
  \end{align*}
  is an estimator of
  $\E[k_h(W_{ij},w) \mid A_i]$.

\end{lemma}

\begin{remark}

  The covariance estimator $\widehat \Sigma_n$
  is not necessarily almost surely positive semi-definite.

\end{remark}

\subsection{Positive semi-definite covariance estimation}
\label{sec:PSD}

In this section we provide a
positive semi-definite estimator $\widehat \Sigma_n^+$
which is uniformly entrywise-consistent for $\Sigma_n$.
Define $\widehat \Sigma_n$ as in
Lemma~\ref{lem:covariance_estimation}
and consider the following optimization problem
over bivariate functions.
\begin{equation}
  \label{eq:SDP}
  \begin{aligned}
    \minimize
    \qquad
     & \sup_{w,w' \in \cW}
    \left|
    \frac{M(w,w') - \widehat\Sigma_n(w,w')}
    {\sqrt{\widehat \Sigma_n(w,w) + \widehat \Sigma_n(w',w')}}
    \right|
    \quad \textup{ over } M: \cW \times \cW \to \R
    \\
    \subjectto
    \qquad
     & M \textup{ is symmetric and positive semi-definite}, \\
     & \big|M(w,w') - M(w, w'')\big|
    \leq \frac{4}{nh^3}
    \Ck \CL
    |w'-w''|
    \textup{ for all }
    w, w', w'' \in \cW.
  \end{aligned}
\end{equation}

\begin{lemma}[Consistency of $\widehat \Sigma_n^+$]
  \label{lem:SDP}

  Suppose that Assumptions~\ref{ass:data}
  and~\ref{ass:kernel_bandwidth} hold and that
  $nh \gtrsim \log n$ and $f_W(w) > 0$ on $\cW$.
  Then the optimization problem (\ref{eq:SDP})
  has an approximately optimal solution $\widehat\Sigma_n^+$
  which is uniformly entrywise-consistent
  for $\Sigma_n$ in the sense that
  \begin{align*}
    \sup_{w,w' \in \cW}
    \left|
    \frac{\widehat \Sigma_n^+(w,w') - \Sigma_n(w,w')}
    {\sqrt{\Sigma_n(w,w) + \Sigma_n(w',w')}}
    \right|
    &\lesssim_\P
    \frac{\sqrt{\log n}}{n}.
  \end{align*}

\end{lemma}

\begin{remark}

  The optimization problem (\ref{eq:SDP}) is stated for functions
  rather than for matrices so is infinite-dimensional.
  However, when restricting to finite-size matrices,
  Lemma~\ref{lem:SDP} still holds and does not depend on
  the size of the matrices.
  Furthermore, the problem
  then becomes a semi-definite program and so can be solved to
  arbitrary precision in polynomial time
  in the size of the matrices
  \citep{laurent2005semidefinite}.

  The Lipschitz-type constraint in the
  optimization problem (\ref{eq:SDP})
  ensures that $\widehat \Sigma_n^+$
  is sufficiently smooth
  and is a technicality required by
  some of the later proofs.
  In practice this constraint is readily verified.

\end{remark}

\begin{lemma}[Positive semi-definite variance estimator bounds]
  \label{lem:variance_estimator_bounds}

  Suppose that Assumptions~\ref{ass:data}
  and~\ref{ass:kernel_bandwidth} hold and that
  $nh \gtrsim \log n$ and $f_W(w) > 0$ on $\cW$.
  Then $\widehat \Sigma_n^+(w,w) \geq 0$
  almost surely for all $w \in \cW$ and
  \begin{align*}
    \frac{\Dl^2}{n} + \frac{1}{n^2h}
    &\lesssim_\P
    \inf_{w \in \cW} \widehat \Sigma_n^+(w,w)
    \leq
    \sup_{w \in \cW} \widehat \Sigma_n^+(w,w)
    \lesssim_\P
    \frac{\Du^2}{n} + \frac{1}{n^2h}.
  \end{align*}

\end{lemma}

\subsection{Feasible uniform confidence bands}
\label{sec:feasible_confidence_bands}

Now we use the strong approximation derived in
Section~\ref{sec:strong_approx} and the
positive semi-definite covariance estimator introduced in
Section~\ref{sec:PSD} to construct feasible uniform confidence bands.

\begin{lemma}[Proximity of the standardized and studentized t-statistics]
  \label{lem:studentized_t_statistic}

  Let Assumptions \ref{ass:data} and
  \ref{ass:kernel_bandwidth} hold and suppose that
  $nh \gtrsim \log n$ and $f_W(w) > 0$ on $\cW$.
  Define for $w \in \cW$
  \begin{align*}
    \widehat T_n(w) = \frac{\widehat f_W(w) - f_W(w)}
    {\sqrt{\widehat\Sigma_n^+(w,w)}}.
  \end{align*}
  Then
  \begin{align*}
    \sup_{w \in \cW}
    \left| \widehat T_n(w) - T_n(w) \right|
    &\lesssim_\P
    \sqrt{\frac{\log n}{n}}
    \left(
      \sqrt{\log n} + \frac{\sqrt n h^{p \wedge \beta}}
      {\Dl + 1/\sqrt{nh}}
    \right)
    \frac{1}{\Dl + 1/\sqrt{nh}}.
  \end{align*}

\end{lemma}

\begin{lemma}[Feasible Gaussian approximation
    of the infeasible Gaussian process]
  \label{lem:distributional_approx_feasible_gaussian}

  Let Assumptions \ref{ass:data} and \ref{ass:kernel_bandwidth}
  hold and suppose that
  $nh \gtrsim \log n$ and $f_W(w) > 0$ on $\cW$.
  Define a process $\widehat Z_n^T(w)$ which,
  conditional on the data $\bW_n$,
  is conditionally mean-zero and
  conditionally Gaussian and whose
  conditional covariance structure is
  \begin{align*}
    \E\left[
      \widehat Z_n^T(w) \widehat Z_n^T(w')
      \bigm| \bW_n \right]
    &=
    \frac{\widehat \Sigma_n^+(w,w')}
    {\sqrt{\widehat \Sigma_n^+(w,w) \widehat \Sigma_n^+(w',w')}}
  \end{align*}
  Then the following conditional
  Kolmogorov--Smirnov result holds.
  \begin{align*}
    \sup_{t \in \R}
    \left|
    \P\left(
      \sup_{w \in \cW}
      \left| Z_n^T(w) \right|
      \leq t
    \right)
    - \P\left(
      \sup_{w \in \cW}
      \left| \widehat Z_n^T(w) \right|
      \leq t
      \biggm\vert \bW_n
    \right)
    \right|
    &\lesssim_\P
    \frac{n^{-1/6}(\log n)^{5/6}}
    {\Dl^{1/3} + (nh)^{-1/6}}.
  \end{align*}

\end{lemma}

\begin{lemma}[Feasible Gaussian approximation of the studentized $t$-statistic]
  \label{lem:feasible_gaussian_approx}

  Let Assumptions \ref{ass:data}, \ref{ass:kernel_bandwidth}
  and \ref{ass:rates} hold and suppose that $f_W(w) > 0$ on $\cW$.
  Then
  \begin{align*}
    \sup_{t \in \R}
    \left|
    \P\left(
      \sup_{w \in \cW}
      \left| \widehat T_n(w) \right|
      \leq t
    \right)
    - \P\left(
      \sup_{w \in \cW}
      \left| \widehat Z_n^T(w) \right|
      \leq t
      \Bigm\vert \bW_n
    \right)
    \right|
    &\ll_\P
    1.
  \end{align*}

\end{lemma}

\begin{theorem}[Feasible uniform confidence bands]
  \label{thm:ucb}

  Let Assumptions \ref{ass:data}, \ref{ass:kernel_bandwidth}
  and \ref{ass:rates} hold and suppose that $f_W(w) > 0$ on $\cW$.
  Let $\alpha \in (0,1)$ be a confidence level
  and define $\widehat q_{1-\alpha}$ as the
  conditional quantile satisfying
  \begin{align*}
    \P\left(
      \sup_{w \in \cW}
      \left| \widehat Z_n^T(w) \right|
      \leq \widehat q_{1-\alpha}
      \Bigm\vert \bW_n
    \right)
    &=
    1 - \alpha.
  \end{align*}
  Then
  \begin{align*}
    \P\left(
      f_W(w)
      \in
      \left[
        \widehat f_W(w)
        \pm
        \widehat q_{1-\alpha}
        \sqrt{\widehat\Sigma_n^+(w,w)}
        \,\right]
      \,\textup{for all }
      w \in \cW
    \right)
    \to 1 - \alpha.
  \end{align*}

\end{theorem}

\begin{remark}

  In practice, suprema over $\cW$ can be
  replaced by maxima over a sufficiently fine finite partition of $\cW$.
  The conditional quantile $\widehat q_{1-\alpha}$
  can be estimated by Monte Carlo
  simulation,
  resampling from the Gaussian process defined by the law of
  $\widehat Z_n^T \mid \bW_n$.
\end{remark}

\subsection{Counterfactual
  dyadic density estimation}
\label{subsec:counterfactual}

As an application we provide methodology for estimation and inference
on counterfactual dyadic density functions,
following the reweighting approach of
\citet{dinardo1996distribution}.
We give the counterfactual data
generating process as
Assumption~\ref{ass:counterfactual}.

\begin{assumption}[Counterfactual data generation]
  \label{ass:counterfactual}

  For each $r \in \{0,1\}$,
  let $\bW_n^r$, $\bA_n^r$ and $\bV_n^r$ be as in Assumption~\ref{ass:data}.
  Let $X_i^r$ be finitely-supported variables,
  setting $\bX_n^r = (X_1^r, \ldots, X_n^r)$.
  Suppose that $(A_i^r, X_i^r)$ are
  independent over $1 \leq i \leq n$
  and that $\bX_n^r$ is independent of $\bV_n^r$.
  Assume that $W_{ij}^r \mid X_i^r, X_j^r$
  has a Lebesgue density
  $f_{W \mid XX}^r(\,\cdot \mid x_1, x_2) \in \cH^\beta_{\CH}(\cW)$
  and that $X_i^r$ has positive
  probability mass function
  $p_X^r(x)$ on a common support $\cX$.
  Suppose that
  $(\bA_n^0, \bV_n^0, \bX_n^0)$
  and $(\bA_n^1, \bV_n^1, \bX_n^1)$
  are independent.

\end{assumption}

The counterfactual density of
$W_{ij}$ in population $1$
had $X_i, X_j$
followed the distribution in population $0$ by
\begin{align*}
  f_W^{1 \triangleright 0}(w)
  &=
  \E\left[
    f_{W \mid XX}^1\big(w \mid X_1^0, X_2^0\big)
  \right]
  = \sum_{x_1 \in \cX}
  \sum_{x_2 \in \cX}
  f_{W \mid XX}^{1}(w \mid x_1, x_2)
  \psi(x_1)
  \psi(x_2)
  p_X^{1}(x_1)
  p_X^{1}(x_2)
\end{align*}
where $\psi(x) = p_X^0(x)/p_X^1(x)$ for $x \in \cX$.
This counterfactual density can be estimated by
the counterfactual dyadic kernel density estimator
\begin{align*}
  \widehat f_W^{\,1 \triangleright 0}(w)
  &=
  \frac{2}{n(n-1)}
  \sum_{i=1}^{n-1}
  \sum_{j=i+1}^n
  \widehat \psi(X_i^1)
  \widehat \psi(X_j^1)
  k_h(W_{ij}^1, w)
\end{align*}
with
$\widehat\psi(x) =
\I\{\widehat p_X^{\,1}(x) > 0\}\widehat p_X^{\,0}(x) / \widehat p_X^{\,1}(x)$
and $\widehat p_X^{\,r}(x) = \frac{1}{n}\sum_{i = 1}^n \I\{X_i^r = x\}$.
Note that since $p_X^r(x) > 0$,
\begin{align*}
  \widehat\psi(x) - \psi(x)
  &=
  \frac{\widehat p_X^{\,0}(x) - p_X^0(x)}{p_X^1(x)}
  - \frac{p_X^0(x)}{p_X^1(x)}
  \frac{\widehat p_X^{\,1}(x) - p_X^1(x)}{p_X^1(x)}
  + \frac{\widehat p_X^{\,1}(x) - p_X^1(x)}{p_X^1(x)}
  \frac{\widehat p_X^{\,1}(x) p_X^0(x) - \widehat p_X^{\,0}(x)p_X^1(x)}
  {\widehat p_X^{\,1}(x) p_X^1(x)} \\
  &=
  \frac{1}{n}
  \sum_{r=1}^n \kappa(X_r^0, X_r^1, x)
  + O_\P\left(\frac{1}{n}\right)
\end{align*}
is an asymptotic linear representation where
\begin{align*}
  \kappa(X_i^0, X_i^1, x)]
  &=
  \frac{\I\{X_i^0 = x\} - p_X^0(x)}{p_X^1(x)}
  - \frac{p_X^0(x)}{p_X^1(x)}
  \frac{\I\{X_i^1 = x\} - p_X^1(x)}{p_X^1(x)}
\end{align*}
satisfies
$\E[\kappa(X_i^0, X_i^1, x)] = 0$.
We now establish uniform consistency and feasible strong
approximation results for the counterfactual density estimator.

\begin{lemma}[Bias of $\widehat f_W^{\,1 \triangleright 0}$]
  \label{lem:counterfactual_bias}

  Suppose that Assumptions~\ref{ass:data},
  \ref{ass:kernel_bandwidth} and \ref{ass:counterfactual} hold.
  Then
  \begin{align*}
    \sup_{w \in \cW}
    \big|
    \E\big[\widehat f_W^{\,1 \triangleright 0}(w)\big]
    - f_W^{1 \triangleright 0}(w)
    \big|
    \lesssim
    h^{p \wedge \beta} + \frac{1}{n}.
  \end{align*}

\end{lemma}

\begin{lemma}[Hoeffding-type decomposition for
  $\widehat f_W^{\,1 \triangleright 0}$]
  \label{lem:counterfactual_hoeffding}

  Suppose that Assumptions~\ref{ass:data},
  \ref{ass:kernel_bandwidth} and
  \ref{ass:counterfactual} hold.
  Writing $k_{ij} = k_h(W_{ij}^1, w)$,
  $\kappa_{r i} = \kappa(X_r^0, X_r^1, X_i^1)$
  and $\psi_i = \psi(X_i^1)$, define the projections
  \begin{align*}
    u
    &=
    \E\left[
      k_{ij}
      \psi_i
      \psi_j
    \right], \\
    u_i
    &=
    \frac{2}{3} \psi_i
    \E\left[
      k_{ij}
      \psi_j
      \mid A_i^1 \right]
    +
    \frac{2}{3} \E\left[
      k_{jr}
      \psi_j \kappa_{ir}
      \mid X_i^0, X_i^1 \right]
    - \frac{2}{3} u, \\
    u_{ij}
    &=
    \frac{1}{3}
    \psi_i
    \psi_j
    \E\left[
      k_{ij}
      \mid A_i^1, A_j^1 \right]
    +
    \frac{1}{3}
    \psi_i
    \E\left[
      k_{ir} \psi_r
      \mid A_i^1 \right]
    +
    \frac{1}{3}
    \psi_i
    \E\left[
      k_{ir}  \kappa_{jr}
      \mid A_i^1, X_j^0, X_j^1 \right]
    +
    \frac{1}{3}
    \kappa_{ji}
    \E\left[
      k_{ir} \psi_r
      \mid A_i^1 \right] \\
    &\quad+
    \frac{1}{3}
    \psi_j
    \E\left[
      k_{jr}  \psi_r
      \mid A_j^1 \right]
    +
    \frac{1}{3}
    \psi_j
    \E\left[
      k_{jr}  \kappa_{ir}
      \mid X_i^0, X_i^1, A_j^1 \right]
    +
    \frac{1}{3}
    \kappa_{ij}
    \E\left[
      k_{jr} \psi_r
      \mid A_j^1 \right]
    - u_i - u_j + u, \\
    u_{ijr}
    &=
    \frac{1}{3}
    \psi_i \psi_j
    \E\left[
      k_{ij}
      \mid A_i^1, A_j^1 \right]
    +
    \frac{1}{3}
    \psi_i \kappa_{r j}
    \E\left[
      k_{ij}
      \mid A_i^1, A_j^1 \right]
    +
    \frac{1}{3}
    \psi_j \kappa_{r i}
    \E\left[
      k_{ij}
      \mid A_i^1, A_j^1 \right]
    + \frac{1}{3}
    \psi_i \psi_r
    \E\left[
      k_{ir}
      \mid A_i^1, A_r^1 \right] \\
    &\quad+
    \frac{1}{3}
    \psi_i \kappa_{jr}
    \E\left[
      k_{ir}
      \mid A_i^1, A_r^1 \right]
    +
    \frac{1}{3}
    \psi_r \kappa_{ji}
    \E\left[
      k_{ir}
      \mid A_i^1, A_r^1 \right]
    + \frac{1}{3}
    \psi_j \psi_r
    \E\left[
      k_{jr}
      \mid A_j^1, A_r^1 \right] \\
    &\quad+
    \frac{1}{3}
    \psi_j \kappa_{ir}
    \E\left[
      k_{jr}
      \mid A_j^1, A_r^1 \right]
    +
    \frac{1}{3}
    \psi_r \kappa_{ij}
    \E\left[
      k_{jr}
      \mid A_j^1, A_r^1 \right]
    - u_{ij} - u_{ir} - u_{jr}
    + u_i + u_j + u_r
    - u, \\
    v_{ijr}
    &=
    \frac{1}{3}
    k_{ij} \big(\psi_i \psi_j +\psi_i \kappa_{r j} +\psi_j \kappa_{r i} \big)
    + \frac{1}{3}
    k_{ir} \big(\psi_i \psi_r +\psi_i \kappa_{jr} +\psi_r \kappa_{ji} \big)
    + \frac{1}{3}
    k_{jr} \big(\psi_j \psi_r +\psi_j \kappa_{ir} +\psi_r \kappa_{ij} \big).
  \end{align*}
  With $l_i^{1 \triangleright 0}(w) = u_i$
  and $e_{ijr}^{1 \triangleright 0}(w) = v_{ijr} - u_{ijr}$,
  set
  \begin{align*}
    L_n^{1 \triangleright 0}(w)
    &=
    \frac{3}{n} \sum_{i=1}^n
    l_i^{1 \triangleright 0}(w)
    &\text{and}                 &
    &E_n^{1 \triangleright 0}(w)
    &=
    \frac{6}{n(n-1)(n-2)}
    \sum_{i=1}^{n-2}
    \sum_{j=i+1}^{n-1}
    \sum_{r=i+1}^n
    e_{ijr}^{1 \triangleright 0}(w).
  \end{align*}
  Then the following Hoeffding-type decomposition holds,
  where $O_\P(1/n)$ is uniform in $w \in \cW$.
  \begin{align*}
    \widehat f_W^{\,1 \triangleright 0}(w)
    = \E\big[\widehat f_W^{\,1 \triangleright 0}(w)\big]
    + L_n^{1 \triangleright 0}(w)
    + E_n^{1 \triangleright 0}(w)
    + O_\P\left( \frac{1}{n} \right).
  \end{align*}
  Further,
  the stochastic processes
  $L_n^{1 \triangleright 0}$
  and $E_n^{1 \triangleright 0}$
  are mean-zero
  and mutually orthogonal
  in $L^2(\P)$.
  Define the upper and lower degeneracy constants as
  \begin{align*}
    \left(
      \Du^{1 \triangleright 0}
    \right)^2
    &=
    \limsup_{n \to \infty}
    \sup_{w \in \cW}
    \Var\big[
    l_i^{1 \triangleright 0}(w)
    \big]
    &\text{and}&
    &
    \left(
      \Dl^{1 \triangleright 0}
    \right)^2
    &=
    \liminf_{n \to \infty}
    \inf_{w \in \cW}
    \Var\big[
    l_i^{1 \triangleright 0}(w)
    \big].
  \end{align*}

\end{lemma}

\begin{lemma}[Uniform consistency of $\widehat f_W^{\,1 \triangleright 0}$]
  \label{lem:counterfactual_uniform_consistency}

  Suppose that Assumptions~\ref{ass:data},
  \ref{ass:kernel_bandwidth} and \ref{ass:counterfactual} hold.
  Then
  \begin{align*}
    \E\left[
      \sup_{w \in \cW}
      \big|\widehat f_W^{\,1 \triangleright 0}(w)
      - f_W^{1 \triangleright 0}(w)
    \right]
    &\lesssim
    h^{p \wedge \beta}
    + \frac{\Du^{1 \triangleright 0}}{\sqrt n}
    + \sqrt{\frac{\log n}{n^2h}}.
  \end{align*}

\end{lemma}

\begin{lemma}[Strong approximation of $\widehat f_W^{\,1 \triangleright 0}$]
  \label{lem:counterfactual_sa}

  On an appropriately enlarged probability space
  and for any sequence $R_n \to \infty$,
  there exists a mean-zero Gaussian process
  $Z_n^{f, 1 \triangleright 0}$
  with the same covariance structure as
  $\widehat f_W^{\,1 \triangleright 0}(w)$ satisfying
  \begin{align*}
    \sup_{w \in \cW}
    \left|
    \widehat f_W^{\,1 \triangleright 0}(w)
    - \E\big[\widehat f_W^{\,1 \triangleright 0}(w)\big]
    - Z_n^{f, 1 \triangleright 0}(w)
    \right|
    &\lesssim_\P
    n^{-1} \log n
    + n^{-5/4} h^{-7/8} (\log n)^{3/8} R_n
    + n^{-7/6} h^{-1/2} (\log n)^{2/3}.
  \end{align*}

\end{lemma}

\begin{lemma}[Counterfactual covariance structure]
  \label{lem:counterfactual_covariance_structure}

  Writing $k_{ij}'$ for $k_h(W_{ij}^1, w')$ etc.,
  the counterfactual covariance function is
  \begin{align*}
    \Sigma_n^{1 \triangleright 0}(w,w')
    &=
    \Cov\left[
      \widehat f_W^{\,1 \triangleright 0}(w),
      \widehat f_W^{\,1 \triangleright 0}(w')
    \right] \\
    &=
    \frac{4}{n}
    \E\left[
      \Big(
      \psi_i
      \E\big[
        k_{ij} \psi_j
        \mid A_i^1
        \big]
      + \E\left[
        k_{r j} \psi_r
        \kappa_{ij}
        \mid X_i^0, X_i^1
      \right]
      \Big)
      \Big(
      \psi_i
      \E\big[
        k_{ij}' \psi_j
        \mid A_i^1
        \big]
      + \E\left[
        k_{r j}' \psi_r \kappa_{ij}
        \mid X_i^0, X_i^1
      \right]
      \Big)
    \right] \\
    &\quad+
    \frac{2}{n^2}
    \E\left[
      k_{ij} k_{ij}'
      \psi_i^2 \psi_j^2
    \right]
    - \frac{4}{n}
    \E\left[
      k_{ij} \psi_i \psi_j
    \right]
    \E\left[
      k_{ij}' \psi_i \psi_j
    \right]
    + O\left( \frac{1}{n^{3/2}} + \frac{1}{\sqrt{n^4h}} \right).
  \end{align*}

\end{lemma}

\begin{lemma}[Infeasible Gaussian approximation
    of the standardized counterfactual $t$-statistic]
  \label{lem:counterfactual_infeasible_t_statistic}

  Let Assumptions \ref{ass:data},
  \ref{ass:kernel_bandwidth} and
  \ref{ass:counterfactual}
  hold and suppose that
  $f_W^{1 \triangleright 0}(w) > 0$ on $\cW$.
  Define for $w \in \cW$
  \begin{align*}
    T_n^{1 \triangleright 0}(w)
    &= \frac{\widehat f_W^{\,1 \triangleright 0}(w)
    - f_W^{1 \triangleright 0}(w)}
    {\sqrt{\Sigma_n^{1 \triangleright 0}(w,w)}}
    \quad\text{and}\quad
    Z_n^{T, 1 \triangleright 0}(w)
    = \frac{Z_n^{f, 1 \triangleright 0}(w)}
    {\sqrt{\Sigma_n^{1 \triangleright 0}(w,w)}}.
  \end{align*}
  Then with $R_n \to \infty$ as in Lemma~\ref{lem:counterfactual_sa},
  \begin{align*}
    \sup_{w \in \cW}
    \left|
    T_n^{1 \triangleright 0}(w) - Z_n^{T, 1 \triangleright 0}(w)
    \right|
    &\lesssim_\P
    \frac{
      n^{-1/2} \log n
      + n^{-3/4} h^{-7/8} (\log n)^{3/8} R_n
      + n^{-2/3} h^{-1/2} (\log n)^{2/3}
      + n^{1/2} h^{p \wedge \beta}}
    {\Dl^{1 \triangleright 0} + 1/\sqrt{nh}}.
  \end{align*}

\end{lemma}

\begin{theorem}[Infeasible counterfactual uniform confidence bands]
  \label{thm:counterfactual_infeasible_ucb}

  Let Assumptions \ref{ass:data}, \ref{ass:kernel_bandwidth},
  \ref{ass:rates} and \ref{ass:counterfactual}
  hold and suppose that $f_W^{1 \triangleright 0}(w) > 0$ on $\cW$.
  Let $\alpha \in (0,1)$ be a confidence level
  and define $q^{1 \triangleright 0}_{1-\alpha}$ as the quantile
  satisfying
  \begin{align*}
    \P\left(
      \sup_{w \in \cW}
      \left| Z_n^{T,1 \triangleright 0}(w) \right|
      \leq q^{1 \triangleright 0}_{1-\alpha}
    \right)
    &=
    1 - \alpha.
  \end{align*}
  Then
  \begin{align*}
    \P\left(
      f_W^{1 \triangleright 0}(w)
      \in
      \left[
        \widehat f_W^{1 \triangleright 0}(w)
        \pm
        q^{1 \triangleright 0}_{1-\alpha}
        \sqrt{\Sigma_n^{1 \triangleright 0}(w,w)}
        \, \right]
      \, \textup{for all }
      w \in \cW
    \right)
    \to 1 - \alpha.
  \end{align*}
\end{theorem}
To conclude this section we propose an estimator
for the counterfactual covariance function
$\Sigma_n^{1 \triangleright 0}$.
First let
\begin{align*}
  \widehat\kappa(X_i^0, X_i^1, x)]
  &=
  \frac{\I\{X_i^0 = x\} - \widehat p_X^0(x)}{\widehat p_X^1(x)}
  - \frac{\widehat p_X^0(x)}{\widehat p_X^1(x)}
  \frac{\I\{X_i^1 = x\} - \widehat p_X^1(x)}{\widehat p_X^1(x)}
\end{align*}
and define the leave-out conditional expectation estimators
\begin{align*}
  S_i^{1 \triangleright 0}(w)
  &=
  \widehat\E\left[
    k_h(W_{ij}^1,w) \psi(X_j^1) \mid A_i^1
  \right]
  =
  \frac{1}{n-1}
  \left(
    \sum_{j=1}^{i-1}
    k_h(W_{ji}^1,w) \widehat\psi(X_j^1)
    + \sum_{j=i+1}^n
    k_h(W_{ij}^1,w) \widehat\psi(X_j^1)
  \right), \\
  \widetilde S_i^{1 \triangleright 0}(w)
  &=
  \widehat\E\left[
    k_h(W_{r j}^1,w) \psi(X_r^1)
    \kappa(X_i^0, X_i^1, X_j^1) \mid X_i^0, X_i^1
  \right]
  = \frac{1}{n-1}
  \sum_{j=1}^n
  \I\{j \neq i\}
  \widehat\kappa(X_i^0, X_i^1, X_j^1)
  S_j^{1 \triangleright 0}(w).
\end{align*}
Then set
\begin{align*}
  \widehat\Sigma_n^{1 \triangleright 0}(w,w')
  &=
  \frac{4}{n^2}
  \sum_{i=1}^n
  \left(
    \widehat\psi(X_i^1)
    S_i^{1 \triangleright 0}(w)
    + \widetilde S_i^{1 \triangleright 0}(w)
  \right)
  \left(
    \widehat\psi(X_i^1)
    S_i^{1 \triangleright 0}(w')
    + \widetilde S_i^{1 \triangleright 0}(w')
  \right) \\
  &\quad-
  \frac{4}{n^3(n-1)}
  \sum_{i<j}
  k_h(W_{ij}^1, w)
  k_h(W_{ij}^1, w')
  \widehat\psi(X_i^1)^2
  \widehat\psi(X_j^1)^2
  - \frac{4}{n}
  \widehat f_W^{\,1 \triangleright 0}(w)
  \widehat f_W^{\,1 \triangleright 0}(w').
\end{align*}
We then use a positive semi-definite approximation to
$\widehat\Sigma_n^{1 \triangleright 0}$, denoted by
$\widehat\Sigma_n^{+, 1 \triangleright 0}$,
following the methodology of
Section~\ref{sec:PSD}.
We omit the proof of consistency of these covariance estimators
in the interest of brevity.
To construct feasible uniform confidence bands,
define a process $\widehat Z_n^{T, 1 \triangleright 0}(w)$ which,
conditional on the data $\bW_n^1$, $\bX_n^0$ and $\bX_n^1$
is conditionally mean-zero and conditionally Gaussian and whose
conditional covariance structure is
\begin{align*}
  \E\left[
    \widehat Z_n^{T, 1 \triangleright 0}(w)
    \widehat Z_n^{T, 1 \triangleright 0}(w')
    \bigm| \bW_n^1, \bX_n^0, \bX_n^1 \right]
  &=
  \frac{\widehat \Sigma_n^{+, 1 \triangleright 0}(w,w')}
  {\sqrt{\widehat \Sigma_n^{+, 1 \triangleright 0}(w,w)
      \widehat \Sigma_n^{+, 1 \triangleright 0}(w',w')}}.
\end{align*}
Let $\alpha \in (0,1)$ be a confidence level and define
$\widehat q_{1-\alpha}^{1 \triangleright 0}$
as the conditional quantile satisfying
\begin{align*}
  \P\left(
    \sup_{w \in \cW}
    \left| \widehat Z_n^{T, 1 \triangleright 0}(w) \right|
    \leq \widehat q_{1-\alpha}^{1 \triangleright 0}
    \Bigm\vert \bW_n^1, \bX_n^0, \bX_n^1
  \right)
  &=
  1 - \alpha.
\end{align*}
Then assuming that the covariance estimator is appropriately consistent,
we have that
\begin{align*}
  \P\left(
    f_W^{1 \triangleright 0}(w)
    \in
    \left[
      \widehat f_W^{\,1 \triangleright 0}(w)
      \pm
      \widehat q^{1 \triangleright 0}_{1-\alpha}
      \sqrt{\widehat\Sigma_n^{+, 1 \triangleright 0}(w,w)}
      \,\right]
    \,\textup{for all }
    w \in \cW
  \right)
  \to 1 - \alpha.
\end{align*}

\section{Technical lemmas}

In this section we present some lemmas which provide the
technical foundations for several of our main results.
These lemmas are stated in as much generality
as is reasonably possible,
and we believe that they may be
of some independent interest.

\subsection{Maximal inequalities for i.n.i.d.\ empirical processes}

Firstly we provide a maximal inequality
for empirical processes of
independent but not necessarily identically distributed
(i.n.i.d.)
random variables,
indexed by a class of functions.
This result is an extension
of Theorem~5.2 from \citet{chernozhukov2014gaussian},
which only covers i.i.d.\ random variables,
and is proven in the same manner.
Such a result is useful in the study of dyadic data
because when conditioning on latent variables,
we may encounter
random variables which are conditionally independent
but which do not necessarily follow the same
conditional distribution.

\begin{lemma}[A maximal inequality for i.n.i.d.\ empirical processes]
  \label{lem:maximal_entropy}

  Let $X_1, \dots, X_n$
  be independent but not necessarily identically distributed
  (i.n.i.d.)
  random variables taking values in a
  measurable space $(S,\cS)$.
  Denote the joint distribution of
  $X_1,\ldots,X_n$ by $\P$
  and the marginal distribution of
  $X_i$ by $\P_i$,
  and let $\bar\P = n^{-1} \sum_i \P_i$.
  Let $\cF$ be a class of Borel measurable functions
  from $S$ to $\R$
  which is pointwise measurable
  (i.e.\ it contains a countable subclass which
  is dense under pointwise convergence).
  Let $F$ be a strictly positive
  measurable envelope function for $\cF$
  (i.e.\ $|f(s)| \leq |F(s)|$ for all $f \in \cF$
  and $s \in S$).
  For a distribution $\Q$ and some $q \geq 1$,
  define the $(\Q,q)$-norm of $f \in \cF$ as
  $\|f\|_{\Q,q}^q = \E_{X \sim \Q}[f(X)^q]$
  and suppose that
  $\|F\|_{\bar\P,2} < \infty$.
  For $f \in \cF$
  define the empirical process
  \begin{align*}
    G_n(f)
    &=
    \frac{1}{\sqrt n}
    \sum_{i=1}^n
    \Big(
    f(X_i) - \E[f(X_i)]
    \Big).
  \end{align*}
  Let $\sigma > 0$ satisfy
  $\sup_{f \in \cF}
  \|f\|_{\bar\P,2}
  \leq
  \sigma
  \leq
  \|F\|_{\bar\P,2}$
  and
  $M = \max_{1 \leq i \leq n} F(X_i)$.
  Then with
  $\delta = \sigma / \|F\|_{\bar\P,2} \in (0,1]$,
  \begin{align*}
    \E \left[
      \sup_{f \in \cF}
      \big| G_n(f) \big|
    \right]
    &\lesssim
    \|F\|_{\bar\P,2}
    \, J\big(\delta, \cF, F \big)
    +
    \frac{\|M\|_{\P,2} \, J(\delta, \cF, F)^2}{\delta^2 \sqrt{n}},
  \end{align*}
  where $\lesssim$ is up to a universal constant,
  and $J(\delta, \cF, F)$ is the covering integral
  \begin{align*}
    J\big(\delta, \cF, F\big)
    &=
    \int_0^\delta
    \sqrt{1 +
    \sup_\Q \log N(\cF, \rho_\Q, \varepsilon \|F\|_{\Q,2})}
    \diff{\varepsilon},
  \end{align*}
  with the supremum taken over finite discrete probability
  measures $\Q$ on $(S, \cS)$.

\end{lemma}

\begin{lemma}[A VC class maximal inequality for i.n.i.d.\ empirical processes]
  \label{lem:maximal_vc_inid}

  Assume the same setup as in
  Lemma~\ref{lem:maximal_entropy},
  and suppose further that $\cF$ forms a VC-type class
  in that
  \begin{align*}
    \sup_\Q N(\cF, \rho_\Q, \varepsilon \|F\|_{\Q,2})
    &\leq
    (C_1/\varepsilon)^{C_2}
  \end{align*}
  for all $\varepsilon \in (0,1]$,
  for some constants
  $C_1 \geq e$
  (where $e$ is the standard exponential constant)
  and $C_2 \geq 1$.
  Then for $\delta \in (0,1]$
  we have the covering integral bound
  \begin{align*}
    J\big(\delta, \cF, F\big)
    &\leq
    3 \delta
    \sqrt{C_2 \log (C_1/\delta)},
  \end{align*}
  and so by Lemma~\ref{lem:maximal_entropy},
  \begin{align*}
    \E \left[
      \sup_{f \in \cF}
      \big| G_n(f) \big|
    \right]
    &\lesssim
    \sigma
    \sqrt{C_2 \log (C_1/\delta)}
    +
    \frac{\|M\|_{\P,2} C_2 \log(C_1/\delta)}{\sqrt{n}} \\
    &\lesssim
    \sigma
    \sqrt{C_2 \log \big(C_1 \|F\|_{\bar\P,2}/\sigma\big)}
    +
    \frac{\|M\|_{\P,2} C_2 \log \big(C_1 \|F\|_{\bar\P,2}/\sigma\big)}
    {\sqrt{n}},
  \end{align*}
  where $\lesssim$ is up to a universal constant.

\end{lemma}

\subsection{Strong approximation results}

Next we provide two strong approximation results.
The first is a corollary of the KMT approximation
\citep{komlos1975approximation}
which applies to bounded-variation functions
of i.i.d.\ variables.
The second is an extension of the Yurinskii coupling
\citep{belloni2019conditional}
which applies to Lipschitz functions
of i.n.i.d.\ variables.

\begin{lemma}[A KMT approximation corollary]
  \label{lem:kmt_corollary}

  For $n \geq 1$
  let $X_1, \ldots, X_n$
  be i.i.d.\ real-valued random variables and
  $g_n: \R \times \R \to \R$
  be a function satisfying
  the total variation bound
  $\sup_{x \in \R} \|g_n(\cdot, x)\|_\TV < \infty$.
  Then on some probability space
  there exist independent copies of
  $X_1, \ldots, X_n$,
  denoted
  $X_1', \ldots, X_n'$,
  and a mean-zero Gaussian process $Z_n(x)$
  such that if we define
  the empirical process
  \begin{align*}
    G_n(x)
    = \frac{1}{\sqrt n} \sum_{i=1}^n
    \Big(g_n(X_i',x) - \E\big[g_n(X_i',x)\big]\Big),
  \end{align*}
  then
  for some universal positive constants
  $C_1,$ $C_2$ and $C_3$,
  \begin{align*}
    \P\left(
      \sup_{x \in \R}
      \big|G_n(x) - Z_n(x)\big|
      > \sup_{x \in \R} \|g_n(\cdot, x)\|_\TV
      \, \frac{t + C_1 \log n}{\sqrt n}
    \right)
    \leq C_2 e^{-C_3 t}.
  \end{align*}
  Further, $Z_n$
  has the same covariance structure as $G_n$
  in the sense that for all $x,\, x' \in \R$,
  \begin{align*}
    \E\big[Z_n(x) Z_n(x')\big]
    = \E\big[G_n(x) G_n(x')\big].
  \end{align*}
  By independently sampling from the law of
  $Z_n$ conditional on $X_1', \ldots, X_n'$,
  we can assume that
  $Z_n$ is a function only of $X_1', \ldots, X_n'$
  and some independent random noise.

\end{lemma}

\begin{lemma}[Yurinskii coupling for Lipschitz i.n.i.d.\ empirical processes]
  \label{lem:yurinskii_corollary}

  For $n \geq 2$
  let $X_1, \dots, X_n$
  be independent but not necessarily identically distributed
  (i.n.i.d.) random variables
  taking values in a measurable space $(S, \cS)$
  and let $\cX_n \subseteq \R$
  be a compact interval
  with $\left|\log \Leb(\cX_n)\right| \leq C_1 \log n$
  where $C_1 > 0$ is a constant.
  Let $g_n$ be a measurable function on
  $S \times \cX_n$
  satisfying
  $\sup_{\xi \in S} \sup_{x \in \cX_n} |g_n(\xi, x)| \leq M_n$
  and
  $\sup_{x \in \cX_n} \max_{1 \leq i \leq n} \Var[g_n(X_i, x)] \leq \sigma_n^2$,
  with $\left|\log M_n\right| \leq C_1 \log n$
  and $\left|\log \sigma_n^2\right| \leq C_1 \log n$.
  Suppose that $g_n$ satisfies the following uniform
  Lipschitz condition:
  \begin{align*}
    \sup_{\xi \in S}
    \sup_{x,x' \in \cX_n}
    \left|
    \frac{g_n(\xi, x) - g_n(\xi, x')}
    {x-x'}
    \right|
    \leq
    l_{n,\infty},
  \end{align*}
  and also the following $L^2$
  Lipschitz condition:
  \begin{align*}
    \sup_{x,x' \in \cX_n}
    \E\left[
      \frac{1}{n}
      \sum_{i=1}^n
      \left|
      \frac{g_n(X_i, x) - g_n(X_i, x')}
      {x-x'}
      \right|^2
    \right]^{1/2}
    \leq
    l_{n,2},
  \end{align*}
  where $0 < l_{n,2} \leq l_{n,\infty}$,
  $\left|\log l_{n,2}\right| \leq C_1 \log n$ and
  $\left|\log l_{n,\infty}\right| \leq C_1 \log n$.
  Then for any $t_n > 0$ with
  $\left|\log t_n\right| \leq C_1 \log n$,
  there is a probability space carrying
  independent copies of
  $X_1, \ldots, X_n$
  denoted
  $X_1', \ldots, X_n'$
  and a mean-zero Gaussian process $Z_n(x)$
  such that if we define the empirical process
  \begin{align*}
    G_n(x)
    &=
    \frac{1}{\sqrt n}
    \sum_{i=1}^n
    \Big(
    g_n(X'_i,x)
    - \E[g_n(X'_i,x)]
    \Big),
  \end{align*}
  then
  \begin{align*}
    \P\left(
      \sup_{x \in \cX_n}
      \big|
      G_n(x) - Z_n(x)
      \big|
      > t_n
    \right)
    &\leq
    \frac{
      C_2
      \sigma_n
      \sqrt{\Leb(\cX_n)}
      \sqrt{\log n}
      \sqrt{M_n + \sigma_n\sqrt{\log n}}
    }{n^{1/4} t_n^2}
    \sqrt{
      l_{n,2}
      \sqrt{\log n}
      + \frac{l_{n,\infty}}{\sqrt n}
      \log n}
  \end{align*}
  where $C_2 > 0$ is a constant depending only on $C_1$.
  Further, $Z_n$
  has the same covariance structure as $G_n$
  in the sense that for all $x, x' \in \cX_n$,
  \begin{align*}
    \E\big[Z_n(x) Z_n(x')\big]
    = \E\big[G_n(x) G_n(x')\big].
  \end{align*}

\end{lemma}

\subsection{The Vorob'ev--Berkes--Philipp theorem}

We present a generalization of the Vorob'ev--Berkes--Philipp theorem
\citep{dudley1999uniform}
which allows one to ``glue'' multiple random variables
or stochastic processes onto the same probability space,
while preserving some pairwise distributions.
We begin with some definitions.

\begin{definition}[Tree]
  A \emph{tree} is
  an undirected graph with finitely many vertices
  which is connected and
  contains no cycles
  or self-loops.
\end{definition}

\begin{definition}[Polish Borel probability space]
  A \emph{Polish Borel probability space}
  is a triple $(\cX, \cF, \P)$,
  where $\cX$ is a Polish space
  (a topological space metrizable by a complete separable metric),
  $\cF$ is the Borel $\sigma$-algebra induced on $\cX$ by its topology,
  and $\P$ is a probability measure on $(\cX, \cF)$.
  Important examples of Polish spaces
  include
  $\R^d$ and
  the Skorokhod space
  $\cD[0,1]^d$
  for some $d \geq 1$.
  In particular,
  one can consider vectors of real-valued random variables
  or stochastic processes indexed by
  compact subsets of $\R^d$ which have
  almost surely continuous trajectories.
\end{definition}

\begin{definition}[Projection of a law]
  Let $(\cX_1, \cF_1)$ and $(\cX_2, \cF_2)$
  be measurable spaces, and
  let $\P_{12}$ be a law on the
  product space
  $(\cX_1 \times \cX_2, \cF_1 \otimes \cF_2)$.
  The \emph{projection} of $\P_{12}$
  onto $\cX_1$ is the law
  $\P_1$ defined on $(\cX_1, \cF_1)$
  by $\P_1 = \P_{12} \circ \pi_1^{-1}$
  where $\pi_1(x_1, x_2) = x_1$
  is the first-coordinate projection.
\end{definition}

\begin{lemma}[Vorob'ev--Berkes--Philipp theorem, tree form]
  \label{lem:vbp}

  Let $\cT$ be a tree
  with vertex set
  $\cV = \{1, \ldots, n\}$
  and edge set $\cE$.
  Suppose that attached to each vertex $i$ is a
  Polish Borel probability space
  $(\cX_i, \cF_i, \P_i)$.
  Suppose that attached to each edge $(i,j) \in \cE$
  (where $i<j$ without loss of generality)
  is a law $\P_{ij}$ on
  $(\cX_i \times \cX_j, \cF_i \otimes \cF_j)$.
  Assume that these laws are pairwise-consistent in the sense that
  the projection of $\P_{ij}$ onto
  $\cX_i$ (resp.\ $\cX_j$) is $\P_i$ (resp.\ $\P_j$)
  for each $(i,j) \in \cE$.
  Then there exists a law $\P$ on
  \begin{align*}
    \left(
      \prod_{i=1}^n \cX_i, \
      \bigotimes_{i=1}^n \cF_i
    \right)
  \end{align*}
  such that the projection of $\P$
  onto $\cX_i \times \cX_j$
  is $\P_{ij}$ for each $(i,j) \in \cE$,
  and therefore also
  the projection of $\P$
  onto $\cX_i$
  is $\P_i$ for each $i \in \cV$.

\end{lemma}

\begin{remark}
  The requirement that $\cT$ must contain no cycles
  is necessary in general.
  To see this, consider the Polish Borel probability spaces
  given by
  $\cX_1 = \cX_2 = \cX_3 = \{0,1\}$,
  their respective Borel $\sigma$-algebras,
  and the pairwise-consistent probability measures:
  \begin{align*}
    1/2
    &=
    \P_1(0) = \P_2(0) = \P_3(0) \\
    1/2
    &=
    \P_{12}(0,1) = \P_{12}(1,0) =
    \P_{13}(0,1) = \P_{13}(1,0) =
    \P_{23}(0,1) = \P_{23}(1,0).
  \end{align*}
  That is, each measure $\P_i$
  places equal mass on 0 and 1,
  while $\P_{ij}$
  asserts that each pair of realizations is a.s.\ not equal.
  The graph of these laws forms a triangle,
  which is not a tree.
  Suppose that $(X_1,X_2,X_3)$ has distribution given by $\P$,
  where $X_i \sim \P_i$ and $(X_i,X_j) \sim \P_{ij}$
  for each $i,j$.
  But then by definition of $\P_{ij}$ we have
  $X_1 = 1-X_2 = X_3 = 1-X_1$ a.s.,
  which is a contradiction.

\end{remark}

\begin{remark}

  Two important applications of
  Lemma~\ref{lem:vbp} include
  the embedding of a random vector into a stochastic process
  and the coupling of stochastic processes
  onto the same probability space:
  \begin{enumerate}[label=(\roman*)]

    \item
      Let $X_1$ and $X_2$
      be stochastic processes with
      trajectories in $\cD[0,1]$.
      For $x_1, \ldots, x_n \in [0,1]$
      let $\widetilde X_1 = (X_1(x_1), \ldots, X_1(x_n))$
      be a random vector
      and suppose that $\widetilde X_1'$
      is a copy of $\widetilde X_1$.
      Then there is a law $\P$ on
      $\cD[0,1] \times \R^n \times \cD[0,1]$
      such that restriction of $\P$ to
      $\cD[0,1] \times \R^n$
      is the law of $(X_1, \widetilde X_1)$,
      while the restriction of $\P$ to
      $\R^n \times \cD[0,1]$
      is the law of $(\widetilde X_1',X_2)$.
      In other words,
      we can embed the vector $\widetilde X_1'$
      into a stochastic process $X_1$
      while maintaining the joint distribution of
      $\widetilde X_1'$ and $X_2$.

    \item
      Let $X_1, X_1', \ldots, X_n, X_n'$
      be stochastic processes with
      trajectories in $\cD[0,1]$,
      where $X_i'$ is a copy of $X_i$
      for each $1 \leq i \leq n-1$.
      Suppose that
      $\P\big(\|X_{i+1} - X_i'\| > t)\leq r_i$
      for each $1 \leq i \leq n-1$,
      where $\|\cdot\|$ is a norm on $\cD[0,1]$.
      Then there exist copies of
      $X_1, \ldots, X_n$
      denoted
      $X_1'', \ldots, X_n''$
      satisfying
      $\P\big(\|X_{i+1}'' - X_i''\| > t)\leq r_i$
      for each $1 \leq i \leq n$.
      That is, all of the approximation inequalities
      can be satisfied simultaneously
      on the same probability space.

  \end{enumerate}

\end{remark}

\begin{remark}
  Note that while we discuss trees in
  Lemma~\ref{lem:vbp},
  these refer to the dependency graph of the relevant
  Polish space-valued random variables and do not
  have any direct relation to the
  networks studied throughout the main paper.
\end{remark}

\section{Additional empirical results}%
\label{sec:empirical}

We present some additional empirical results
using the International Monetary Fund's
Direction of Trade Statistics (DOTS) data set,
to complement those given in the main paper.
This data set contains information about
the yearly trade flows among
$n = 207$ economies ($N = 21\,321$ pairs),
and we focus on the years
$1995$, $2000$ and $2005$.

We define the \emph{trade volume} between countries
$i$ and $j$ as the logarithm of the sum of
the trade flow (in billions of US dollars)
from $i$ to $j$ and the trade flow from $j$ to $i$.
In each year several pairs of countries
did not trade directly, yielding trade flows of zero
and hence a trade volume of $-\infty$.
We therefore assume that the distribution of trade volumes
is a mixture of a point mass at $-\infty$ and a
Lebesgue density on $\R$.
The local nature of our estimator means that
observations taking the value of $-\infty$
can simply be removed from the data set.

For counterfactual analysis we use the
gross domestic product (GDP) of each country as a
covariate, using $10\%$-percentiles to group the values into
$10$ different levels for ease of estimation.
This allows for a comparison of the observed distribution of trade at each year
with, for example, the counterfactual distribution of trade
had the GDP distribution remained as it was in $1995$.
As such we can measure how much of the change in trade
distribution is attributable
to a shift in the GDP distribution.

\begin{figure}[H]
  \centering
  \begin{subfigure}{0.32\textwidth}
    \centering
    \includegraphics[scale=0.52]{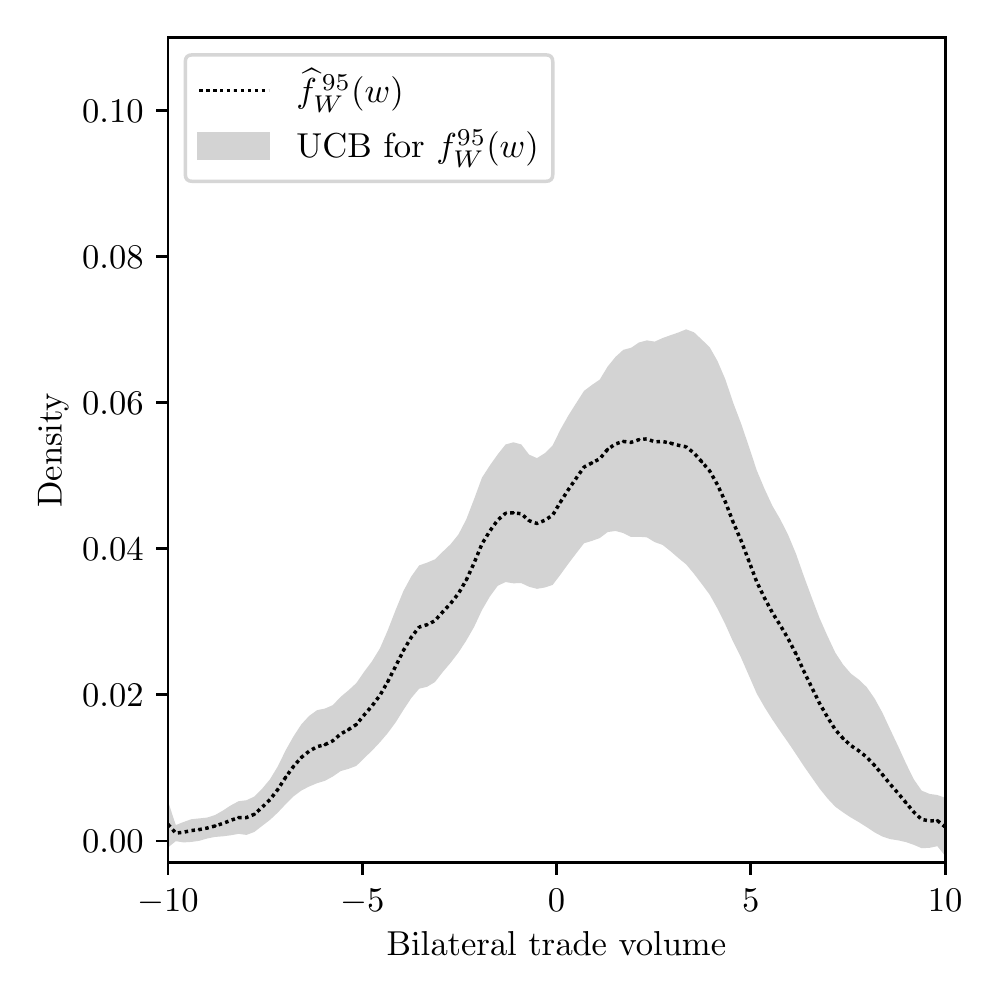}
    \caption{Year 1995, $\widehat h_{\ROT} = 1.27$}
    \label{subfig:trade_1995}
  \end{subfigure}
  \begin{subfigure}{0.32\textwidth}
    \centering
    \includegraphics[scale=0.52]{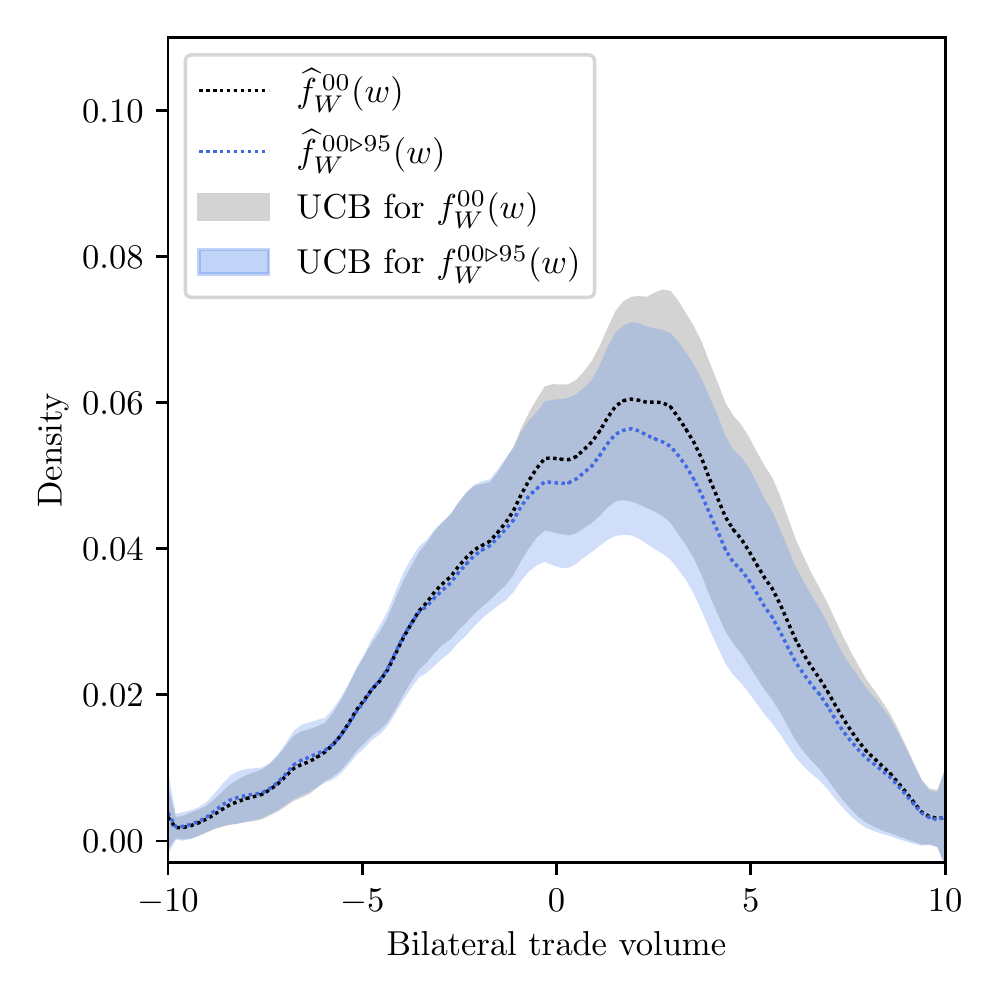}
    \caption{Year 2000, $\widehat h_{\ROT} = 1.31$}
    \label{subfig:trade_2000}
  \end{subfigure}
  \begin{subfigure}{0.32\textwidth}
    \centering
    \includegraphics[scale=0.52]{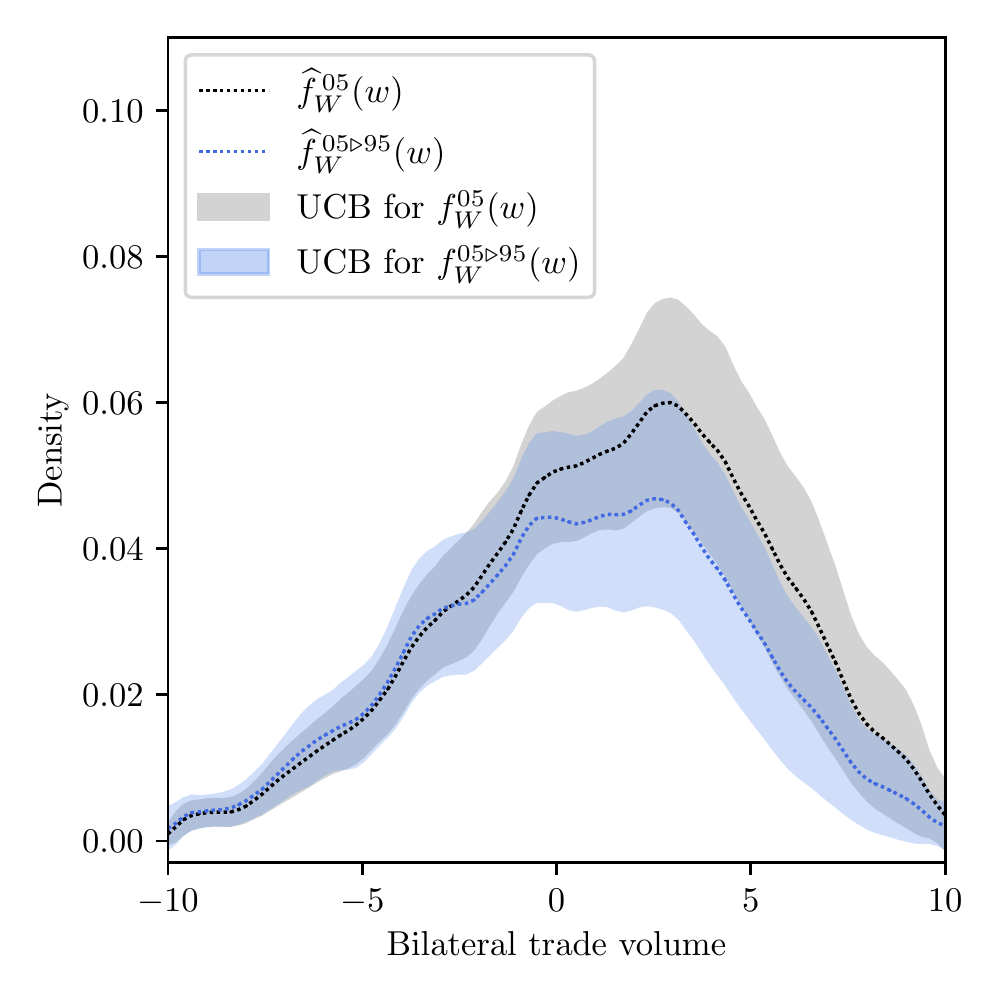}
    \caption{Year 2005, $\widehat h_{\ROT} = 1.37$}
    \label{subfig:trade_2005}
  \end{subfigure}
  \caption{Real and counterfactual density estimates
    and confidence bands for the DOTS data with histogram-based
    covariate estimation.}
  \label{fig:trade}
\end{figure}

To estimate the trade volume density function
we use the counterfactual dyadic kernel density estimator
from Section~\ref{subsec:counterfactual}
with $d=100$ equally-spaced evaluation points in $[-10,10]$,
using the rule-of-thumb bandwidth selector $\widehat h_{\ROT}$
described in the main paper with $p=2$ and $C(K) = 2.435$.
For inference we use an
Epanechnikov kernel of order $p=4$ and
resample the Gaussian process
$B = 10\,000$ times.
We also estimate the counterfactual trade distributions
in 2000 and 2005 respectively,
replacing the GDP distribution with that from 1995.
For each year, Figure~\ref{fig:trade}
plots the real and counterfactual density estimates
along with their respective uniform confidence bands (UCB)
at the nominal coverage rate of $95\%$.
Our empirical results show that
the counterfactual distribution drifts further from the truth
in 2005 compared with 2000, indicating a more significant shift
in the GDP distribution.

\begin{figure}[H]
  \centering
  \begin{subfigure}{0.32\textwidth}
    \centering
    \includegraphics[scale=0.52]{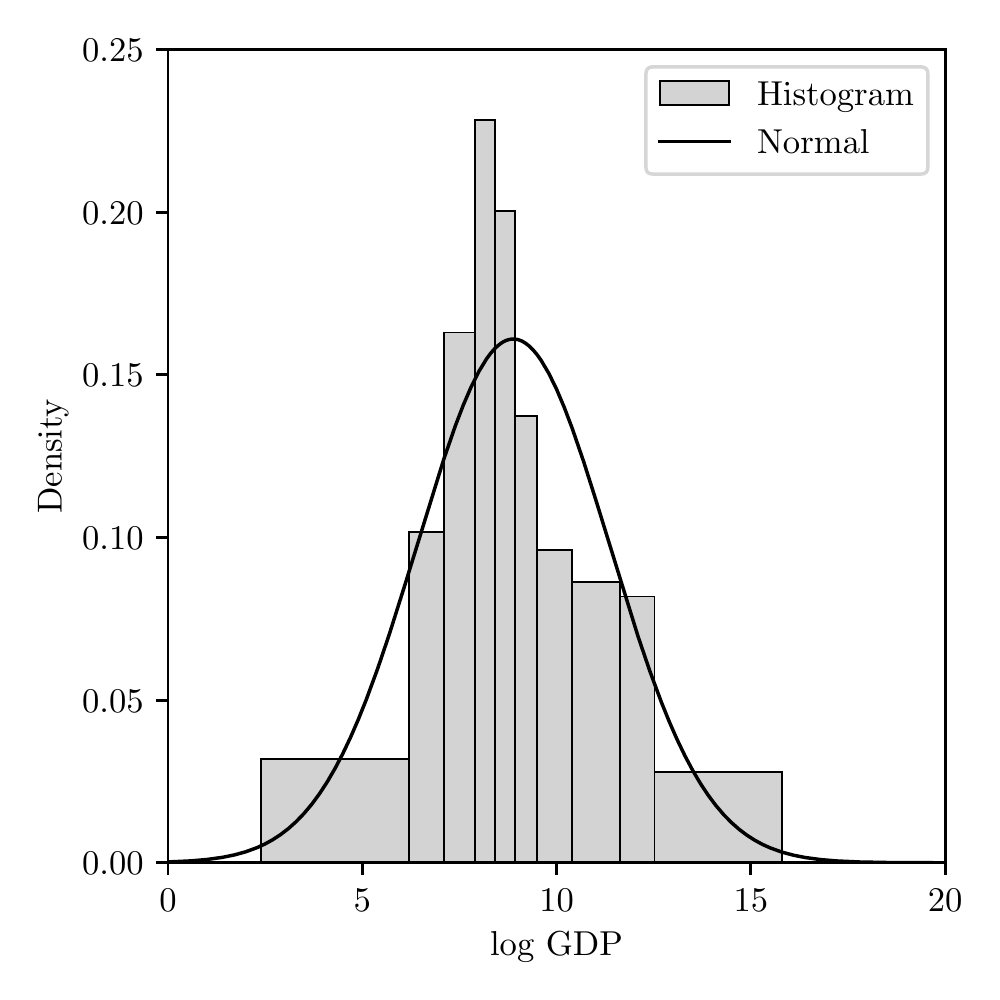}
    \caption{Year 1995}
    \label{subfig:gdp_1995}
  \end{subfigure}
  \begin{subfigure}{0.32\textwidth}
    \centering
    \includegraphics[scale=0.52]{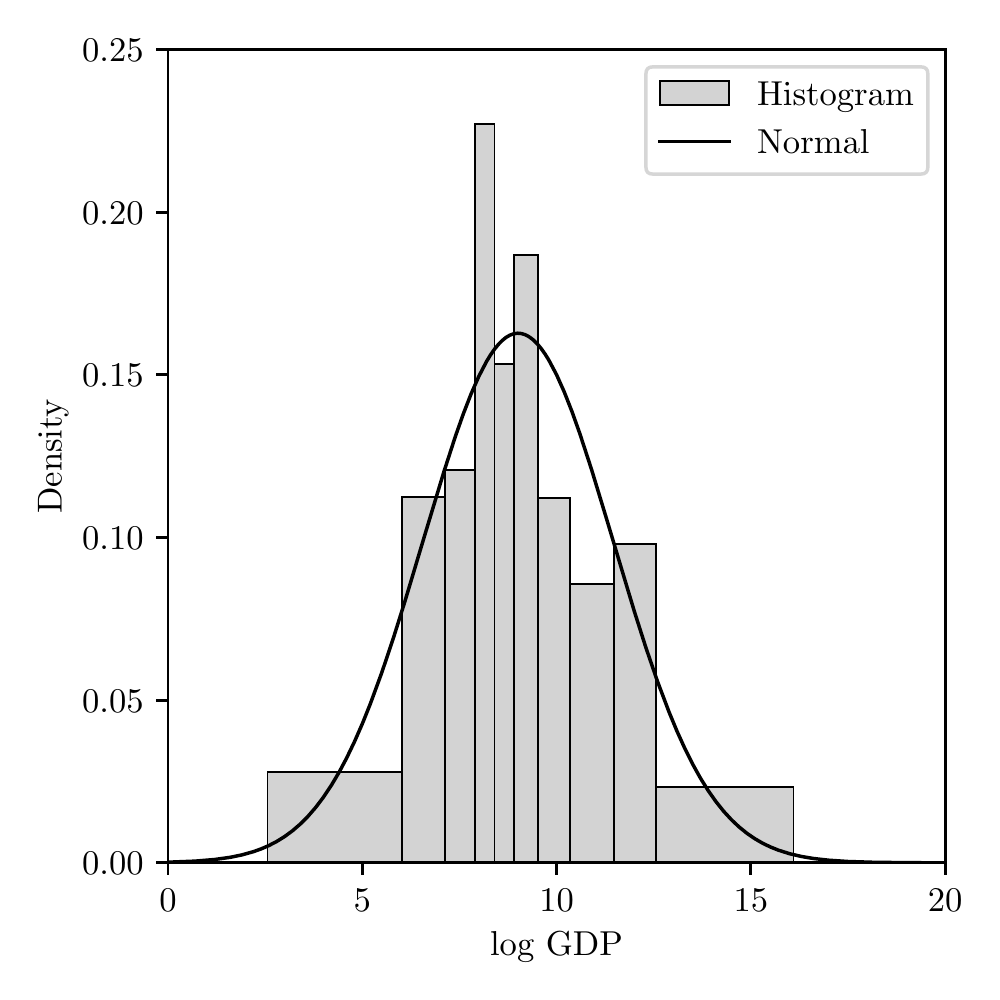}
    \caption{Year 2000}
    \label{subfig:gdp_2000}
  \end{subfigure}
  \begin{subfigure}{0.32\textwidth}
    \centering
    \includegraphics[scale=0.52]{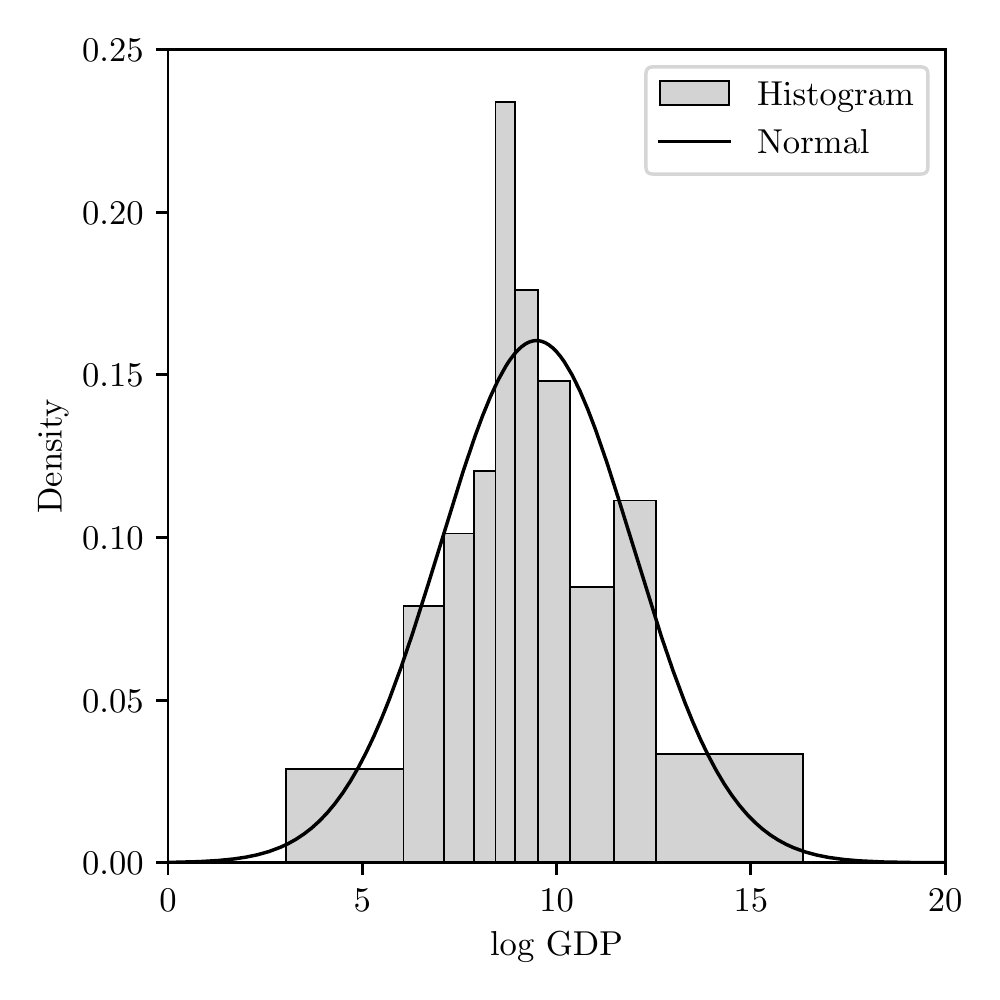}
    \caption{Year 2005}
    \label{subfig:gdp_2005}
  \end{subfigure}
  \caption{Estimated GDP distributions for the DOTS data.}
  \label{fig:gdp}
\end{figure}

In Figure~\ref{fig:gdp} we illustrate how, in the preliminary
step of the counterfactual analysis, the distribution of
log GDP is approximated using the histogram estimators
$\widehat p_X^{\,0}$ and $\widehat p_X^{\,1}$
defined in Section~\ref{subsec:counterfactual}.
We also plot the density function of a normal distribution,
fitted using maximum likelihood estimation,
and this seems to capture the distribution of log GDP reasonably well.
Such a parametric approach to the preliminary step
may be favored in cases where a choice of model
is clear or where the histogram estimators perform poorly.

\begin{figure}[H]
  \centering
  \begin{subfigure}{0.32\textwidth}
    \centering
    \includegraphics[scale=0.52]{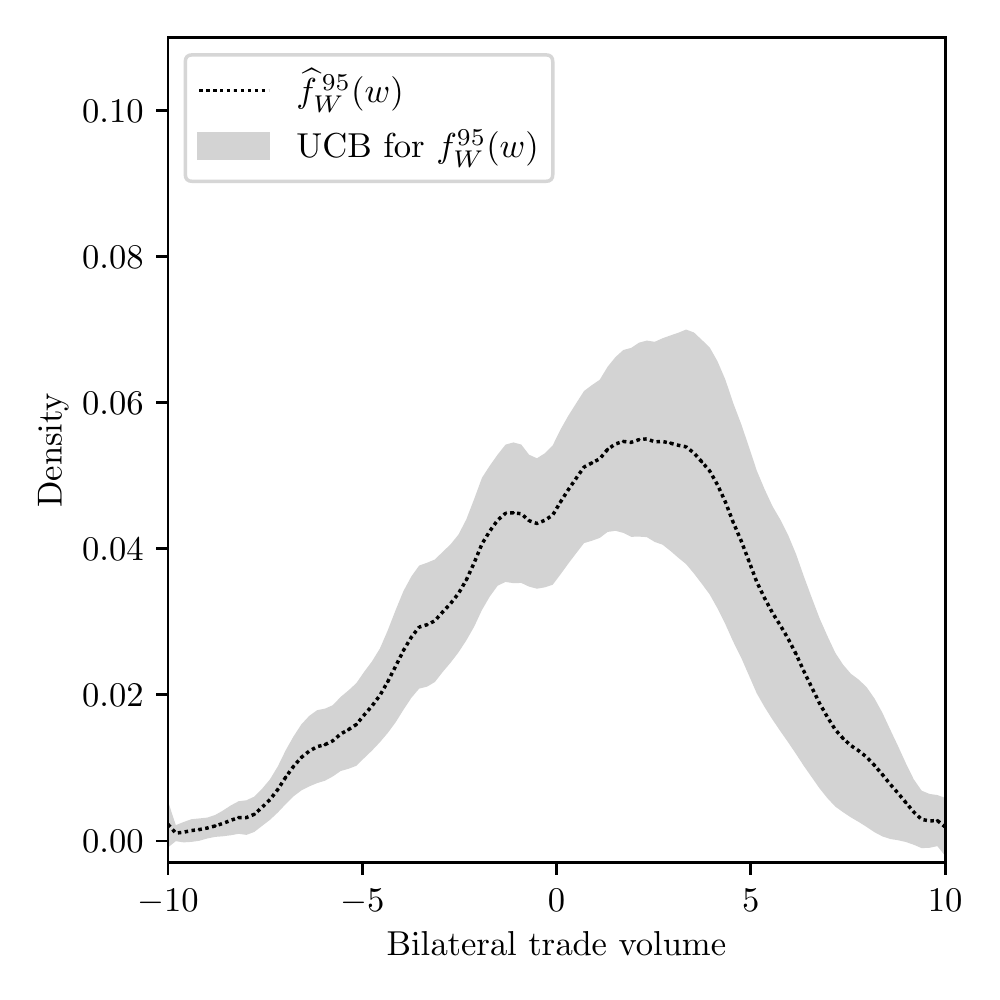}
    \caption{Year 1995, $\widehat h_{\ROT} = 1.27$}
    \label{subfig:trade_para_1995}
  \end{subfigure}
  \begin{subfigure}{0.32\textwidth}
    \centering
    \includegraphics[scale=0.52]{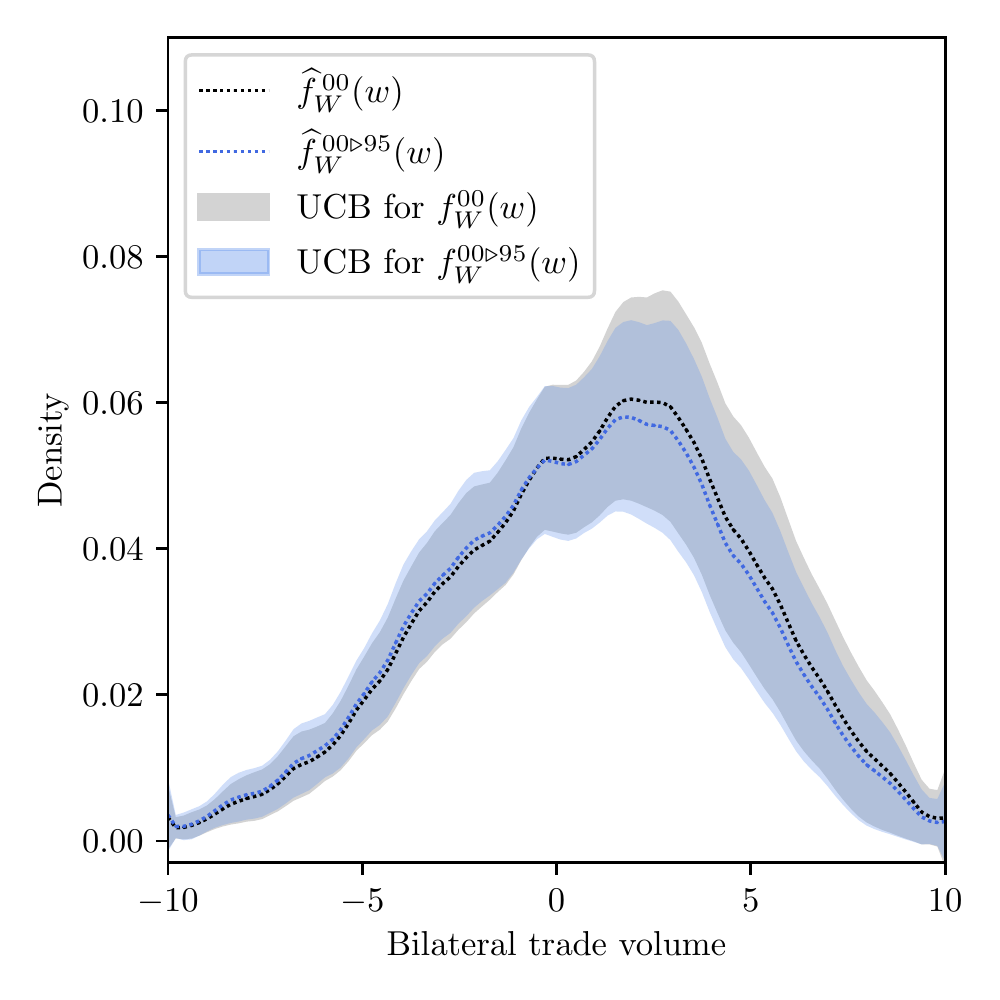}
    \caption{Year 2000, $\widehat h_{\ROT} = 1.31$}
    \label{subfig:trade_para_2000}
  \end{subfigure}
  \begin{subfigure}{0.32\textwidth}
    \centering
    \includegraphics[scale=0.52]{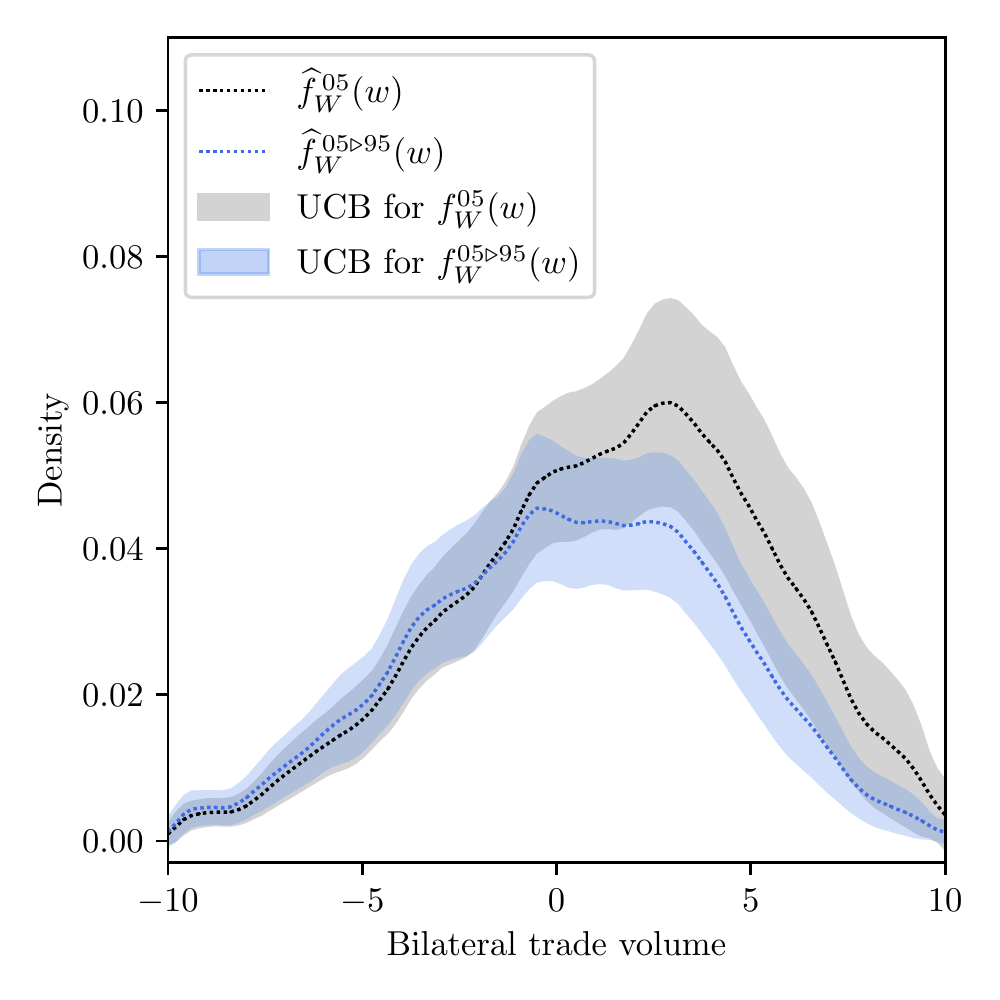}
    \caption{Year 2005, $\widehat h_{\ROT} = 1.37$}
    \label{subfig:trade_para_2005}
  \end{subfigure}
  \caption{Real and counterfactual density estimates
    and confidence bands for the DOTS data
    with parametric covariate estimation.}
  \label{fig:trade_para}
\end{figure}

To demonstrate the relative robustness
of our counterfactual analysis to the choice of preliminary
estimation step, we also provide results using a parametric
estimator of the distribution of GDP.
Figure~\ref{fig:trade_para} repeats the procedure used for
Figure~\ref{fig:trade}, but this time
replacing the histogram estimators by parametric estimators of the log GDP
based on normal likelihood maximization.
The point estimates are qualitatively similar, with the counterfactual
distribution drifting in the same direction over time.
The confidence bands are also similar, with the band based on the
parametric fit being slightly narrower in general.
This could be due to the more stringent model specification
leading to less estimated variance in the fitted values.

\section{Proofs}

\subsection{Preliminary lemmas}

In this section we list some results
in probability and U-statistic theory
which are used in proofs of this paper's main results.
Other auxiliary lemmas will be introduced when
they are needed.

\subsubsection{Standard probabilistic results}

\begin{lemma}[Bernstein's inequality for independent random variables]
  \label{lem:bernstein}

  Let $X_1, \ldots, X_n$ be independent real-valued
  random variables with
  $\E[X_i] = 0$ and $|X_i| \leq M$ and
  $\E[X_i^2] \leq \sigma^2$,
  where $M$ and $\sigma$ are non-random.
  Then for all $t>0$,
  \begin{align*}
    \P \left(
      \left| \frac{1}{n} \sum_{i=1}^n X_i \right| \geq t
    \right)
    \leq 2 \exp \left( -
      \frac{t^2 n}
      {2 \sigma^2 + \frac{2}{3} M t}
    \right).
  \end{align*}

\end{lemma}

\begin{proof}[Lemma~\ref{lem:bernstein}]

  See for example
  Lemma~2.2.9 in~\citet{van1996weak}.
\end{proof}

\begin{lemma}[The matrix Bernstein inequality]
  \label{lem:matrix_bernstein}

  For $1 \leq i \leq n$
  let $X_i$ be independent symmetric $d \times d$
  real random matrices
  with expected values $\mu_i = \E[X_i]$.
  Suppose that
  $\|X_i - \mu_i\|_2 \leq M$
  almost surely for all $1 \leq i \leq n$
  where $M$ is non-random, and define
  $\sigma^2 = \big\| \sum_i \E[(X_i - \mu_i)^2] \big\|_2$.
  Then there exists a universal constant $C > 0$
  such that
  for any $t > 0$ and $q \geq 1$,
  \begin{align*}
    \P\left(
      \left\|
      \sum_{i=1}^n
      \left(
        X_i - \mu_i
      \right)
      \right\|_2
      \geq
      2 \sigma \sqrt{t}
      + \frac{4}{3} M t
    \right)
    &\leq
    2 d e^{-t}, \\
    \E\left[
      \left\|
      \sum_{i=1}^n
      \left(
        X_i - \mu_i
      \right)
      \right\|_2^q
    \right]^{1/q}
    &\leq
    C \sigma \sqrt{q + \log 2d}
    + C M (q + \log 2d).
  \end{align*}
  Another simplified version of this is as follows:
  suppose that
  $\|X_i\|_2 \leq M$ almost surely,
  so that
  $\|X_i - \mu_i\|_2 \leq 2M$.
  Then since
  $\sigma^2 \leq n M^2$,
  we have
  \begin{align*}
    \P\left(
      \left\|
      \sum_{i=1}^n
      \left(
        X_i - \mu_i
      \right)
      \right\|_2
      \geq
      4M \big(t + \sqrt{n t}\big)
    \right)
    &\leq
    2 d e^{-t}, \\
    \E\left[
      \left\|
      \sum_{i=1}^n
      \left(
        X_i - \mu_i
      \right)
      \right\|_2^q
    \right]^{1/q}
    &\leq
    C M
    \big(q + \log 2d + \sqrt{n(q + \log 2d)}\big).
  \end{align*}

\end{lemma}

\begin{proof}[Lemma~\ref{lem:matrix_bernstein}]

  See Lemma~3.2 in \citet{minsker2019moment}.
\end{proof}

\begin{lemma}[A maximal inequality for Gaussian vectors]
  \label{lem:gaussian_vector_maximal}

  Take $n \geq 2$.
  Let $X_i \sim \cN(0, \sigma_i^2)$
  for $1 \leq i \leq n$
  (not necessarily independent),
  with $\sigma_i^2 \leq \sigma^2$.
  Then
  \begin{align}
    \label{eq:gaussian_vector_maximal}
    \E\left[
      \max_{1 \leq i \leq n}
      X_i
    \right]
    &\leq
    \sigma \sqrt{2 \log n}, \\
    \label{eq:gaussian_vector_maximal_abs}
    \E\left[
      \max_{1 \leq i \leq n}
      |X_i|
    \right]
    &\leq
    2 \sigma \sqrt{\log n}.
  \end{align}
  If $\Sigma_1$ and $\Sigma_2$ are constant
  positive semi-definite $n \times n$ matrices
  and $N \sim \cN(0,I_n)$,
  then
  \begin{align}
    \label{eq:gaussian_difference_psd}
    \E\Big[
      \big\|
      \Sigma_1^{1/2} N
      - \Sigma_2^{1/2} N
      \big\|_\infty
      \Big]
    &\leq
    2 \sqrt{\log n} \,
    \big\|
    \Sigma_1 - \Sigma_2
    \big\|_2^{1/2}.
  \end{align}
  If further $\Sigma_1$ is
  positive definite,
  then
  \begin{align}
    \label{eq:gaussian_difference_pd}
    \E\Big[
      \big\|
      \Sigma_1^{1/2} N
      - \Sigma_2^{1/2} N
      \big\|_\infty
      \Big]
    &\leq
    \sqrt{\log n} \,
    \lambda_\mini(\Sigma_1)^{-1/2} \,
    \big\|
    \Sigma_1 - \Sigma_2
    \big\|_2.
  \end{align}

\end{lemma}

\begin{proof}[Lemma~\ref{lem:gaussian_vector_maximal}]

  For $t > 0$,
  Jensen's inequality on the concave logarithm function
  gives
  \begin{align*}
    \E\left[
      \max_{1 \leq i \leq n}
      X_i
    \right]
    &=
    \frac{1}{t}
    \E\left[
      \log
      \exp
      \max_{1 \leq i \leq n}
      t X_i
    \right]
    \leq
    \frac{1}{t}
    \log
    \E\left[
      \exp
      \max_{1 \leq i \leq n}
      t X_i
    \right]
    \leq
    \frac{1}{t}
    \log
    \sum_{i=1}^n
    \E\left[
      \exp
      t X_i
    \right] \\
    &=
    \frac{1}{t}
    \log
    \sum_{i=1}^n
    \exp
    \left(
      \frac{t^2 \sigma_i^2}{2}
    \right)
    \leq
    \frac{1}{t}
    \log n
    + \frac{t \sigma^2}{2},
  \end{align*}
  where we use the Gaussian moment generating function.
  Minimizing this upper bound over $t$
  by setting $t = \sqrt{2 \log n} / \sigma$
  yields Equation~\ref{eq:gaussian_vector_maximal}:
  \begin{align*}
    \E\left[
      \max_{1 \leq i \leq n}
      X_i
    \right]
    &\leq
    \sigma \sqrt{2 \log n}.
  \end{align*}
  For Equation~\ref{eq:gaussian_vector_maximal_abs},
  we use the symmetry of the Gaussian distribution:
  \begin{align*}
    \E\left[
      \max_{1 \leq i \leq n}
      |X_i|
    \right]
    &=
    \E\left[
      \max_{1 \leq i \leq n}
      \{X_i, -X_i\}
    \right]
    \leq
    \sigma \sqrt{2 \log 2n}
    \leq
    2 \sigma \sqrt{\log n}.
  \end{align*}
  For
  Equations~\ref{eq:gaussian_difference_psd}
  and~\ref{eq:gaussian_difference_pd},
  note that
  $\Sigma_1^{1/2} N - \Sigma_2^{1/2} N$
  is a Gaussian vector with covariance matrix
  $\big(\Sigma_1^{1/2} - \Sigma_2^{1/2}\big)^2$.
  The variances of of its components are the diagonal
  elements of this matrix, namely
  \begin{align*}
    \sigma_i^2
    &=
    \Var\big[
      \big(\Sigma_1^{1/2} N - \Sigma_2^{1/2} N\big)_i
      \big]
    =
    \Big(\big(
    \Sigma_1^{1/2} - \Sigma_2^{1/2}
    \big)^2\Big)_{ii}.
  \end{align*}
  Note that if $e_i$ is the
  $i$th standard unit basis vector,
  then for any real symmetric matrix $A$,
  we have
  $e_i^\T A^2 e_i = (A^2)_{ii}$,
  so in particular
  $(A^2)_{ii} \leq \|A\|_2^2$.
  Therefore
  \begin{align*}
    \sigma_i^2
    &\leq
    \big\|
    \Sigma_1^{1/2} - \Sigma_2^{1/2}
    \big\|_2^2
    =\vcentcolon
    \sigma^2.
  \end{align*}
  Applying
  Equation~\ref{eq:gaussian_vector_maximal_abs}
  then gives
  \begin{align*}
    \E\Big[
      \big\|
      \Sigma_1^{1/2} N
      - \Sigma_2^{1/2} N
      \big\|_\infty
      \Big]
    &\leq
    2 \sqrt{\log n} \,
    \big\|
    \Sigma_1^{1/2} - \Sigma_2^{1/2}
    \big\|_2.
  \end{align*}
  By Theorem~X.1.1
  in \citet{bhatia1997matrix},
  we can deduce
  \begin{align*}
    \big\|
    \Sigma_1^{1/2} - \Sigma_2^{1/2}
    \big\|_2
    &\leq
    \big\|
    \Sigma_1 - \Sigma_2
    \big\|_2^{1/2},
  \end{align*}
  giving
  Equation~\ref{eq:gaussian_difference_psd}.
  If further $\Sigma_1$
  is positive definite,
  then by
  Theorem~X.3.8 in
  \citet{bhatia1997matrix},
  \begin{align*}
    \big\|
    \Sigma_1^{1/2} - \Sigma_2^{1/2}
    \big\|_2
    &\leq
    \frac{1}{2}
    \lambda_\mini(\Sigma_1)^{-1/2} \,
    \big\|
    \Sigma_1 - \Sigma_2
    \big\|_2,
  \end{align*}
  giving
  Equation~\ref{eq:gaussian_difference_pd}.
\end{proof}

\begin{lemma}[Maximal inequalities for Gaussian processes]
  \label{lem:gaussian_process_maximal}

  Let $Z$ be a separable
  mean-zero Gaussian process indexed
  by $x \in \cX$.
  Recall that $Z$ is separable for example if
  $\cX$ is Polish
  and $Z$ has
  continuous trajectories.
  Define its covariance structure on $\cX \times \cX$
  by $\Sigma(x, x') = \E[Z(x) Z(x')]$,
  and the corresponding semimetric on $\cX$
  by
  \begin{align*}
    \rho(x,x')
    &=
    \E\big[\big(Z(x) - Z(x')\big)^2\big]^{1/2}
    = \big(\Sigma(x,x)
    - 2 \Sigma(x,x')
    + \Sigma(x',x')\big)^{1/2}.
  \end{align*}
  Let $N(\varepsilon, \cX, \rho)$
  denote the $\varepsilon$-covering number of $\cX$
  with respect to the semimetric $\rho$.
  Define $\sigma = \sup_x \Sigma(x,x)^{1/2}$.
  Then there exists a universal constant $C > 0$
  such that for any $\delta > 0$,
  \begin{align*}
    \E\left[
      \sup_{x \in \cX}
      |Z(x)|
    \right]
    &\leq
    C \sigma
    + C \int_0^{2\sigma}
    \sqrt{\log N(\varepsilon, \cX, \rho)}
    \diff{\varepsilon}, \\
    \E\left[
      \sup_{\rho(x,x') \leq \delta}
      |Z(x) - Z(x')|
    \right]
    &\leq
    C \int_0^{\delta}
    \sqrt{\log N(\varepsilon, \cX, \rho)}
    \diff{\varepsilon}.
  \end{align*}

\end{lemma}

\begin{proof}[Lemma~\ref{lem:gaussian_process_maximal}]

  See Corollary~2.2.8 in \citet{van1996weak},
  noting that for any $x,x' \in \cX$, we have
  $\E[|Z(x)|] \lesssim \sigma$ and
  $\rho(x,x') \leq 2\sigma$,
  implying that
  $\log N(\varepsilon, \cX, \rho) = 0$
  for all
  $\varepsilon > 2 \sigma$.
\end{proof}

\begin{lemma}[Anti-concentration for Gaussian process absolute suprema]
  \label{lem:anti-concentration}

  Let $Z$ be a separable mean-zero Gaussian process
  indexed by a semimetric space $\cX$
  satisfying
  $\E[Z(x)^2] = 1$
  for all $x \in \cX$.
  Then for any $\varepsilon > 0$,
  \begin{align*}
    \sup_{t \in \R}
    \P\left(
      \left|
      \sup_{x \in \cX}
      \big| Z(x) \big|
      - t
      \right|
      \leq \varepsilon
    \right)
    &\leq
    4 \varepsilon
    \left(
      1 + \E\left[
        \sup_{x \in \cX}
        \big| Z(x) \big|
      \right]
    \right).
  \end{align*}

\end{lemma}

\begin{proof}[Lemma~\ref{lem:anti-concentration}]

  See Corollary~2.1
  in \citet{chernozhukov2014anti}.
\end{proof}

\begin{lemma}[No slowest rate of convergence in probability]
  \label{lem:slow_convergence}

  Let $X_n$ be a sequence of real-valued random
  variables with
  $X_n = o_\P(1)$.
  Then there exists a deterministic sequence
  $\varepsilon_n \to 0$
  such that
  $\P\big(|X_n| > \varepsilon_n\big) \leq \varepsilon_n$
  for all $n \geq 1$.

\end{lemma}

\begin{proof}[Lemma~\ref{lem:slow_convergence}]

  Define the following deterministic sequence
  for $k \geq 1$.
  \begin{align*}
    \tau_k
    &=
    \sup
    \big\{
    n \geq 1:
    \P\big(|X_n| > 1/k\big)
    > 1/k
    \big\}
    \vee
    (\tau_{k-1} +1)
  \end{align*}
  with $\tau_0 = 0$.
  Since $X_n = o_\P(1)$,
  each $\tau_k$ is finite
  and so we can define
  $\varepsilon_n = \frac{1}{k}$
  where $\tau_k < n \leq \tau_{k+1}$.
  Then, noting that $\varepsilon_n \to 0$,
  we have
  $\P\big(|X_n| > \varepsilon_n\big)
  = \P\big(|X_n| > 1/k\big) \leq 1/k = \varepsilon_n$.
\end{proof}

\subsubsection{U-statistics}

\begin{lemma}[General Hoeffding-type decomposition]
  \label{lem:general_hoeffding}

  Let $\cU$ be a vector space.
  Let $u_{ij} \in \cU$ be defined for
  $1 \leq i, j \leq n$
  and
  $i \neq j$.
  Suppose that $u_{ij} = u_{ji}$
  for all $i,j$.
  Then for any $u_i \in \cU$
  (for $1 \leq i \leq n$)
  and any $u \in \cU$,
  the following decomposition holds:
  \begin{align*}
    \sum_{i=1}^n
    \sum_{\substack{j=1 \\ j \neq i}}^n
    \big(u_{ij} - u\big)
    &=
    2(n-1)
    \sum_{i=1}^n
    \big(u_i - u\big)
    +
    \sum_{i=1}^n
    \sum_{\substack{j=1 \\ j \neq i}}^n
    \big(u_{ij} - u_i - u_j + u\big).
  \end{align*}

\end{lemma}

\begin{proof}[Lemma~\ref{lem:general_hoeffding}]

  We compute the left hand side minus the right hand side,
  beginning by observing that all of the
  $u_{ij}$ and $u$ terms clearly cancel.
  \begin{align*}
    &\sum_{i=1}^n
    \sum_{\substack{j=1 \\ j \neq i}}^n
    \big(u_{ij} - u\big)
    - 2(n-1)
    \sum_{i=1}^n
    \big(u_i - u\big)
    -
    \sum_{i=1}^n
    \sum_{j \neq i}
    \big(u_{ij} - u_i - u_j + u\big) \\
    &\qquad=
    - 2(n-1)
    \sum_{i=1}^n
    u_i
    -
    \sum_{i=1}^n
    \sum_{\substack{j=1 \\ j \neq i}}^n
    \big(- u_i - u_j\big)
    =
    - 2(n-1)
    \sum_{i=1}^n
    u_i
    +
    \sum_{i=1}^n
    \sum_{\substack{j=1 \\ j \neq i}}^n
    u_i
    +
    \sum_{j=1}^n
    \sum_{\substack{i=1 \\ i \neq j}}^n
    u_j \\
    &\qquad=
    - 2(n-1)
    \sum_{i=1}^n
    u_i
    +
    (n-1)
    \sum_{i=1}^n
    u_i
    +
    (n-1)
    \sum_{j=1}^n
    u_j
    = 0.
  \end{align*}
\end{proof}

\begin{lemma}[A U-statistic concentration inequality]
  \label{lem:ustat_concentration}

  Let $(S,\cS)$ be a measurable space and
  $X_1, \ldots, X_n$ be i.i.d.\ $S$-valued random variables.
  Let $H: S^m \to \R$ be a function of $m$ variables
  satisfying the symmetry property
  $H(x_1, \ldots, x_m) = H(x_{\tau (1)}, \ldots, x_{\tau (m)})$
  for any $m$-permutation $\tau$.
  Suppose also that
  $\E[H(X_1, \ldots, X_m)] = 0$.
  Let
  $M = \|H\|_\infty$
  and
  $\sigma^2 = \E\big[\E[H(X_1, \ldots, X_m) \mid X_1]^2\big]$.
  Define the (not necessarily degenerate)
  U-statistic
  \begin{align*}
    U_n
    &=
    \frac{m!(n-m)!}{n!}
    \sum_{1 \leq i_1 < \cdots < i_m \leq n}
    H(X_1, \ldots, X_n).
  \end{align*}
  Then for any $t > 0$,
  \begin{align*}
    \P\left(
      |U_n| > t
    \right)
    &\leq
    4 \exp \left(
      - \frac{n t^2}{C_1(m) \sigma^2 + C_2(m) M t}
    \right),
  \end{align*}
  where $C_1(m)$, $C_2(m)$
  are positive constants depending only on $m$.

\end{lemma}

\begin{proof}[Lemma~\ref{lem:ustat_concentration}]

  See Theorem~2 in \citet{arcones1995bernstein}.
\end{proof}

\begin{lemma}[A second-order U-process maximal inequality]
  \label{lem:uprocess_maximal}

  Let $X_1 \ldots, X_n$
  be i.i.d.\ random variables taking values
  in a measurable space $(S, \cS)$
  with distribution $\P$.
  Let $\cF$ be a class of measurable functions from
  $S \times S$ to $\R$ which is also pointwise measurable.
  Define the degenerate second-order U-process
  \begin{align*}
    U_n(f)
    &=
    \frac{2}{n(n-1)}
    \sum_{i<j}
    \Big(
    f(X_i, X_j)
    - \E\big[f(X_i,X_j) \mid X_i\big]
    - \E\big[f(X_i,X_j) \mid X_j\big]
    + \E\big[f(X_i,X_j)\big]
    \Big)
  \end{align*}
  for $f \in \cF$.
  Suppose that each $f \in \cF$ is symmetric in the sense that
  $f(s_1,s_2) = f(s_2,s_1)$
  for all $s_1, s_2 \in S$.
  Let $F$ be a measurable envelope function for $\cF$
  satisfying $|f(s_1,s_2)| \leq F(s_1,s_2)$
  for all $s_1,s_2 \in S$.
  For a law $\Q$ on
  $(S \times S, \, \cS \otimes \cS)$,
  define the $(\Q,q)$-norm of $f \in \cF$ by
  $\|f\|_{\Q,q}^q = \E_\Q[|f|^q]$.
  Assume that $\cF$ is VC-type in the following manner.
  \begin{align*}
    \sup_\Q
    N\big(
    \cF, \|\cdot\|_{\Q,2}, \varepsilon \|F\|_{\Q,2}
    \big)
    &\leq
    (C_1/\varepsilon)^{C_2}
  \end{align*}
  for some constants
  $C_1 \geq e$
  and
  $C_2 \geq 1$,
  and for all $\varepsilon \in (0,1]$,
  where $\Q$ ranges over all finite discrete laws
  on
  $S \times S$.
  Let $\sigma > 0$ be any deterministic value satisfying
  $\sup_{f \in \cF} \|f\|_{\P,2} \leq \sigma \leq \|F\|_{\P,2}$,
  and define the random variable $M = \max_{i,j} |F(X_i, X_j)|$.
  Then there exists a universal constant $C_3 > 0$
  satisfying
  \begin{align*}
    n
    \E\left[
      \sup_{f \in \cF}
      \big| U_n(f) \big|
    \right]
    &\leq
    C_3 \sigma
    \Big(
    C_2 \log\big(C_1 \|F\|_{\P,2} / \sigma \big)
    \Big)
    + \frac{C_3 \|M\|_{\P,2}}{\sqrt{n}}
    \Big(
    C_2 \log\big(C_1 \|F\|_{\P,2} / \sigma \big)
    \Big)^2.
  \end{align*}

\end{lemma}

\begin{proof}[Lemma~\ref{lem:uprocess_maximal}]

  Apply Corollary~5.3
  from \citet{chen2020jackknife}
  with the order of the U-statistic fixed at
  $r=2$,
  and with $k=2$.
\end{proof}

\begin{lemma}[A U-statistic matrix concentration inequality]
  \label{lem:ustat_matrix_concentration}

  Let $X_1, \ldots, X_n$ be i.i.d.\ random variables
  taking values in a measurable space $(S, \cS)$.
  Suppose
  $H: S^2 \to \R^{d \times d}$
  is a measurable matrix-valued function
  of two variables
  satisfying the following assumptions:
  \begin{enumerate}[label=(\roman*)]

    \item
      $H(X_1, X_2)$ is an almost surely symmetric matrix.

    \item
      $\|H(X_1, X_2)\|_2 \leq M$ almost surely.

    \item
      $H$ is a symmetric function in its arguments in that
      $H(X_1, X_2) = H(X_2, X_1)$ almost surely.

    \item
      $H$ is degenerate in the sense that
      $\E[H(X_1, x_2)] = 0$ for all $x_2 \in S$.

  \end{enumerate}
  Let $U_n = \sum_i \sum_{j \neq i} H(X_i, X_j)$
  be a U-statistic,
  and define the variance-type constant
  \begin{align*}
    \sigma^2
    &=
    \E\left[
      \left\|
      \E\left[
        H(X_i, X_j)^2
        \mid X_j
      \right]
      \right\|_2
    \right].
  \end{align*}
  Then for a universal constant $C > 0$
  and for all $t > 0$,
  \begin{align*}
    \P\left(
      \|U_n\|_2
      \geq
      C \sigma n (t + \log d)
      + C M \sqrt{n} (t + \log d)^{3/2}
    \right)
    &\leq
    C e^{-t}.
  \end{align*}
  We remark here that clearly
  by Jensen's inequality,
  $\sigma^2 \leq \E[ \| H(X_i, X_j)^2 \|_2 ]
  = \E[ \| H(X_i, X_j) \|_2^2 ]
  \leq M^2$,
  giving the weaker but simpler concentration inequality
  \begin{align*}
    \P\left(
      \|U_n\|_2
      \geq
      2 C M n
      (t + \log d)^{3/2}
    \right)
    &\leq
    C e^{-t}.
  \end{align*}
  From this last inequality we can deduce the following moment bound
  by integration of tail probabilities.
  \begin{align*}
    \E\left[
      \|U_n\|_2
    \right]
    &\lesssim
    M n (\log d)^{3/2}.
  \end{align*}

\end{lemma}

\begin{proof}[Lemma~\ref{lem:ustat_matrix_concentration}]

  We apply results from \citet{minsker2019moment}.

  \proofparagraph{decoupling}

  Let $\bar U_n = \sum_{i=1}^n \sum_{j=1}^n H(X_i^{(1)}, X_j^{(2)})$
  be a decoupled matrix U-statistic,
  where $X^{(1)}$ and $X^{(2)}$
  are i.i.d.\ copies of the sequence $X_1, \ldots, X_n$.
  By Lemma~5.2 in \citet{minsker2019moment},
  since we are only stating this result for
  degenerate U-statistics of order 2,
  there exists a universal constant $D_2$
  such that for any $t > 0$,
  we have
  \begin{align*}
    \P\left(
      \|U_n\|_2 \geq t
    \right)
    &\leq
    D_2
    \P\left(
      \|\bar U_n\|_2 \geq t / D_2
    \right).
  \end{align*}

  \proofparagraph{concentration of the decoupled U-statistic}

  By Equation~11 in \citet{minsker2019moment},
  we have the following concentration inequality
  for decoupled degenerate U-statistics.
  For some universal constant $C_1$
  and for any $t > 0$,
  \begin{align*}
    \P\left(
      \|\bar U_n\|_2
      \geq
      C_1 \sigma n (t + \log d)
      + C_1 M \sqrt{n} (t + \log d)^{3/2}
    \right)
    &\leq
    e^{-t}.
  \end{align*}

  \proofparagraph{concentration of the original U-statistic}

  Hence we have
  \begin{align*}
    &\P\left(
      \|U_n\|_2
      \geq
      C_1 D_2 \sigma n (t + \log d)
      + C_1 D_2 M \sqrt{n} (t + \log d)^{3/2}
    \right) \\
    &\quad\leq
    D_2 \P\left(
      \|\bar U_n\|_2
      \geq
      C_1 \sigma n (t + \log d)
      + C_1 M \sqrt{n} (t + \log d)^{3/2}
    \right) \\
    &\quad\leq
    D_2 e^{-t}.
  \end{align*}
  The main result follows by setting
  $C = C_1 + C_1 D_2$.

  \proofparagraph{moment bound}

  The final equation,
  giving a moment bound for the simplified version,
  can be deduced as follows.
  We already have that
  \begin{align*}
    \P\left(
      \|U_n\|_2
      \geq
      2 C M n
      (t + \log d)^{3/2}
    \right)
    &\leq
    C e^{-t}.
  \end{align*}
  This implies that for any $t \geq \log d$,
  we have
  \begin{align*}
    \P\left(
      \|U_n\|_2
      \geq
      8 C M n
      t^{3/2}
    \right)
    &\leq
    C e^{-t}.
  \end{align*}
  Defining
  $s = 8 C M n t^{3/2}$,
  or equivalently
  $t = \left( \frac{s}{8C M n} \right)^{2/3}$,
  shows that
  for any $s \geq 8C M n(\log d)^{3/2}$,
  \begin{align*}
    \P\left(
      \|U_n\|_2
      \geq
      s
    \right)
    &\leq
    C e^{-\left( \frac{s}{8C M n} \right)^{2/3}}.
  \end{align*}
  Hence the moment bound is obtained:
  \begin{align*}
    \E\left[
      \|U_n\|_2
    \right]
    &=
    \int_0^\infty
    \P\left(
      \|U_n\|_2
      \geq
      s
    \right)
    \diff{s}
    = \int_0^{8C M n(\log d)^{3/2}}
    \P\left(
      \|U_n\|_2
      \geq
      s
    \right)
    \diff{s}
    +
    \int_{8C M n(\log d)^{3/2}}^\infty
    \P\left(
      \|U_n\|_2
      \geq
      s
    \right)
    \diff{s} \\
    &\leq
    8C M n(\log d)^{3/2}
    +
    \int_0^\infty
    C e^{-\left( \frac{s}{8C M n} \right)^{2/3}}
    \diff{s}
    =
    8C M n(\log d)^{3/2}
    +
    8C M n
    \int_0^\infty
    e^{s^{-2/3}}
    \diff{s} \\
    &\lesssim
    Mn(\log d)^{3/2}.
  \end{align*}
\end{proof}

\subsection{Technical results}

\subsubsection{Maximal inequalities for i.n.i.d.\ empirical processes}

Before presenting the proof of
Lemma~\ref{lem:maximal_entropy},
we give some auxiliary lemmas;
namely a symmetrization inequality
(Lemma~\ref{lem:symmetrization}),
a Rademacher contraction principle
(Lemma~\ref{lem:contraction})
and a Hoffman--J{\o}rgensen inequality
(Lemma~\ref{lem:hoffmann}).
Recall that the Rademacher distribution
places probability mass of $1/2$
on each of the points $-1$ and $1$.

\begin{lemma}[A symmetrization inequality for i.n.i.d.\ variables]
  \label{lem:symmetrization}

  Let $(S, \cS)$ be a measurable space and
  $\cF$ a class of Borel-measurable functions
  from $S$ to $\R$ which is pointwise measurable
  (i.e.\ it contains a countable dense subset
  under pointwise convergence).
  Let $X_1, \ldots, X_n$
  be independent
  but not necessarily identically distributed
  $S$-valued random variables.
  Let $a_1, \ldots, a_n$ be arbitrary points in $S$
  and $\phi$ a non-negative non-decreasing convex function
  from $\R$ to $\R$.
  Define $\varepsilon_1, \ldots, \varepsilon_n$
  as independent Rademacher
  random variables,
  independent of $X_1, \ldots, X_n$.
  Then
  \begin{align*}
    \E \left[
      \phi \left(
        \sup_{f \in \cF}
        \left|
        \sum_{i=1}^n
        \Big(
        f(X_i)
        - \E[f(X_i)]
        \Big)
        \right|
      \right)
    \right]
    &\leq
    \E \left[
      \phi \left(
        2
        \sup_{f \in \cF}
        \left|
        \sum_{i=1}^n
        \varepsilon_i
        \Big(
        f(X_i)
        - a_i
        \Big)
        \right|
      \right)
    \right].
  \end{align*}
  Note that in particular this holds with $a_i = 0$
  and also holds with $\phi(t) = t \vee 0$.

\end{lemma}

\begin{proof}[Lemma~\ref{lem:symmetrization}]

  See Lemma~2.3.6 in
  \citet{van1996weak}.
\end{proof}

\begin{lemma}[A Rademacher contraction principle]
  \label{lem:contraction}

  Let $\varepsilon_1, \ldots, \varepsilon_n$
  be independent Rademacher random variables
  and $\cT$ be a bounded subset of $\R^n$.
  Define
  $M = \sup_{t \in \cT} \max_{1 \leq i \leq n} |t_i|$.
  Then, noting that the supremum is measurable
  because $\cT$ is a subset of a separable metric space
  and is therefore itself separable,
  \begin{align*}
    \E
    \left[
      \sup_{t \in \cT}
      \left|
      \sum_{i=1}^n
      \varepsilon_i
      t_i^2
      \right|
    \right]
    &\leq
    4M \,
    \E
    \left[
      \sup_{t \in \cT}
      \left|
      \sum_{i=1}^n
      \varepsilon_i
      t_i
      \right|
    \right].
  \end{align*}
  This gives the following corollary.
  Let $X_1, \ldots, X_n$ be mutually independent
  and also independent of $\varepsilon_1, \ldots, \varepsilon_n$.
  Let $\cF$ be a pointwise measurable class of functions
  from a measurable space $(S, \cS)$ to $\R$,
  with measurable envelope $F$.
  Define $M = \max_i F(X_i)$.
  Then we obtain that
  \begin{align*}
    \E
    \left[
      \sup_{f \in \cF}
      \left|
      \sum_{i=1}^n
      \varepsilon_i
      f(X_i)^2
      \right|
    \right]
    &\leq
    4
    \E
    \left[
      M
      \sup_{f \in \cF}
      \left|
      \sum_{i=1}^n
      \varepsilon_i
      f(X_i)
      \right|
    \right].
  \end{align*}

\end{lemma}

\begin{proof}[Lemma~\ref{lem:contraction}]

  We apply Theorem~4.12
  from \citet{ledoux1991probability}
  with $F$ as the identity function
  and
  \begin{align*}
    \psi_i(s)
    = \psi(s)
    &=
    \min
    \left(
      \frac{s^2}{2M},
      \frac{M}{2}
    \right).
  \end{align*}
  This is a weak contraction
  (i.e.\ 1-Lipschitz)
  because it is continuous,
  differentiable on $(-M,M)$
  with derivative bounded by
  $|\psi'(s)| \leq |s|/M \leq 1$,
  and constant outside $(-M,M)$.
  Note that since $|t_i| \leq M$
  by definition,
  we have $\psi_i(t_i) = t_i^2 / (2M)$.
  Hence
  by Theorem~4.12
  from \citet{ledoux1991probability},
  \begin{align*}
    \E
    \left[
      F
      \left(
        \frac{1}{2}
        \sup_{t \in \cT}
        \left|
        \sum_{i=1}^n
        \varepsilon_i
        \psi_i(t_i)
        \right|
      \right)
    \right]
    &\leq
    \E
    \left[
      F
      \left(
        \sup_{t \in \cT}
        \left|
        \sum_{i=1}^n
        \varepsilon_i
        t_i
        \right|
      \right)
    \right],
    &\E
    \left[
      \frac{1}{2}
      \sup_{t \in \cT}
      \left|
      \sum_{i=1}^n
      \varepsilon_i
      \frac{t_i^2}{2M}
      \right|
    \right]
    &\leq
    \E
    \left[
      \sup_{t \in \cT}
      \left|
      \sum_{i=1}^n
      \varepsilon_i
      t_i
      \right|
    \right], \\
    \E
    \left[
      \sup_{t \in \cT}
      \left|
      \sum_{i=1}^n
      \varepsilon_i
      t_i^2
      \right|
    \right]
    &\leq
    4M \,
    \E
    \left[
      \sup_{t \in \cT}
      \left|
      \sum_{i=1}^n
      \varepsilon_i
      t_i
      \right|
    \right].
  \end{align*}
  To see the corollary,
  set
  $\cT = \left\{\big(f(X_1), \ldots, f(X_n)\big) : f \in \cF\right\}$
  and note that for fixed realization
  $X_1, \ldots, X_n$,
  \begin{align*}
    \E_\varepsilon
    \left[
      \sup_{f \in \cF}
      \left|
      \sum_{i=1}^n
      \varepsilon_i
      f(X_i)^2
      \right|
    \right]
    &=
    \E_\varepsilon
    \left[
      \sup_{t \in \cT}
      \left|
      \sum_{i=1}^n
      \varepsilon_i
      t_i^2
      \right|
    \right]
    \leq 4
    \E_\varepsilon
    \left[
      M
      \sup_{t \in \cT}
      \left|
      \sum_{i=1}^n
      \varepsilon_i
      t_i
      \right|
    \right]
    = 4 \E_\varepsilon
    \left[
      M
      \sup_{f \in \cF}
      \left|
      \sum_{i=1}^n
      \varepsilon_i
      f(X_i)
      \right|
    \right].
  \end{align*}
  Taking an expectation over $X_1, \ldots, X_n$
  and applying Fubini's theorem yields the result.
\end{proof}

\begin{lemma}[A Hoffmann--J{\o}rgensen inequality]
  \label{lem:hoffmann}

  Let $(S, \cS)$ be a measurable space
  and $X_1, \ldots, X_n$
  be $S$-valued random variables.
  Suppose that
  $\cF$ is a pointwise measurable class of functions from $S$ to $\R$
  with finite envelope $F$.
  Let $\varepsilon_1, \ldots, \varepsilon_n$
  be independent Rademacher random variables
  which are independent of $X_1, \ldots, X_n$.
  Then for any $q \in (1, \infty)$,
  \begin{align*}
    \E \left[
      \sup_{f \in \cF}
      \left|
      \sum_{i=1}^n
      \varepsilon_i
      f(X_i)
      \right|
      ^q
    \right]
    ^{1/q}
    &\leq
    C_q
    \left(
      \E \left[
        \sup_{f \in \cF}
        \left|
        \sum_{i=1}^n
        \varepsilon_i
        f(X_i)
        \right|
      \right]
      +
      \E \left[
        \max_{1 \leq i \leq n}
        \sup_{f \in \cF}
        \big| f(X_i) \big|^q
      \right]^{1/q}
    \right),
  \end{align*}
  where $C_q$ is a positive constant depending only on $q$.

\end{lemma}

\begin{proof}[Lemma~\ref{lem:hoffmann}]

  We use Talagrand's formulation of
  a Hoffmann--J{\o}rgensen inequality.
  Consider the
  independent
  $\ell^\infty(\cF)$-valued
  random functionals $u_i$ defined by
  $u_i(f) = \varepsilon_i f(X_i)$,
  where $\ell^\infty(\cF)$
  is the Banach space of bounded functions from
  $\cF$ to $\R$,
  equipped with the norm
  $\|u\|_\cF = \sup_{f \in \cF} |u(f)|$.
  Then Remark~3.4 in \citet{kwapien1991hypercontraction} gives
  \begin{align*}
    \E \left[
      \sup_{f \in \cF}
      \left|
      \sum_{i=1}^n
      u_i(f)
      \right|
      ^q
    \right]
    ^{1/q}
    &\leq
    C_q
    \left(
      \E \left[
        \sup_{f \in \cF}
        \left|
        \sum_{i=1}^n
        u_i(f)
        \right|
      \right]
      +
      \E \left[
        \max_{1 \leq i \leq n}
        \sup_{f \in \cF}
        \left|
        u_i(f)
        \right|^q
      \right]^{1/q}
    \right) \\
    \E \left[
      \sup_{f \in \cF}
      \left|
      \sum_{i=1}^n
      \varepsilon_i
      f(X_i)
      \right|
      ^q
    \right]
    ^{1/q}
    &\leq
    C_q
    \left(
      \E \left[
        \sup_{f \in \cF}
        \left|
        \sum_{i=1}^n
        \varepsilon_i
        f(X_i)
        \right|
      \right]
      +
      \E \left[
        \max_{1 \leq i \leq n}
        \sup_{f \in \cF}
        \big| f(X_i) \big|^q
      \right]^{1/q}
    \right).
  \end{align*}
\end{proof}

\begin{proof}[Lemma~\ref{lem:maximal_entropy}]

  We follow the proof of Theorem~5.2
  from \citet{chernozhukov2014gaussian},
  using our i.n.i.d.\ versions of the symmetrization inequality
  (Lemma~\ref{lem:symmetrization}),
  Rademacher contraction principle
  (Lemma~\ref{lem:contraction})
  and Hoffmann--J{\o}rgensen inequality
  (Lemma~\ref{lem:hoffmann}).

  Without loss of generality,
  we may assume that $J(1, \cF, F) < \infty$
  as otherwise there is nothing to prove,
  and that $F > 0$ everywhere on $S$.
  Let $\P_n = n^{-1} \sum_i \delta_{X_i}$
  be the empirical distribution
  of $X_i$,
  and define the empirical variance bound
  $\sigma_n^2 = \sup_\cF n^{-1} \sum_i f(X_i)^2$.
  By the i.n.i.d.\ symmetrization inequality
  (Lemma~\ref{lem:symmetrization}),
  \begin{align*}
    \E \left[
      \sup_{f \in \cF}
      \big| G_n(f) \big|
    \right]
    &=
    \frac{1}{\sqrt n}
    \E \left[
      \sup_{f \in \cF}
      \left|
      \sum_{i=1}^n
      \Big(
      f(X_i)
      - \E[f(X_i)]
      \Big)
      \right|
    \right]
    \leq
    \frac{2}{\sqrt n}
    \E \left[
      \sup_{f \in \cF}
      \left|
      \sum_{i=1}^n
      \varepsilon_i
      f(X_i)
      \right|
    \right],
  \end{align*}
  where $\varepsilon_1, \ldots, \varepsilon_n$
  are independent Rademacher random variables
  which are independent of $X_1, \ldots, X_n$.
  Then the standard entropy integral inequality
  from the proof of Theorem~5.2 in
  the supplemental materials for
  \citet{chernozhukov2014gaussian}
  gives for a universal constant $C_1 > 0$,
  \begin{align*}
    \frac{1}{\sqrt n}
    \E \left[
      \sup_{f \in \cF}
      \left|
      \sum_{i=1}^n
      \varepsilon_i
      f(X_i)
      \right|
      \Bigm\vert
      X_1, \ldots, X_n
    \right]
    &\leq
    C_1 \|F\|_{\P_n,2}
    J(\sigma_n / \|F\|_{\P_n,2}, \cF, F).
  \end{align*}
  Taking marginal expectations
  and applying Jensen's inequality along with
  a convexity result for the covering integral,
  as in Lemma~A.2 in \citet{chernozhukov2014gaussian},
  gives
  \begin{align*}
    Z
    &\vcentcolon=
    \frac{1}{\sqrt n}
    \E \left[
      \sup_{f \in \cF}
      \left|
      \sum_{i=1}^n
      \varepsilon_i
      f(X_i)
      \right|
    \right]
    \leq
    C_1 \|F\|_{\bar\P,2}
    J(\E[\sigma_n^2]^{1/2} / \|F\|_{\bar\P,2}, \cF, F).
  \end{align*}
  Next we use the symmetrization inequality
  (Lemma~\ref{lem:symmetrization}),
  the contraction principle
  (Lemma~\ref{lem:contraction}),
  the Cauchy--Schwarz inequality
  and the Hoffmann--J{\o}rgensen inequality
  (Lemma~\ref{lem:hoffmann})
  to deduce that
  \begin{align*}
    \E[\sigma_n^2]
    &=
    \E\left[
      \sup_{f \in \cF}
      \frac{1}{n}
      \sum_{i=1}^n
      f(X_i)^2
    \right]
    \leq
    \sup_{f \in \cF}
    \E_{\bar\P} \left[
      f(X_i)^2
    \right]
    + \frac{1}{n}
    \E\left[
      \sup_{f \in \cF}
      \left|
      \sum_{i=1}^n
      f(X_i)^2
      - \E \left[
        f(X_i)^2
      \right]
      \right|
    \right] \\
    &\leq
    \sigma^2
    + \frac{2}{n}
    \E\left[
      \sup_{f \in \cF}
      \left|
      \sum_{i=1}^n
      \varepsilon_i
      f(X_i)^2
      \right|
    \right]
    \leq
    \sigma^2
    + \frac{8}{n}
    \E\left[
      M
      \sup_{f \in \cF}
      \left|
      \sum_{i=1}^n
      \varepsilon_i
      f(X_i)
      \right|
    \right] \\
    &\leq
    \sigma^2
    + \frac{8}{n}
    \E\left[
      M^2
    \right]^{1/2}
    \E\left[
      \sup_{f \in \cF}
      \left|
      \sum_{i=1}^n
      \varepsilon_i
      f(X_i)
      \right|^2
    \right]^{1/2} \\
    &\leq
    \sigma^2
    + \frac{8}{n}
    \|M\|_{\P,2} \,
    C_2
    \left(
      \E \left[
        \sup_{f \in \cF}
        \left|
        \sum_{i=1}^n
        \varepsilon_i
        f(X_i)
        \right|
      \right]
      +
      \E \left[
        \max_{1 \leq i \leq n}
        \sup_{f \in \cF}
        \big| f(X_i) \big|^2
      \right]^{1/2}
    \right) \\
    &=
    \sigma^2
    + \frac{8C_2}{n}
    \|M\|_{\P,2} \,
    \left(
      \sqrt{n} Z
      +
      \|M\|_{\P,2}
    \right)
    \lesssim
    \sigma^2
    +
    \frac{\|M\|_{\P,2} Z}{\sqrt n}
    +
    \frac{\|M\|_{\P,2}^2}{n},
  \end{align*}
  where $\lesssim$ indicates a bound up to a universal constant.
  Hence taking a square root we see that,
  following the notation from the proof of Theorem~5.2
  in the supplemental materials to
  \citet{chernozhukov2014gaussian},
  \begin{align*}
    \sqrt{\E[\sigma_n^2]}
    &\lesssim
    \sigma
    +
    \|M\|_{\P,2}^{1/2} Z^{1/2} n^{-1/4}
    +
    \|M\|_{\P,2} n^{-1/2}
    \lesssim
    \|F\|_{\bar\P,2}
    \left( \Delta \vee \sqrt{DZ} \right),
  \end{align*}
  where
  $\Delta^2 = \|F\|_{\bar\P,2}^{-2}
  \big(\sigma^2 \vee (\|M\|_{\P,2}^2 / n) \big) \geq \delta^2$
  and
  $D = \|M\|_{\P,2} n^{-1/2} \|F\|_{\bar\P,2}^{-2}$.
  Thus returning to our bound on $Z$,
  we now have
  \begin{align*}
    Z
    &\lesssim
    \|F\|_{\bar\P,2}
    J(\Delta \vee \sqrt{DZ}, \cF, F).
  \end{align*}
  The final steps proceed as
  in the proof of Theorem~5.2
  from \citet{chernozhukov2014gaussian},
  considering cases separately for
  $\Delta \geq \sqrt{DZ}$
  and
  $\Delta < \sqrt{DZ}$,
  and applying convexity properties of
  the entropy integral $J$.
\end{proof}

\begin{proof}[Lemma~\ref{lem:maximal_vc_inid}]

  We are assuming the VC-type condition that
  \begin{align*}
    \sup_\Q N(\cF, \rho_\Q, \varepsilon \|F\|_{\Q,2})
    &\leq
    (C_1/\varepsilon)^{C_2}
  \end{align*}
  for all $\varepsilon \in (0,1]$,
  for some constants
  $C_1 \geq e$
  and $C_2 \geq 1$.
  Hence for $\delta \in (0,1]$,
  the entropy integral can be bounded as follows.
  \begin{align*}
    J\big(\delta, \cF, F\big)
    &=
    \int_0^\delta
    \sqrt{1 +
    \sup_\Q \log N(\cF, \rho_\Q, \varepsilon \|F\|_{\Q,2})}
    \diff{\varepsilon}
    \leq
    \int_0^\delta
    \sqrt{1 +
      C_2 \log (C_1/\varepsilon)}
    \diff{\varepsilon} \\
    &\leq
    \int_0^\delta
    \left(
      1 +
      \sqrt{C_2 \log (C_1/\varepsilon)}
    \right)
    \diff{\varepsilon}
    =
    \delta
    + \sqrt{C_2}
    \int_0^\delta
    \sqrt{\log (C_1/\varepsilon)}
    \diff{\varepsilon} \\
    &\leq
    \delta
    + \sqrt{\frac{C_2}{\log (C_1/\delta)}}
    \int_0^\delta
    \log (C_1/\varepsilon)
    \diff{\varepsilon}
    =
    \delta
    + \sqrt{\frac{C_2}{\log (C_1/\delta)}}
    \big(
    \delta
    + \delta \log (C_1/\delta)
    \big)
    \leq
    3 \delta
    \sqrt{C_2 \log (C_1/\delta)}.
  \end{align*}
  The remaining equations now follow
  by Lemma~\ref{lem:maximal_entropy}.
\end{proof}

\subsubsection{Strong approximation results}

Before proving Lemma~\ref{lem:kmt_corollary},
we require the elementary characterization
of bounded-variation functions given in
Lemma~\ref{lem:bv_characterization}.

\begin{lemma}[A characterization of bounded-variation functions]
  \label{lem:bv_characterization}

  Let $\cV_1$ be
  the class of real-valued functions on $[0,1]$
  which are 0 at 1 and have total variation bounded by 1.
  Also define the class of
  half-interval indicator functions $\cI = \{\I[0,t]: t \in [0,1]\}$.
  For any topological vector space $\cX$,
  define the symmetric convex hull of a subset $\cY \subseteq \cX$
  as
  \begin{align*}
    \symconv \cY
    &=
    \left\{
    \sum_{i=1}^n
    \lambda_i
    y_i :
    \sum_{i=1}^n
    \lambda_i
    = 1, \
    \lambda_i
    \geq 0, \
    y_i \in \cY \cup -\cY, \
    n \in \bbN
    \right\}.
  \end{align*}
  Denote its topological closure by
  $\overline\symconv \ \cY$.
  Then under the topology induced
  by pointwise convergence,
  \begin{align*}
    \cV_1 \subseteq \overline\symconv \ \cI.
  \end{align*}

\end{lemma}

\begin{proof}[Lemma~\ref{lem:bv_characterization}]

  Firstly, let $\cD \subseteq \cV_1$
  be the class of real-valued functions
  on $[0,1]$
  which are
  0 at 1,
  have total variation exactly 1,
  and are weakly monotone decreasing.
  Therefore for $g \in \cD$,
  we have
  $\|g\|_\TV = g(0) = 1$.
  Let $S = \{s_1, s_2, \dots\} \subseteq [0,1]$
  be the countable set of discontinuity points of $g$.
  We want to find a sequence of
  convex combinations of elements of
  $\cI$ which converges pointwise to $g$.
  To do this, first define the sequence of meshes
  \begin{align*}
    A_n =
    \{s_k : 1 \leq k \leq n\}
    \cup
    \{k/n : 0 \leq k \leq n\}
  \end{align*}
  which satisfies
  $\bigcup_n A_n = S \cup ([0,1] \cap \Q)$.
  Endow $A_n$ with the ordering
  induced by the canonical order on $\R$,
  giving $A_n = \{a_1, a_2, \ldots\}$,
  and define the sequence of functions
  \begin{align*}
    g_n(x)
    = \sum_{k = 1}^{|A_n|-1}
    \I[0,a_k]
    \big( g(a_k) - g(a_{k+1}) \big),
  \end{align*}
  where clearly
  $\I[0, a_k] \in \cI$
  and
  $g(a_k) - g(a_{k+1}) \geq 0$
  and
  $\sum_{k = 1}^{|A_n|-1}
  \big(
  g(a_k) - g(a_{k+1})
  \big)
  = g(0) - g(1) = 1$.
  Therefore
  $g_n$ is a convex combination of elements of
  $\cI$.
  Further, note that for
  $a_k \in A_n$,
  we have
  \begin{align*}
    g_n(a_k)
    = \sum_{j = k}^{|A_n|-1}
    \big( g(a_j) - g(a_{j+1}) \big)
    = g(a_k) - g(a_{|A_n|})
    = g(a_k) - g(1)
    = g(a_k).
  \end{align*}
  Hence if $x \in S$,
  then eventually $x \in A_n$
  so $g_n(x) \to g(x)$.
  Alternatively in $x \not\in S$,
  then $g$ is continuous at $x$.
  But $g_n \to g$ on the dense set $\bigcup_n A_n$,
  so also $g_n(x) \to g(x)$.
  Hence $g_n \to g$
  pointwise on $[0,1]$.

  Now take $f \in \cV_1$.
  By the Jordan decomposition for
  total variation functions
  \citep{royden1988real},
  we can write
  $f = f^+ - f^-$,
  with
  $f^+$ and $f^-$ weakly decreasing,
  $f^+(1) = f^-(1) = 0$,
  and
  $\|f^+\|_\TV + \|f^-\|_\TV = \|f\|_\TV$.
  Supposing that both
  $\|f^+\|_\TV$ and $\|f^-\|_\TV$
  are strictly positive, let
  $g_n^+$ approximate
  the unit-variation function
  $f^+/\|f^+\|_\TV$
  and
  $g_n^-$ approximate $f^-/\|f^-\|_\TV$
  as above.
  Then since trivially
  \begin{align*}
    f =
    \|f^+\|_\TV f^+ / \|f^+\|_\TV
    - \|f^-\|_\TV f^- / \|f^-\|_\TV
    + \big(1 - \|f^+\|_\TV - \|f^-\|_\TV) \cdot 0,
  \end{align*}
  we have that
  the convex combination
  \begin{align*}
    g_n^+ \|f^+\|_\TV
    - g_n^- \|f^-\|_\TV
    + \big(1 - \|f^+\|_\TV - \|f^-\|_\TV) \cdot 0
  \end{align*}
  converges pointwise to $f$.
  This also holds if either of the total variations
  $\|f^\pm\|_\TV$
  are zero,
  since then the corresponding sequence $g_n^\pm$
  need not be defined.
  Now note that each of
  $g_n^+$, $\,-g_n^-$ and $0$
  are in $\symconv \cI$, so
  $f \in \overline\symconv \ \cI$
  under pointwise convergence.
\end{proof}

\begin{proof}[Lemma~\ref{lem:kmt_corollary}]

  We follow the Gaussian approximation method given in
  Section~2 of \citet{gine2004kernel}.
  The KMT approximation theorem \citep{komlos1975approximation}
  asserts the existence
  of a probability space
  carrying $n$ i.i.d.\ uniform random variables
  $\xi_1, \ldots, \xi_n \sim \cU[0,1]$
  and a standard Brownian motion
  $B_n(s): s \in [0,1]$
  such that if
  \begin{align*}
    \alpha_n(s)
    &\vcentcolon=
    \frac{1}{\sqrt{n}}
    \sum_{i=1}^n
    \big(
    \I\{\xi_i \leq s\} - s
    \big),
    &\beta_n(s)
    &\vcentcolon=
    B_n(s) - s B_n(1),
  \end{align*}
  then
  for some universal positive constants
  $C_1$, $C_2$, $C_3$
  and for all $t > 0$,
  \begin{align*}
    \P\left(
      \sup_{s \in [0,1]}
      \big| \alpha_n(s) - \beta_n(s) \big|
      > \frac{t + C_1\log n}{\sqrt{n}}
    \right)
    \leq C_2 e^{-C_3 t}.
  \end{align*}
  We can
  view $\alpha_n$ and $\beta_n$ as random functionals
  defined on the class of
  half-interval indicator functions
  $\cI = \big\{\I[0,s]: s \in [0,1]\big\}$
  in the following way.
  \begin{align*}
    \alpha_n(\I[0,s])
    &= \frac{1}{\sqrt{n}}
    \sum_{i=1}^n
    \big( \I[0,s](\xi_i) - \E[\I[0,s](\xi_i)]), \\
    \beta_n(\I[0,s])
    &= \int_0^1 \I[0,s](u) \diff{B_n(u)}
    - B_n(1) \int_0^1 \I[0,s](u) \diff{u},
  \end{align*}
  where the integrals are defined as It\^o and
  Riemann--Stieltjes integrals in
  the usual way for stochastic integration against semimartingales
  \citep[Chapter~5]{legall2016brownian}.
  Now we extend their definitions to the class
  $\cV_1$
  of functions on $[0,1]$
  which are 0 at 1 and have total variation bounded by 1.
  This is achieved by
  noting that by Lemma~\ref{lem:bv_characterization},
  we have
  $\cV_1 \subseteq \overline\symconv \ \cI$
  where $\overline{\symconv} \ \cI$ is the
  smallest
  symmetric convex class containing $\cI$
  which is closed under pointwise convergence.
  Thus by the dominated convergence theorem,
  every function in $\cV_1$ is approximated in $L^2$ by finite convex
  combinations of functions in $\pm\cI$,
  and the extension to $g \in \cV_1$ follows
  by linearity and $L^2$ convergence of (stochastic) integrals:
  \begin{align*}
    \alpha_n(g)
    &=
    \frac{1}{\sqrt{n}}
    \sum_{i=1}^n
    \big( g(\xi_i) - \E[g(\xi_i)]),
    &\beta_n(g)
    &= \int_0^1 g(s) \diff{B_n(s)}
    - B_n(1) \int_0^1 g(s) \diff{s}.
  \end{align*}
  Now we show that the norm induced on
  $(\alpha_n - \beta_n)$
  by the function class $\cV_1$ is a.s.\ identical to the
  supremum norm.
  Writing the sums as integrals and using integration by parts
  for finite-variation Lebesgue--Stieltjes and It\^o integrals,
  and recalling that $g(1) = \alpha_n(0) = B_n(0) = 0$,
  we see
  \begin{align*}
    \sup_{g \in \cV_1}
    \big|\alpha_n(g) - \beta_n(g)\big|
    &=
    \sup_{g \in \cV_1}
    \left|
    \int_0^1 g(s) \diff{\alpha_n(s)}
    - \int_0^1 g(s) \diff{B_n(s)}
    + B_n(1) \int_0^1 g(s) \diff{s}
    \right| \\
    &=
    \sup_{g \in \cV_1}
    \left|
    \int_0^1 \alpha_n(s) \diff{g(s)}
    - \int_0^1 B_n(s) \diff{g(s)}
    + B_n(1) \int_0^1 s \diff{g(s)}
    \right| \\
    &=
    \sup_{g \in \cV_1}
    \left|
    \int_0^1 \big(\alpha_n(s) - \beta_n(s)\big)
    \diff{g(s)}
    \right|
    = \sup_{s \in [0,1]}
    \big|
    \alpha_n(s) - \beta_n(s)
    \big|,
  \end{align*}
  where in the last line
  the upper bound is because $\|g\|_\TV \leq 1$
  and the lower bound is by taking
  $g_\varepsilon = \pm \I[0,s_\varepsilon]$ where
  $|\alpha_n(s_\varepsilon) - \beta_n(s_\varepsilon)|
  \geq \sup_s |\alpha_n(s) - \beta_n(s)| -
  \varepsilon$.
  Hence we obtain
  \begin{align}
    \label{eq:kmt_concentration}
    \P\left(
      \sup_{g \in \cV_1}
      \big|\alpha_n(g) - \beta_n(g)\big|
      > \frac{t + C_1\log n}{\sqrt{n}}
    \right)
    \leq C_2 e^{-C_3 t}.
  \end{align}
  Now define $V_n = \sup_{x \in \R} \|g_n(\cdot, x)\|_\TV$,
  noting that if $V_n = 0$ then the result is trivially true
  by setting $Z_n = 0$.
  Let $F_X$ be the common c.d.f.\ of $X_i$,
  and define the quantile function
  $F_X^{-1}(s) = \inf \{u: F_X(u) \geq s\}$ for $s \in [0,1]$,
  writing $\inf \emptyset = \infty$
  and $\inf \R = -\infty$.
  Consider the function class
  \begin{align*}
    \cG_n = \big\{
    V_n^{-1} g_n\big(F_X^{-1}(\cdot), x\big)
    - V_n^{-1} g_n\big(F_X^{-1}(1), x\big)
    : x \in \R \big\},
  \end{align*}
  noting that $g_n(\cdot,x)$
  is finite-variation so
  $g_n(\pm \infty, x)$
  can be interpreted as
  the relevant limit.
  By monotonicity of $F_X$ and the definition of $V_n$,
  the members of $\cG_n$ have total variation of at most $1$
  and are 0 at 1, implying that
  $\cG_n \subseteq \cV_1$.
  Noting that $\alpha_n$ and $\beta_n$ are random
  linear operators which a.s.\ annihilate
  constant functions,
  define
  \begin{align*}
    Z_n(x)
    &=
    \beta_n \Big(g_n\big(F_X^{-1}(\cdot), x\big)\Big)
    = V_n \beta_n \Big(
    V_n^{-1} g_n\big(F_X^{-1}(\cdot), x\big)
    - V_n^{-1} g_n\big(F_X^{-1}(1), x\big)
    \Big),
  \end{align*}
  which is a mean-zero Gaussian process with continuous trajectories.
  Its covariance structure is
  \begin{align*}
    &\E[Z_n(x) Z_n(x')] \\
    &=
    \E\bigg[
      \left(
        \int_0^1 g_n\big(F_X^{-1}(s),x\big) \diff{B_n(s)}
        - B_n(1) \int_0^1 g_n\big(F_X^{-1}(s),x\big) \diff{s}
      \right) \\
      &\quad\times
      \left(
        \int_0^1 g_n\big(F_X^{-1}(s),x'\big) \diff{B_n(s)}
        - B_n(1) \int_0^1 g_n\big(F_X^{-1}(s),x'\big) \diff{s}
      \right)
      \bigg] \\
    &=
    \E\left[
      \int_0^1 g_n\big(F_X^{-1}(s),x\big) \diff{B_n(s)}
      \int_0^1 g_n\big(F_X^{-1}(s),x'\big) \diff{B_n(s)}
    \right] \\
    &\quad- \int_0^1 g_n\big(F_X^{-1}(s),x\big) \diff{s} \
    \E\left[
      B_n(1) \int_0^1 g_n\big(F_X^{-1}(s),x'\big) \diff{B_n(s)}
    \right] \\
    &\quad-
    \int_0^1 g_n\big(F_X^{-1}(s),x'\big) \diff{s} \
    \E\left[
      B_n(1) \int_0^1 g_n\big(F_X^{-1}(s),x\big) \diff{B_n(s)}
    \right] \\
    &\quad+
    \int_0^1 g_n\big(F_X^{-1}(s),x\big) \diff{s}
    \int_0^1 g_n\big(F_X^{-1}(s),x'\big) \diff{s} \
    \E\left[
      B_n(1)^2
    \right] \\
    &=
    \int_0^1 g_n\big(F_X^{-1}(s),x\big)
    g_n\big(F_X^{-1}(s),x'\big) \diff{s}
    - \int_0^1 g_n\big(F_X^{-1}(s),x\big) \diff{s}
    \int_0^1 g_n\big(F_X^{-1}(s),x'\big) \diff{s} \\
    &=
    \E\Big[
    g_n\big(F_X^{-1}(\xi_i), x\big)
    g_n\big(F_X^{-1}(\xi_i), x'\big)
    \Big]
    - \E\Big[
    g_n\big(F_X^{-1}(\xi_i), x\big)
    \Big]
    \E\Big[
    g_n\big(F_X^{-1}(\xi_i), x'\big)
    \Big] \\
    &=
    \E\Big[
      g_n\big(X_i, x\big)
      g_n\big(X_i, x'\big)
      \Big]
    - \E\Big[
      g_n\big(X_i, x\big)
      \Big]
    \E\Big[
      g_n\big(X_i, x'\big)
      \Big]
    =
    \E\big[
      G_n(x)
      G_n(x')
      \big]
  \end{align*}
  as desired, where we used the It\^o isometry for stochastic integrals,
  writing
  $B_n(1) = \int_0^1 \diff{B_n(s)}$;
  and noting that $F_X^{-1}(\xi_i)$
  has the same distribution as $X_i$.
  Finally, note that
  \begin{align*}
    G_n(x)
    &=
    \alpha_n \Big(g_n\big(F_X^{-1}(\cdot), x\big)\Big)
    = V_n \alpha_n \Big(
    V_n^{-1} g_n\big(F_X^{-1}(\cdot), x\big)
    - V_n^{-1} g_n\big(F_X^{-1}(1), x\big)
    \Big),
  \end{align*}
  and so by Equation~\ref{eq:kmt_concentration}
  \begin{align*}
    \P\left(
      \sup_{x \in \R}
      \Big|G_n(x) - Z_n(x)\Big|
      > V_n \frac{t + C_1 \log n}{\sqrt n}
    \right)
    \leq
    \P\left(
      \sup_{g \in \cV_1}
      \big|\alpha_n(g) - \beta_n(g)\big|
      > \frac{t + C_1\log n}{\sqrt{n}}
    \right)
    \leq C_2 e^{-C_3 t}.
  \end{align*}
\end{proof}

\begin{proof}[Lemma~\ref{lem:yurinskii_corollary}]

  Take $0 < \delta_n \leq \Leb(\cX_n)$ and let
  $\cX_n^\delta = \big\{ x_1, \dots, x_{|\cX_n^\delta|}\big\}$
  be a $\delta_n$-covering of $\cX_n$ with cardinality
  $|\cX_n^\delta| \leq \Leb(\cX_n)/\delta_n$.
  Suppose that $\left|\log \delta_n\right| \lesssim C_1 \log n$
  up to a universal constant.
  We first use the Yurinskii coupling to
  construct a Gaussian process
  $Z_n$
  which is close to $G_n$
  on this finite cover.
  Then we bound the fluctuations in $G_n$
  and in $Z_n$
  using entropy methods.

  \proofparagraph{Yurinskii coupling}

  Define the i.n.i.d.\
  and mean-zero variables
  \begin{align*}
    h_i(x)
    &=
    \frac{1}{\sqrt n}
    \Big(
    g_n(X_i', x)
    - \E[g_n(X_i', x)]
    \Big),
  \end{align*}
  where $X_1', \ldots, X_n'$
  are independent copies of $X_1, \ldots, X_n$
  on some new probability space,
  so that we have
  $G_n(x) = \sum_{i=1}^n h_i(x)$
  in distribution.
  Also define the length-$|\cX_n^\delta|$ random vector
  \begin{align*}
    h_i^\delta
    &=
    \big(
    h_i(x): x \in \cX_n^\delta
    \big).
  \end{align*}
  By an extension of
  Yurinskii's coupling
  to general norms
  \citep[supplemental materials, Lemma~38]{belloni2019conditional},
  there exists on the new probability space a
  Gaussian length-$|\cX_n^\delta|$ vector $Z_n^\delta$
  which is mean-zero
  and with the same covariance structure as
  $
  \sum_{i=1}^n
  h_i^\delta
  $
  satisfying
  \begin{align*}
    \P\left(
      \bigg\|
      \sum_{i=1}^n
      h_i^\delta
      - Z_n^\delta
      \bigg\|_\infty
      > 3 t_n
    \right)
    \leq
    \min_{s > 0}
    \left(
      2 \P\big( \|N\|_\infty > s)
      + \frac{\beta s^2}{t_n^3}
    \right),
  \end{align*}
  where
  \begin{align*}
    \beta
    = \sum_{i=1}^n
    \Big(
    \E\big[\|h_i^\delta\|_2^2 \,
      \|h_i^\delta\|_\infty
      \big]
    + \E\big[\|z_i\|_2^2 \,
      \|z_i\|_\infty
      \big]
    \Big),
  \end{align*}
  with
  $z_i \sim \cN(0, \Var[h_i^\delta])$
  independent and
  $N \sim \cN(0, I_{|\cX_n^\delta|})$.
  By the almost sure
  and variance bounds
  on $g_n$,
  \begin{align*}
    \E\big[\|h_i^\delta\|_2^2 \,
      \|h_i^\delta\|_\infty \,
      \big]
    \leq
    \frac{M_n}{\sqrt n}
    \E\big[\|h_i^\delta\|_2^2 \,
      \big]
    =
    \frac{M_n}{\sqrt n}
    \sum_{x \in \cX_n^\delta}
    \E\big[h_i(x)^2 \,
      \big]
    \leq
    \frac{M_n}{\sqrt n}
    \frac{|\cX_n^\delta| \sigma_n^2}{n}
    \leq
    \frac{M_n \sigma_n^2 \Leb(\cX_n)}{n^{3/2}\delta_n}.
  \end{align*}
  By the fourth moment bound for Gaussian variables,
  \begin{align*}
    \E\big[
      \|z_i\|_2^4 \,
      \big]
    &\leq
    |\cX_n^\delta| \,
    \E\big[
      \|z_i\|_4^4
      \big]
    \leq
    |\cX_n^\delta|^2 \,
    \max_j
    \E\big[
    (z_i^{(j)})^4
    \big]
    \leq
    3
    |\cX_n^\delta|^2 \,
    \max_j
    \E\big[
    (z_i^{(j)})^2
    \big]^2 \\
    &=
    3
    |\cX_n^\delta|^2 \,
    \max_{x \in \cX_n^\delta}
    \E\big[
      h_i(x)^2
      \big]^2
    \leq
    \frac{3\sigma_n^4 \Leb(\cX_n)^2}{n^2\delta_n^2} .
  \end{align*}
  Also by Jensen's inequality
  and for $|\cX_n^\delta| \geq 2$,
  assuming $C_1 > 1$ without loss of generality,
  \begin{align*}
    \E\big[
      \|z_i\|_\infty^2
      \big]
    &\leq
    \frac{4 \sigma_n^2}{n}
    \log
    \E\big[
      e^{\|z_i\|_\infty^2 / (4\sigma_n^2/n)}
      \big]
    \leq
    \frac{4 \sigma_n^2}{n}
    \log
    \E\left[
      \sum_{j=1}^{|\cX_n^\delta|}
      e^{(z_i^{(j)})^2 / (4\sigma_n^2/n)}
    \right]
    \leq
    \frac{4\sigma_n^2}{n}
    \log \big(2|\cX_n^\delta|\big) \\
    &\leq
    \frac{4\sigma_n^2}{n}
    \left(
      \log 2 + \log \Leb(\cX_n) - \log \delta_n
    \right)
    \leq
    \frac{12 C_1 \sigma_n^2 \log n}{n},
  \end{align*}
  where we used the moment
  generating function of a $\chi_1^2$ random variable.
  Therefore we can apply the Cauchy--Schwarz inequality
  to obtain
  \begin{align*}
    \E\big[\|z_i\|_2^2 \,
      \|z_i\|_\infty
      \big]
    &\leq
    \sqrt{
      \E\big[\|z_i\|_2^4
        \big]}
    \sqrt{
      \E\big[
        \|z_i\|_\infty^2
        \big]}
    \leq
    \sqrt{
      \frac{3\sigma_n^4 \Leb(\cX_n)^2}{n^2\delta_n^2}}
    \sqrt{
      \frac{12 C_1 \sigma_n^2 \log n}{n}
    }
    \leq
    \frac{6\sigma_n^3
      \Leb(\cX_n)
      \sqrt{C_1 \log n}}{n^{3/2} \delta_n}.
  \end{align*}
  Now summing over the $n$ samples gives
  \begin{align*}
    \beta
    \leq
    \frac{M_n \sigma_n^2 \Leb(\cX_n)}{\sqrt n \delta_n}
    + \frac{6\sigma_n^3 \Leb(\cX_n) \sqrt{C_1 \log n}}
    {\sqrt n \delta_n}
    =
    \frac{\sigma_n^2 \Leb(\cX_n)}{\sqrt n \delta_n}
    \Big(M_n + 6\sigma_n \sqrt{C_1 \log n}\Big).
  \end{align*}
  By a union bound
  and Gaussian tail probabilities,
  we have that
  $\P\big( \|N\|_\infty > s)
  \leq 2|\cX_n^\delta| e^{-s^2/2}$.
  Thus we get the following Yurinskii coupling inequality
  for all $s > 0$:
  \begin{align*}
    \P\left(
      \bigg\|
      \sum_{i=1}^n
      h_i^\delta
      - Z_n^\delta
      \bigg\|_\infty
      > t_n
    \right)
    &\leq
    \frac{4 \Leb(\cX_n)}{\delta_n}
    e^{-s^2/2}
    + \frac{\sigma_n^2 \Leb(\cX_n) s^2}{\sqrt n \delta_n t_n^3}
    \Big(M_n + 6 \sigma_n \sqrt{C_1 \log n}\Big).
  \end{align*}
  Note that
  $Z_n^\delta$
  now extends
  by the Vorob'ev--Berkes--Philipp Theorem
  (Lemma~\ref{lem:vbp})
  to a mean-zero Gaussian
  process
  $Z_n$ on the compact interval $\cX_n$
  with covariance structure given by
  \begin{align*}
    \E\big[
      Z_n(x)
      Z_n(x')
      \big]
    =
    \E\big[
      G_n(x)
      G_n(x')
      \big],
  \end{align*}
  satisfying for any $s' > 0$
  \begin{align*}
    \P\left(
      \sup_{x \in \cX_n^\delta}
      \big|
      G_n(x) - Z_n(x)
      \big|
      > t_n
    \right)
    &\leq
    \frac{4 \Leb(\cX_n)}{\delta_n}
    e^{-s^2/2}
    + \frac{\sigma_n^2 \Leb(\cX_n) s^2}{\sqrt n \delta_n t_n^3}
    \Big(M_n + 6 \sigma_n \sqrt{C_1 \log n}\Big).
  \end{align*}

  \proofparagraph{Regularity of $G_n$}

  Next we bound the fluctuations in
  the empirical process $G_n$.
  Consider the following classes of functions on $S$
  and their associated (constant) envelope functions.
  By continuity of $g_n$,
  each class is pointwise measurable
  (to see this, restrict the index sets to rationals).
  \begin{align*}
    \cG_n
    &=
    \big\{
    g_n(\cdot, x):
    x \in \cX_n
    \big\}, \\
    \Env(\cG_n)
    &=
    M_n, \\
    \cG_n^\delta
    &=
    \big\{
    g_n(\cdot, x)
    - g_n(\cdot, x'):
    x, x' \in \cX_n,
    |x-x'| \leq \delta_n
    \big\}, \\
    \Env(\cG_n^\delta)
    &=
    l_{n,\infty} \delta_n.
  \end{align*}
  We first show that
  for each $n$
  these are VC-type classes.
  To see this,
  note that by the uniform Lipschitz assumption
  we have that
  \begin{align*}
    \big\|
    g_n(\cdot, x)
    - g_n(\cdot, x')
    \big\|_\infty
    &\leq l_{n,\infty} |x-x'|
  \end{align*}
  for all $x,x' \in \cX_n$.
  Therefore with $\Q$ ranging over the
  finitely-supported distributions
  on $(S, \cS)$,
  noting that any $\|\cdot\|_\infty$-cover
  is a $\rho_\Q$-cover,
  \begin{align*}
    \sup_\Q
    N\big(\cG_n, \rho_\Q, \varepsilon l_{n,\infty} \Leb(\cX_n)\big)
    &\leq
    N\big(\cG_n, \|\cdot\|_\infty,
    \varepsilon l_{n,\infty} \Leb(\cX_n)\big)
    \leq
    N\big(\cX_n, |\cdot|, \varepsilon \Leb(\cX_n)\big)
    \leq
    1/\varepsilon.
  \end{align*}
  Replacing $\varepsilon$ by
  $\varepsilon M_n/(l_{n,\infty} \Leb(\cX_n))$
  gives
  \begin{align*}
    \sup_\Q
    N\big(\cG_n, \rho_\Q, \varepsilon M_n \big)
    &\leq
    \frac{l_{n,\infty} \Leb(\cX_n)}{\varepsilon M_n},
  \end{align*}
  and so $\cG_n$
  is a VC-type class.
  To see that $\cG_n^\delta$
  is also a VC-type class,
  we construct a cover in the following way.
  Let $\cF_n$ be an $\varepsilon$-cover
  for $(\cG_n, \|\cdot\|_\infty)$.
  Then by the triangle inequality,
  $\cF_n - \cF_n$ is a $2\varepsilon$-cover
  for $(\cG_n - \cG_n, \|\cdot\|_\infty)$
  of cardinality at most $|\cF_n|^2$,
  where the subtractions are defined as set subtractions.
  Since $\cG_n^\delta \subseteq \cG_n - \cG_n$,
  we see that
  $\cF_n - \cF_n$ is a $2\varepsilon$-external cover
  for $\cG_n^\delta$.
  Thus
  \begin{align*}
    \sup_\Q
    N\big(\cG_n^\delta, \rho_\Q, \varepsilon l_{n,\infty} \Leb(\cX_n)\big)
    &\leq
    N\big(\cG_n^\delta, \|\cdot\|_\infty,
    \varepsilon l_{n,\infty} \Leb(\cX_n)\big)
    \leq
    N\big(\cG_n, \|\cdot\|_\infty,
    \varepsilon l_{n,\infty} \Leb(\cX_n)\big)^2
    \leq
    1/\varepsilon^2.
  \end{align*}
  Replacing $\varepsilon$ by
  $\varepsilon \delta_n/\Leb(\cX_n)$
  gives
  \begin{align*}
    \sup_\Q
    N\big(\cG_n^\delta, \rho_\Q, \varepsilon l_{n,\infty} \delta_n \big)
    &\leq
    \frac{\Leb(\cX_n)^2}{\varepsilon^2 \delta_n^2}
    \leq
    (C_{1,n}/\varepsilon)^{2}
  \end{align*}
  with $C_{1,n} = \Leb(\cX_n) / \delta_n$,
  demonstrating that $\cG_n^\delta$
  forms a VC-type class.
  We now apply the maximal inequality
  for i.n.i.d.\ data
  given in
  Lemma~\ref{lem:maximal_vc_inid}.
  To do this,
  note that
  $\sup_{\cG_n^\delta} \|g\|_{\bar\P,2}
  \leq l_{n,2} \delta_n$
  by the $L^2$ Lipschitz condition,
  and recall
  $\Env(\cG_n^\delta) = l_{n,\infty} \delta_n$.
  Therefore
  Lemma~\ref{lem:maximal_vc_inid}
  with
  $\|F\|_{\bar\P,2} = l_{n,\infty} \delta_n$,
  $\|M\|_{\P,2} = l_{n,\infty} \delta_n$
  and $\sigma = l_{n,2} \delta_n$
  gives,
  up to universal constants
  \begin{align*}
    \E\left[
      \sup_{g \in \cG_n^\delta}
      \left|
      \frac{1}{\sqrt{n}}
      \sum_{i=1}^n
      \Big(
      g(X_i)
      - \E[g(X_i)]
      \Big)
      \right|
    \right]
    &\lesssim
    \sigma
    \sqrt{2 \log \big(C_{1,n} \|F\|_{\bar\P,2}/\sigma\big)}
    +
    \frac{\|M\|_{\P,2} 2 \log \big(C_{1,n} \|F\|_{\bar\P,2}/\sigma\big)}
    {\sqrt{n}} \\
    &\lesssim
    l_{n,2} \delta_n
    \sqrt{C_1 \log n}
    +
    \frac{l_{n,\infty} \delta_n}{\sqrt n}
    C_1 \log n,
  \end{align*}
  and hence by Markov's inequality,
  \begin{align*}
    &\P\left(
      \sup_{|x-x'| \leq \delta_n}
      \big|
      G_n(x) - G_n(x')
      \big|
      > t_n
    \right) \\
    &=
    \P\left(
      \sup_{|x-x'| \leq \delta_n}
      \frac{1}{\sqrt{n}}
      \left|
      \sum_{i=1}^n
      \Big(
      g_n(X_i, x) - \E[g_n(X_i, x)]
      - g_n(X_i, x') + \E[g_n(X_i, x')]
      \Big)
      \right|
      > t_n
    \right) \\
    &=
    \P\left(
      \sup_{g \in \cG_n^\delta}
      \left|
      \frac{1}{\sqrt{n}}
      \sum_{i=1}^n
      \Big(
      g(X_i) - \E[g(X_i)]
      \Big)
      \right|
      > t_n
    \right)
    \leq
    \frac{1}{t}
    \E\left[
      \sup_{g \in \cG_n^\delta}
      \left|
      \frac{1}{\sqrt{n}}
      \sum_{i=1}^n
      \Big(
      g(X_i) - \E[g(X_i)]
      \Big)
      \right|
    \right] \\
    &\lesssim
    \frac{l_{n,2} \delta_n}{t_n}
    \sqrt{C_1 \log n}
    + \frac{l_{n,\infty} \delta_n}{t_n \sqrt n} C_1 \log n.
  \end{align*}

  \proofparagraph{regularity of $Z_n$}

  Next we bound the fluctuations in the Gaussian process
  $Z_n$.
  Let $\rho$ be the following semimetric:
  \begin{align*}
    \rho(x, x')^2
    &=
    \E\big[\big( Z_n(x) - Z_n(x') \big)^2\big]
    =
    \E\big[\big( G_n(x) - G_n(x') \big)^2\big] \\
    &=
    \frac{1}{n}
    \sum_{i=1}^n
    \E\big[\big( h_i(x) - h_i(x') \big)^2\big]
    \leq
    l_{n,2}^2 \, |x - x'|^2.
  \end{align*}
  Hence
  $\rho(x, x')
  \leq
  l_{n,2} \, |x - x'|$.
  By
  the Gaussian process maximal inequality from
  Lemma~\ref{lem:gaussian_process_maximal},
  we obtain that
  \begin{align*}
    &\E\bigg[
      \sup_{|x - x'| \leq \delta_n}
      \big|
      Z_n(x) - Z_n(x')
      \big|
      \bigg]
    \lesssim
    \E\bigg[
      \sup_{\rho(x,x') \leq l_{n,2} \delta_n}
      \big|
      Z_n(x) - Z_n(x')
      \big|
      \bigg]
    \leq
    \int_0^{l_{n,2} \delta_n}
    \sqrt{\log N(\varepsilon, \cX_n, \rho)}
    \diff{\varepsilon} \\
    &\quad\leq
    \int_0^{l_{n,2} \delta_n}
    \sqrt{\log N(\varepsilon / l_{n,2}, \cX_n, |\cdot|)}
    \diff{\varepsilon}
    \leq
    \int_0^{l_{n,2} \delta_n}
    \sqrt{\log \left( 1 + \frac{\Leb(\cX_n) l_{n,2}}{\varepsilon} \right)}
    \diff{\varepsilon} \\
    &\quad\leq
    \int_0^{l_{n,2} \delta_n}
    \sqrt{\log \left( \frac{2\Leb(\cX_n) l_{n,2}}{\varepsilon} \right)}
    \diff{\varepsilon}
    \leq
    \log \left(\frac{2\Leb(\cX_n)}{\delta_n} \right)^{-1/2}
    \int_0^{l_{n,2} \delta_n}
    \log \left( \frac{2\Leb(\cX_n) l_{n,2}}{\varepsilon} \right)
    \diff{\varepsilon} \\
    &\quad=
    \log \left(\frac{2\Leb(\cX_n)}{\delta_n} \right)^{-1/2}
    \left(
      l_{n,2} \delta_n \log \left( 2 \Leb(\cX_n) l_{n,2} \right)
      + l_{n,2} \delta_n
      + l_{n,2} \delta_n \log \left( \frac{1}{l_{n,2} \delta_n} \right)
    \right) \\
    &\quad=
    \log \left(\frac{2\Leb(\cX_n)}{\delta_n} \right)^{-1/2}
    l_{n,2} \delta_n
    \left(
      1 +
      \log \left( \frac{2\Leb(\cX_n)}{\delta_n} \right)
    \right)
    \lesssim
    l_{n,2} \delta_n
    \sqrt{\log \left( \frac{\Leb(\cX_n)}{\delta_n} \right)}
    \lesssim
    l_{n,2} \delta_n
    \sqrt{C_1 \log n},
  \end{align*}
  where we used that $\delta_n \leq \Leb(\cX_n)$.
  So by Markov's inequality,
  \begin{align*}
    \P\left(
      \sup_{|x - x'| \leq \delta_n}
      \big|
      Z_n(x) - Z_n(x')
      \big|
      > t_n
    \right)
    &\lesssim
    t_n^{-1}
    l_{n,2} \delta_n
    \sqrt{C_1 \log n}.
  \end{align*}

  \proofparagraph{conclusion}

  By the results of the previous parts,
  we have up to universal constants that
  \begin{align*}
    &\P\left(
      \sup_{x \in \cX_n}
      \big|
      G_n(x) - Z_n(x)
      \big|
      > t_n
    \right) \\
    &\quad\leq
    \P\left(
      \sup_{x \in \cX_n^\delta}
      \big|
      G_n(x) - Z_n(x)
      \big|
      > t_n / 3
    \right)
    + \P\left(
      \sup_{|x-x'| \leq \delta_n}
      \big|
      G_n(x) - G_n(x')
      \big|
      > t_n / 3
    \right) \\
    &\qquad+
    \P\left(
      \sup_{|x - x'| \leq \delta_n}
      \big|
      Z_n(x) - Z_n(x')
      \big|
      > t_n / 3
    \right) \\
    &\quad\lesssim
    \frac{4 \Leb(\cX_n)}{\delta_n}
    e^{-s^2/2}
    + \frac{\sigma_n^2 \Leb(\cX_n) s^2}{\sqrt n \delta_n t_n^3}
    \Big(M_n + 6 \sigma_n \sqrt{C_1 \log n}\Big)
    + \frac{l_{n,2} \delta_n}{t_n}
    \sqrt{C_1 \log n}
    + \frac{l_{n,\infty} \delta_n}{t_n \sqrt n} C_1 \log n.
  \end{align*}
  Choosing an approximately optimal mesh size of
  \begin{align*}
    \delta_n
    &=
    \sqrt{
      \frac{\sigma_n^2 \Leb(\cX_n) \log n}{\sqrt n t_n^3}
      \Big(M_n + \sigma_n \sqrt{\log n}\Big)
    } \Bigg/
    \sqrt{
    t_n^{-1}
    l_{n,2}
    \sqrt{\log n}
    \left(
      1 + \frac{l_{n,\infty} \sqrt{\log n}}{l_{n,2} \sqrt{n}}
    \right)
    }
  \end{align*}
  satisfies $\log |\delta_n| \lesssim C_1 \log n$
  up to a universal constant,
  so taking $s$ a large enough multiple of $\sqrt{\log n}$ gives
  \begin{align*}
    &\P\left(
      \sup_{x \in \cX_n}
      \big|
      G_n(x) - Z_n(x)
      \big|
      > t_n
    \right) \\
    &\quad\lesssim
    \frac{4 \Leb(\cX_n)}{\delta_n}
    e^{-s^2/2}
    + \frac{\sigma_n^2 \Leb(\cX_n) s^2}{\sqrt n \delta_n t_n^3}
    \Big(M_n + 6 \sigma_n \sqrt{C_1 \log n}\Big)
    + \frac{l_{n,2} \delta_n}{t_n}
    \sqrt{C_1 \log n}
    + \frac{l_{n,\infty} \delta_n}{t_n \sqrt n} C_1 \log n \\
    &\quad\lesssim
    \delta_n
    \frac{l_{n,2} \sqrt {\log n}}{t_n}
    \left( 1 + \frac{l_{n,\infty} \sqrt{\log n}}{l_{n,2} \sqrt n} \right) \\
    &\quad\lesssim
    \frac{\sigma_n \sqrt{\Leb(\cX_n)} \sqrt{\log n}
      \sqrt{M_n + \sigma_n \sqrt{\log n}}}
    {n^{1/4} t_n^2}
    \sqrt{l_{n,2} \sqrt {\log n}
      + \frac{l_{n,\infty}}{\sqrt n} \log n}
  \end{align*}
  up to constants depending only on $C_1$.
\end{proof}

\subsubsection{The Vorob'ev--Berkes--Philipp theorem}

\begin{proof}[Lemma~\ref{lem:vbp}]

  The proof is by induction on the number of vertices in the tree.
  Let $\cT$ have $n$ vertices,
  and suppose
  that vertex $n$ is a leaf
  connected to vertex $n-1$ by an edge,
  relabelling the vertices if necessary.
  By the induction hypothesis we assume that there is a
  probability measure $\P^{(n-1)}$
  on $\prod_{i=1}^{n-1} \cX_i$
  whose projections onto $\cX_i$ are $\P_i$
  and whose projections onto $\cX_i \times \cX_j$ are $\P_{ij}$,
  for $i,j \leq n-1$.
  Now apply the original
  Vorob'ev--Berkes--Philipp theorem,
  which can be found as Theorem~1.1.10 in
  \citet{dudley1999uniform},
  to the spaces
  $\prod_{i=1}^{n-2} \cX_i$,\,
  $\cX_{n-1}$ and
  $\cX_n$;
  and to the laws
  $\P^{(n-1)}$
  and
  $\P_{n-1, n}$.
  This gives a law $\P^{(n)}$
  which agrees with $\P_i$
  at every vertex by definition,
  and agrees with
  $\P_{ij}$ for all $i,j \leq n-1$.
  It also agrees with $\P_{n-1,n}$,
  and this is the only edge touching vertex $n$.
  Hence $\P^{(n)}$ satisfies the desired properties.
\end{proof}

\subsection{Main results}

We give our main results on
consistency, minimax optimality, strong approximation,
covariance estimation, feasible inference
and counterfactual estimation.

\subsubsection{Bias}

We begin with a basic fact about Lipschitz functions.

\begin{lemma}[Lipschitz kernels are bounded]
  \label{lem:lipschitz_kernels_bounded}

  Let $\cX \subseteq \R$ be a connected set.
  Let $f: \cX \to \R$ satisfy the Lipschitz condition
  $|f(x) - f(x')| \leq C |x-x'|$ for some $C > 0$
  and all $x, x' \in \cX$.
  Suppose also that $f$ is a kernel in the sense that
  $\int_\cX f(x) \diff{x} = 1$.
  Then we have
  \begin{align*}
    \sup_{x \in \cX} |f(x)|
    &\leq
    C \Leb(\cX) + \frac{1}{\Leb(\cX)}.
  \end{align*}
  Now let $g: \cX \to [0,\infty)$ satisfy
    $|g(x) - g(x')| \leq C |x-x'|$ for some $C > 0$
    and all $x, x' \in \cX$.
    Suppose also that $g$ is a sub-kernel in the sense that
    $\int_\cX g(x) \diff{x} \leq 1$.
    Then for any $M \in \big(0, \Leb(\cX)\big]$,
  we have
  \begin{align*}
    \sup_{x \in \cX} f(x)
    &\leq
    C M + \frac{1}{M}.
  \end{align*}

\end{lemma}

\begin{remark}

  Applying Lemma~\ref{lem:lipschitz_kernels_bounded}
  to the density and kernel functions defined in
  Assumptions~\ref{ass:data} and~\ref{ass:kernel_bandwidth}
  yields the following.
  Firstly, since $k_h(\cdot, w)$ is $\CL/h^2$-Lipschitz
  on $[w \pm h] \cap \cW$ and integrates to one,
  we have by the first inequality in
  Lemma~\ref{lem:lipschitz_kernels_bounded} that
  \begin{align*}
    |k_h(s,w)|
    &\leq \frac{2 \CL + 1}{h} + \frac{1}{\Leb(\cW)}.
  \end{align*}
  Since each of
  $f_{W \mid AA}(\cdot \mid a,a')$,
  $f_{W \mid A}(\cdot \mid a)$ and
  $f_W$ is non-negative and $\CH$-Lipschitz on $\cW$
  and integrates to at most one over $\cW$,
  taking $M = \frac{1}{\sqrt \CH} \wedge \Leb(\cW)$
  in the second inequality in
  Lemma~\ref{lem:lipschitz_kernels_bounded}
  gives
  \begin{align*}
    f_{W \mid AA}(w \mid a,a')
    &\leq 2 \sqrt{\CH} + \frac{1}{\Leb(\cW)},
    &f_{W \mid A}(w \mid a)
    &\leq 2 \sqrt{\CH} + \frac{1}{\Leb(\cW)},
    &f_W(w)
    &\leq 2 \sqrt{\CH} + \frac{1}{\Leb(\cW)}.
  \end{align*}

\end{remark}

\begin{proof}[Lemma~\ref{lem:lipschitz_kernels_bounded}]

  We begin with the first inequality.
  Note that if $\Leb(\cX) = \infty$ there is nothing to prove.
  Suppose for contradiction that
  $|f(x)| > C \Leb(\cX) + \frac{1}{\Leb(\cX)}$
  for some $x \in \cX$.
  If $f(x) \geq 0$
  then by the Lipschitz property, for any $y \in \cX$,
  \begin{align*}
    f(y)
    \geq f(x) - C|y-x|
    > C \Leb(\cX) + \frac{1}{\Leb(\cX)} - C\Leb(\cX)
    = \frac{1}{\Leb(\cX)}.
  \end{align*}
  Similarly if $f(x) \leq 0$ then
  \begin{align*}
    f(y)
    \leq f(x) + C|y-x|
    < - C \Leb(\cX) - \frac{1}{\Leb(\cX)} + C\Leb(\cX)
    = -\frac{1}{\Leb(\cX)}.
  \end{align*}
  But then either
  $\int_\cX f(x) \diff{x} > \int_\cX 1/\Leb(\cX) \diff{x} = 1$
  or
  $\int_\cX f(x) \diff{x} < \int_\cX -1/\Leb(\cX) \diff{x} = -1 < 1$,
  giving a contradiction.

  For the second inequality,
  assume that $f$ is non-negative on $\cX$
  and take $M \in \big(0, \Leb(\cX)\big]$.
  Suppose for contradiction that
  $f(x) > C M + \frac{1}{M}$
  for some $x \in \cX$.
  Then again by the Lipschitz property,
  we have $f(y) \geq 1/M$
  for all $y$ such that $|y - x| \leq M$.
  Since $\cX$ is connected, we have
  $\Leb(\cX \cap [x \pm M]) \geq M$
  and so we deduce that
  $\int_\cX f(x) \diff{x} > M/M = 1$
  which is a contradiction.
\end{proof}

\begin{proof}[Lemma~\ref{lem:bias}]

  Begin by defining
  \begin{align*}
    P_p(s,w)
    &=
    \sum_{r = 0}^p
    \frac{f_W^{(r)}(w)}{r!}
    {(s-w)^r}
  \end{align*}
  for $s, w \in \cW$
  as the degree-$p$ Taylor polynomial of $f_W$,
  centered at $w$ and evaluated at $s$.
  Note that
  for $p \leq \flbeta-1$,
  by Taylor's theorem with Lagrange remainder,
  \begin{align*}
    f_W(s) - P_p(s,w)
    &=
    \frac{f_W^{(p+1)}(w')}{(p+1)!}
    (s-w)^{p+1}
  \end{align*}
  for some $w'$ between $w$ and $s$.
  Also note that for any $p$,
  \begin{align*}
    \int_{\cW}
    k_h(s,w)
    \big(
    P_p(s,w)
    - P_{p-1}(s,w)
    \big)
    \diff{s}
    &=
    \int_{\cW}
    k_h(s,w)
    \frac{f_W^{(p)}(w)}{p!}
    (s-w)^p
    \diff{s}
    = h^p B_p(w).
  \end{align*}
  Further, by the order of the kernel,
  \begin{align*}
    \E\big[\widehat f_W(w)\big]
    - f_W(w)
    &=
    \int_{\cW}
    k_h(s,w)
    f_W(s)
    \diff{s}
    - f_W(w)
    =
    \int_{\cW}
    k_h(s,w)
    \big(f_W(s) - f_W(w)\big)
    \diff{s} \\
    &=
    \int_{\cW}
    k_h(s,w)
    \big(f_W(s) - P_{p-1}(s,w)\big)
    \diff{s}.
  \end{align*}

  \proofparagraph{low-order kernel}
  Suppose that $p \leq \flbeta - 1$.
  Then
  \begin{align*}
    &\sup_{w \in \cW}
    \big|
    \E[\widehat f_W(w)]
    - f_W(w)
    - h^p B_p(w)
    \big|
    =
    \sup_{w \in \cW}
    \left|
    \int_{\cW}
    k_h(s,w)
    \big(f_W(s) - P_{p-1}(s,w)\big)
    \diff{s}
    - h^p B_p(w)
    \right| \\
    &\quad=
    \sup_{w \in \cW}
    \left|
    \int_{\cW}
    k_h(s,w)
    \big(
    f_W(s) - P_{p}(s,w)
    + P_{p}(s,w) - P_{p-1}(s,w)
    \big)
    \diff{s}
    - h^p B_p(w)
    \right| \\
    &\quad=
    \sup_{w \in \cW}
    \left|
    \int_{\cW}
    k_h(s,w)
    \big(
    f_W(s) - P_{p}(s,w)
    \big)
    \diff{s}
    \right|
    =
    \sup_{w \in \cW}
    \left|
    \int_{\cW}
    k_h(s,w)
    \frac{f_W^{(p+1)}(w')}{(p+1)!}
    (s-w)^{p+1}
    \diff{s}
    \right| \\
    &\quad\leq
    \sup_{w \in \cW}
    \left|
    \int_{[w \pm h]}
    \frac{\Ck}{h}
    \frac{\CH}{(p+1)!}
    h^{p+1}
    \diff{s}
    \right|
    \leq
    \frac{2\Ck \CH}{(p+1)!}
    h^{p+1}.
  \end{align*}

  \proofparagraph{order of kernel matches smoothness}
  Suppose that $p = \flbeta$.
  Then
  \begin{align*}
    &\sup_{w \in \cW}
    \big|
    \E[\widehat f_W(w)]
    - f_W(w)
    - h^p B_p(w)
    \big|
    =
    \sup_{w \in \cW}
    \left|
    \int_{\cW}
    k_h(s,w)
    \big(f_W(s) - P_{\flbeta - 1}(s,w)\big)
    \diff{s}
    - h^p B_p(w)
    \right| \\
    &\quad=
    \sup_{w \in \cW}
    \left|
    \int_{\cW}
    k_h(s,w)
    \big(
    f_W(s) - P_{\flbeta}(s,w)
    + P_{\flbeta}(s,w) - P_{\flbeta - 1}(s,w)
    \big)
    \diff{s}
    - h^{\flbeta} B_{\flbeta}(w)
    \right| \\
    &\quad=
    \sup_{w \in \cW}
    \left|
    \int_{\cW}
    k_h(s,w)
    \big(
    f_W(s) - P_{\flbeta}(s,w)
    \big)
    \diff{s}
    \right|
    =
    \sup_{w \in \cW}
    \left|
    \int_{\cW}
    k_h(s,w)
    \frac{f_W^{(\flbeta)}(w') - f_W^{(\flbeta)}(w)}{\flbeta!}
    (s-w)^{\flbeta}
    \diff{s}
    \right| \\
    &\quad\leq
    \sup_{w \in \cW}
    \left|
    \int_{[w \pm h]}
    \frac{\Ck}{h}
    \frac{\CH h^{\beta - \flbeta}}{\flbeta !}
    h^{\flbeta}
    \diff{s}
    \right|
    \leq
    \frac{2 \Ck \CH}{\flbeta !}
    h^\beta.
  \end{align*}

  \proofparagraph{high-order kernel}
  Suppose that $p \geq \flbeta+1$.
  Then as in the previous part
  \begin{align*}
    \sup_{w \in \cW}
    \big|
    \E[\widehat f_W(w)]
    - f_W(w)
    \big|
    &=
    \sup_{w \in \cW}
    \left|
    \int_{[w \pm h] \cap \cW}
    k_h(s,w)
    \big(
    f_W(s) - P_{\flbeta}(s,w)
    \big)
    \diff{s}
    \right|
    \leq
    \frac{2 \Ck \CH}{\flbeta !}
    h^\beta.
  \end{align*}
\end{proof}

\subsubsection{Uniform consistency}

\begin{proof}[Lemma~\ref{lem:hoeffding}]

  \proofparagraph{Hoeffding-type decomposition}

  Note that
  \begin{align*}
    \widehat f_W(w)
    - E_n(w)
    - \E[\widehat f_W(w)]
    &=
    \frac{2}{n(n-1)}
    \sum_{i=1}^{n-1}
    \sum_{j=i+1}^{n}
    \Big(
    \E[k_h(W_{ij},w) \mid A_i, A_j]
    - \E[k_h(W_{ij},w)]
    \Big) \\
    &=
    \frac{1}{n(n-1)}
    \sum_{i=1}^{n-1}
    \sum_{j \neq i}
    \Big(
    \E[k_h(W_{ij},w) \mid A_i, A_j]
    - \E[k_h(W_{ij},w)]
    \Big)
  \end{align*}
  and apply Lemma~\ref{lem:general_hoeffding}
  with
  \begin{align*}
    u_{ij}
    &=
    \frac{1}{n(n-1)}
    \E\big[k_h(W_{ij},w) \mid A_i, A_j\big],
    &u_i
    &=
    \frac{1}{n(n-1)}
    \E\big[k_h(W_{ij},w) \mid A_i\big],
    &u
    &=
    \frac{1}{n(n-1)}
    \E\big[k_h(W_{ij},w)\big],
  \end{align*}
  to see
  \begin{align*}
    \widehat f_W(w)
    - E_n(w)
    - \E[\widehat f_W(w)]
    &=
    \frac{2}{n}
    \sum_{i=1}^n
    \big(u_i - u\big)
    + \frac{1}{n(n-1)}
    \sum_{i=1}^n
    \sum_{j \neq i}
    \big(
    u_{ij} - u_i - u_j + u
    \big) \\
    &=
    \frac{2}{n}
    \sum_{i=1}^n
    l_i(w)
    + \frac{2}{n(n-1)}
    \sum_{i=1}^n
    \sum_{j = i+1}^n
    q_{ij}(w)
    =
    L_n + Q_n.
  \end{align*}

  \proofparagraph{Expectation and covariance of $L_n$, $Q_n$ and $E_n$}

  Observe that
  $L_n$, $Q_n$ and $E_n$
  are clearly mean-zero.
  For orthogonality,
  note that their summands
  have the following properties,
  for any $1 \leq i < j \leq n$
  and $1 \leq r < s \leq n$,
  and for any $w, w' \in \cW$:
  \begin{align*}
    \E\big[
      l_i(w)
      q_{rs}(w')
      \big]
    &=
    \E\big[
      l_i(w)
      \E\big[
        q_{rs}(w') \mid A_i
        \big]
      \big]
    = 0, \\
    \E\big[
      l_i(w)
      e_{rs}(w')
      \big]
    &=
    \begin{cases}
      \E\big[
        l_i(w)
        \big]
      \E\big[
        e_{rs}(w')
        \big],
      \text{ if } i \notin \{r,s\} \\
      \E\big[
        l_i(w)
        \E\big[
          e_{rs}(w') \mid A_r, A_s
          \big]
        \big],
      \text{ if } i \in \{r,s\}
    \end{cases} \\
    &=
    0, \\
    \E\big[
      q_{ij}(w)
      e_{rs}(w')
      \big]
    &=
    \begin{cases}
      \E\big[
        q_{ij}(w)
        \big]
      \E\big[
        e_{rs}(w')
        \big],
      \text{ if } \{i,j\} \cap \{r,s\} = \emptyset \\
      \E\big[
        \E\big[
          q_{ij}(w) \mid A_i
          \big]
        \E\big[
          e_{rs}(w') \mid A_i
          \big]
        \big],
      \text{ if } \{i,j\} \cap \{r,s\} = \{i\}     \\
      \E\big[
        \E\big[
          q_{ij}(w) \mid A_j
          \big]
        \E\big[
          e_{rs}(w') \mid A_j
          \big]
        \big],
      \text{ if } \{i,j\} \cap \{r,s\} = \{j\}     \\
      \E\big[
        q_{ij}(w)
        \E\big[
          e_{rs}(w') \mid A_r, A_s
          \big]
        \big],
      \text{ if } \{i,j\} = \{r,s\}
    \end{cases} \\
    &=
    0,
  \end{align*}
  where we used mutual independence of
  $\bA_n$ and $\bV_n$
  and also
  $\E[q_{rs}(w) \mid A_i] = 0$
  and
  $\E[e_{ij}(w) \mid A_i, A_j] = 0$.
\end{proof}

\begin{proof}[Lemma~\ref{lem:trichotomy}]

  \proofparagraph{total degeneracy}

  Suppose
  $\Dl = 0$, so
  $\Var[f_{W \mid A}(w \mid A_i)] = 0$
  for all $w \in \cW$.
  Therefore for all $w \in \cW$,
  we have $f_{W \mid A}(w) = f_W(w)$ almost surely.
  By taking a union over $\cW \cap \Q$
  and by continuity of $f_{W \mid A}$ and $f_W$,
  this implies that $f_{W \mid A}(w) = f_W(w)$
  for all $w \in \cW$
  almost surely. Thus
  \begin{align*}
    \E\left[
      k_h(W_{ij},w) \mid A_i
    \right]
    &=
    \int_{\cW}
    k_h(s,w)
    f_{W \mid A}(s \mid A_i)
    \diff{s}
    =
    \int_{\cW}
    k_h(s,w)
    f_W(s)
    \diff{s}
    =
    \E\left[
      k_h(W_{ij},w)
    \right]
  \end{align*}
  for all $w \in \cW$
  almost surely.
  Hence $l_i(w) = 0$
  and therefore $L_n(w) = 0$
  for all $w \in \cW$ almost surely.

  \proofparagraph{no degeneracy}

  Suppose that $\Dl > 0$.
  Now since $f_{W|A}(\cdot \mid a)$ is $\CH$-Lipschitz
  for all $a \in \cA$
  and since $|k_h| \leq \Ck/h$,
  \begin{align*}
    \sup_{w \in \cW}
    \left|
    \E[k_h(W_{ij},w) \mid A_i]
    - f_{W \mid A}(w \mid A_i)
    \right|
    &=
    \sup_{w \in \cW}
    \left|
    \int_{\cW}
    k_h(s,w)
    f_{W \mid A}(s \mid A_i)
    \diff{s}
    - f_{W \mid A}(w \mid A_i)
    \right| \\
    &=
    \sup_{w \in \cW}
    \left|
    \int_{\cW \cap [w \pm h]}
    k_h(s,w)
    \left(
      f_{W \mid A}(s \mid A_i)
      - f_{W \mid A}(w \mid A_i)
    \right)
    \diff{s}
    \right| \\
    &\leq
    2h
    \frac{\Ck}{h}
    \CH h
    \leq
    2 \Ck \CH h
  \end{align*}
  almost surely.
  Therefore since $f_{W \mid A}(w \mid a) \leq \Cd$,
  we have
  \begin{align*}
    \sup_{w \in \cW}
    \left|
    \Var\big[
      \E[k_h(W_{ij},w) \mid A_i]
      \big]
    - \Var\left[
      f_{W \mid A}(w \mid A_i)
    \right]
    \right|
    &\leq
    16 \Ck \CH \Cd h
  \end{align*}
  whenever $h$ is small enough that
  $2 \Ck \CH h \leq \Cd$. Thus
  \begin{align*}
    \inf_{w \in \cW} \Var\big[\E[k_h(W_{ij},w) \mid A_i]\big]
    &\geq
    \inf_{w \in \cW}\Var[f_{W \mid A}(w \mid A_i)]
    - 16 \Ck \CH \Cd h.
  \end{align*}
  Therefore if
  $\Dl > 0$
  then eventually
  $\inf_{w \in \cW} \Var\big[\E[k_h(W_{ij},w) \mid A_i]\big] \geq \Dl/2$.
  Finally note that
  \begin{align*}
    \inf_{w \in \cW}\Var[L_n(w)]
    &=
    \frac{4}{n}
    \inf_{w \in \cW}
    \Var\big[\E[k_h(W_{ij},w) \mid A_i]\big]
    \geq
    \frac{2 \Dl}{n}.
  \end{align*}

  \proofparagraph{partial degeneracy}

  Since $f_{W \mid A}(w \mid A_i)$
  is bounded by $\Cd$ and $\CH$-Lipschitz in $w$,
  we must have that
  $\Var[f_{W \mid A}(w \mid A_i)]$
  is a continuous function on $\cW$.
  Thus if $\Dl = 0$,
  there must be at least one point $w \in \cW$
  for which
  $\Var[f_{W \mid A}(w \mid A_i)] = 0$
  by compactness.
  Let $w$ be any such degenerate point.
  Then by the previous part,
  \begin{align*}
    \Var[L_n(w)] =
    \frac{4}{n} \Var\big[\E[k_h(W_{ij},w) \mid A_i]\big]
    &\leq
    64 \Ck \CH \Cd \frac{h}{n}.
  \end{align*}
  If conversely $w$ is not a degenerate point
  then
  $\Var[f_{W \mid A}(w \mid A_i)] > 0$
  so eventually
  \begin{align*}
    \Var[L_n(w)]
    = \frac{4}{n}
    \Var\big[\E[k_h(W_{ij},w) \mid A_i]\big]
    &\geq
    \frac{2}{n}
    \Var[f_{W \mid A}(w \mid A_i)].
  \end{align*}
\end{proof}

\begin{proof}[Lemma~\ref{lem:uniform_concentration}]

  We establish VC-type properties of some function
  classes and apply results from empirical process theory.

  \proofparagraph{establishing VC-type classes}

  Consider the following function classes:
  \begin{align*}
    \cF_1
    &=
    \Big\{
    W_{ij} \mapsto
    k_h(W_{ij},w)
    : w \in \cW
    \Big\}, \\
    \cF_2
    &=
    \Big\{
    (A_i, A_j) \mapsto
    \E\big[ k_h(W_{ij},w) \mid A_i, A_j \big]
    : w \in \cW
    \Big\}, \\
    \cF_3
    &=
    \Big\{
    A_i \mapsto
    \E\big[ k_h(W_{ij},w) \mid A_i \big]
    : w \in \cW
    \Big\}.
  \end{align*}
  For $\cF_1$,
  take $0 < \varepsilon \leq \Leb(\cW)$
  and let $\cW_\varepsilon$ be an
  $\varepsilon$-cover of $\cW$
  of cardinality at most $\Leb(\cW)/\varepsilon$.
  Since
  \begin{align*}
    \sup_{s, w, w' \in \cW}
    \left|
    \frac{k_h(s,w) - k_h(s,w')}
    {w-w'}
    \right|
    &\leq
    \frac{C_\mathrm{L}}{h^2}
  \end{align*}
  almost surely,
  we see that
  \begin{align*}
    \sup_\Q
    N\left(\cF_1, \rho_\Q,
      \frac{C_\mathrm{L}}{h^2} \varepsilon \right)
    &\leq
    N\left(\cF_1, \|\cdot\|_\infty,
      \frac{C_\mathrm{L}}{h^2} \varepsilon \right)
    \leq
    \frac{\Leb(\cW)}{\varepsilon},
  \end{align*}
  where $\Q$ ranges over Borel
  probability measures on $\cW$.
  Since
  $\frac{\Ck}{h}$
  is an envelope function for $\cF_1$,
  \begin{align*}
    \sup_\Q
    N\left(\cF_1, \rho_\Q,
      \frac{\Ck}{h} \varepsilon \right)
    &\leq
    \frac{\CL}{\Ck}
    \frac{\Leb(\cW)}{h \varepsilon}.
  \end{align*}
  Thus for all $\varepsilon \in (0,1]$,
  \begin{align*}
    \sup_\Q
    N\left(\cF_1, \rho_\Q,
      \frac{\Ck}{h} \varepsilon \right)
    &\leq
    \frac{\CL}{\Ck}
    \frac{\Leb(\cW) \vee 1}{h \varepsilon}
    \leq
    (C_1/(h\varepsilon))^{C_2},
  \end{align*}
  where
  $C_1 = \frac{\CL}{\Ck} (\Leb(\cW) \vee 1)$
  and $C_2 = 1$.
  Next, $\cF_2$ forms a smoothly parametrized class of functions
  since for $w,w' \in \cW$ we have
  by the uniform Lipschitz properties of
  $f_{W \mid AA}(\cdot \mid A_i, A_j)$ and
  $k_h(s, \cdot)$,
  with $|w-w'| \leq h$,
  \begin{align*}
    &\left|
    \E\big[ k_h(W_{ij},w) \mid A_i, A_j \big]
    - \E\big[ k_h(W_{ij},w') \mid A_i, A_j \big]
    \right| \\
    &\quad=
    \left|
    \int_{[w \pm h] \cap \cW}
    k_h(s,w)
    f_{W \mid AA}(s \mid A_i, A_j)
    \diff{s}
    - \int_{[w' \pm h] \cap \cW}
    k_h(s,w')
    f_{W \mid AA}(s \mid A_i, A_j)
    \diff{s}
    \right| \\
    &\quad=
    \left|
    \int_{[w \pm 2h] \cap \cW}
    \big(
    k_h(s,w)
    - k_h(s,w')
    \big)
    f_{W \mid AA}(s \mid A_i, A_j)
    \diff{s}
    \right| \\
    &\quad=
    \left|
    \int_{[w \pm 2h] \cap \cW}
    \big(
    k_h(s,w)
    - k_h(s,w')
    \big)
    \big(
    f_{W \mid AA}(s \mid A_i, A_j)
    - f_{W \mid AA}(w \mid A_i, A_j)
    \big)
    \diff{s}
    \right| \\
    &\quad\leq
    4h
    \frac{\CL}{h^2}
    |w-w'|
    2 \CH h
    \leq
    8 \CL \CH
    |w-w'|
    \leq
    C_3
    |w-w'|,
  \end{align*}
  where $C_3 = 8 \CL \CH$.
  The same holds for $|w-w'| > h$
  because Lipschitzness is a local property.
  By taking $\E[\, \cdot \mid A_i]$,
  it can be seen
  by the contraction property of conditional expectation that
  the same holds for the
  singly-conditioned terms:
  \begin{align*}
    \left|
    \E\big[ k_h(W_{ij},w) \mid A_i \big]
    - \E\big[ k_h(W_{ij},w') \mid A_i \big]
    \right|
    &\leq
    C_3
    |w-w'|.
  \end{align*}
  Therefore $\cF_3$ is also smoothly parametrized
  in exactly the same manner.
  Let
  \begin{align*}
    C_4
    &=
    \sup_{w \in \cW}
    \esssup_{A_i, A_j}
    \big|
    \E\big[ k_h(W_{ij},w) \mid A_i, A_j \big]
    \big|
    =
    \sup_{w \in \cW}
    \esssup_{A_i, A_j}
    \left|
    \int_{[w \pm h] \cap \cW}
    k_h(s,w)
    f_{W \mid AA}(s \mid A_i, A_j)
    \diff{s}
    \right| \\
    &\leq 2h \frac{\Ck}{h} \Cd
    \leq 2 \Ck \Cd.
  \end{align*}
  For any $\varepsilon \in (0,1]$,
  take an $(\varepsilon C_4/C_3)$-cover of $\cW$
  of cardinality at most $C_3 \Leb(\cW) / (\varepsilon C_4)$.
  By the above smooth parametrization properties,
  this cover induces an
  $\varepsilon C_4$-cover for both $\cF_2$ and $\cF_3$:
  \begin{align*}
    \sup_\Q
    N\big(\cF_2, \rho_\Q, \varepsilon C_4 \big)
    &\leq
    N\big(\cF_2, \|\cdot\|_\infty, \varepsilon C_4 \big)
    \leq
    C_3 \Leb(\cW) / (\varepsilon C_4), \\
    \sup_\Q
    N\big(\cF_3, \rho_\Q, \varepsilon C_4 \big)
    &\leq
    N\big(\cF_3, \|\cdot\|_\infty, \varepsilon C_4 \big)
    \leq
    C_3 \Leb(\cW) / (\varepsilon C_4).
  \end{align*}
  Hence
  $\cF_1$,
  $\cF_2$ and
  $\cF_3$
  form VC-type classes with envelopes
  $F_1 = \Ck / h$ and
  $F_2 = F_3 = C_4$ respectively:
  \begin{align*}
    \sup_\Q
    N\left(\cF_1, \rho_\Q,
      \varepsilon \Ck / h \right)
    &\leq
    (C_1/(h\varepsilon))^{C_2},
    &\sup_\Q
    N\big(\cF_2, \rho_\Q, \varepsilon C_4 \big)
    &\leq
    (C_1/\varepsilon)^{C_2}, \\
    \sup_\Q
    N\big(\cF_3, \rho_\Q, \varepsilon C_4 \big)
    &\leq
    (C_1/\varepsilon)^{C_2},
  \end{align*}
  for some constants $C_1 \geq e$ and $C_2 \geq 1$,
  where we augment the constants if necessary.

  \proofparagraph{controlling $L_n$}

  Observe that
  $\sqrt{n}L_n$
  is the empirical process of the i.i.d.\ variables $A_i$
  indexed by $\cF_3$.
  We apply Lemma~\ref{lem:maximal_vc_inid}
  with $\sigma = C_4$:
  \begin{align*}
    \E \left[
      \sup_{w \in \cW}
      \big| \sqrt{n} L_
      n(w) \big|
    \right]
    &\lesssim
    C_4
    \sqrt{C_2 \log C_1}
    +
    \frac{C_4 C_2 \log C_1}
    {\sqrt{n}}
    \lesssim 1.
  \end{align*}
  By Lemma~\ref{lem:trichotomy},
  the left hand side is zero whenever
  $\Du = 0$,
  so we can also write
  \begin{align*}
    \E \left[
      \sup_{w \in \cW}
      \big| \sqrt{n} L_n(w) \big|
    \right]
    &\lesssim
    \Du.
  \end{align*}

  \proofparagraph{controlling $Q_n$}

  Observe that $n Q_n$
  is the completely degenerate second-order U-process
  of the i.i.d.\ variables $A_i$
  indexed by $\cF_2$.
  This function class is again uniformly bounded
  and VC-type,
  so applying the U-process maximal inequality from
  Lemma~\ref{lem:uprocess_maximal}
  yields with $\sigma = C_4$
  \begin{align*}
    \E \left[
      \sup_{w \in \cW}
      \big| n Q_n(w) \big|
    \right]
    &\lesssim
    C_4
    C_2 \log C_1
    +
    \frac{C_4 (C_2 \log C_1)^2}
    {\sqrt{n}}
    \lesssim 1.
  \end{align*}

  \proofparagraph{controlling $E_n$}

  Conditional on $\bA_n$,
  note that $n E_n$
  is the empirical process of the conditionally
  i.n.i.d.\ variables $W_{ij}$
  indexed by $\cF_1$.
  We apply Lemma~\ref{lem:maximal_vc_inid}
  conditionally with
  \begin{align*}
    \sigma^2
    &=
    \sup_{w \in \cW}
    \E\Big[
      \big(
      k_h(W_{ij},w)
      - \E[k_h(W_{ij},w) \mid A_i, A_j]
      \big)^2
      \mid A_i, A_j
      \Big]
    \leq
    \sup_{w \in \cW}
    \E\Big[
      k_h(W_{ij},w)^2
      \mid A_i, A_j
      \Big] \\
    &\leq
    \sup_{w \in \cW}
    \int_{[w \pm h] \cap \cW}
    k_h(s,w)^2
    f_{W \mid AA}(s \mid A_i, A_j)
    \diff{s}
    \leq 2h \frac{\Ck^2}{h^2}
    \lesssim 1/h
  \end{align*}
  and noting that we have
  a sample size of
  $\frac{1}{2}n(n-1)$,
  giving
  \begin{align*}
    \E \left[
      \sup_{w \in \cW}
      \big| n E_n(w) \big|
    \right]
    &\lesssim
    \sigma
    \sqrt{C_2 \log \big((C_1/h) F_1 / \sigma \big)}
    +
    \frac{F_1 C_2 \log \big((C_1/h) F_1 / \sigma\big)}
    {n} \\
    &\lesssim
    \frac{1}{\sqrt h}
    \sqrt{C_2 \log \big((C_1/h) (\Ck/h) \sqrt h \big)}
    +
    \frac{(\Ck/h) C_2 \log \big((C_1/h) (\Ck/h) \sqrt h \big)}
    {n} \\
    &\lesssim
    \sqrt{\frac{\log 1/h}{h}}
    +
    \frac{\log \big(1/h\big)}
    {nh}
    \lesssim
    \sqrt{\frac{\log n}{h}},
  \end{align*}
  where the last line follows by the bandwidth assumption
  of $\frac{\log n}{n^2h} \to 0$.
\end{proof}

\begin{proof}[Theorem~\ref{thm:uniform_consistency}]

  This follows from Lemma~\ref{lem:bias}
  and Lemma~\ref{lem:uniform_concentration}.
\end{proof}

\subsubsection{Minimax optimality}

Before proving Theorem~\ref{thm:minimax}
we first give a lower bound result
for parametric point estimation in
Lemma~\ref{lem:neyman_pearson_bernoulli}.

\begin{lemma}[A Neyman--Pearson result for Bernoulli random variables]
  \label{lem:neyman_pearson_bernoulli}

  Recall that the Bernoulli distribution
  $\Ber(\theta)$
  places mass $\theta$ at $1$ and mass
  $1-\theta$ at $0$.
  Define $\P_\theta^n$ as the law of
  $(A_1, A_2, \ldots, A_n, V)$,
  where $A_1, \ldots, A_n$
  are i.i.d.\ $\Ber(\theta)$,
  and $V$ is an $\R^d$-valued random variable
  for some $d \geq 1$
  which is independent of the $A$ variables
  and with a fixed distribution that does not depend on $\theta$.
  Let $\theta_0 = \frac{1}{2}$
  and $\theta_{1,n} = \frac{1}{2} + \frac{1}{\sqrt{8n}}$.
  Then for any estimator $\widetilde \theta_n$
  which is a function of
  $(A_1, A_2, \ldots, A_n, V)$ only,
  \begin{align*}
    \P_{\theta_0}^n \left(
      \big| \widetilde \theta_n - \theta_0 \big|
      \geq \frac{1}{\sqrt{32n}}
    \right)
    + \P_{\theta_{1,n}}^n \left(
      \big| \widetilde \theta_n - \theta_{1,n} \big|
      \geq \frac{1}{\sqrt{32n}}
    \right)
    \geq \frac{1}{2}.
  \end{align*}

\end{lemma}

\begin{proof}[Lemma~\ref{lem:neyman_pearson_bernoulli}]

  Let $f: \{0,1\}^n \to \{0,1\}$
  be any function.
  Considering this function as a statistical test,
  the Neyman--Pearson lemma
  and Pinsker's inequality
  \citep{gine2021mathematical}
  give
  \begin{align*}
    \P_{\theta_0}^n \big(
    f=1
    \big)
    +\P_{\theta_{1,n}}^n \big(
    f=0
    \big)
    &\geq
    1-
    \TV\left(
      \P_{\theta_0}^n,
      \P_{\theta_{1,n}}^n
    \right)
    \geq
    1-
    \sqrt{
      \frac{1}{2}
      \KL \left(
        \P_{\theta_0}^n
        \bigm\|
        \P_{\theta_{1,n}}^n
      \right)} \\
    &=
    1-
    \sqrt{
      \frac{n}{2}
      \KL \left(
        \Ber(\theta_0)
        \bigm\|
        \Ber(\theta_{1,n})
      \right)
      + \frac{n}{2}
      \KL \left(
        V
        \bigm\|
        V
      \right)} \\
    &=
    1-
    \sqrt{
      \frac{n}{2}
      \KL \left(
        \Ber(\theta_0)
        \bigm\|
        \Ber(\theta_{1,n})
      \right)},
  \end{align*}
  where $\TV$ is the total variation distance
  and $\KL$ is the Kullback--Leibler divergence.
  In the penultimate line
  we used the tensorization of Kullback--Leibler divergence
  \citep{gine2021mathematical},
  noting that the law of $V$ is fixed and hence does not contribute.
  We now evaluate this Kullback--Leibler divergence at the specified
  parameter values.
  \begin{align*}
    \P_{\theta_0}^n \big(
    f=1
    \big)
    +\P_{\theta_{1,n}}^n \big(
    f=0
    \big)
    &\geq
    1-
    \sqrt{
      \frac{n}{2}
      \KL \left(
        \Ber(\theta_0)
        \bigm\|
        \Ber(\theta_{1,n})
      \right)} \\
    &=
    1-
    \sqrt{\frac{n}{2}}
    \sqrt{
      \theta_0 \log \frac{\theta_0}{\theta_{1,n}}
      + (1 - \theta_0) \log \frac{1 - \theta_0}{1 - \theta_{1,n}}} \\
    &=
    1-
    \sqrt{\frac{n}{2}}
    \sqrt{
      \frac{1}{2} \log \frac{1/2}{1/2 + 1/\sqrt{8n}}
      + \frac{1}{2} \log \frac{1/2}{1/2 - 1/\sqrt{8n}}} \\
    &=
    1-
    \frac{\sqrt n}{2}
    \sqrt{\log \frac{1}{1 - 1/(2n)}}
    \geq
    1-
    \frac{\sqrt n}{2}
    \sqrt{\frac{1}{n}}
    =
    \frac{1}{2},
  \end{align*}
  where in the penultimate line we used that
  $\log \frac{1}{1-x} \leq 2x$
  for $x \in [0,1/2]$.
  Now define a test $f$ by
  $f = 1$ if $\widetilde \theta_n > \frac{1}{2} + \frac{1}{\sqrt{32n}}$
  and $f=0$ otherwise,
  to see
  \begin{align*}
    \P_{\theta_0}^n \left(
      \widetilde \theta_n > \frac{1}{2} + \frac{1}{\sqrt{32n}}
    \right)
    + \P_{\theta_{1,n}}^n \left(
      \widetilde \theta_n \leq \frac{1}{2} + \frac{1}{\sqrt{32n}}
    \right)
    \geq \frac{1}{2}.
  \end{align*}
  By the triangle inequality,
  recalling that
  $\theta_0 = \frac{1}{2}$
  and $\theta_{1,n} = \frac{1}{2} + \frac{1}{\sqrt{8n}}$,
  we have the event inclusions
  \begin{align*}
    \left\{
    \widetilde \theta_n > \frac{1}{2} + \frac{1}{\sqrt{32n}}
    \right\}
    &\subseteq
    \left\{
    \left| \widetilde \theta_n - \theta_0 \right|
    \geq \frac{1}{\sqrt{32n}}
    \right\} \\
    \left\{
    \widetilde \theta_n \leq \frac{1}{2} + \frac{1}{\sqrt{32n}}
    \right\}
    &\subseteq
    \left\{
    \left| \widetilde \theta_n - \theta_{1,n} \right|
    \geq \frac{1}{\sqrt{32n}}
    \right\}.
  \end{align*}
  Thus by the monotonicity of measures,
  \begin{align*}
    \P_{\theta_0}^n \left(
      \big| \widetilde \theta_n - \theta_0 \big|
      \geq \frac{1}{\sqrt{32n}}
    \right)
    + \P_{\theta_{1,n}}^n \left(
      \big| \widetilde \theta_n - \theta_{1,n} \big|
      \geq \frac{1}{\sqrt{32n}}
    \right)
    \geq \frac{1}{2}.
  \end{align*}
\end{proof}

\begin{proof}[Theorem~\ref{thm:minimax}]

  \proofparagraph{lower bound for $\cP$}

  By translation and scaling of the data,
  we may assume without loss of generality that $\cW = [-1,1]$.
  We may also assume that $\CH \leq 1/2$,
  since reducing $\CH$ can only shrink the class of distributions.
  Define the dyadic distribution $\P_\theta$
  with parameter $\theta \in [1/2, 1]$
  as follows:
  $A_1, \ldots, A_n$ are i.i.d.\ $\Ber(\theta)$, while
  $V_{ij}$ for $1 \leq i < j \leq n$ are i.i.d.\
  and independent of $\bA_n$.
  The distribution of $V_{ij}$ is given by its density function
  $f_V(v) = \frac{1}{2} + \CH v$ on $[-1,1]$.
  Finally generate
  $W_{ij} = W(A_i, A_j, V_{ij}) \vcentcolon=
  (2 A_i A_j - 1) V_{ij}$.
  Note that the function $W$ does not depend on $\theta$.
  The conditional and marginal densities of $W_{ij}$ are
  for $w \in [-1,1]$
  \begin{align*}
    f_{W \mid AA}(w \mid A_i, A_j)
    &=
    \begin{cases}
      \frac{1}{2} + \CH w & \text{if } A_i = A_j = 1                \\
      \frac{1}{2} - \CH w & \text{if } A_i = 0 \text{ or } A_j = 0, \\
    \end{cases} \\
    f_{W \mid A}(w \mid A_i)
    &=
    \begin{cases}
      \frac{1}{2} + (2 \theta - 1) \CH w
      & \text{if } A_i = 1   \\
      \frac{1}{2} - \CH w & \text{if } A_i = 0 , \\
    \end{cases} \\
    f_W(w)&= \frac{1}{2} + (2\theta^2 - 1) \CH w.
  \end{align*}
  Clearly
  $f_W \in \cH^\beta_{\CH}(\cW)$ and
  $f_{W \mid AA}(\cdot \mid a, a') \in \cH^1_{\CH}(\cW)$.
  Also
  $\sup_{w \in \cW} \|f_{W \mid A}(w \mid \cdot\,)\|_\TV \leq 1$.
  Therefore
  $\P_\theta$ satisfies Assumption~\ref{ass:data}
  and so
  $\big\{\P_\theta : \theta \in [1/2, 1] \big\} \subseteq \cP$.

  Note that $f_W(1) = \frac{1}{2} + (2\theta^2 - 1) \CH $,
  so $\theta^2 = \frac{1}{2 \CH}(f_W(1) - 1/2 + \CH)$.
  Thus if $\widetilde f_W$ is some
  density estimator
  depending only on the data $\bW_n$,
  we can naturally define
  the non-negative parameter point estimator
  \begin{align*}
    \widetilde \theta_n^2
    &\vcentcolon=
    \frac{1}{2 \CH}\left(
      \widetilde f_W(1) - \frac{1}{2} + \CH
    \right)
    \vee 0.
  \end{align*}
  This gives the inequality
  \begin{align*}
    \big|
    \widetilde \theta_n^2 - \theta^2
    \big|
    &=
    \left|
    \frac{1}{2 \CH}\left(
      \widetilde f_W(1) - \frac{1}{2} + \CH
    \right)
    \vee 0
    -
    \frac{1}{2 \CH}\left(
      f_W(1) - \frac{1}{2} + \CH
    \right)
    \right|
    \leq
    \frac{1}{2 \CH}
    \sup_{w \in \cW}
    \left|
    \widetilde f_W(w)  - f_W(w)
    \right|.
  \end{align*}
  Therefore since also $\widetilde \theta \geq 0$
  and
  $\theta \geq \frac{1}{2}$,
  \begin{align*}
    \big|
    \widetilde \theta_n - \theta
    \big|
    &=
    \frac{\big|\widetilde \theta_n^2 - \theta^2\big|}
    {\widetilde \theta_n + \theta}
    \leq
    \frac{1}{\CH}
    \sup_{w \in \cW}
    \left|
    \widetilde f_W(w)  - f_W(w)
    \right|.
  \end{align*}
  Now we apply the point estimation lower bound from
  Lemma~\ref{lem:neyman_pearson_bernoulli},
  setting $\theta_0 = \frac{1}{2}$
  and $\theta_{1,n} = \frac{1}{2} + \frac{1}{\sqrt{8n}}$,
  noting that the estimator
  $\widetilde \theta_n$
  is a function of $\bW_n$ only,
  thus is a function of $\bA_n$ and
  $\bV_n$ only and so satisfies the conditions.
  \begin{align*}
    &\P_{\theta_0} \left(
      \sup_{w \in \cW} \big| \widetilde f_W(w) - f^{(0)}_W(w) \big|
      \geq \frac{1}{C\sqrt{n}}
    \right)
    + \P_{\theta_{1,n}} \left(
      \sup_{w \in \cW} \big| \widetilde f_W(w) - f^{(1)}_W(w) \big|
      \geq \frac{1}{C\sqrt{n}}
    \right) \\
    &\quad\geq
    \P_{\theta_0} \left(
      \big| \widetilde \theta_n - \theta_0 \big|
      \geq \frac{1}{C \CH \sqrt{n}}
    \right)
    + \P_{\theta_{1,n}} \left(
      \big| \widetilde \theta_n - \theta_{1,n} \big|
      \geq \frac{1}{C \CH \sqrt{n}}
    \right) \\
    &\quad\geq
    \P_{\theta_0} \left(
      \big| \widetilde \theta_n - \theta_0 \big|
      \geq \frac{1}{\sqrt{32n}}
    \right)
    + \P_{\theta_{1,n}} \left(
      \big| \widetilde \theta_n - \theta_{1,n} \big|
      \geq \frac{1}{\sqrt{32n}}
    \right)
    \geq
    \frac{1}{2},
  \end{align*}
  where we set $C \geq \frac{\sqrt{32}}{\CH}$.
  Therefore we deduce that
  \begin{align*}
    \inf_{\widetilde f_W}
    \sup_{\P \in \cP}
    \P\left(
      \sup_{w \in \cW}
      \big|
      \widetilde f_W(w) - f_W(w)
      \big|
      \geq
      \frac{1}{C \sqrt n}
    \right)
    \geq \frac{1}{4}
  \end{align*}
  and so
  \begin{align*}
    \inf_{\widetilde f_W}
    \sup_{\P \in \cP}
    \E_\P\left[
      \sup_{w \in \cW}
      \big|
      \widetilde f_W(w) - f_W(w)
      \big|
    \right]
    \geq \frac{1}{4 C \sqrt{n}}.
  \end{align*}

  \proofparagraph{lower bound for $\cPd$}

  For the subclass of totally degenerate distributions,
  we rely on the main theorem
  from \citet{khasminskii1978lower}.
  Let $\cP_0$ be the subclass of $\cPd$
  consisting of the distributions which satisfy
  $A_1 = \cdots = A_n = 0$
  and $W_{ij} \vcentcolon= A_i + A_j + V_{ij} = V_{ij}$,
  so that $W_{ij}$ are i.i.d.\ with common density $f_W = f_V$.
  Define the class
  \begin{align*}
    \cF
    &=
    \left\{
    f \text{ density function on } \R, \
    f \in \cH^\beta_{\CH}(\cW)
    \right\}.
  \end{align*}
  Write $\E_f$ for the expectation under $W_{ij}$ having density $f$.
  Then by the main theorem in \citet{khasminskii1978lower},
  \begin{align*}
    \liminf_{n \to \infty}
    \inf_{\widetilde f_W}
    \sup_{f \in \cF}
    \E_f\left[
      \left( \frac{n^2}{\log n} \right)^{\frac{\beta}{2\beta + 1}}
      \sup_{w \in \cW}
      \big| \widetilde f_W(w) - f_W(w) \big|
    \right]
    > 0,
  \end{align*}
  where $\widetilde f_W$ is any
  density estimator
  depending only on the $\frac{1}{2}n(n-1)$ i.i.d.\ data samples $\bW_n$.
  Now every density function in
  $\cH^\beta_{\CH}(\cW)$
  corresponds to a distribution in
  $\cP_0$ and therefore to a distribution in $\cPd$.
  Thus for large enough $n$ and
  some positive constant $C$,
  \begin{align*}
    \inf_{\widetilde f_W}
    \sup_{\P \in \cPd}
    \E_\P\left[
      \sup_{w \in \cW}
      \big| \widetilde f_W(w) - f_W(w) \big|
    \right]
    \geq
    \frac{1}{C}
    \left( \frac{\log n}{n^2} \right)^{\frac{\beta}{2\beta + 1}}.
  \end{align*}

  \proofparagraph{upper bounds}

  The corresponding upper bounds follow by
  using a dyadic kernel density estimator $\widehat f_W$
  with a boundary bias-corrected
  Lipschitz kernel of order $p \geq \beta$
  and using a bandwidth of $h$.
  Firstly Lemma~\ref{lem:bias} gives
  \begin{align*}
    \sup_{\P \in \cP}
    \sup_{w \in \cW}
    \big|
    \E_\P\big[\widehat f_W(w)\big]
    - f_W(w)
    \big|
    \leq
    \frac{4\Ck \CH}{\flbeta !}
    h^\beta.
  \end{align*}
  Then,
  treating the degenerate and non-degenerate cases separately
  and noting that all inequalities hold uniformly over
  $\cP$ and $\cPd$,
  the proof of Lemma~\ref{lem:uniform_concentration}
  shows that
  \begin{align*}
    \sup_{\P \in \cP}
    \E_\P\left[
      \sup_{w \in \cW}
      \big|\widehat f_W(w) - \E_\P[\widehat f_W(w)]\big|
    \right]
    &\lesssim
    \frac{1}{\sqrt n}
    + \sqrt{\frac{\log n}{n^2h}}, \\
    \sup_{\P \in \cPd}
    \E_\P\left[
      \sup_{w \in \cW}
      \big|\widehat f_W(w) - \E_\P[\widehat f_W(w)]\big|
    \right]
    &\lesssim
    \sqrt{\frac{\log n}{n^2h}}.
  \end{align*}
  Thus combining these yields that
  \begin{align*}
    \sup_{\P \in \cP}
    \E_\P\left[
      \sup_{w \in \cW}
      \big|\widehat f_W(w) - f_W(w)\big|
    \right]
    &\lesssim
    h^\beta
    + \frac{1}{\sqrt n}
    + \sqrt{\frac{\log n}{n^2h}}, \\
    \sup_{\P \in \cPd}
    \E_\P\left[
      \sup_{w \in \cW}
      \big|\widehat f_W(w) - f_W(w)\big|
    \right]
    &\lesssim
    h^\beta
    + \sqrt{\frac{\log n}{n^2h}}.
  \end{align*}
  Set $h = \left( \frac{\log n}{n^2} \right)^{\frac{1}{2\beta+1}}$
  and note that $\beta \geq 1$ implies that
  $\left(\frac{\log n}{n^2} \right)^{\frac{\beta}{2\beta+1}}
  \ll \frac{1}{\sqrt n}$.
  Therefore for some constant $C > 0$,
  \begin{align*}
    \sup_{\P \in \cP}
    \E_\P\left[
      \sup_{w \in \cW}
      \big|\widehat f_W(w) - f_W(w)\big|
    \right]
    &\lesssim
    \frac{1}{\sqrt n}
    + \left(
      \frac{\log n}{n^2}
    \right)^{\frac{\beta}{2\beta+1}}
    \leq
    \frac{C}{\sqrt n}, \\
    \sup_{\P \in \cPd}
    \E_\P\left[
      \sup_{w \in \cW}
      \big|\widehat f_W(w) - f_W(w)\big|
    \right]
    &\leq
    C\left(
      \frac{\log n}{n^2}
    \right)^{\frac{\beta}{2\beta+1}}.
  \end{align*}
\end{proof}

\subsubsection{Covariance structure}

\begin{proof}[Lemma~\ref{lem:covariance_structure}]

  Throughout this proof we will write
  $k_{ij}$ for $k_h(W_{ij},w)$
  and
  $k_{ij}'$ for $k_h(W_{ij},w')$,
  in the interest of brevity.
  \begin{align*}
    \Sigma_n(w,w')
    &=
    \E\Big[
      \big(
      \widehat f_W(w)
      - \E[\widehat f_W(w)]
      \big)
      \big(
      \widehat f_W(w')
      - \E[\widehat f_W(w')]
      \big)
      \Big] \\
    &=
    \E\left[
      \left(
        \frac{2}{n(n-1)}
        \sum_{i<j}
        \big(
        k_{ij} - \E k_{ij}
        \big)
      \right)
      \left(
        \frac{2}{n(n-1)}
        \sum_{r<s}
        \big(
        k_{rs}' - \E k_{rs}'
        \big)
      \right)
    \right] \\
    &=
    \frac{4}{n^2(n-1)^2}
    \sum_{i<j}
    \sum_{r<s}
    \E\left[
      \big(
      k_{ij} - \E k_{ij}
      \big)
      \big(
      k_{rs}' - \E k_{rs}'
      \big)
    \right]
    =
    \frac{4}{n^2(n-1)^2}
    \sum_{i<j}
    \sum_{r<s}
    \Cov\left[
      k_{ij},
      k_{rs}'
    \right].
  \end{align*}
  Note first that
  for $i,j,r,s$ all distinct,
  $k_{ij}$ is independent of $k_{rs}'$
  and so the covariance is zero.
  By a counting argument,
  it can be seen that
  there are
  $n(n-1)/2$
  summands where
  $|\{i,j,r,s\}| = 2$,
  and
  $n(n-1)(n-2)$
  summands where
  $|\{i,j,r,s\}| = 3$.
  Therefore since the samples
  are identically distributed,
  the value of the summands
  depends only on the number of distinct indices
  and we have the decomposition
  \begin{align*}
    \Sigma_n(w,w')
    &=
    \frac{4}{n^2(n-1)^2}
    \bigg(
    \frac{n(n-1)}{2}
    \Cov[k_{ij}, k_{ij}']
    + n(n-1)(n-2)
    \Cov[k_{ij}, k_{ir}']
    \bigg) \\
    &=
    \frac{2}{n(n-1)}
    \Cov[k_{ij}, k_{ij}']
    + \frac{4(n-2)}{n(n-1)}
    \Cov[k_{ij}, k_{ir}'],
  \end{align*}
  giving the first representation.
  To obtain the second representation,
  note that since
  $W_{ij}$ and $W_{ir}$
  are independent conditional
  on $A_i$,
  \begin{align*}
    \Cov\big[
      k_{ij}
      k_{ir}'
      \big]
    &=
    \E\big[
      k_{ij}
      k_{ir}'
      \big]
    -
    \E[k_{ij}]
    \E[k_{ir}']
    =
    \E\big[
      \E\big[
        k_{ij}
        k_{ir}'
        \mid A_i
        \big]
      \big]
    -
    \E[k_{ij}]
    \E[k_{ir}'] \\
    &=
    \E\big[
      \E[k_{ij} \mid A_i]
      \E[k_{ir}' \mid A_i]
      \big]
    -
    \E[k_{ij}]
    \E[k_{ir}']
    =
    \Cov\big[
      \E[k_{ij} \mid A_i],
      \E[k_{ir}' \mid A_i]
      \big].
  \end{align*}
\end{proof}

\begin{proof}[Lemma~\ref{lem:variance_bounds}]

  By Lemma~\ref{lem:covariance_structure},
  the diagonal elements of $\Sigma_n$ are
  \begin{align*}
    \Sigma_n(w,w)
    &=
    \frac{2}{n(n-1)}
    \Var\big[
      k_h(W_{ij},w)
      \big]
    +
    \frac{4(n-2)}{n(n-1)}
    \Var\big[
      \E[k_h(W_{ij},w) \mid A_i]
      \big].
  \end{align*}
  We bound each of the two terms separately.
  Firstly, note that since $k_h$ is bounded by $\Ck/h$,
  \begin{align*}
    \Var\big[
      k_h(W_{ij},w)
      \big]
    \leq
    \E\big[
      k_h(W_{ij},w)^2
      \big]
    =
    \int_{\cW \cap [w \pm h]}
    k_h(s,w)^2
    f_W(s)
    \diff{s}
    \leq 2 \Ck^2 / h.
  \end{align*}
  Conversely, since
  $\big|\E[k_h(W_{ij},w)]\big|
  = \big|\int_{[w \pm h] \cap \cW} k_h(s,w) f_W(s) \diff{s}\big|
  \leq 2 \Ck \Cd$,
  Jensen's integral inequality shows
  \begin{align*}
    \Var\big[
      k_h(W_{ij},w)
      \big]
    &\geq
    \int_{\cW \cap [w \pm h]}
    k_h(s,w)^2
    f_W(s)
    \diff{s}
    - 4 \Ck^2 \Cd^2 \\
    &\geq
    \inf_{w \in \cW} f_W(w)
    \frac{1}{2h}
    \left(
      \int_{\cW \cap [w \pm h]}
      k_h(s,w)
      \diff{s}
    \right)^2
    - 4 \Ck^2 \Cd^2 \\
    &\geq
    \frac{1}{2h}
    \inf_{w \in \cW} f_W(w)
    - 4 \Ck^2 \Cd^2
    \geq
    \frac{1}{4h}
    \inf_{w \in \cW} f_W(w)
  \end{align*}
  for small enough $h$, noting that this is trivially true
  if the infimum is zero.
  For the other term, we have
  \begin{align*}
    \Var\big[
      \E[k_h(W_{ij},w) \mid A_i]
      \big]
    &\leq
    \Var\big[
      f_{W \mid A}(w \mid A_i)
      \big]
    + 16 \CH \Ck \Cd h
    \leq
    2 \Du^2
  \end{align*}
  for small enough $h$, by a result from
  the proof of Lemma~\ref{lem:trichotomy}.
  Also
  \begin{align*}
    \Var\big[
      \E[k_h(W_{ij},w) \mid A_i]
      \big]
    &\geq
    \Var\big[
      f_{W \mid A}(w \mid A_i)
      \big]
    - 16 \CH \Ck \Cd h
    \geq
    \frac{\Dl^2}{2}
  \end{align*}
  for small enough $h$.
  Combining these four inequalities yields
  that for all large enough $n$,
  \begin{align*}
    &\frac{2}{n(n-1)}
    \frac{1}{4h}
    \inf_{w \in \cW} f_W(w)
    + \frac{4(n-2)}{n(n-1)}
    \frac{\Dl^2}{2}
    \leq
    \inf_{w \in \cW} \Sigma_n(w,w) \\
    &\qquad\leq
    \sup_{w \in \cW} \Sigma_n(w,w)
    \leq
    \frac{2}{n(n-1)}
    \frac{2 \Ck^2}{h}
    + \frac{4(n-2)}{n(n-1)}
    2 \Du^2,
  \end{align*}
  so that
  \begin{align*}
    \frac{\Dl^2}{n}
    + \frac{1}{n^2h}
    \inf_{w \in \cW} f_W(w)
    &\lesssim
    \inf_{w \in \cW} \Sigma_n(w,w)
    \leq
    \sup_{w \in \cW} \Sigma_n(w,w)
    \lesssim
    \frac{\Du^2}{n}
    + \frac{1}{n^2h}.
  \end{align*}
\end{proof}

\subsubsection{Strong approximation}

\begin{proof}[Lemma~\ref{lem:strong_approx_Ln}]

  To obtain the strong approximation,
  we apply the KMT corollary from
  Lemma~\ref{lem:kmt_corollary}.
  Define the functions
  \begin{align*}
    k_h^A(a, w) = 2\E[k_h(W_{ij},w) \mid A_i = a],
  \end{align*}
  which are of bounded variation in $a$ uniformly over $w$ since
  \begin{align*}
    \sup_{w \in \cW} \|k_h^A(\cdot,w)\|_\TV
    &=
    2\sup_{w \in \cW}
    \sup_{m \in \bbN}
    \sup_{a_0 \leq \cdots \leq a_m}
    \sum_{i=1}^m
    \big|k_h^A(a_i,w) - k_h^A(a_{i-1},w)\big| \\
    &=
    2\sup_{w \in \cW}
    \sup_{m \in \bbN}
    \sup_{a_0 \leq \cdots \leq a_m}
    \sum_{i=1}^m
    \left|
    \int_{[w \pm h] \cap \cW}
    k_h(s,w)
    \big(
    f_{W \mid A}(s \mid  a_i)
    - f_{W \mid A}(s \mid  a_{i-1})
    \big)
    \diff{s}
    \right| \\
    &\leq
    2 \sup_{w \in \cW}
    \int_{[w \pm h] \cap \cW}
    |k_h(s,w)|
    \sup_{m \in \bbN}
    \sup_{a_0 \leq \cdots \leq a_m}
    \sum_{i=1}^m
    \big|
    f_{W \mid A}(s \mid  a_i)
    - f_{W \mid A}(s \mid  a_{i-1})
    \big|
    \diff{s} \\
    &\leq
    2 \sup_{w \in \cW}
    \int_{[w \pm h] \cap \cW}
    |k_h(s,w)|
    \,
    \big\|
    f_{W \mid A}(w \mid  \cdot)
    \big\|_\TV
    \diff{s} \\
    &\leq
    4 \Ck \sup_{w \in \cW}
    \big\|
    f_{W \mid A}(w \mid  \cdot)
    \big\|_\TV
    \lesssim
    \Du,
  \end{align*}
  where the last line is by observing that the total variation
  is zero whenever $\Du = 0$.
  Hence by Lemma~\ref{lem:kmt_corollary}
  there exist (on some probability space)
  $n$ independent copies of $A_i$,
  denoted $A_i'$
  and a centered Gaussian process $Z_n^{L\prime}$
  such that if we define
  \begin{align*}
    L_n'(w)
    &=
    \frac{1}{n}
    \sum_{i=1}^n
    \big(k_h^A(A_i',w) -
    \E[k_h^A(A_i',w)]\big),
  \end{align*}
  then for some positive constants
  $C_1, C_2, C_3$,
  by defining the processes as zero outside $\cW$
  we have
  \begin{align*}
    \P\left(
      \sup_{w \in \cW}
      \Big|\sqrt{n} L_n'(w) - Z_n^{L\prime}(w)\Big|
      > \Du \frac{t + C_1 \log n}{\sqrt n}
    \right)
    \leq C_2 e^{-C_3 t}.
  \end{align*}
  Integrating tail probabilities shows that
  \begin{align*}
    \E\left[
      \sup_{w \in \cW}
      \Big|\sqrt{n} L_n'(w) - Z_n^{L\prime}(w)\Big|
    \right]
    &\leq
    \Du \frac{C_1 \log n}{\sqrt n}
    + \int_0^\infty
    \frac{\Du}{\sqrt n}
    C_2 e^{-C_3 t}
    \diff{t}
    \lesssim
    \frac{\Du \log n}{\sqrt n}.
  \end{align*}
  Further,
  $Z_n^{L\prime}$ has the
  same covariance structure as $G_n^{L\prime}$ in the
  sense that for all $w, w' \in \cW$,
  \begin{align*}
    \E\big[Z_n^{L\prime}(w) Z_n^{L\prime}(w')\big]
    = \E\big[G_n^{L\prime}(w) G_n^{L\prime}(w')\big],
  \end{align*}
  and clearly
  $L_n'$
  is equal in distribution to $L_n$.
  To obtain the trajectory regularity property of
  $Z_n^{L\prime}$,
  note that it was shown in the proof of
  Lemma~\ref{lem:uniform_concentration}
  that for all $w,w' \in \cW$,
  \begin{align*}
    \left|
    k_h^A(A_i,w)
    - k_h^A(A_i,w')
    \right|
    &\leq
    C
    |w-w'|
  \end{align*}
  for some constant $C > 0$.
  Therefore
  since the $A_i$ are i.i.d.,
  \begin{align*}
    &\E\left[
      \big|
      Z_n^{L\prime}(w)
      - Z_n^{L\prime}(w')
      \big|^2
    \right]^{1/2}
    =
    \sqrt{n}
    \E\left[
      \big|
      L_n(w)
      - L_n(w')
      \big|^2
    \right]^{1/2} \\
    &\quad=
    \sqrt{n}
    \E\left[
      \left|
      \frac{1}{n}
      \sum_{i=1}^n
      \Big(
      k_h^A(A_i,w)
      - k_h^A(A_i,w')
      - \E\big[k_h^A(A_i,w)]
      + \E\big[k_h^A(A_i,w')]
      \Big)
      \right|^2
    \right]^{1/2} \\
    &\quad=
    \E\left[
      \Big|
      k_h^A(A_i,w)
      - k_h^A(A_i,w')
      - \E\big[k_h^A(A_i,w)]
      + \E\big[k_h^A(A_i,w')]
      \Big|^2
    \right]^{1/2}
    \lesssim
    |w-w'|.
  \end{align*}
  Therefore by
  the regularity result for Gaussian processes in
  Lemma~\ref{lem:gaussian_process_maximal},
  with $\delta_n \in (0, 1/2]$:
  \begin{align*}
    \E\left[
      \sup_{|w-w'| \leq \delta_n}
      \big|
      Z_n^{L\prime}(w)
      - Z_n^{L\prime}(w')
      \big|
    \right]
    &\lesssim
    \int_0^{\delta_n}
    \sqrt{\log 1/\varepsilon}
    \diff{\varepsilon}
    \lesssim
    \delta_n \sqrt{\log 1/\delta_n}
    \lesssim
    \Du
    \delta_n \sqrt{\log 1/\delta_n},
  \end{align*}
  where the last inequality is because
  $Z_n^{L\prime} \equiv 0$ whenever $\Du = 0$.
  There is a modification of $Z_n^{L\prime}$
  with continuous trajectories
  by Kolmogorov's continuity criterion
  \citep[Theorem~2.9]{legall2016brownian}.
  Note that $L_n'$ is $\bA_n'$-measurable
  and so by Lemma~\ref{lem:kmt_corollary}
  we can assume that $Z_n^{L\prime}$
  depends only on $\bA_n'$ and some
  random noise which is independent of
  $(\bA_n', \bV_n')$.
  Finally in order to have
  $\bA_n', \bV_n', L_n'$ and $Z_n^{L\prime}$
  all defined on the same probability space,
  we note that $\bA_n$ and $\bV_n$ are random vectors
  while $L_n'$ and $Z_n^{L\prime}$
  are stochastic processes
  with continuous sample paths
  indexed on
  the compact interval $\cW$.
  Hence the Vorob'ev--Berkes--Philipp theorem
  (Lemma~\ref{lem:vbp})
  allows us to ``glue'' them together
  in the desired way
  on another new probability space,
  giving
  $\big(\bA_n', \bV_n', L_n', Z_n^{L\prime}\big)$,
  where we retain the single prime notation for clarity.
\end{proof}

\begin{proof}[Lemma~\ref{lem:conditional_strong_approx_En}]

  We apply Lemma~\ref{lem:yurinskii_corollary} conditional on
  $\bA_n$. While this lemma is not in its current form
  stated for conditional distributions,
  the Yurinskii coupling on which it depends can be readily extended
  by following the proof of \citet[Lemma~38]{belloni2019conditional},
  using a conditional version of Strassen's theorem
  \cite[Theorem~B.2]{chen2020jackknife}.
  Care must similarly be taken in embedding the conditionally Gaussian vectors
  into a conditionally Gaussian process, using the
  Vorob'ev--Berkes--Philipp theorem (Lemma~\ref{lem:vbp}).

  By the mutual independence of $A_i$ and $V_{ij}$,
  we have that the observations
  $W_{ij}$ are independent
  (but not necessarily identically distributed)
  conditionally on $\bA_n$.
  Note that
  $\sup_{s,w \in \cW} |k_h(s,w)| \lesssim M_n = h^{-1}$
  and
  $\E[k_h(W_{ij},w)^2 \mid \bA_n] \lesssim \sigma_n^2 = h^{-1}$.
  The following uniform Lipschitz condition holds
  with $l_{n,\infty} = \CL h^{-2}$,
  by the Lipschitz property of the kernels:
  \begin{align*}
    \sup_{s,w,w' \in \cW}
    \left|
    \frac{k_h(s, w) - k_h(s, w')}
    {w-w'}
    \right|
    \leq
    l_{n,\infty}.
  \end{align*}
  Also, the following $L^2$ Lipschitz condition holds
  uniformly with $l_{n,2} = 2 \CL \sqrt{\Cd} h^{-3/2}$:
  \begin{align*}
    &\E\big[
      \big|
      k_h(W_{ij}, w) - k_h(W_{ij}, w')
      \big|^2
      \mid \bA_n
      \big]^{1/2} \\
    &\quad\leq
    \frac{\CL}{h^2}
    |w-w'|
    \left(
      \int_{([w \pm h] \cup [w' \pm h]) \cap \cW}
      f_{W \mid AA}(s \mid  \bA_n)
      \diff{s}
    \right)^{1/2} \\
    &\quad\leq
    \frac{\CL}{h^2}
    |w-w'|
    \sqrt{4h \Cd}
    \leq
    l_{n,2}
    |w-w'|.
  \end{align*}
  So we can apply
  Lemma~\ref{lem:yurinskii_corollary}
  conditionally on $\bA_n$
  to the $\frac{1}{2}n(n-1)$ observations,
  noting that
  \begin{align*}
    \sqrt{n^2h} E_n(w)
    =
    \sqrt{\frac{2nh}{n-1}}
    \sqrt{\frac{2}{n(n-1)}}
    \sum_{i=1}^{n-1}
    \sum_{j=i+1}^{n}
    \Big(
    k_h(W_{ij},w)
    - \E[k_h(W_{ij},w)  \mid  A_i, A_j]
    \Big),
  \end{align*}
  to deduce that for $t_n > 0$ there exist
  (an enlarged probability space)
  conditionally mean-zero
  and conditionally Gaussian processes
  $\widetilde Z_n^{E\prime}(w)$
  with the same conditional covariance structure as
  $\sqrt{n^2 h} E_n(w)$ and
  satisfying
  \begin{align*}
    &\P\left(
      \sup_{w \in \cW}
      \big|
      \sqrt{n^2h} E_n(w) - \widetilde Z_n^{E\prime}(w)
      \big|
      > t_n
      \Bigm\vert \bA_n'
    \right) \\
    &\quad=
    \P\left(
      \sup_{w \in \cW}
      \left|
      \sqrt{\frac{n(n-1)}{2}} E_n(w)
      - \sqrt{\frac{n-1}{2nh}} \widetilde Z_n^{E\prime}(w)
      \right|
      > \sqrt{\frac{n-1}{2nh}}
      t_n
      \Bigm\vert \bA_n'
    \right) \\
    &\quad\lesssim
    \frac{
      \sigma_n
      \sqrt{\Leb(\cW)}
      \sqrt{\log n}
      \sqrt{M_n + \sigma_n\sqrt{\log n}}
    }{n^{1/2} t_n^2 / h}
    \sqrt{
      l_{n,2}
      \sqrt{\log n}
      + \frac{l_{n,\infty}}{n}
      \log n} \\
    &\quad\lesssim
    \frac{
      h^{-1/2}
      \sqrt{\log n}
      \sqrt{h^{-1} + h^{-1/2} \sqrt{\log n}}
    }{n^{1/2} t_n^2 / h}
    \sqrt{
      h^{-3/2}
      \sqrt{\log n}
      + \frac{h^{-2}}{n}
      \log n} \\
    &\quad\lesssim
    \sqrt{\frac{\log n}{n}}
    \frac{
      \sqrt{1 + \sqrt{h \log n}}
    }{t_n^2}
    \sqrt{
      \sqrt{\frac{\log n}{h^3}}
      \left( 1 + \sqrt{\frac{\log n}{n^2 h}} \right)
    } \\
    &\quad\lesssim
    \sqrt{\frac{\log n}{n}}
    \frac{ 1 }{t_n^2}
    \left(
      \frac{\log n}{h^3}
    \right)^{1/4}
    \lesssim
    t_n^{-2}
    n^{-1/2}
    h^{-3/4}
    (\log n)^{3/4},
  \end{align*}
  where we used
  $h \lesssim 1 / \log n$
  and $\frac{\log n}{n^2 h} \lesssim 1$.
  To obtain the trajectory regularity property of
  $\widetilde Z_n^{E\prime}$,
  note that
  for $w, w' \in \cW$,
  by conditional independence,
  \begin{align*}
    &\E\left[
      \big|
      \widetilde Z_n^{E\prime}(w)
      - \widetilde Z_n^{E\prime}(w')
      \big|^2
      \mid \bA_n'
    \right]^{1/2} \\
    &\quad=
    \sqrt{n^2h} \,
    \E\left[
      \big|
      E_n(w)
      - E_n(w')
      \big|^2
      \mid \bA_n
    \right]^{1/2} \\
    &\quad\lesssim
    \sqrt{n^2h} \,
    \E\left[
      \left|
      \frac{2}{n(n-1)}
      \sum_{i=1}^{n-1}
      \sum_{j=i+1}^{n}
      \Big(
      k_h(W_{ij},w)
      - k_h(W_{ij},w')
      \Big)
      \right|^2
      \Bigm\vert \bA_n
    \right]^{1/2} \\
    &\quad\lesssim
    \sqrt{h} \,
    \E\left[
      \big|
      k_h(W_{ij},w)
      - k_h(W_{ij},w')
      \big|^2
      \bigm\vert \bA_n
    \right]^{1/2}
    \lesssim
    h^{-1} |w-w'|.
  \end{align*}
  Therefore by
  the regularity result for Gaussian processes in
  Lemma~\ref{lem:gaussian_process_maximal},
  with $\delta_n \in (0, 1/(2h)]$:
  \begin{align*}
    \E\left[
      \sup_{|w-w'| \leq \delta_n}
      \big|
      \widetilde Z_n^{E\prime}(w)
      - \widetilde Z_n^{E\prime}(w')
      \big|
      \mid \bA_n'
    \right]
    &\lesssim
    \int_0^{\delta_n/h}
    \sqrt{\log (\varepsilon^{-1} h^{-1})}
    \diff{\varepsilon}
    \lesssim
    \frac{\delta_n}{h}
    \sqrt{\log \frac{1}{h\delta_n}},
  \end{align*}
  and there exists a modification
  with continuous trajectories.
  Finally in order to have
  $\bA_n', \bV_n', E_n'$ and $\widetilde Z_n^{E\prime}$
  all defined on the same probability space,
  we note that $\bA_n$ and $\bV_n$ are random vectors
  while $E_n'$ and $\widetilde Z_n^{E\prime}$
  are stochastic processes
  with continuous sample paths
  indexed on
  the compact interval $\cW$.
  Hence the Vorob'ev--Berkes--Philipp theorem
  (Lemma~\ref{lem:vbp})
  allows us to ``glue together''
  $\big(\bA_n, \bV_n, E_n\big)$
  and
  $\big(E_n', \widetilde Z_n^{E\prime}\big)$
  in the desired way
  on another new probability space,
  giving
  $\big(\bA_n', \bV_n', E_n', \widetilde Z_n^{E\prime}\big)$,
  where we retain the single prime notation for clarity.

  Note that the trajectories of the conditionally Gaussian processes
  $\widetilde Z_n^{E\prime}$ depend on the choice of $t_n$,
  necessitating the use of a divergent sequence $R_n$ to establish
  bounds in probability on the coupling error.
\end{proof}

\begin{proof}[Lemma~\ref{lem:unconditional_strong_approx_En}]

  \proofparagraph{defining $Z_n^{E\dprime}$}

  Pick $\delta_n \to 0$
  with $\log 1/\delta_n \lesssim \log n$.
  Let $\cW_\delta$ be a $\delta_n$-covering of $\cW$
  with cardinality $\Leb(\cW)/\delta_n$
  which is also a $\delta_n$-packing.
  Let $\widetilde Z_{n,\delta}^{E\prime}$
  be the restriction of $\widetilde Z_n^{E\prime}$
  to $\cW_\delta$.
  Let
  $\widetilde \Sigma_n^E(w, w') =
  \E\big[\widetilde Z_n^{E\prime}(w) \widetilde Z_n^{E\prime}(w')
  \mid \bA_n' \big]$
  be the conditional covariance function of $\widetilde Z_n^{E\prime}$,
  and define
  $\Sigma_n^E(w,w') = \E\big[\widetilde \Sigma_n^E(w,w')\big]$.
  Let $\widetilde \Sigma^E_{n,\delta}$
  and $\Sigma^E_{n,\delta}$
  be the restriction matrices of
  $\widetilde \Sigma^E_n$
  and $\Sigma^E_n$
  respectively
  to $\cW_\delta \times \cW_\delta$,
  noting that, as (conditional) covariance matrices,
  these are
  (almost surely)
  positive semi-definite.

  Let $N \sim \cN(0, I_{|\cW_\delta|})$
  be independent of $\bA_n'$,
  and define using the matrix square root
  $\widetilde Z_{n,\delta}^{E\dprime}
  = \big(\widetilde \Sigma^E_{n,\delta})^{1/2} N$,
  which has the same distribution as
  $\widetilde Z_{n,\delta}^{E\prime}$,
  conditional on $\bA_n'$.
  Extend it using
  the Vorob'ev--Berkes--Philipp theorem
  (Lemma~\ref{lem:vbp})
  to the compact interval $\cW$,
  giving a conditionally Gaussian process
  $\widetilde Z_n^{E\dprime}$
  which has the same distribution as
  $\widetilde Z_{n}^{E\prime}$,
  conditional on $\bA_n'$.
  Define
  $Z_{n,\delta}^{E\dprime} = \big(\Sigma^E_{n,\delta})^{1/2} N$,
  noting that this is independent of $\bA_n'$,
  and extend it using
  the Vorob'ev--Berkes--Philipp theorem
  (Lemma~\ref{lem:vbp})
  to a Gaussian process
  $Z_n^{E\dprime}$ on the compact interval $\cW$,
  which is independent of $\bA_n'$
  and has covariance structure given by
  $\Sigma_n^E$.

  \proofparagraph{closeness of $Z_n^{E\dprime}$ and
  $\widetilde Z_n^{E\dprime}$ on the mesh}

  Note that conditionally on $\bA_n'$,
  $\widetilde Z_{n,\delta}^{E\dprime} - Z_{n,\delta}^{E\dprime}$
  is a length-$|\cW_\delta|$
  Gaussian random vector with covariance matrix
  $\big(
  \big(\widetilde \Sigma^E_{n,\delta}\big)^{1/2}
  - \big(\Sigma^E_{n,\delta}\big)^{1/2}
  \big)^2$.
  So by the Gaussian maximal inequality in
  Lemma~\ref{lem:gaussian_vector_maximal}
  applied conditionally on $\bA_n'$,
  we have
  \begin{align*}
    \E\left[
      \max_{w \in \cW_\delta}
      \big|\widetilde Z_n^{E\dprime}(w) - Z_n^{E\dprime}(w)\big|
      \Bigm| \bA_n'
    \right]
    &\lesssim
    \sqrt{\log n}
    \left\|
    \widetilde\Sigma^E_{n,\delta}
    - \Sigma^E_{n,\delta}
    \right\|_2^{1/2},
  \end{align*}
  since $\log |\cW_\delta| \lesssim \log n$.
  Next, we apply some U-statistic theory to
  $\widetilde\Sigma^E_{n,\delta} - \Sigma^E_{n,\delta}$,
  with the aim of applying the
  matrix concentration result
  for second-order U-statistics
  presented in Lemma~\ref{lem:ustat_matrix_concentration}.
  Firstly we note that
  since
  the conditional covariance structures of
  $\widetilde Z_n^{E\prime}$ and $\sqrt{n^2h} E_n$
  are equal in distribution,
  we have,
  writing $E_n(\cW_\delta)$
  for the vector $\big(E_n(w) : w \in \cW_\delta\big)$
  and similarly for $k_h(W_{ij}, \cW_\delta)$,
  \begin{align*}
    \widetilde\Sigma^E_{n,\delta}
    &=
    n^2h \E[E_n(\cW_\delta) E_n(\cW_\delta)^\T \mid \bA_n] \\
    &=
    n^2h
    \frac{4}{n^2(n-1)^2}
    \sum_{i=1}^{n-1}
    \sum_{j=i+1}^{n} \\
    &\qquad\E\left[
      \Big(
      k_h(W_{ij}, \cW_\delta)
      - \E\left[
        k_h(W_{ij}, \cW_\delta)
        \mid \bA_n
      \right]
      \Big)
      \Big(
      k_h(W_{ij}, \cW_\delta)
      - \E\left[
        k_h(W_{ij}, \cW_\delta)
        \mid \bA_n
      \right]
      \Big)^\T
      \bigm\vert \bA_n
    \right] \\
    &=
    \frac{4h}{(n-1)^2}
    \sum_{i=1}^{n-1}
    \sum_{j=i+1}^{n}
    u(A_i, A_j),
  \end{align*}
  where we
  define the random
  $|\cW_\delta| \times |\cW_\delta|$
  matrices
  \begin{align*}
    u(A_i, A_j)
    &=
    \E\left[
      k_h(W_{ij}, \cW_\delta)
      k_h(W_{ij}, \cW_\delta)^\T
      \mid \bA_n
    \right]
    -
    \E\left[
      k_h(W_{ij}, \cW_\delta)
      \mid \bA_n
    \right]
    \E\left[
      k_h(W_{ij}, \cW_\delta)
      \mid \bA_n
    \right]^\T.
  \end{align*}
  Let $u(A_i) = \E[u(A_i, A_j) \mid A_i]$
  and
  $u = \E[u(A_i, A_j)]$.
  The following Hoeffding decomposition holds,
  by Lemma~\ref{lem:general_hoeffding}:
  \begin{align*}
    \widetilde \Sigma^E_{n,\delta} - \Sigma^E_{n,\delta}
    &=
    \widetilde L +\widetilde Q,
  \end{align*}
  where
  \begin{align*}
    \widetilde L
    &=
    \frac{4h}{n-1}
    \sum_{i=1}^n
    \big(
    u(A_i) - u
    \big), \\
    \widetilde Q
    &=
    \frac{4h}{(n-1)^2}
    \sum_{i=1}^{n-1}
    \sum_{j=i+1}^{n}
    \big(
    u(A_i, A_j) - u(A_i) - u(A_j) + u
    \big).
  \end{align*}
  Next, we seek an almost sure upper bound on
  $\|u(A_i, A_j)\|_2$.
  Since this is a symmetric matrix,
  we have by H\"older's inequality
  \begin{align*}
    \|u(A_i, A_j)\|_2
    &\leq
    \|u(A_i, A_j)\|_1^{1/2}
    \|u(A_i, A_j)\|_\infty^{1/2}
    =
    \max_{1 \leq k \leq |\cW_\delta|}
    \sum_{l=1}^{|\cW_\delta|}
    |u(A_i, A_j)_{kl}|.
  \end{align*}
  The terms on the right hand side can be bounded as follows,
  writing $w, w'$ for the $k$th and $l$th
  points in $\cW_\delta$ respectively:
  \begin{align*}
    |u(A_i, A_j)_{kl}|
    &=
    \big|
    \E\left[
      k_h(W_{ij}, w)
      k_h(W_{ij}, w')
      \mid \bA_n
    \right]
    -
    \E\left[
      k_h(W_{ij}, w)
      \mid \bA_n
    \right]
    \E\left[
      k_h(W_{ij}, w')
      \mid \bA_n
    \right]
    \big| \\
    &\lesssim
    \E\left[
      |
      k_h(W_{ij}, w)
      k_h(W_{ij}, w')
      |
      \mid \bA_n
    \right]
    +
    \E\left[
      |
      k_h(W_{ij}, w)
      |
      \mid \bA_n
    \right]
    \E\left[
      |
      k_h(W_{ij}, w')
      |
      \mid \bA_n
    \right] \\
    &\lesssim
    h^{-1}
    \I\big\{ |w-w'| \leq 2h \big\}
    + 1
    \lesssim
    h^{-1}
    \I\big\{ |k-l| \leq 2h/\delta_n \big\}
    + 1,
  \end{align*}
  where we used that
  $|w-w'| \geq |k-l| \delta_n$
  because $\cW_\delta$
  is a $\delta_n$-packing.
  Hence
  \begin{align*}
    \|u(A_i, A_j)\|_2
    &\leq
    \max_{1 \leq k \leq |\cW_\delta|}
    \sum_{l=1}^{|\cW_\delta|}
    |u(A_i, A_j)_{kl}|
    \lesssim
    \max_{1 \leq k \leq |\cW_\delta|}
    \sum_{l=1}^{|\cW_\delta|}
    \Big(
    h^{-1}
    \I\big\{ |k-l| \leq 2h/\delta_n \big\}
    + 1
    \Big) \\
    &\lesssim
    1/\delta_n
    + 1/h
    + |\cW_\delta|
    \lesssim
    1/\delta_n
    + 1/h.
  \end{align*}
  Clearly the same bound holds for
  $\|u(A_i)\|_2$
  and $\|u\|_2$, by Jensen's inequality.
  Therefore applying the matrix Bernstein inequality
  (Lemma~\ref{lem:matrix_bernstein})
  to the zero-mean matrix $\widetilde L$ gives
  \begin{align*}
    \E\left[
      \left\|
      \widetilde L
      \right\|_2
    \right]
    &\lesssim
    \frac{h}{n}
    \left(\frac{1}{\delta_n} + \frac{1}{h} \right)
    \left(
      \log |\cW_\delta| + \sqrt{n \log |\cW_\delta|}
    \right)
    \lesssim
    \left(\frac{h}{\delta_n} + 1 \right)
    \sqrt{\frac{\log n}{n}}.
  \end{align*}
  Applying the matrix U-statistic concentration inequality
  (Lemma~\ref{lem:ustat_matrix_concentration})
  to the zero-mean matrix $\widetilde Q$ gives
  \begin{align*}
    \E\left[
      \left\|
      \widetilde Q
      \right\|_2
    \right]
    &\lesssim
    \frac{h}{n^2}
    n
    \left(\frac{1}{\delta_n} + \frac{1}{h} \right)
    \left(
      \log |\cW_\delta|
    \right)^{3/2}
    \lesssim
    \left(\frac{h}{\delta_n} + 1 \right)
    \frac{(\log n)^{3/2}}{n}.
  \end{align*}
  Hence putting everything together,
  taking a marginal expectation
  and applying Jensen's inequality,
  \begin{align*}
    \E\left[
      \max_{w \in \cW_\delta}
      \big|\widetilde Z_n^{E\dprime}(w) - Z_n^{E\dprime}(w)\big|
    \right]
    &\lesssim
    \sqrt{\log n} \
    \E\left[
      \left\|
      \widetilde\Sigma^E_{n,\delta} - \Sigma^E_{n,\delta}
      \right\|_2^{1/2}
    \right]
    \lesssim
    \sqrt{\log n} \
    \E\left[
      \left\|
      \widetilde\Sigma^E_{n,\delta} - \Sigma^E_{n,\delta}
      \right\|_2
    \right]^{1/2} \\
    &\lesssim
    \sqrt{\log n} \
    \E\left[
      \left\|
      \widetilde L
      + \widetilde Q
      \right\|_2
    \right]^{1/2}
    \lesssim
    \sqrt{\log n} \
    \E\left[
      \left\|
      \widetilde L
      \right\|_2
      + \left\|
      \widetilde Q
      \right\|_2
    \right]^{1/2} \\
    &\lesssim
    \sqrt{\log n}
    \left(
      \left(\frac{h}{\delta_n} + 1 \right)
      \sqrt{\frac{\log n}{n}}
      + \left(\frac{h}{\delta_n} + 1 \right)
      \frac{(\log n)^{3/2}}{n}
    \right)^{1/2} \\
    &\lesssim
    \sqrt{\frac{h}{\delta_n} + 1}
    \frac{(\log n)^{3/4}}{n^{1/4}}.
  \end{align*}

  \proofparagraph{regularity of $Z_n^E$ and $\widetilde Z_n^{E\prime}$}

  Define the semimetrics
  \begin{align*}
    \rho(w, w')^2
    &=
    \E\left[
      \big|Z_n^{E\dprime}(w) - Z_n^{E\dprime}(w')\big|^2
    \right],
    &\widetilde\rho(w, w')^2
    &=
    \E\left[
      \big|\widetilde Z_n^{E\dprime}(w) - \widetilde Z_n^{E\dprime}(w')\big|^2
      \mid \bA_n
    \right].
  \end{align*}
  We can bound $\widetilde \rho$ as follows,
  since $\widetilde Z_n^{E\dprime}$ and $\sqrt{n^2h} E_n$
  have the same conditional covariance structure:
  \begin{align*}
    \widetilde\rho(w, w')
    &=
    \E\left[
      \big|\widetilde Z_n^{E\dprime}(w) - \widetilde Z_n^{E\dprime}(w')\big|^2
      \mid \bA_n'
    \right]^{1/2}
    =
    \sqrt{n^2 h} \,
    \E\left[
      \big|E_n(w) - E_n(w')\big|^2
      \mid \bA_n'
    \right]^{1/2}
    \lesssim
    h^{-1}
    |w-w'|,
  \end{align*}
  uniformly in $\bA_n'$,
  where the last line was shown in
  the proof of Lemma~\ref{lem:conditional_strong_approx_En}.
  Note that also
  \begin{align*}
    \rho(w, w')
    &=
    \sqrt{\E[\widetilde \rho(w,w')^2]}
    \lesssim
    h^{-1}
    |w-w'|.
  \end{align*}
  Thus Lemma~\ref{lem:gaussian_process_maximal}
  applies directly to $Z_n^E$
  and conditionally to $\widetilde Z_n^{E\prime}$,
  with $\delta_n \in (0, 1/(2h)]$,
  yielding
  \begin{align*}
    \E\left[
      \sup_{|w-w'| \leq \delta_n}
      \big|\widetilde Z_n^{E\dprime}(w) - \widetilde Z_n^{E\dprime}(w')\big|
      \bigm\vert \bA_n'
    \right]
    &\lesssim
    \int_0^{\delta_n / h}
    \sqrt{\log (1 / (\varepsilon h))}
    \diff{\varepsilon}
    \lesssim
    \frac{\delta_n}{h}
    \sqrt{\log \frac{1}{h \delta_n}}, \\
    \E\left[
      \sup_{|w-w'| \leq \delta_n}
      |Z_n^{E\dprime}(w) - Z_n^{E\dprime}(w')|
    \right]
    &\lesssim
    \int_0^{\delta_n / h}
    \sqrt{\log (1 / (\varepsilon h))}
    \diff{\varepsilon}
    \lesssim
    \frac{\delta_n}{h}
    \sqrt{\log \frac{1}{h \delta_n}}.
  \end{align*}
  Continuity of trajectories follows from this.

  \proofparagraph{conclusion}

  We use the previous parts
  to deduce that
  \begin{align*}
    \E\left[
      \sup_{w \in \cW}
      \big|\widetilde Z_n^{E\dprime}(w) - Z_n^{E\dprime}(w)\big|
    \right]
    &\lesssim
    \E\left[
      \max_{w \in \cW_\delta}
      \big|\widetilde Z_n^{E\dprime}(w) - Z_n^{E\dprime}(w)\big|
    \right]
    + \E\left[
      \sup_{|w-w'| \leq \delta_n}
      \big|\widetilde Z_n^{E\dprime}(w) - \widetilde Z_n^{E\dprime}(w')\big|
    \right] \\
    &\qquad+
    \E\left[
      \sup_{|w-w'| \leq \delta_n}
      \big|Z_n^{E\dprime}(w) - Z_n^{E\dprime}(w')\big|
    \right] \\
    &\lesssim
    \sqrt{\frac{h}{\delta_n} + 1}
    \frac{(\log n)^{3/4}}{n^{1/4}}
    + \frac{\delta_n \sqrt{\log n}}{h}.
  \end{align*}
  Setting
  $\delta_n = h \left( \frac{\log n}{n} \right)^{1/6}$
  gives
  \begin{align*}
    \E\left[
      \sup_{w \in \cW}
      \big|\widetilde Z_n^{E\dprime}(w) - Z_n^{E\dprime}(w)\big|
    \right]
    &\lesssim
    n^{-1/6} (\log n)^{2/3}.
  \end{align*}
  Finally, independence of $Z_n^{E\dprime}$ and $\bA_n''$
  follows by another application of
  the Vorob'ev--Berkes--Philipp theorem
  from Lemma~\ref{lem:vbp},
  this time conditionally on $\bA_n'$,
  to the random variables
  $\big(\bA_n', \widetilde Z_n^{E\prime}\big)$
  and
  $\big(\widetilde Z_n^{E\dprime}, Z_n^{E\dprime}\big)$.
\end{proof}

\begin{proof}[Theorem~\ref{thm:strong_approx_fW}]

  We add together the strong approximations
  for the $L_n$ and $E_n$ terms,
  and then add an independent Gaussian process
  to account for the variance of $Q_n$.

  \proofparagraph{gluing together the strong approximations}

  Let $\big(\bA_n', \bV_n', L_n', Z_n^{L\prime}\big)$
  be the strong approximation for $L_n$
  from Lemma~\ref{lem:strong_approx_Ln}.
  Likewise let
  $\big(\bA_n'', \bV_n'', E_n'', \widetilde Z_n^{E\dprime}\big)$
  and
  $\big(\bA_n''', \bV_n''', \widetilde Z_n^{E\tprime}, Z_n^{E\tprime}\big)$
  be the conditional and unconditional strong approximations for $E_n$
  given in Lemmas~\ref{lem:conditional_strong_approx_En}
  and \ref{lem:unconditional_strong_approx_En}
  respectively.
  The first step is to define copies of all of these variables
  and processes on the same probability space.
  This is achieved by applying the
  Vorob'ev--Berkes--Philipp Theorem
  (Lemma~\ref{lem:vbp}).
  In particular, dropping the prime notation for clarity,
  we construct
  $\big(\bA_n, \bV_n, L_n, Z_n^L, E_n, \widetilde Z_n^E, Z_n^E\big)$
  with the following properties:
  \begin{enumerate}[label=(\roman*)]

    \item
      $\sup_{w \in \cW}
      \big| \sqrt{n} L_n(w) -  Z_n^L(w)\big|
      \lesssim_\P n^{-1/2} \log n$,

    \item
      $\sup_{w \in \cW}
      \big|\sqrt{n^2h} E_n(w) - \widetilde Z^E_n(w) \big|
      \lesssim_\P n^{-1/4}  h^{-3/8} (\log n)^{3/8} R_n$,

    \item
      $\sup_{w \in \cW}
      \big| \widetilde Z^E_n(w) - Z^E_n(w) \big|
      \lesssim_\P n^{-1/6} (\log n)^{2/3}$,

    \item
      $Z_n^L$ is independent of $Z_n^E$.

  \end{enumerate}
  Note that the independence of
  $Z_n^L$ and $Z_n^E$
  follows since $Z_n^L$
  depends only on $\bA_n$ and some independent random noise,
  while $Z_n^E$ is independent of $\bA_n$.
  Therefore $(Z_n^L, Z_n^E)$ are jointly Gaussian.
  To obtain the strong approximation result
  for $\widehat f_W$,
  define the mean-zero Gaussian process
  \begin{align*}
    Z_n^f(w)
    &=
    \frac{1}{\sqrt n} Z_n^L(w)
    + \frac{1}{n} Z_n^Q(w)
    + \frac{1}{\sqrt{n^2h}} Z_n^E(w),
  \end{align*}
  where $Z_n^Q(w)$
  is a mean-zero Gaussian process
  independent of everything else and
  with covariance
  \begin{align*}
    \E\big[
      Z_n^Q(w)
      Z_n^Q(w')
      \big]
    &=
    n^2 \E\big[
      Q_n(w)
      Q_n(w')
      \big].
  \end{align*}
  As shown in the proof of
  Lemma~\ref{lem:uniform_concentration},
  the process
  $Q_n(w)$ is uniformly Lipschitz
  and uniformly bounded in $w$.
  Thus by Lemma~\ref{lem:gaussian_process_maximal},
  we have
  $\E\big[\sup_{w \in \cW}
    |Z_n^Q(w)|\big]
  \lesssim 1$.
  Therefore the uniform approximation error is given by
  \begin{align*}
    &
    \sup_{w \in \cW}
    \big|
    \widehat f_W(w) - \E[\widehat f_W(w)]
    - Z_n^f(w)
    \big|
    \\
    &\quad=
    \sup_{w \in \cW}
    \left|
    \frac{1}{\sqrt n} Z_n^L(w)
    + \frac{1}{n} Z_n^Q(w)
    + \frac{1}{\sqrt{n^2h}} Z_n^E(w)
    - \Big(
    L_n(w) + Q_n(w) + E_n(w)
    \Big)
    \right| \\
    &\quad\leq
    \sup_{w \in \cW}
    \bigg(
    \frac{1}{\sqrt n}
    \left|
    Z_n^L(w) - \sqrt{n} L_n(w)
    \right|
    + \frac{1}{\sqrt{n^2h}}
    \left|
    \widetilde Z_n^E(w) - \sqrt{n^2h} E_n(w)
    \right|
    + \frac{1}{\sqrt{n^2h}}
    \left|
    Z_n^E(w) - \widetilde Z_n^E(w)
    \right| \\
    &\qquad+
    \big| Q_n(w) \big|
    + \frac{1}{n}
    \big| Z_n^Q(w) \big|
    \bigg) \\
    &\quad\lesssim_\P
    n^{-1} \log n
    + n^{-5/4} h^{-7/8} (\log n)^{3/8} R_n
    + n^{-7/6} h^{-1/2} (\log n)^{2/3}.
  \end{align*}

  \proofparagraph{covariance structure}

  Since $L_n$, $Q_n$ and $E_n$
  are mutually orthogonal in $L^2$
  (as shown in Lemma~\ref{lem:hoeffding}),
  we have the following covariance
  structure:
  \begin{align*}
    \E\big[Z_n^f(w) Z_n^f(w')\big]
    &=
    \frac{1}{n} \E\big[ Z_n^L(w) Z_n^L(w') \big]
    + \frac{1}{n^2} \E\big[ Z_n^Q(w) Z_n^Q(w') \big]
    + \frac{1}{n^2h} \E\big[ Z_n^E(w) Z_n^E(w') \big] \\
    &=
    \E\big[ L_n(w) L_n(w') \big]
    + \E\big[ Q_n(w) Q_n(w') \big]
    + \E\big[ E_n(w) E_n(w') \big] \\
    &=
    \E\big[
      \big(\widehat f_W(w) - \E[\widehat f_W(w)]\big)
      \big(\widehat f_W(w') - \E[\widehat f_W(w')]\big)
      \big].
  \end{align*}

  \proofparagraph{trajectory regularity}

  The trajectory regularity of the process
  $Z_n^f$ follows directly
  by adding the regularities
  of the processes $\frac{1}{\sqrt n} Z_n^L$
  and
  $\frac{1}{n} Z_n^Q$
  and
  $\frac{1}{\sqrt{n^2h}} Z_n^E$.
  Similarly, $Z_n^f$ has continuous trajectories.
\end{proof}

\subsubsection{Infeasible uniform confidence bands}

\begin{proof}[Lemma~\ref{lem:infeasible_t_statistic}]

  Note that
  \begin{align*}
    \left| T_n(w) - Z_n^T(w) \right|
    &= \frac{\big| \widehat f_W(w) - f_W(w) - Z_n^f(w) \big|}
    {\sqrt{\Sigma_n(w,w)}}.
  \end{align*}
  By Theorem~\ref{thm:strong_approx_fW}
  and Lemma~\ref{lem:bias},
  the numerator can be bounded above by
  \begin{align*}
    &\sup_{w \in \cW}
    \left|
    \widehat f_W(w) - f_W(w)
    -
    Z_n^f(w)
    \right|
    \leq
    \sup_{w \in \cW}
    \left|
    \widehat f_W(w)
    - \E\big[\widehat f_W(w)\big]
    -
    Z_n^f(w)
    \right|
    + \sup_{w \in \cW}
    \left|
    \E\big[\widehat f_W(w)\big]
    - f_W(w)
    \right| \\
    &\quad\lesssim_\P
    n^{-1} \log n
    + n^{-5/4} h^{-7/8} (\log n)^{3/8} R_n
    + n^{-7/6} h^{-1/2} (\log n)^{2/3}
    + h^{p \wedge \beta}.
  \end{align*}
  By Lemma~\ref{lem:variance_bounds}
  with $\inf_\cW f_W(w) > 0$,
  the denominator is bounded below by
  \begin{align*}
    \inf_{w \in \cW}
    \sqrt{\Sigma_n(w,w)}
    &\gtrsim
    \frac{\Dl}{\sqrt n} + \frac{1}{\sqrt{n^2h}},
  \end{align*}
  and the result follows.
\end{proof}

\begin{proof}[Theorem~\ref{thm:infeasible_ucb}]

  Note that the covariance structure of $Z_n^T$ is given by
  \begin{align*}
    \Cov\big[
      Z_n^T(w),
      Z_n^T(w')
      \big]
    &=
    \frac{\Sigma_n(w,w')}
    {\sqrt{\Sigma_n(w,w) \Sigma_n(w',w')}}.
  \end{align*}
  We apply an anti-concentration result
  to establish that all quantiles of the random variable
  $\sup_{w \in \cW} \big|Z_n^T(w)\big|$ exist.
  To do this, we must first establish regularity
  properties of $Z_n^T$.

  \proofparagraph{$L^2$ regularity of $Z_n^T$}

  Writing $k_{ij}'$
  for $k_h(W_{ij},w')$ etc.,
  note that by Lemma~\ref{lem:covariance_structure},
  \begin{align*}
    &\big|
    \Sigma_n(w,w')
    -
    \Sigma_n(w, w'')
    \big| \\
    &\quad=
    \left|
    \frac{2}{n(n-1)}
    \Cov\big[
      k_{ij},
      k_{ij}'
      \big]
    +
    \frac{4(n-2)}{n(n-1)}
    \Cov\big[
      k_{ij},
      k_{ir}'
      \big]
    -
    \frac{2}{n(n-1)}
    \Cov\big[
      k_{ij},
      k_{ij}''
      \big]
    -
    \frac{4(n-2)}{n(n-1)}
    \Cov\big[
      k_{ij},
      k_{ir}''
      \big]
    \right| \\
    &\quad\leq
    \frac{2}{n(n-1)}
    \Big|
    \Cov\big[
      k_{ij},
      k_{ij}' - k_{ij}''
      \big]
    \Big|
    +
    \frac{4(n-2)}{n(n-1)}
    \Big|
    \Cov\big[
      k_{ij},
      k_{ir}' - k_{ir}''
      \big]
    \Big| \\
    &\quad\leq
    \frac{2}{n(n-1)}
    \|k_{ij}\|_\infty
    \|k_{ij}' - k_{ij}''\|_\infty
    +
    \frac{4(n-2)}{n(n-1)}
    \|k_{ij}\|_\infty
    \|k_{ir}' - k_{ir}''\|_\infty \\
    &\quad\leq
    \frac{4}{nh^3}
    \Ck \CL
    |w'-w''|
    \lesssim
    n^{-1}h^{-3} |w'-w''|
  \end{align*}
  uniformly in $w, w', w'' \in \cW$.
  Therefore by Lemma~\ref{lem:variance_bounds},
  with $\delta_n \leq n^{-2} h^2$,
  we have
  \begin{align*}
    \inf_{|w-w'| \leq \delta_n}
    \Sigma_n(w,w')
    &\gtrsim
    \frac{\Dl^2}{n}
    + \frac{1}{n^2h}
    - n^{-1} h^{-3} \delta_n
    \gtrsim
    \frac{\Dl^2}{n}
    + \frac{1}{n^2h}
    - \frac{1}{n^3h}
    \gtrsim
    \frac{\Dl^2}{n}
    + \frac{1}{n^2h}, \\
    \sup_{|w-w'| \leq \delta_n}
    \Sigma_n(w,w')
    &\lesssim
    \frac{\Du^2}{n}
    + \frac{1}{n^2h}
    + n^{-1} h^{-3} \delta_n
    \lesssim
    \frac{\Du^2}{n}
    + \frac{1}{n^2h}
    + \frac{1}{n^3h}
    \lesssim
    \frac{\Du^2}{n}
    + \frac{1}{n^2h}.
  \end{align*}
  The $L^2$
  regularity of $Z_n^T$
  is
  \begin{align*}
    \E\left[
      \big(
      Z_n^T(w) - Z_n^T(w')
      \big)^2
    \right]
    &=
    2 - 2
    \frac{\Sigma_n(w,w')}
    {\sqrt{\Sigma_n(w,w) \Sigma_n(w',w')}}.
  \end{align*}
  Applying the elementary result
  that for $a,b,c > 0$,
  \begin{align*}
    1 - \frac{a}{\sqrt{b c}}
    &=
    \frac{b(c-a) + a(b-a)}
    {\sqrt{b c}\big(\sqrt{b c} + a\big)},
  \end{align*}
  with
  \begin{align*}
    a&= \Sigma_n(w,w'),
    &b                &= \Sigma_n(w,w),
    &c                &= \Sigma_n(w',w')
  \end{align*}
  and noting that
  $|c-a| \lesssim n^{-1} h^{-3} |w-w'|$
  and
  $|b-a| \lesssim n^{-1} h^{-3} |w-w'|$
  and
  $\frac{\Dl^2}{n} + \frac{1}{n^2h}
  \lesssim a,b,c
  \lesssim \frac{\Du^2}{n} + \frac{1}{n^2h}$
  yields
  \begin{align*}
    \E\left[
      \big(
      Z_n^T(w) - Z_n^T(w')
      \big)^2
    \right]
    &\lesssim
    \frac{(\Du^2/n + 1/(n^2h))n^{-1}h^{-3}|w-w'|}
    {(\Dl^2/n + 1/(n^2h))^2}
    \lesssim
    \frac{n^{2} h^{-4}|w-w'|}
    {n^{-4}h^{-2}}
    \lesssim
    n^2 h^{-2} |w-w'|.
  \end{align*}
  Thus the semimetric
  induced by $Z_n^T$ on $\cW$ is
  \begin{align*}
    \rho(w,w')
    &\vcentcolon=
    \E\left[
      \big(
      Z_n^T(w) - Z_n^T(w')
      \big)^2
    \right]^{1/2}
    \lesssim
    n h^{-1} \sqrt{|w-w'|}.
  \end{align*}

  \proofparagraph{trajectory regularity of $Z_n^T$}

  By the bound on $\rho$ established in the previous part,
  we can deduce the following covering number bound:
  \begin{align*}
    N(\varepsilon, \cW, \rho)
    &\lesssim
    N\big(
    \varepsilon,
    \cW,
    n h^{-1} \sqrt{|\cdot|}
    \big)
    \lesssim
    N\big(
    n^{-1} h \varepsilon,
    \cW,
    \sqrt{|\cdot|}
    \big)
    \lesssim
    N\big(
    n^{-2} h^2 \varepsilon^2,
    \cW,
    |\cdot|
    \big)
    \lesssim
    n^2 h^{-2} \varepsilon^{-2}.
  \end{align*}
  Now apply the Gaussian process regularity result from
  Lemma~\ref{lem:gaussian_process_maximal}.
  \begin{align*}
    \E\left[
      \sup_{\rho(w,w') \leq \delta}
      \big| Z_n^T(w) - Z_n^T(w') \big|
    \right]
    &\lesssim
    \int_0^{\delta}
    \sqrt{\log N(\varepsilon, \cW, \rho)}
    \diff{\varepsilon}
    \lesssim
    \int_0^{\delta}
    \sqrt{\log (n^2 h^{-2} \varepsilon^{-2})}
    \diff{\varepsilon} \\
    &\lesssim
    \int_0^{\delta}
    \left(
      \sqrt{\log n}
      + \sqrt{\log 1/\varepsilon}
    \right)
    \diff{\varepsilon}
    \lesssim
    \delta
    \left(
      \sqrt{\log n}
      + \sqrt{\log 1/\delta}
    \right),
  \end{align*}
  and so
  \begin{align*}
    \E\left[
      \sup_{|w-w'| \leq \delta_n}
      \big| Z_n^T(w) - Z_n^T(w') \big|
    \right]
    &\lesssim
    \E\left[
      \sup_{\rho(w,w') \leq n h^{-1} \delta_n^{1/2}}
      \big| Z_n^T(w) - Z_n^T(w') \big|
    \right]
    \lesssim
    n h^{-1}
    \sqrt{\delta_n \log n},
  \end{align*}
  whenever $1/\delta_n$
  is at most polynomial in $n$.

  \proofparagraph{existence of the quantile}

  Apply the Gaussian anti-concentration
  result from Lemma~\ref{lem:anti-concentration},
  noting that $Z_n^T$ is separable,
  mean-zero and has unit variance:
  \begin{align*}
    \sup_{t \in \R}
    \P\left(
      \left|
      \sup_{w \in \cW}
      \big| Z_n^T(w) \big|
      - t
      \right|
      \leq 2\varepsilon_n
    \right)
    &\leq
    8 \varepsilon_n
    \left(
      1 + \E\left[
        \sup_{w \in \cW}
        \big| Z_n^T(w) \big|
      \right]
    \right).
  \end{align*}
  To bound the supremum on the right hand side,
  apply the Gaussian process maximal inequality from
  Lemma~\ref{lem:gaussian_process_maximal}
  with
  $\sigma \leq 1$
  and
  $N(\varepsilon, \cW, \rho) \lesssim n^2 h^{-2} \varepsilon^{-2}$:
  \begin{align*}
    \E\left[
      \sup_{w \in \cW}
      \big|Z_n^T(w)\big|
    \right]
    &\lesssim
    1
    + \int_0^{2}
    \sqrt{\log (n^2 h^{-2} \varepsilon^{-2})}
    \diff{\varepsilon}
    \lesssim
    \sqrt{\log n}.
  \end{align*}
  Therefore
  \begin{align*}
    \sup_{t \in \R}
    \P\left(
      \left|
      \sup_{w \in \cW}
      \big| Z_n^T(w) \big|
      - t
      \right|
      \leq \varepsilon
    \right)
    &\lesssim
    \varepsilon
    \sqrt{\log n}.
  \end{align*}
  Letting $\varepsilon \to 0$
  shows that the distribution function of
  $\sup_{w \in \cW} \big|Z_n^T(w)\big|$
  is continuous,
  and therefore all of its quantiles exist.

  \proofparagraph{validity of the infeasible uniform confidence band}

  Under Assumption~\ref{ass:rates} and with a
  sufficiently slowly diverging sequence $R_n$,
  the strong approximation rate established in
  Lemma~\ref{lem:infeasible_t_statistic} is
  \begin{align*}
    \sup_{w \in \cW} \left| T_n(w) - Z_n^T(w) \right|
    &\lesssim_\P
    \frac{
      n^{-1/2} \log n
      + n^{-3/4} h^{-7/8} (\log n)^{3/8} R_n
      + n^{-2/3} h^{-1/2} (\log n)^{2/3}
      + n^{1/2} h^{p \wedge \beta}}
    {\Dl + 1/\sqrt{nh}} \\
    &\ll
    \frac{1}{\sqrt{\log n}}.
  \end{align*}
  So by Lemma~\ref{lem:slow_convergence},
  take $\varepsilon_n$ such that
  \begin{align*}
    \P \left(
      \sup_{w \in \cW} \left| T_n(w) - Z_n^T(w) \right|
      > \varepsilon_n
    \right)
    &\leq
    \varepsilon_n \sqrt{\log n}
  \end{align*}
  and $\varepsilon_n \sqrt{\log n} \to 0$.
  So by the previously established anti-concentration result,
  \begin{align*}
    &\P\left(
      \left|
      \widehat f_W(w) - f_W(w)
      \right|
      \leq
      q_{1-\alpha}
      \sqrt{\Sigma_n(w,w)}
      \textup{ for all }
      w \in \cW
    \right) \\
    &\quad=
    \P\left(
      \sup_{w \in \cW}
      \left| T_n(w) \right|
      \leq
      q_{1-\alpha}
    \right) \\
    &\quad\leq
    \P\left(
      \sup_{w \in \cW}
      \left| Z_n^T(w) \right|
      \leq
      q_{1-\alpha}
      + \varepsilon_n
    \right)
    + \P \left(
      \sup_{w \in \cW} \left| T_n(w) - Z_n^T(w) \right|
      > \varepsilon_n
    \right) \\
    &\quad\leq
    \P\left(
      \sup_{w \in \cW}
      \left|
      Z_n^T(w)
      \right|
      \leq
      q_{1-\alpha}
    \right)
    + \P\left(
      \left|
      \sup_{w \in \cW}
      \big| Z_n^T(w) \big|
      - q_{1-\alpha}
      \right|
      \leq \varepsilon_n
    \right)
    + \varepsilon_n \sqrt{\log n} \\
    &\quad\leq
    1 - \alpha
    + 2 \varepsilon_n \sqrt{\log n}.
  \end{align*}
  The lower bound follows analogously:
  \begin{align*}
    &\P\left(
      \left|
      \widehat f_W(w) - f_W(w)
      \right|
      \leq
      q_{1-\alpha}
      \sqrt{\Sigma_n(w,w)}
      \textup{ for all }
      w \in \cW
    \right) \\
    &\quad\geq
    \P\left(
      \sup_{w \in \cW}
      \left| Z_n^T(w) \right|
      \leq
      q_{1-\alpha}
      - \varepsilon_n
    \right)
    - \varepsilon_n \sqrt{\log n} \\
    &\quad\geq
    \P\left(
      \sup_{w \in \cW}
      \left|
      Z_n^T(w)
      \right|
      \leq
      q_{1-\alpha}
    \right)
    - \P\left(
      \left|
      \sup_{w \in \cW}
      \big| Z_n^T(w) \big|
      - q_{1-\alpha}
      \right|
      \leq \varepsilon_n
    \right)
    - \varepsilon_n \sqrt{\log n} \\
    &\quad\leq
    1 - \alpha
    - 2 \varepsilon_n \sqrt{\log n}.
  \end{align*}
  Finally we apply $\varepsilon_n \sqrt{\log n} \to 0$
  to see
  \begin{align*}
    \left|
    \P\left(
      \left|
      \widehat f_W(w) - f_W(w)
      \right|
      \leq
      q_{1-\alpha}
      \sqrt{\Sigma_n(w,w)}
      \textup{ for all }
      w \in \cW
    \right)
    - (1 - \alpha)
    \right|
    &\to 0.
  \end{align*}
\end{proof}

\subsubsection{Covariance estimation}

Before proving
Lemma~\ref{lem:covariance_estimation},
we provide the following useful
concentration inequality.
This is essentially a corollary of the
U-statistic concentration inequality given in
Theorem~3.3 in \citet{gine2000exponential}.

\begin{lemma}[A concentration inequality]
  \label{lem:dyadic_concentration}

  Let $X_{ij}$ be mutually independent random variables
  for $1 \leq i < j \leq n$
  taking values in a measurable space $\cX$.
  Let $h_1$, $h_2$
  be measurable functions from $\cX$ to $\R$
  satisfying the following
  for all $i$ and $j$.
  \begin{align*}
    \E\big[h_1(X_{ij})\big]
    &= 0,
    &\E\big[h_2(X_{ij})\big]
    &=0, \\
    \E\big[h_1(X_{ij})^2\big]
    &\leq \sigma^2,
    &\E\big[h_2(X_{ij})^2\big]
    &\leq \sigma^2, \\
    \big|h_1(X_{ij})\big|
    &\leq M,
    &\big|h_2(X_{ij})\big|
    &\leq M.
  \end{align*}
  Consider the sum
  \begin{align*}
    S_n
    &=
    \sum_{1 \leq i < j < r \leq n}
    h_1(X_{ij})
    h_2(X_{ir}).
  \end{align*}
  Then $S_n$ satisfies the concentration inequality
  \begin{align*}
    \P\big(
    |S_n| \geq t
    \big)
    &\leq
    C \exp\left(
      -\frac{1}{C}
      \min \left\{
      \frac{t^2}{n^3 \sigma^4},
      \frac{t}{\sqrt{n^3 \sigma^4}},
      \frac{t^{2/3}}{(n M \sigma)^{2/3}},
      \frac{t^{1/2}}{M}
      \right\}
    \right)
  \end{align*}
  for some universal constant
  $C > 0$
  and for all $t>0$.

\end{lemma}

\begin{proof}[Lemma~\ref{lem:dyadic_concentration}]

  We proceed in three main steps.
  Firstly we write $S_n$ as a second-order U-statistic
  where we use double indices instead of single indices.
  Then we use a decoupling result to introduce extra independence.
  Finally a concentration result is applied
  to the decoupled U-statistic.

  \proofparagraph{writing $S_n$ as a second-order U-statistic}

  Note that we can write $S_n$ as
  the second-order U-statistic
  \begin{align*}
    S_n
    &=
    \sum_{1 \leq i < j \leq n}
    \sum_{1 \leq q < r \leq n}
    h_{ijqr}
    (X_{ij}, X_{qr}),
  \end{align*}
  where
  \begin{align*}
    h_{ijqr}
    (a,b)
    &=
    h_1(a) h_2(b) \,
    \I\{j<r,\, q=i\}.
  \end{align*}
  Although this may look like a fourth-order
  U-statistic,
  it is in fact second-order.
  This is due to independence of the variables
  $X_{ij}$,
  and by treating $(i,j)$ as a single index.

  \proofparagraph{decoupling}

  By the decoupling result of Theorem~1
  from \citet{de1995decoupling}, there exists a universal
  constant $C_1 > 0$ satisfying
  \begin{align*}
    \P\big(
    |S_n| \geq t
    \big)
    &\leq
    C_1 \P\left(
      C_1 |\widetilde S_n| \geq t
    \right),
  \end{align*}
  where
  \begin{align*}
    \widetilde S_n
    &=
    \sum_{1 \leq i < j \leq n}
    \sum_{1 \leq q < r \leq n}
    h_{ijqr}
    (X_{ij}, X'_{qr}),
  \end{align*}
  with $(X'_{ij})$
  an independent copy of $(X_{ij})$.

  \proofparagraph{U-statistic concentration}

  The U-statistic kernel $h_{ijqr}(X_{ij}, X'_{qr})$
  is totally degenerate in the sense that
  \begin{align*}
    \E[h_{ijqr}(X_{ij}, X'_{qr}) \mid X_{ij}]
    &=
    \E[h_{ijqr}(X_{ij}, X'_{qr}) \mid X'_{qr}]
    = 0.
  \end{align*}
  Define and bound the following quantities:
  \begingroup
  \allowdisplaybreaks
  \begin{align*}
    A
    &=
    \max_{ijqr}
    \|h_{ijqr}(X_{ij}, X'_{qr})\|_\infty
    \leq M^2, \\
    B
    &=
    \max
    \left\{
    \left\|
    \sum_{1 \leq i < j \leq n}
    \E\Big[
    h_{ijqr}(X_{ij}, X'_{qr})^2
    \mid X_{ij}
    \Big]
    \right\|_\infty, \
    \left\|
    \sum_{1 \leq q < r \leq n}
    \E\Big[
    h_{ijqr}(X_{ij}, X'_{qr})^2
    \mid X'_{qr}
    \Big]
    \right\|_\infty
    \right\}^{1/2} \\
    &=
    \max
    \Bigg\{
    \left\|
    \sum_{1 \leq i < j \leq n}
    h_1(X_{ij})^2
    \E\big[
      h_2(X_{qr}')^2
      \big]
    \I\{j<r, q=i\}
    \right\|_\infty, \\
    &\qquad\qquad\quad
    \left\|
    \sum_{1 \leq q < r \leq n}
    \E\big[
      h_1(X_{ij})^2
      \big]
    h_2(X_{qr}')^2
    \I\{j<r, q=i\}
    \right\|_\infty
    \Bigg\}^{1/2} \\
    &\leq
    \max
    \left\{
    n^2 M^2 \sigma^2,
    n M^2 \sigma^2
    \right\}^{1/2}
    =
    n M \sigma, \\
    C
    &=
    \left(
      \sum_{1 \leq i < j \leq n}
      \sum_{1 \leq q < r \leq n}
      \E\big[
      h_{ijqr}(X_{ij}, X'_{qr})^2
      \big]
    \right)^{1/2} \\
    &=
    \left(
      \sum_{1 \leq i < j < r \leq n}
      \E\big[
        h_1(X_{ij})^2
        h_2(X_{ir}')^2
        \big]
    \right)^{1/2}
    \leq
    \sqrt{n^3 \sigma^4}, \\
    D
    &=
    \sup_{f,g} \Bigg\{
    \sum_{1 \leq i < j \leq n}
    \sum_{1 \leq q < r \leq n}
    \E\big[
    h_{ijqr}(X_{ij}, X'_{qr})
    f_{ij}(X_{ij})
    g_{qr}(X'_{qr})
    \big]
    \ : \\
    &\qquad\qquad\quad
    \sum_{1 \leq i < j \leq n}
    \E\big[f_{ij}(X_{ij})^2\big]
    \leq 1, \
    \sum_{1 \leq q < r \leq n}
    \E\big[g_{qr}(X'_{qr})^2\big]
    \leq 1
    \Bigg\} \\
    &=
    \sup_{f,g} \Bigg\{
    \sum_{1 \leq i < j < r \leq n}
    \E\big[
      h_1(X_{ij})
      f_{ij}(X_{ij})
      \big]
    \E\big[
    h_2(X'_{ir})
    g_{ir}(X'_{ir})
    \big]
    \ : \\
    &\qquad\qquad\quad
    \sum_{1 \leq i < j \leq n}
    \E\big[f_{ij}(X_{ij})^2\big]
    \leq 1, \
    \sum_{1 \leq q < r \leq n}
    \E\big[g_{qr}(X'_{qr})^2\big]
    \leq 1
    \Bigg\} \\
    &\leq
    \sup_{f,g} \Bigg\{
    \sum_{1 \leq i < j < r \leq n}
    \E\big[ h_1(X_{ij})^2 \big]^{1/2}
    \E\big[ f_{ij}(X_{ij})^2 \big]^{1/2}
    \E\big[ h_2(X_{ir}')^2 \big]^{1/2}
    \E\big[ g_{ir}(X'_{ir})^2 \big]^{1/2}
    \ : \\
    &\qquad\qquad\quad
    \sum_{1 \leq i < j \leq n}
    \E\big[f_{ij}(X_{ij})^2\big]
    \leq 1, \
    \sum_{1 \leq q < r \leq n}
    \E\big[g_{qr}(X'_{qr})^2\big]
    \leq 1
    \Bigg\} \\
    &\leq
    \sigma^2
    \sup_{f,g} \Bigg\{
    \sum_{1 \leq i < j \leq n}
    \E\big[ f_{ij}(X_{ij})^2 \big]^{1/2}
    \sum_{1 \leq r \leq n }
    \E\big[ g_{ir}(X'_{ir})^2 \big]^{1/2}
    \ : \\
    &\qquad\qquad\quad
    \sum_{1 \leq i < j \leq n}
    \E\big[f_{ij}(X_{ij})^2\big]
    \leq 1, \
    \sum_{1 \leq q < r \leq n}
    \E\big[g_{qr}(X'_{qr})^2\big]
    \leq 1
    \Bigg\} \\
    &\leq
    \sigma^2
    \sup_{f,g} \Bigg\{
    \Bigg(
    n^2
    \sum_{1 \leq i < j \leq n}
    \E\big[ f_{ij}(X_{ij})^2 \big]
    \Bigg)^{1/2}
    \Bigg(
    n
    \sum_{1 \leq r \leq n }
    \E\big[ g_{ir}(X'_{ir})^2 \big]
    \Bigg)^{1/2}
    \ : \\
    &\qquad\qquad\quad
    \sum_{1 \leq i < j \leq n}
    \E\big[f_{ij}(X_{ij})^2\big]
    \leq 1, \
    \sum_{1 \leq q < r \leq n}
    \E\big[g_{qr}(X'_{qr})^2\big]
    \leq 1
    \Bigg\}
    \leq
    \sqrt{n^3 \sigma^4}.
  \end{align*}
  \endgroup
  By Theorem~3.3 in
  \citet{gine2000exponential},
  we have that for some universal constant
  $C_2 > 0$
  and for all $t > 0$,
  \begin{align*}
    \P\left(
      |\widetilde S_n| \geq t
    \right)
    &\leq
    C_2 \exp\left(
      -\frac{1}{C_2}
      \min \left\{
      \frac{t^2}{C^2},
      \frac{t}{D},
      \frac{t^{2/3}}{B^{2/3}},
      \frac{t^{1/2}}{A^{1/2}}
      \right\}
    \right) \\
    &\leq
    C_2 \exp\left(
      -\frac{1}{C_2}
      \min \left\{
      \frac{t^2}{n^3 \sigma^4},
      \frac{t}{\sqrt{n^3 \sigma^4}},
      \frac{t^{2/3}}{(n M \sigma)^{2/3}},
      \frac{t^{1/2}}{M}
      \right\}
    \right)
  \end{align*}

  \proofparagraph{Conclusion}

  By the previous parts
  and absorbing constants into a new constant $C > 0$,
  we therefore have
  \begin{align*}
    \P\left(
      |S_n| \geq t
    \right)
    &\leq
    C_1 \P\left(
      C_1 |\widetilde S_n| \geq t
    \right) \\
    &\leq
    C_1 C_2 \exp\left(
      -\frac{1}{C_2}
      \min \left\{
      \frac{t^2}{n^3 \sigma^4 C_1^2},
      \frac{t}{\sqrt{n^3 \sigma^4 C_1}},
      \frac{t^{2/3}}{(n M \sigma C_1)^{2/3}},
      \frac{t^{1/2}}{M C_1^{1/2}}
      \right\}
    \right) \\
    &\leq
    C \exp\left(
      -\frac{1}{C}
      \min \left\{
      \frac{t^2}{n^3 \sigma^4},
      \frac{t}{\sqrt{n^3 \sigma^4}},
      \frac{t^{2/3}}{(n M \sigma)^{2/3}},
      \frac{t^{1/2}}{M}
      \right\}
    \right).
  \end{align*}
\end{proof}

\begin{proof}[Lemma~\ref{lem:covariance_estimation}]

  Throughout this proof we will write
  $k_{ij}$ for $k_h(W_{ij},w)$ and
  $k_{ij}'$ for $k_h(W_{ij},w')$,
  in the interest of brevity.
  Similarly we write $S_{ijr}$ to denote $S_{ijr}(w,w')$.
  The estimand and estimator are reproduced below for clarity.
  \begin{align*}
    \Sigma_n(w,w')
    &=
    \frac{2}{n(n-1)}
    \E[k_{ij} k_{ij}']
    + \frac{4(n-2)}{n(n-1)}
    \E[k_{ij} k_{ir}']
    - \frac{4n-6}{n(n-1)}
    \E[k_{ij}]
    \E[k_{ij}'] \\
    \widehat \Sigma_n(w,w')
    &=
    \frac{2}{n(n-1)}
    \frac{2}{n(n-1)}
    \sum_{i<j}
    k_{ij}
    k_{ij}'
    +
    \frac{4(n-2)}{n(n-1)}
    \frac{6}{n(n-1)(n-2)}
    \sum_{i<j<r}
    S_{ijr}
    -
    \frac{4n-6}{n(n-1)}
    \widehat f_W(w)
    \widehat f_W(w'),
  \end{align*}
  where
  \begin{align*}
    S_{ijr}
    &=
    \frac{1}{6}
    \Big(
    k_{ij}
    k_{ir}'
    + k_{ij}
    k_{jr}'
    + k_{ir}
    k_{ij}'
    + k_{ir}
    k_{jr}'
    + k_{jr}
    k_{ij}'
    + k_{jr}
    k_{ir}'
    \Big).
  \end{align*}
  We will prove uniform consistency of each of the three terms separately.

  \proofparagraph{uniform consistency of the
    $\widehat f_W(w) \widehat f_W(w')$ term}

  By boundedness of $f_W$ and
  Theorem~\ref{thm:uniform_consistency},
  $\widehat f_W$ is uniformly bounded in probability.
  Noting that
  $\E[\widehat f_W(w)] = \E[k_{ij}]$
  and by Lemma~\ref{lem:variance_bounds},
  \begin{align*}
    &\sup_{w,w' \in \cW}
    \left|
    \frac{
      \widehat f_W(w) \widehat f_W(w')
      - \E\big[k_{ij}\big] \E\big[k_{ij'}\big]}
    {\sqrt{\Sigma_n(w,w) + \Sigma_n(w',w')}}
    \right| \\
    &\quad=
    \sup_{w,w' \in \cW}
    \left|
    \frac{
      \widehat f_W(w) \widehat f_W(w')
      - \E\big[\widehat f_W(w)\big] \E\big[\widehat f_W(w')\big]}
    {\sqrt{\Sigma_n(w,w) + \Sigma_n(w',w')}}
    \right| \\
    &\quad\leq
    \sup_{w,w' \in \cW}
    \left|
    \frac{\widehat f_W(w) - \E\big[\widehat f_W(w)\big]}
    {\sqrt{\Sigma_n(w,w)}}
    \widehat f_W(w')
    + \frac{\widehat f_W(w') - \E\big[\widehat f_W(w')\big]}
    {\sqrt{\Sigma_n(w',w')}}
    \E\big[\widehat f_W(w)]
    \right| \\
    &\quad\lesssim_\P
    \sup_{w \in \cW}
    \left|
    \frac{\widehat f_W(w) - \E\big[\widehat f_W(w)\big]}
    {\sqrt{\Sigma_n(w,w)}}
    \right|
    \lesssim_\P
    \sup_{w \in \cW}
    \left|
    \frac{L_n(w)}
    {\sqrt{\Sigma_n(w,w)}}
    \right|
    + \sqrt{n^2h} \sup_{w \in \cW} \left| Q_n(w) \right|
    + \sqrt{n^2h} \sup_{w \in \cW} \left| E_n(w) \right| \\
    &\quad\lesssim_\P
    \sup_{w \in \cW}
    \left|
    \frac{L_n(w)}
    {\sqrt{\Sigma_n(w,w)}}
    \right|
    + \sqrt{n^2h} \frac{1}{n}
    + \sqrt{n^2h} \sqrt{\frac{\log n}{n^2h}}
    \lesssim_\P
    \sup_{w \in \cW}
    \left|
    \frac{L_n(w)}
    {\sqrt{\Sigma_n(w,w)}}
    \right|
    + \sqrt{\log n}.
  \end{align*}
  Now consider the function class
  \begin{align*}
    \cF
    &=
    \left\{
    a \mapsto
    \frac{
      \E\big[k_h(W_{ij},w) \mid A_i = a \big]
      - \E\big[k_h(W_{ij},w) \big]}
    {\sqrt{n \Sigma_n(w,w)}}:
    w \in \cW
    \right\},
  \end{align*}
  noting that
  \begin{align*}
    \frac{L_n(w)}
    {\Sigma_n(w,w)^{1/2}}
    &=
    \frac{1}{\sqrt n}
    \sum_{i=1}^n
    g_w(A_i)
  \end{align*}
  is an empirical process evaluated at
  $g_w \in \cF$.
  By the lower bound on $\Sigma_n(w,w)$
  from Lemma~\ref{lem:variance_bounds}
  with $\inf_\cW f_W(w) > 0$ and since $nh \gtrsim \log n$,
  the class $\cF$ has a constant envelope function
  given by $F(a) \lesssim \sqrt{nh}$.
  Clearly $M = \sup_a F(a) \lesssim \sqrt{nh}$.
  Also by definition of $\Sigma_n$
  and orthogonality of $L_n$, $Q_n$ and $E_n$,
  we have
  $\sup_{f \in \cF} \E[f(A_i)^2] \leq \sigma^2 = 1$.
  To verify a VC-type condition on $\cF$
  we need to establish the regularity of the process.
  By Lipschitz properties
  of $L_n$ and $\Sigma_n$
  derived in the proofs of Lemma~\ref{lem:uniform_concentration}
  and Theorem~\ref{thm:infeasible_ucb}
  respectively,
  we have
  \begin{align*}
    \left|
    \frac{L_n(w)}
    {\sqrt{\Sigma_n(w,w)}}
    - \frac{L_n(w')}
    {\sqrt{\Sigma_n(w',w')}}
    \right|
    &\lesssim
    \frac{\big|L_n(w) - L_n(w')\big|}
    {\sqrt{\Sigma_n(w,w)}}
    +
    \left| L_n(w') \right|
    \left|
    \frac{1}
    {\sqrt{\Sigma_n(w,w)}}
    - \frac{1}
    {\sqrt{\Sigma_n(w',w')}}
    \right| \\
    &\lesssim
    \sqrt{n^2h}
    |w-w'|
    +\left|
    \frac{\Sigma_n(w,w) - \Sigma_n(w',w')}
    {\Sigma_n(w,w)\sqrt{\Sigma_n(w',w')}}
    \right| \\
    &\lesssim
    \sqrt{n^2h}
    |w-w'|
    + (n^2h)^{3/2}
    \left|
    \Sigma_n(w,w) - \Sigma_n(w',w')
    \right| \\
    &\lesssim
    \sqrt{n^2h}
    |w-w'|
    + (n^2h)^{3/2}
    n^{-1} h^{-3}
    |w-w'|
    \lesssim
    n^4 |w-w'|
  \end{align*}
  uniformly over $w,w' \in \cW$.
  Therefore by compactness of $\cW$
  we have the covering number bound
  \begin{align*}
    N(\cF, \|\cdot\|_\infty, \varepsilon)
    &\lesssim
    N(\cW, |\cdot|, n^{-4} \varepsilon)
    \lesssim n^4 \varepsilon^{-1}.
  \end{align*}
  Thus by Lemma~\ref{lem:maximal_vc_inid},
  \begin{align*}
    \E \left[
      \sup_{w \in \cW}
      \left|
      \frac{L_n(w)}
      {\sqrt{\Sigma_n(w,w)}}
      \right|
    \right]
    &\lesssim
    \sqrt{\log n}
    + \frac{\sqrt{nh} \log n}{\sqrt{n}}
    \lesssim
    \sqrt{\log n}.
  \end{align*}
  Thus
  \begin{align*}
    \sup_{w,w' \in \cW}
    \left|
    \frac{
      \widehat f_W(w) \widehat f_W(w')
      - \E\big[k_{ij}\big] \E\big[k_{ij'}\big]}
    {\sqrt{\Sigma_n(w,w) + \Sigma_n(w',w')}}
    \right|
    &\lesssim_\P
    \sqrt{\log n}.
  \end{align*}

  \proofparagraph{decomposition of the $S_{ijr}$ term}

  We first decompose the $S_{ijr}$ term into two parts,
  and obtain a pointwise concentration result for each.
  This is extended to a uniform concentration result
  by considering the regularity of the covariance estimator process.
  Note that
  $\E[S_{ijr}] = \E[k_{ij} k_{ir}']$,
  and hence
  \begin{align*}
    \frac{6}{n(n-1)(n-2)}
    \sum_{i<j<r}
    \big(
    S_{ijr}
    - \E[k_{ij} k_{ir}']
    \big)
    &=
    \frac{6}{n(n-1)(n-2)}
    \sum_{i<j<r}
    S_{ijr}^{(1)}
    + \frac{6}{n(n-1)(n-2)}
    \sum_{i<j<r}
    S_{ijr}^{(2)}
  \end{align*}
  where
  \begin{align*}
    S_{ijr}^{(1)}
    &=
    S_{ijr} -
    \E[S_{ijr} \mid \bA_n],
    &S_{ijr}^{(2)}
    &=
    \E[S_{ijr} \mid \bA_n]
    - \E[S_{ijr}].
  \end{align*}

  \proofparagraph{pointwise concentration of the $S_{ijr}^{(1)}$ term}

  By symmetry in $i, j$ and $r$
  it is sufficient to consider only the first summand
  in the definition of $S_{ijr}$.
  By conditional independence properties,
  we have the decomposition
  \begin{align}
    \nonumber
    &\frac{6}{n(n-1)(n-2)}
    \sum_{i<j<r}
    \Big(
    k_{ij}k_{ir}'
    - \E[k_{ij}k_{ir}' \mid \bA_n]
    \Big) \\
    \nonumber
    &\quad=
    \frac{6}{n(n-1)(n-2)}
    \sum_{i<j<r}
    \Big(
    k_{ij}k_{ir}'
    - \E[k_{ij} \mid \bA_n]
    \E[k_{ir}' \mid \bA_n]
    \Big) \\
    \nonumber
    &\quad=
    \frac{6}{n(n-1)(n-2)}
    \sum_{i<j<r}
    \Big(
    \big(
    k_{ij}
    - \E[k_{ij} \mid \bA_n]
    \big)
    \big(
    k_{ir}'
    - \E[k_{ir}' \mid \bA_n]
    \big) \\
    \nonumber
    &\qquad+
    \big(
    k_{ij}
    - \E[k_{ij} \mid \bA_n]
    \big)
    \E[k_{ir}' \mid \bA_n]
    + \big(
    k_{ir}'
    - \E[k_{ir}' \mid \bA_n]
    \big)
    \E[k_{ij} \mid \bA_n]
    \Big) \\
    \label{eq:Sijr1_decomp1}
    &\quad=
    \frac{6}{n(n-1)(n-2)}
    \sum_{i<j<r}
    \Big(
    k_{ij}
    - \E[k_{ij} \mid \bA_n]
    \Big)
    \Big(
    k_{ir}'
    - \E[k_{ir}' \mid \bA_n]
    \Big) \\
    \label{eq:Sijr1_decomp2}
    &\qquad+
    \frac{2}{(n-1)(n-2)}
    \sum_{i=1}^{n-2}
    \sum_{j=i+1}^{n-1}
    \Big(
    k_{ij}
    - \E[k_{ij} \mid \bA_n]
    \Big)
    \cdot \frac{3}{n}
    \sum_{r=j+1}^n
    \E[k_{ir}' \mid \bA_n] \\
    \label{eq:Sijr1_decomp3}
    &\qquad+
    \frac{2}{(n-1)(n-2)}
    \sum_{i=1}^{n-2}
    \sum_{r=i+2}^n
    \Big(
    k_{ir}'
    - \E[k_{ir}' \mid \bA_n]
    \Big)
    \cdot \frac{3}{n}
    \sum_{j=i+1}^{r-1}
    \E[k_{ij} \mid \bA_n].
  \end{align}
  For the term in~(\ref{eq:Sijr1_decomp1}),
  note that conditional on $\bA_n$,
  we have that
  $k_{ij} - \E[k_{ij} \mid \bA_n]$
  are conditionally mean-zero
  and conditionally independent,
  as the only randomness is from $\bV_n$.
  Also
  $\Var[k_{ij} \mid \bA_n] \lesssim \sigma^2 := 1/h$
  and
  $|k_{ij}| \lesssim M := 1/h$
  uniformly.
  The same is true for $k_{ij}'$.
  Thus by Lemma~\ref{lem:dyadic_concentration}
  for some universal constant $C_1 > 0$:
  \begin{align*}
    &\P\left(
      \left|
      \sum_{i<j<r}
      \Big(
      k_{ij}
      - \E[k_{ij} \mid \bA_n]
      \Big)
      \Big(
      k_{ir}'
      - \E[k_{ir}' \mid \bA_n]
      \Big)
      \right|
      > t
      \biggm\vert \bA_n
    \right) \\
    &\quad\leq
    C_1 \exp\left(
      -\frac{1}{C_1}
      \min \left\{
      \frac{t^2}{n^3 \sigma^4},
      \frac{t}{\sqrt{n^3 \sigma^4}},
      \frac{t^{2/3}}{(n M \sigma)^{2/3}},
      \frac{t^{1/2}}{M}
      \right\}
    \right) \\
    &\quad\leq
    C_1 \exp\left(
      -\frac{1}{C_1}
      \min \left\{
      \frac{t^2 h^2}{n^3},
      \frac{t h}{\sqrt{n^3}},
      \frac{t^{2/3} h}{n^{2/3}},
      t^{1/2} h
      \right\}
    \right),
  \end{align*}
  and therefore
  with $t \geq 1$
  and since
  $nh \gtrsim \log n$,
  introducing and adjusting a new
  constant $C_2$ where necessary,
  \begin{align*}
    &\P\left(
      \left|
      \frac{6}{n(n-1)(n-2)}
      \sum_{i<j<r}
      \Big(
      k_{ij}
      - \E[k_{ij} \mid \bA_n]
      \Big)
      \Big(
      k_{ir}'
      - \E[k_{ir}' \mid \bA_n]
      \Big)
      \right|
      > t
      \frac{\log n}{\sqrt{n^3 h^2}}
      \Bigm\vert \bA_n
    \right) \\
    &\quad\leq
    \P\left(
      \left|
      \sum_{i<j<r}
      \Big(
      k_{ij}
      - \E[k_{ij} \mid \bA_n]
      \Big)
      \Big(
      k_{ir}'
      - \E[k_{ir}' \mid \bA_n]
      \Big)
      \right|
      > t
      n^{3/2} h^{-1} \log n / 24
      \Bigm\vert \bA_n
    \right) \\
    &\quad\leq
    C_2 \exp\left(
      -\frac{1}{C_2}
      \min \left\{
      (t \log n)^2,
      t \log n,
      (t \log n)^{2/3} (nh)^{1/3},
      (t nh \log n)^{1/2} n^{1/4}
      \right\}
    \right) \\
    &\quad\leq
    C_2 \exp\left(
      -\frac{1}{C_2}
      \min \left\{
      t \log n,
      t \log n,
      t^{2/3} \log n,
      t^{1/2} n^{1/4} \log n
      \right\}
    \right) \\
    &\quad=
    C_2 \exp\left(
      -\frac{t^{2/3} \log n}{C_2}
    \right)
    =
    C_2
    n^{-t^{2/3} / C_2}.
  \end{align*}
  Now for the term
  in~(\ref{eq:Sijr1_decomp2}),
  note that
  $\frac{3}{n} \sum_{r=j+1}^n \E[k_{ir}' \mid \bA_n]$
  is $\bA_n$-measurable and bounded uniformly in $i,j$.
  Also, using the previously established conditional variance
  and almost sure bounds on $k_{ij}$,
  Bernstein's inequality
  (Lemma~\ref{lem:bernstein})
  applied conditionally
  gives for some constant $C_3 > 0$
  \begin{align*}
    &\P\left(
      \Bigg|
      \frac{2}{(n-1)(n-2)}
      \sum_{i=1}^{n-2}
      \sum_{j=i+1}^{n-1}
      \Big(
      k_{ij}
      - \E[k_{ij} \mid \bA_n]
      \Big)
      \cdot \frac{3}{n}
      \sum_{r=j+1}^n
      \E[k_{ir}' \mid \bA_n]
      \Bigg|
      > t
      \sqrt{\frac{\log n}{n^2h}}
      \Bigm\vert \bA_n
    \right) \\
    &\qquad\leq
    2 \exp \left( -
      \frac{t^2 n^2 \log n / (n^2h)}
      {C_3/(2h) + C_3 t \sqrt{\log n / (n^2h)} / (2h)}
    \right) \\
    &\qquad=
    2 \exp \left( -
      \frac{t^2 \log n}
      {C_3/2 + C_3 t \sqrt{\log n / (n^2h)} / 2}
    \right)
    \leq
    2 \exp \left( -
      \frac{t^2 \log n}{C_3}
    \right)
    =
    2 n^{-t^2 / C_3}.
  \end{align*}
  The term in~(\ref{eq:Sijr1_decomp3})
  is controlled in exactly the same way.
  Putting these together, noting the symmetry in $i,j,r$
  and taking a marginal expectation,
  we obtain the unconditional pointwise concentration inequality
  \begin{align*}
    \P\left(
      \Bigg|
      \frac{6}{n(n-1)(n-2)}
      \sum_{i<j<r}
      S_{ijr}^{(1)}
      \Bigg|
      > t
      \frac{\log n}{\sqrt{n^3h^2}}
      + t \sqrt{\frac{\log n}{n^2h}}
    \right)
    &\leq
    C_2 n^{-t^{2/3} / C_2}
    + 4 n^{-t^2 / (4C_3)}.
  \end{align*}
  Multiplying by
  $\big(\Sigma_n(w,w) + \Sigma_n(w',w')\big)^{-1/2} \lesssim \sqrt{n^2h}$
  gives (adjusting constants if necessary)
  \begin{align*}
    \P\left(
      \Bigg|
      \frac{6}{n(n-1)(n-2)}
      \sum_{i<j<r}
      \frac{S_{ijr}^{(1)}}
      {\sqrt{\Sigma_n(w,w) + \Sigma_n(w',w')}}
      \Bigg|
      > t \frac{\log n}{\sqrt{nh}}
      + t \sqrt{\log n}
    \right)
    &\leq
    C_2 n^{-t^{2/3} / C_2}
    + 4 n^{-t^2 / (4C_3)}.
  \end{align*}

  \proofparagraph{pointwise concentration of the $S_{ijr}^{(2)}$ term}

  We apply the U-statistic concentration inequality from
  Lemma~\ref{lem:ustat_concentration}.
  Note that the terms
  $\E[S_{ijr} \mid \bA_n]$
  are permutation-symmetric functions of
  the random variables
  $A_i, A_j$ and $A_r$ only,
  making $S_{ijr}^{(2)}$ the summands of
  a (non-degenerate) mean-zero third-order U-statistic.
  While we could apply a third-order Hoeffding decomposition
  here to achieve degeneracy,
  it is unnecessary as Lemma~\ref{lem:ustat_concentration}
  is general enough to deal with the non-degenerate case directly.
  The quantity of interest here is
  \begin{align*}
    \frac{6}{n(n-1)(n-2)}
    \sum_{i<j<r}
    S_{ijr}^{(2)}
    &=
    \frac{6}{n(n-1)(n-2)}
    \sum_{i<j<r}
    \Big(
    \E[S_{ijr} \mid \bA_n]
    - \E[S_{ijr}]
    \Big).
  \end{align*}
  Note that by conditional independence,
  \begin{align*}
    \big|
    \E\big[
    k_{ij}k_{ir} \mid \bA_n
    \big]
    \big|
    &=
    \big|
    \E\big[
      k_{ij} \mid \bA_n
      \big]
    \E\big[
      k_{ir} \mid \bA_n
      \big]
    \big|
    \lesssim 1,
  \end{align*}
  and similarly for the other summands in $S_{ijr}$,
  giving the almost-sure bound
  $|S_{ijr}^{(2)}| \lesssim 1$.
  We also have
  \begin{align*}
    \Var\big[ \E[k_{ij} \mid A_i] \E[k_{ir}' \mid A_i] \big]
    &\lesssim
    \Var\big[\E[k_{ij} \mid A_i]\big]
    + \Var\big[\E[k_{ir}' \mid A_i]\big]
    \lesssim
    n \Var[L_n(w)] + n \Var[L_n(w')] \\
    &\lesssim
    n \Sigma_n(w,w) + n \Sigma_n(w',w')
  \end{align*}
  and similarly for the other summands in $S_{ijr}$,
  giving the conditional variance bound
  \begin{align*}
    \E[\E[S_{ijr}^{(2)} \mid A_i]^2] \lesssim
    n \Sigma_n(w,w) + n \Sigma_n(w',w').
  \end{align*}
  So Lemma~\ref{lem:ustat_concentration}
  and Lemma~\ref{lem:variance_bounds}
  give the pointwise concentration inequality
  \begin{align*}
    &\P\left(
      \Bigg|
      \frac{6}{n(n-1)(n-2)}
      \sum_{i<j<r}
      S_{ijr}^{(2)}
      \Bigg|
      > t \sqrt{\log n} \sqrt{\Sigma_n(w,w) + \Sigma_n(w',w')}
    \right) \\
    &\quad\leq
    4 \exp \left(
      - \frac{n t^2 (\Sigma_n(w,w) + \Sigma_n(w',w')) \log n}
      {C_4 (n\Sigma_n(w,w) + n\Sigma_n(w',w'))
        + C_4 t \sqrt{\Sigma_n(w,w) + \Sigma_n(w',w')}\sqrt{\log n}}
    \right) \\
    &\quad\leq
    4 \exp \left(
      - \frac{t^2 \log n}
      {C_4
        + C_4 t (\Sigma_n(w,w) + \Sigma_n(w',w'))^{-1/2} \sqrt{\log n} / n}
    \right)
    \leq
    4 \exp \left(
      - \frac{t^2 \log n}
      {C_4
        + C_4 t \sqrt{h}}
    \right)
    \leq
    4 n^{-t^2 / C_4}
  \end{align*}
  for some universal constant $C_4 > 0$
  (which may change from line to line),
  since the order of this U-statistic is fixed at three.

  \proofparagraph{concentration of the $S_{ijr}$ term on a mesh}

  Pick $\delta_n \to 0$
  with $\log 1/\delta_n \lesssim \log n$.
  Let $\cW_\delta$ be a $\delta_n$-covering of $\cW$
  with cardinality $O(1/\delta_n)$.
  Then $\cW_\delta \times \cW_\delta$
  is a $2\delta_n$-covering of $\cW \times \cW$
  with cardinality $O(1/\delta_n^2)$,
  under the Manhattan metric
  $d\big((w_1, w_1'), (w_2, w_2')\big)
  = |w_1 - w_2| + |w_1' - w_2'|$.
  By the previous parts,
  we have that for fixed $w$ and $w'$:
  \begin{align*}
    &\P\Bigg(
    \Bigg|
    \frac{6}{n(n-1)(n-2)}
    \sum_{i<j<r}
    \frac{S_{ijr}(w,w') - \E[S_{ijr}(w,w')]}
      {\sqrt{\Sigma_n(w,w) + \Sigma_n(w',w')}}
    \Bigg|
    > t \frac{\log n}{\sqrt{nh}}
    + 2t \sqrt{\log n}
    \Bigg) \\
    &\quad\leq
    C_2 n^{-t^{2/3} / C_2}
    + 4 n^{-t^2 / (4C_3)}
    + 4 n^{-t^2 / C_4}.
  \end{align*}
  Taking a union bound over $\cW_\delta \times \cW_\delta$,
  noting that $nh \gtrsim \log n$
  and adjusting constants gives
  \begin{align*}
    &\P\Bigg(
    \sup_{w, w' \in \cW_\delta}
    \Bigg|
    \frac{6}{n(n-1)(n-2)}
    \sum_{i<j<r}
    \frac{S_{ijr}(w,w') - \E[S_{ijr}(w,w')]}
      {\sqrt{\Sigma_n(w,w) + \Sigma_n(w',w')}}
    \Bigg|
    > t \sqrt{\log n}
    \Bigg) \\
    &\quad\lesssim
    \delta_n^{-2}
    \Big(
    C_2 n^{-t^{2/3} / C_2}
    + 4 n^{-t^2 / (4C_3)}
    + 4 n^{-t^2 / C_4}
    \Big)
    \lesssim
    \delta_n^{-2}
    n^{-t^{2/3} / C_5},
  \end{align*}
  for some constant $C_5 > 0$.

  \proofparagraph{regularity of the $S_{ijr}$ term}

  Next we bound the fluctuations in $S_{ijr}(w,w')$.
  Writing $k_{ij}(w)$ for $k_h(W_{ij},w)$,
  note that
  \begin{align*}
    \big|
    k_{ij}(w_1)
    k_{ir}(w_1')
    - k_{ij}(w_2)
    k_{ir}(w_2')
    \big|
    &\lesssim
    \frac{1}{h}
    \big| k_{ij}(w_1) - k_{ij}(w_2) \big|
    +
    \frac{1}{h}
    \big| k_{ir}(w_1') - k_{ir}(w_2') \big| \\
    &\lesssim
    \frac{1}{h^3}
    \Big(
    |w_1 - w_2|
    + |w_1' - w_2'|
    \Big),
  \end{align*}
  using the Lipschitz property of the kernel.
  Similarly for the other summands in $S_{ijr}$.
  Therefore
  \begin{align*}
    \sup_{|w_1-w_2| \leq \delta_n}
    \sup_{|w_1'-w_2'| \leq \delta_n}
    \big|
    S_{ijr}(w_1, w_1')
    - S_{ijr}(w_2, w_2')
    \big|
    &\lesssim
    \delta_n h^{-3}.
  \end{align*}
  Also as noted in the proof of Theorem~\ref{thm:infeasible_ucb},
  \begin{align*}
    \sup_{|w_1-w_2| \leq \delta_n}
    \sup_{|w_1'-w_2'| \leq \delta_n}
    \big|
    \Sigma_n(w_1,w_1')
    -
    \Sigma_n(w_2, w_2')
    \big|
    &\lesssim
    \delta_n n^{-1}h^{-3}.
  \end{align*}
  Therefore since $\sqrt{\Sigma_n(w,w)} \gtrsim \sqrt{n^2h}$
  and $|S_{ijr}| \lesssim h^{-2}$,
  using the elementary fact
  $\frac{a}{\sqrt b} - \frac{c}{\sqrt d}
  = \frac{a-c}{\sqrt b} + c \frac{d-b}{\sqrt{b d} \sqrt{b+d}}$,
  \begin{align*}
    &\sup_{|w_1-w_2| \leq \delta_n}
    \sup_{|w_1'-w_2'| \leq \delta_n}
    \left|
    \frac{S_{ijr}(w_1, w_1')}
    {\sqrt{\Sigma_n(w_1,w_1) + \Sigma_n(w_1',w_1')}}
    - \frac{S_{ijr}(w_2, w_2')}
    {\sqrt{\Sigma_n(w_2,w_2) + \Sigma_n(w_2',w_2')}}
    \right| \\
    &\quad\lesssim
    \delta_n h^{-3} \sqrt{n^2h}
    + h^{-2} \delta_n n^{-1} h^{-3} (n^2h)^{3/2}
    \lesssim
    \delta_n n h^{-5/2}
    + \delta_n n^{2} h^{-7/2}
    \lesssim
    \delta_n n^{6},
  \end{align*}
  where in the last line we use that
  $1/h \lesssim n$.

  \proofparagraph{uniform concentration of the $S_{ijr}$ term}

  By setting
  $\delta_n = n^{-6} \sqrt{\log n}$,
  the fluctuations can be at most $\sqrt{\log n}$,
  so we have for $t \geq 1$
  \begin{align*}
    &\P\Bigg(
    \sup_{w, w' \in \cW}
    \Bigg|
    \frac{6}{n(n-1)(n-2)}
    \sum_{i<j<r}
    \frac{S_{ijr}(w,w') - \E[S_{ijr}(w,w')]}
      {\sqrt{\Sigma_n(w,w) + \Sigma_n(w',w')}}
    \Bigg|
    > 2t \sqrt{\log n}
    \Bigg)
    \lesssim
    \delta_n^{-2}
    n^{-t^{2/3} / C_5}
    \lesssim
    n^{12-t^{2/3} / C_5}.
  \end{align*}
  This converges to zero for any sufficiently large $t$, so
  \begin{align*}
    \sup_{w, w' \in \cW}
    \Bigg|
    \frac{6}{n(n-1)(n-2)}
    \sum_{i<j<r}
    \frac{S_{ijr}(w,w') - \E[S_{ijr}(w,w')]}
    {\sqrt{\Sigma_n(w,w) + \Sigma_n(w',w')}}
    \Bigg|
    &\lesssim_\P
    \sqrt{\log n}.
  \end{align*}

  \proofparagraph{decomposition of the $k_{ij}k_{ij}'$ term}

  We move on to the final term in
  the covariance estimator.
  We have the decomposition
  \begin{align*}
    \frac{2}{n(n-1)}
    \sum_{i<j}
    \Big(
    k_{ij} k_{ij}'
    - \E\big[k_{ij} k_{ij}']
    \Big)
    &=
    \frac{2}{n(n-1)}
    \sum_{i<j}
    S_{ij}^{(1)}
    +
    \frac{2}{n(n-1)}
    \sum_{i<j}
    S_{ij}^{(2)},
  \end{align*}
  where
  \begin{align*}
    S_{ij}^{(1)}
    &=
    k_{ij} k_{ij}'
    - \E\big[ k_{ij} k_{ij}' \mid \bA_n \big],
    &S_{ij}^{(2)}
    &=
    \E\big[ k_{ij} k_{ij}' \mid \bA_n \big]
    - \E\big[ k_{ij} k_{ij}' \big].
  \end{align*}

  \proofparagraph{pointwise concentration of the $S_{ij}^{(1)}$ term}

  Conditioning on $\bA_n$,
  the variables $S_{ij}^{(1)}$
  are conditionally independent
  and conditionally mean-zero.
  They further satisfy
  $|S_{ij}^{(1)}| \lesssim h^{-2}$
  and the conditional variance bound
  $\E\big[\big( S_{ij}^{(1)} \big)^2 \mid \bA_n \big] \lesssim h^{-3}$.
  Therefore applying Bernstein's inequality
  (Lemma~\ref{lem:bernstein})
  conditional on $\bA_n$,
  we obtain the pointwise in $w,w'$
  concentration inequality
  \begin{align*}
    \P\left(
      \Bigg|
      \frac{2}{n(n-1)}
      \sum_{i<j}
      S_{ij}^{(1)}
      \Bigg|
      > t
      \sqrt{\frac{\log n}{n^2h^3}}
      \Bigm\vert \bA_n
    \right)
    &\leq
    2 \exp\left(
      - \frac{t^2 n^2 \log n / (n^2h^3)}
      {C_6 h^{-3} / 2 + C_6 t h^{-2} \sqrt{\log n / (n^2h^3)} / 2}
    \right) \\
    &\leq
    2 \exp\left(
      - \frac{t^2 \log n}
      {C_6 / 2 + C_6 t \sqrt{\log n / (n^2h)} / 2}
    \right) \\
    &\leq
    2 \exp\left(
      - \frac{t^2 \log n}{C_6}
    \right)
    =
    2 n^{-t^2 / C_6},
  \end{align*}
  where $C_6$ is a universal positive constant.

  \proofparagraph{pointwise concentration of the $S_{ij}^{(2)}$ term}

  We apply the U-statistic concentration inequality from
  Lemma~\ref{lem:ustat_concentration}.
  Note that $S_{ij}^{(2)}$
  are permutation-symmetric functions of
  the random variables
  $A_i$ and $A_j$ only,
  making them the summands of
  a (non-degenerate) mean-zero second-order U-statistic.
  Note that
  $\big|S_{ij}^{(2)}\big| \lesssim h^{-1}$
  and so trivially
  $\E\big[\E[S_{ij}^{(2)} \mid A_i ]^2 \big] \lesssim h^{-2}$.
  Thus by Lemma~\ref{lem:ustat_concentration},
  since the order of this U-statistic is fixed at two,
  for some universal positive constant $C_7$ we have
  \begin{align*}
    \P\left(
      \Bigg|
      \frac{2}{n(n-1)}
      \sum_{i<j}
      S_{ij}^{(2)}
      \Bigg|
      > t
      \sqrt{\frac{\log n}{nh^2}}
    \right)
    &\leq
    2 \exp\left(
      - \frac{t^2 n \log n / (nh^2)}
      {C_7 h^{-2} / 2 + C_7 t h^{-1} \sqrt{\log n / (nh^2)} / 2}
    \right) \\
    &\leq
    2 \exp\left(
      - \frac{t^2 \log n}
      {C_7 / 2 + C_7 t \sqrt{\log n / n} / 2}
    \right) \\
    &\leq
    2 \exp\left(
      - \frac{t^2 \log n}{C_7}
    \right)
    =
    2 n^{-t^2 / C_7}.
  \end{align*}

  \proofparagraph{concentration of the $k_{ij}k_{ij}'$ term on a mesh}

  As before, use a union bound
  on the mesh $\cW_\delta \times \cW_\delta$.
  \begin{align*}
    &\P\left(
      \sup_{w,w' \in \cW_\delta}
      \left|
      \frac{2}{n(n-1)}
      \sum_{i<j}
      \Big(
      k_{ij} k_{ij}'
      - \E\big[k_{ij} k_{ij}']
      \Big)
      \right|
      > t \sqrt{\frac{\log n}{n^2h^3}}
      + t \sqrt{\frac{\log n}{nh^2}}
    \right) \\
    &\quad\leq
    \P\left(
      \sup_{w,w' \in \cW_\delta}
      \Bigg|
      \frac{2}{n(n-1)}
      \sum_{i<j}
      S_{ij}^{(1)}
      \Bigg|
      > t
      \sqrt{\frac{\log n}{n^2h^3}}
    \right)
    + \P\left(
      \sup_{w,w' \in \cW_\delta}
      \Bigg|
      \frac{2}{n(n-1)}
      \sum_{i<j}
      S_{ij}^{(2)}
      \Bigg|
      > t
      \sqrt{\frac{\log n}{nh^2}}
    \right) \\
    &\quad\lesssim
    \delta_n^{-2} n^{-t^2 / C_6}
    + \delta_n^{-2} n^{-t^2 / C_7}.
  \end{align*}

  \proofparagraph{regularity of the $k_{ij}k_{ij}'$ term}

  Just as for the $S_{ijr}$ term,
  we have
  \begin{align*}
    \big|
    k_{ij}(w_1)
    k_{ij}(w_1')
    - k_{ij}(w_2)
    k_{ij}(w_2')
    \big|
    &\lesssim
    \frac{1}{h^3}
    \Big(
    |w_1 - w_2|
    + |w_1' - w_2'|
    \Big).
  \end{align*}

  \proofparagraph{uniform concentration of the $k_{ij}k_{ij}'$ term}

  By setting
  $\delta_n = h^3\sqrt{\log n / (nh^2)}$,
  the fluctuations can be at most $\sqrt{\log n / (nh^2)}$,
  so we have for $t \geq 1$
  \begin{align*}
    &\P\left(
      \sup_{w,w' \in \cW}
      \left|
      \frac{2}{n(n-1)}
      \sum_{i<j}
      \Big(
      k_{ij} k_{ij}'
      - \E\big[k_{ij} k_{ij}']
      \Big)
      \right|
      > t \sqrt{\frac{\log n}{n^2h^3}}
      + 2t \sqrt{\frac{\log n}{nh^2}}
    \right) \\
    &\quad\leq
    \P\left(
      \sup_{w,w' \in \cW_\delta}
      \left|
      \frac{2}{n(n-1)}
      \sum_{i<j}
      \Big(
      k_{ij} k_{ij}'
      - \E\big[k_{ij} k_{ij}']
      \Big)
      \right|
      > t \sqrt{\frac{\log n}{n^2h^3}}
      + t \sqrt{\frac{\log n}{nh^2}}
    \right) \\
    &\qquad+
    \P\left(
      \sup_{|w_1-w_2| \leq \delta_n}
      \sup_{|w_1'-w_2'| \leq \delta_n}
      \big|
      k_{ij}(w_1)
      k_{ij}(w_1')
      - k_{ij}(w_2)
      k_{ij}(w_2')
      \big|
      > t \sqrt{\frac{\log n}{nh^2}}
    \right) \\
    &\quad\lesssim
    \delta_n^{-2} n^{-t^2 / C_6}
    + \delta_n^{-2} n^{-t^2 / C_7}
    \lesssim
    n^{1-t^2 / C_6} h^{-4}
    + n^{1-t^2 / C_7} h^{-4}
    \lesssim
    n^{5-t^2 / C_8},
  \end{align*}
  where $C_8 > 0$ is a constant and
  in the last line we use $1/h \lesssim n$.
  This converges to zero for any sufficiently large $t$,
  so by Lemma~\ref{lem:variance_bounds} we have
  \begin{align*}
    \sup_{w,w' \in \cW}
    \left|
    \frac{2}{n(n-1)}
    \sum_{i<j}
    \frac{k_{ij} k_{ij}' - \E\big[k_{ij} k_{ij}']}
    {\sqrt{\Sigma_n(w,w) + \Sigma_n(w',w')}}
    \Big)
    \right|
    &\lesssim_\P
    \left(
      \sqrt{\frac{\log n}{n^2h^3}}
      + \sqrt{\frac{\log n}{nh^2}}
    \right)
    \sqrt{n^2h}
    \lesssim_\P
    \sqrt{\frac{n \log n}{h}}.
  \end{align*}

  \proofparagraph{conclusion}

  By the uniform bounds in probability
  derived in the previous parts,
  and with $nh \gtrsim \log n$,
  we conclude that
  \begin{align*}
    \sup_{w,w' \in \cW}
    \left|
    \frac{\widehat \Sigma_n(w,w') - \Sigma_n(w,w')}
    {\sqrt{\Sigma_n(w,w) + \Sigma_n(w',w')}}
    \right|
    &\leq
    \frac{2}{n(n-1)}
    \sup_{w,w' \in \cW}
    \left|
    \frac{2}{n(n-1)}
    \sum_{i<j}
    \frac{k_{ij} k_{ij}' - \E\big[k_{ij} k_{ij}']}
    {\sqrt{\Sigma_n(w,w) + \Sigma_n(w',w')}}
    \right| \\
    &\quad+
    \frac{4(n-2)}{n(n-1)}
    \sup_{w,w' \in \cW}
    \left|
    \frac{6}{n(n-1)(n-2)}
    \sum_{i<j<r}
    \frac{S_{ijr} - \E\big[k_{ij} k_{ir}']}
    {\sqrt{\Sigma_n(w,w) + \Sigma_n(w',w')}}
    \right| \\
    &\quad+
    \frac{4n-6}{n(n-1)}
    \sup_{w,w' \in \cW}
    \left|
    \frac{\widehat f_W(w) \widehat f_W(w') - \E[k_{ij}] \E[k_{ij}']}
    {\sqrt{\Sigma_n(w,w) + \Sigma_n(w',w')}}
    \right| \\
    &\lesssim_\P
    \sqrt{\frac{\log n}{n^3h}}
    + \frac{\sqrt{\log n}}{n}
    + \frac{\sqrt{\log n}}{n}
    \lesssim_\P
    \frac{\sqrt{\log n}}{n}.
  \end{align*}
\end{proof}

\begin{proof}[Lemma~\ref{lem:alternative_covariance_estimator}]

  Since there is no ambiguity,
  we may understand
  $k_{ij}$
  to mean $k_h(W_{ij},w)$
  when $i<j$
  and
  to mean $k_h(W_{ji},w)$
  when $j<i$.
  We use a prime
  to denote evaluation at $w'$
  rather than $w$.
  In this notation we may write
  \begin{align*}
    S_i(w)
    &=
    \frac{1}{n-1}
    \sum_{j \neq i}
    k_{ij}
  \end{align*}
  Let $\sum_{i \neq j \neq r}$
  indicate that all the indices are distinct.
  Then
  \begin{align*}
    \frac{4}{n^2}
    \sum_{i=1}^n
    S_i(w) S_i(w')
    &=
    \frac{4}{n^2}
    \sum_{i=1}^n
    \frac{1}{n-1}
    \sum_{j \neq i}
    k_{ij}
    \frac{1}{n-1}
    \sum_{r \neq i}
    k_{ir}'
    =
    \frac{4}{n^2(n-1)^2}
    \sum_{i=1}^n
    \sum_{j \neq i}
    \sum_{r \neq i}
    k_{ij}
    k_{ir}' \\
    &=
    \frac{4}{n^2(n-1)^2}
    \sum_{i=1}^n
    \sum_{j \neq i}
    \left(
      \sum_{r \neq i, r \neq j}
      k_{ij}
      k_{ir}'
      + k_{ij}
      k_{ij}'
    \right) \\
    &=
    \frac{4}{n^2(n-1)^2}
    \sum_{i \neq j \neq r}
    k_{ij}
    k_{ir}'
    + \frac{4}{n^2(n-1)^2}
    \sum_{i \neq j}
    k_{ij}
    k_{ij}' \\
    &=
    \frac{24}{n^2(n-1)^2}
    \sum_{i < j < r}
    S_{ijr}(w,w')
    + \frac{8}{n^2(n-1)^2}
    \sum_{i < j}
    k_{ij}
    k_{ij}' \\
    &=
    \widehat \Sigma_n(w,w')
    + \frac{4}{n^2(n-1)^2}
    \sum_{i < j}
    k_{ij}
    k_{ir}'
    + \frac{4n-6}{n(n-1)}
    \widehat f
    \widehat f',
  \end{align*}
  and the result follows.
\end{proof}

\subsubsection{Positive semi-definite covariance estimation}

\begin{proof}[Lemma~\ref{lem:SDP}]

  Firstly we prove that the true covariance function
  $\Sigma_n$
  is feasible for the optimization problem
  (\ref{eq:SDP}) in the sense that it satisfies the constraints.
  Clearly as a covariance function,
  $\Sigma_n$ is symmetric and positive semi-definite.
  The Lipschitz constraint is satisfied because
  as established in
  the proof of
  Theorem~\ref{thm:infeasible_ucb},
  \begin{align*}
    \big|
    \Sigma_n(w,w')
    -
    \Sigma_n(w, w'')
    \big|
    &\leq
    \frac{4}{nh^3}
    \Ck
    \CL
    |w'-w''|
  \end{align*}
  for all $w,w',w'' \in \cW$.
  Denote the (random) objective function
  in the optimization problem (\ref{eq:SDP}) by
  \begin{align*}
    \objective(M) = \sup_{w,w' \in \cW}
    \left|
    \frac{M(w,w') - \widehat\Sigma_n(w,w')}
    {\sqrt{\widehat \Sigma_n(w,w) + \widehat \Sigma_n(w',w')}}
    \right|.
  \end{align*}
  By Lemma~\ref{lem:covariance_estimation}
  with $w = w'$ we deduce that
  $\sup_{w \in \cW}
  \left|\frac{\widehat \Sigma_n(w,w)}{\Sigma_n(w,w)} - 1\right|
  \lesssim_\P \sqrt{h \log n}$
  and so
  \begin{align*}
    \objective(\Sigma_n)
    &= \sup_{w,w' \in \cW}
    \left|
    \frac{\widehat\Sigma_n(w,w') - \Sigma_n}
    {\sqrt{\Sigma_n(w,w) + \Sigma_n(w',w')}}
    \right|
    \sqrt{\frac{\Sigma_n(w,w) + \Sigma_n(w',w')}
      {\widehat \Sigma_n(w,w) + \widehat \Sigma_n(w',w')}} \\
    &\lesssim_\P
    \frac{\sqrt{\log n}}{n}
    \left(
      1 - \frac{\big|\widehat \Sigma_n(w,w) - \Sigma_n(w,w)\big|}
      {\Sigma_n(w,w)}
      - \frac{\big|\widehat \Sigma_n(w',w') - \Sigma_n(w',w')\big|}
      {\Sigma_n(w',w')}
    \right)^{-1/2} \\
    &\lesssim_\P
    \frac{\sqrt{\log n}}{n}
    \left(
      1 - \sqrt{h \log n}
    \right)^{-1/2}
    \lesssim_\P
    \frac{\sqrt{\log n}}{n}.
  \end{align*}
  Since the objective function
  is non-negative and because we have established
  at least one feasible function $M$ with
  an almost surely finite objective value,
  we can conclude the following.
  Let $\objective^* = \inf_M \objective(M)$,
  where the infimum is over feasible functions $M$.
  Then for all $\varepsilon > 0$
  there exists a feasible function $M_\varepsilon$ with
  $\objective(M_\varepsilon) \leq \objective^* + \varepsilon$,
  and we call such a solution $\varepsilon$-optimal.
  Let $\widehat \Sigma_n^+$ be an $n^{-1}$-optimal solution.
  Then
  \begin{align*}
    \objective(\widehat \Sigma_n^+)
    &\leq \objective^* + n^{-1}
    \leq \objective(\Sigma_n) + n^{-1}.
  \end{align*}
  Thus by the triangle inequality,
  \begin{align*}
    \sup_{w,w' \in \cW}
    \left|
    \frac{\widehat \Sigma_n^+(w,w') - \Sigma_n(w,w')}
    {\sqrt{\Sigma_n(w,w) + \Sigma_n(w',w')}}
    \right|
    &\leq
    \objective(\widehat \Sigma_n^+)
    + \objective(\Sigma_n)
    \leq 2 \, \objective(\Sigma_n) + n^{-1}
    \lesssim_\P
    \frac{\sqrt{\log n}}{n}.
  \end{align*}
\end{proof}

\begin{proof}[Lemma~\ref{lem:variance_estimator_bounds}]

  Since $\widehat \Sigma_n^+$ is positive semi-definite,
  we must have $\widehat \Sigma_n^+(w,w) \geq 0$.
  Now Lemma~\ref{lem:SDP}
  implies that for all $\varepsilon \in (0,1)$
  there exists a $C_\varepsilon$ such that
  \begin{align*}
    &\P\left(
      \Sigma_n(w,w) - C_\varepsilon \frac{\sqrt{\log n}}{n} \sqrt{\Sigma_n(w,w)}
      \leq
      \widehat \Sigma_n^+(w,w)
      \leq
      \Sigma_n(w,w) + C_\varepsilon \frac{\sqrt{\log n}}{n} \sqrt{\Sigma_n(w,w)}
      \text{ for all } w \in \cW
    \right) \\
    &\quad\geq 1-\varepsilon.
  \end{align*}
  Consider the function
  $g_a(t) = t - a \sqrt{t}$
  and note that it is increasing on $\{t \geq a^2/4\}$.
  Applying this with $t = \Sigma_n(w,w)$
  and $a = \frac{\sqrt{\log n}}{n}$,
  noting that by Lemma~\ref{lem:variance_bounds} we have
  $t = \Sigma_n(w,w) \gtrsim \frac{1}{n^2h}
  \gg \frac{\log n}{4n^2} = a^2/4$,
  shows that for $n$ large enough,
  \begin{align*}
    \inf_{w \in \cW} \Sigma_n(w,w)
    - \frac{\sqrt{\log n}}{n} \sqrt{\inf_{w \in \cW} \Sigma_n(w,w)}
    \lesssim_\P
    \inf_{w \in \cW}\widehat \Sigma_n^+(w,w), \\
    \sup_{w \in \cW}\widehat \Sigma_n^+(w,w)
    \lesssim_\P
    \sup_{w \in \cW} \Sigma_n(w,w)
    + \frac{\sqrt{\log n}}{n} \sqrt{\sup_{w \in \cW} \Sigma_n(w,w)}.
  \end{align*}
  Applying the bounds from Lemma~\ref{lem:variance_bounds}
  yields
  \begin{align*}
    \frac{\Dl^2}{n} + \frac{1}{n^2h}
    - \frac{\sqrt{\log n}}{n}
    \left( \frac{\Dl}{\sqrt n} + \frac{1}{\sqrt{n^2h}} \right)
    \lesssim_\P
    \inf_{w \in \cW}\widehat \Sigma_n^+(w,w), \\
    \sup_{w \in \cW}\widehat \Sigma_n^+(w,w)
    \lesssim_\P
    \frac{\Du^2}{n} + \frac{1}{n^2h}
    + \frac{\sqrt{\log n}}{n}
    \left( \frac{\Du}{\sqrt n} + \frac{1}{\sqrt{n^2h}} \right)
  \end{align*}
  and so
  \begin{align*}
    \frac{\Dl^2}{n} + \frac{1}{n^2h}
    \lesssim_\P
    \inf_{w \in \cW}\widehat \Sigma_n^+(w,w)
    \leq
    \sup_{w \in \cW}\widehat \Sigma_n^+(w,w)
    \lesssim_\P
    \frac{\Du^2}{n} + \frac{1}{n^2h}.
  \end{align*}
\end{proof}

\subsubsection{Feasible uniform confidence bands}

\begin{proof}[Lemma~\ref{lem:studentized_t_statistic}]

  \begin{align*}
    \sup_{w \in \cW}
    \left| \widehat T_n(w) - T_n(w) \right|
    &=
    \sup_{w \in \cW}
    \left\{
    \left|
    \widehat f_W(w) - f_W(w)
    \right|
    \cdot
    \left|
    \frac{1}
    {\sqrt{\widehat\Sigma_n^+(w,w)}}
    -
    \frac{1}{\sqrt{\Sigma_n(w,w)}}
    \right|
    \right\} \\
    &\leq
    \sup_{w \in \cW}
    \left|
    \frac{\widehat f_W(w) - \E\big[\widehat f_W(w)\big]}
    {\sqrt{\Sigma_n(w,w)}}
    + \frac{\E\big[\widehat f_W(w)\big] - f_W(w)}
    {\sqrt{\Sigma_n(w,w)}}
    \right|
    \cdot
    \sup_{w \in \cW}
    \left|
    \frac{\widehat\Sigma_n^+(w,w) - \Sigma_n(w,w)}
    {\sqrt{\Sigma_n(w,w) \widehat\Sigma_n^+(w,w)}}
    \right|
  \end{align*}
  Now from the proof of Lemma~\ref{lem:covariance_estimation} we have that
  $\sup_{w \in \cW} \left|
  \frac{\widehat f_W(w) - \E\big[\widehat f_W(w)\big]}
  {\sqrt{\Sigma_n(w,w)}} \right|
  \lesssim_\P \sqrt{\log n}$,
  while Lemma~\ref{lem:bias} gives
  $\sup_{w \in \cW} \big| \E\big[\widehat f_W(w)\big] - f_W(w) \big|
  \lesssim h^{p \wedge \beta}$.
  By Lemma~\ref{lem:variance_bounds} we have
  $\sup_{w \in \cW} \Sigma_n(w,w)^{-1/2}
  \lesssim \frac{1}{\Dl/\sqrt{n} + 1/\sqrt{n^2h}}$,
  and similarly Lemma~\ref{lem:variance_estimator_bounds} gives
  $\sup_{w \in \cW} \widehat \Sigma_n^+(w,w)^{-1/2}
  \lesssim_\P \frac{1}{\Dl/\sqrt{n} + 1/\sqrt{n^2h}}$,
  Thus, applying Lemma~\ref{lem:SDP} to control the
  covariance estimation error,
  \begin{align*}
    \sup_{w \in \cW}
    \left| \widehat T_n(w) - T_n(w) \right|
    &\lesssim_\P
    \left(
      \sqrt{\log n} + \frac{h^{p \wedge \beta}}{\Dl/\sqrt{n} + 1/\sqrt{n^2h}}
    \right)
    \frac{\sqrt{\log n}}{n}
    \frac{1}{\Dl/\sqrt{n} + 1/\sqrt{n^2h}} \\
    &\lesssim_\P
    \sqrt{\frac{\log n}{n}}
    \left(
      \sqrt{\log n} + \frac{\sqrt n h^{p \wedge \beta}}
      {\Dl + 1/\sqrt{nh}}
    \right)
    \frac{1}{\Dl + 1/\sqrt{nh}}.
  \end{align*}
\end{proof}

\begin{proof}[Lemma~\ref{lem:distributional_approx_feasible_gaussian}]

  Firstly note that the process $\widehat Z_n^T$ exists
  by noting that $\widehat \Sigma_n^+(w,w')$ and therefore also
  $\frac{\widehat \Sigma_n^+(w,w')}
  {\sqrt{\widehat \Sigma_n^+(w,w) \widehat \Sigma_n^+(w',w')}}$
  are positive semi-definite
  functions and appealing to the
  Kolmogorov consistency theorem \citep{gine2021mathematical}.
  To obtain the desired Kolmogorov--Smirnov result we discretize and
  use the Gaussian--Gaussian comparison result found in
  Lemma~3.1 in \citet{chernozhukov2013gaussian}.

  \proofparagraph{bounding the covariance discrepancy}

  Define the maximum discrepancy in the (conditional) covariances
  of $\widehat Z_n^T$ and $Z_n^T$ by
  \begin{align*}
    \Delta
    &\vcentcolon=
    \sup_{w, w' \in \cW}
    \left|
    \frac{\widehat \Sigma_n^+(w,w')}
    {\sqrt{\widehat \Sigma_n^+(w,w) \widehat \Sigma_n^+(w',w')}}
    - \frac{\Sigma_n(w,w')}
    {\sqrt{\Sigma_n(w,w) \Sigma_n(w',w')}}
    \right|.
  \end{align*}
  This random variable can be bounded in probability
  in the following manner.
  First note that by the Cauchy--Schwarz inequality
  for covariances,
  $|\Sigma_n(w,w')| \leq
  \sqrt{\Sigma_n(w,w) \Sigma_n(w',w')}$.
  Hence
  \begin{align*}
    \Delta
    &\leq
    \sup_{w, w' \in \cW}
    \left|
    \frac{\widehat \Sigma_n^+(w,w') - \Sigma_n(w,w')}
    {\sqrt{\widehat \Sigma_n^+(w,w) \widehat \Sigma_n^+(w',w')}}
    \right|
    + \sup_{w, w' \in \cW}
    \left|
    \frac{\sqrt{\widehat \Sigma_n^+(w,w) \widehat \Sigma_n^+(w',w')}
      - \sqrt{\Sigma_n(w,w) \Sigma_n(w',w')}}
    {\sqrt{\widehat \Sigma_n^+(w,w) \widehat \Sigma_n^+(w',w')}}
    \right| \\
    &\leq
    \sup_{w, w' \in \cW}
    \left\{
    \sqrt{\frac{\Sigma_n(w,w) + \Sigma_n(w',w')}
      {\widehat \Sigma_n^+(w,w) \widehat \Sigma_n^+(w',w')}}
    \left|
    \frac{\widehat \Sigma_n^+(w,w') - \Sigma_n(w,w')}
    {\sqrt{\Sigma_n(w,w) + \Sigma_n(w',w')}}
    \right|
    \right\} \\
    &\quad+
    \sup_{w, w' \in \cW}
    \left|
    \frac{\widehat \Sigma_n^+(w,w)\widehat \Sigma_n^+(w',w')
      - \Sigma_n(w,w) \Sigma_n(w',w')}
    {\sqrt{\widehat \Sigma_n^+(w,w) \widehat \Sigma_n^+(w',w')
        \Sigma_n(w,w) \Sigma_n(w',w')}}
    \right|.
  \end{align*}
  For the first term note that
  $\inf_{w \in \cW} \widehat \Sigma_n^+(w,w)
  \gtrsim \frac{\Dl^2}{n} + \frac{1}{n^2h}$
  by Lemma~\ref{lem:variance_estimator_bounds} and
  $\sup_{w \in \cW}
  \left|\frac{\widehat \Sigma_n(w,w)}{\Sigma_n(w,w)} - 1\right|
  \lesssim_\P \sqrt{h \log n}$
  by the proof of Lemma~\ref{lem:SDP}.
  Thus by Lemma~\ref{lem:SDP},
  \begin{align*}
    \sup_{w, w' \in \cW}
    \left\{
    \sqrt{\frac{\Sigma_n(w,w) + \Sigma_n(w',w')}
      {\widehat \Sigma_n^+(w,w) \widehat \Sigma_n^+(w',w')}}
    \left|
    \frac{\widehat \Sigma_n^+(w,w') - \Sigma_n(w,w')}
    {\sqrt{\Sigma_n(w,w) + \Sigma_n(w',w')}}
    \right|
    \right\}
    &\lesssim_\P
    \frac{\sqrt{\log n}}{n}
    \frac{1}{\Dl/\sqrt{n} + 1/\sqrt{n^2h}} \\
    &\lesssim_\P
    \sqrt{\frac{\log n}{n}}
    \frac{1}{\Dl + 1/\sqrt{nh}}.
  \end{align*}
  For the second term, we have by the same bounds
  \begin{align*}
    &\sup_{w, w' \in \cW}
    \left|
    \frac{\widehat \Sigma_n^+(w,w) \widehat \Sigma_n^+(w',w')
      - \Sigma_n(w,w) \Sigma_n(w',w')}
    {\sqrt{\widehat \Sigma_n^+(w,w) \widehat \Sigma_n^+(w',w')
        \Sigma_n(w,w) \Sigma_n(w',w')}}
    \right| \\
    &\quad\leq
    \sup_{w, w' \in \cW}
    \left\{
    \frac{\big| \widehat \Sigma_n^+(w,w) - \Sigma_n(w,w)\big|
      \widehat \Sigma_n^+(w',w')}
    {\sqrt{\widehat \Sigma_n^+(w,w) \widehat \Sigma_n^+(w',w')
        \Sigma_n(w,w) \Sigma_n(w',w')}}
    + \frac{\big| \widehat \Sigma_n^+(w',w') - \Sigma_n(w',w')\big|
      \Sigma_n(w,w)}
    {\sqrt{\widehat \Sigma_n^+(w,w) \widehat \Sigma_n^+(w',w')
        \Sigma_n(w,w) \Sigma_n(w',w')}}
    \right\} \\
    &\quad\leq
    \sup_{w, w' \in \cW}
    \left\{
    \frac{\big| \widehat \Sigma_n^+(w,w) - \Sigma_n(w,w)\big|}
    {\sqrt{\Sigma_n(w,w)}}
    \frac{\sqrt{\widehat \Sigma_n^+(w',w')}}
    {\sqrt{\widehat \Sigma_n^+(w,w) \Sigma_n(w',w')}}
    \right\} \\
    &\qquad+
    \sup_{w, w' \in \cW}
    \left\{
    \frac{\big| \widehat \Sigma_n^+(w',w') - \Sigma_n(w',w')\big|}
    {\sqrt{\Sigma_n(w',w')}}
    \frac{\sqrt{\Sigma_n(w,w)}}
    {\sqrt{\widehat \Sigma_n^+(w,w) \widehat \Sigma_n^+(w',w')}}
    \right\} \\
    &\quad\lesssim_\P
    \sqrt{\frac{\log n}{n}}
    \frac{1}{\Dl + 1/\sqrt{nh}}.
  \end{align*}
  Therefore
  \begin{align*}
    \Delta \lesssim_\P
    \sqrt{\frac{\log n}{n}}
    \frac{1}{\Dl + 1/\sqrt{nh}}.
  \end{align*}

  \proofparagraph{Gaussian comparison on a mesh}

  Let $\cW_\delta$ be a $\delta_n$-covering of $\cW$
  with cardinality $O(1/\delta_n)$,
  where $1/\delta_n$ is at most polynomial in $n$.
  The scaled (conditionally) Gaussian
  processes $Z_n^T$ and $\widehat Z_n^T$
  both have pointwise (conditional) variances of 1.
  Therefore by Lemma~3.1 in \citet{chernozhukov2013gaussian},
  \begin{align*}
    \sup_{t \in \R}
    \left|
    \P\left(
      \sup_{w \in \cW_\delta}
      Z_n^T(w)
      \leq t
    \right)
    - \P\left(
      \sup_{w \in \cW_\delta}
      \widehat Z_n^T(w)
      \leq t
      \Bigm\vert \bW_n
    \right)
    \right|
    &\lesssim
    \Delta^{1/3}
    \Big(
    1 \vee \log \frac{1}{\Delta \delta_n}
    \Big)^{2/3}
  \end{align*}
  uniformly in the data. By the previous part and
  since $x (\log 1/x)^2$ is increasing on $\big(0, e^{-2}\big)$,
  \begin{align*}
    \sup_{t \in \R}
    \left|
    \P\left(
      \sup_{w \in \cW_\delta}
      Z_n^T(w)
      \leq t
    \right)
    - \P\left(
      \sup_{w \in \cW_\delta}
      \widehat Z_n^T(w)
      \leq t
      \Bigm\vert \bW_n
    \right)
    \right|
    &\lesssim_\P
    \left(
      \sqrt{\frac{\log n}{n}}
      \frac{1}{\Dl + 1/\sqrt{nh}}
    \right)^{1/3}
    (\log n)^{2/3} \\
    &\lesssim_\P
    \frac{n^{-1/6}(\log n)^{5/6}}
    {\Dl^{1/3} + (nh)^{-1/6}}.
  \end{align*}

  \proofparagraph{trajectory regularity of $Z_n^T$}

  During the proof of
  Theorem~\ref{thm:infeasible_ucb}
  we established that
  $Z_n^T$ satisfies the following
  trajectory regularity property:
  \begin{align*}
    \E\left[
      \sup_{|w-w'| \leq \delta_n}
      \big| Z_n^T(w) - Z_n^T(w') \big|
    \right]
    &\lesssim
    n h^{-1}
    \sqrt{\delta_n \log n},
  \end{align*}
  whenever $1/\delta_n$
  is at most polynomial in $n$.

  \proofparagraph{conditional $L^2$ regularity of $\widehat Z_n^T$}

  By Lemma~\ref{lem:SDP},
  with $nh \gtrsim \log n$,
  we have
  uniformly in $w,w'$,
  \begin{align*}
    \big|
    \widehat \Sigma_n^+(w,w')
    - \widehat \Sigma_n^+(w,w)
    \big|
    &\lesssim
    n^{-1} h^{-3} |w-w'|.
  \end{align*}
  Taking
  $\delta_n \leq n^{-2} h^2$,
  Lemma~\ref{lem:variance_estimator_bounds}
  gives
  \begin{align*}
    \inf_{|w-w'| \leq \delta_n}
    \widehat \Sigma_n^+(w,w')
    \gtrsim
    \frac{\Dl^2}{n}
    + \frac{1}{n^2h}
    - n^{-1} h^{-3} \delta_n
    \gtrsim
    \frac{\Dl^2}{n}
    + \frac{1}{n^2h}
    - \frac{1}{n^3h}
    \gtrsim
    \frac{\Dl^2}{n}
    + \frac{1}{n^2h}.
  \end{align*}
  The conditional $L^2$
  regularity of $\widehat Z_n^T$ is
  \begin{align*}
    \E\left[
      \big(
      \widehat Z_n^T(w) - \widehat Z_n^T(w')
      \big)^2
      \bigm\vert \bW_n
    \right]
    &=
    2 - 2
    \frac{\widehat \Sigma_n^+(w,w')}
    {\sqrt{\widehat \Sigma_n^+(w,w) \widehat \Sigma_n^+(w',w')}}.
  \end{align*}
  Applying the same elementary result as for $Z_n^T$
  in the proof of Theorem~\ref{thm:infeasible_ucb} yields
  \begin{align*}
    \E\left[
      \big(
      \widehat Z_n^T(w) - \widehat Z_n^T(w')
      \big)^2
      \bigm\vert \bW_n
    \right]
    &\lesssim_\P
    n^2 h^{-2} |w-w'|.
  \end{align*}
  Thus the conditional semimetric
  induced by $\widehat Z_n^T$ on $\cW$ is
  \begin{align*}
    \widehat\rho(w,w')
    &\vcentcolon=
    \E\left[
      \big(
      \widehat Z_n^T(w) - \widehat Z_n^T(w')
      \big)^2
      \bigm\vert \bW_n
    \right]^{1/2}
    \lesssim_\P
    n h^{-1} \sqrt{|w-w'|}.
  \end{align*}

  \proofparagraph{conditional trajectory regularity of $\widehat Z_n^T$}

  Just as for $Z_n^T$
  in the proof of Theorem~\ref{thm:infeasible_ucb}
  we apply Lemma~\ref{lem:gaussian_process_maximal},
  this time conditionally,
  to obtain that also
  \begin{align*}
    \E\left[
      \sup_{|w-w'| \leq \delta_n}
      \left| \widehat Z_n^T(w) - \widehat Z_n^T(w') \right|
      \Bigm\vert \bW_n
    \right]
    &\lesssim_\P
    n h^{-1}
    \sqrt{\delta_n \log n},
  \end{align*}
  whenever $1/\delta_n$
  is at most polynomial in $n$.

  \proofparagraph{uniform Gaussian comparison}

  Now we use the trajectory regularity properties to
  extend the Gaussian--Gaussian comparison result from a finite mesh
  to all of $\cW$.
  Write the previously established
  approximation rate as
  \begin{align*}
    r_n
    &=
    \frac{n^{-1/6}(\log n)^{5/6}}
    {\Dl^{1/3} + (nh)^{-1/6}}.
  \end{align*}
  Take $\varepsilon_n > 0$ and observe that
  uniformly in $t \in \R$,
  \begin{align*}
    &\P\left(
      \sup_{w \in \cW}
      \big| \widehat Z_n^T(w) \big|
      \leq t
      \Bigm\vert \bW_n
    \right) \\
    &\quad\leq
    \P\left(
      \sup_{w \in \cW_\delta}
      \big| \widehat Z_n^T(w) \big|
      \leq t + \varepsilon_n
      \Bigm\vert \bW_n
    \right)
    + \P\left(
      \sup_{|w-w'| \leq \delta_n}
      \left|
      \widehat Z_n^T(w)
      - \widehat Z_n^T(w')
      \right|
      \geq \varepsilon_n
      \Bigm\vert \bW_n
    \right) \\
    &\quad\leq
    \P\left(
      \sup_{w \in \cW_\delta}
      \big| Z_n^T(w) \big|
      \leq t + \varepsilon_n
    \right)
    + O_\P(r_n)
    + \P\left(
      \sup_{|w-w'| \leq \delta_n}
      \left|
      \widehat Z_n^T(w)
      - \widehat Z_n^T(w')
      \right|
      \geq \varepsilon_n
      \Bigm\vert \bW_n
    \right) \\
    &\quad\leq
    \P\left(
      \sup_{w \in \cW}
      \big| Z_n^T(w) \big|
      \leq t + 2\varepsilon_n
    \right)
    + O_\P(r_n) \\
    &\qquad+
    \P\left(
      \sup_{|w-w'| \leq \delta_n}
      \left|
      Z_n^T(w)
      - Z_n^T(w')
      \right|
      \geq \varepsilon_n
    \right)
    + \P\left(
      \sup_{|w-w'| \leq \delta_n}
      \left|
      \widehat Z_n^T(w)
      - \widehat Z_n^T(w')
      \right|
      \geq \varepsilon_n
      \Bigm\vert \bW_n
    \right) \\
    &\quad\leq
    \P\left(
      \sup_{w \in \cW}
      \big| Z_n^T(w) \big|
      \leq t + 2\varepsilon_n
    \right)
    + O_\P(r_n)
    + O_\P(\varepsilon_n^{-1} nh^{-1} \sqrt{\delta_n \log n}) \\
    &\quad\leq
    \P\left(
      \sup_{w \in \cW}
      \big| Z_n^T(w) \big|
      \leq t
    \right)
    + \P\left(
      \left|
      \sup_{w \in \cW}
      \big| Z_n^T(w) \big|
      - t
      \right|
      \leq 2\varepsilon_n
    \right)
    + O_\P(r_n)
    + O_\P(\varepsilon_n^{-1} nh^{-1} \sqrt{\delta_n \log n}).
  \end{align*}
  The converse inequality is obtained
  analogously as follows:
  \begin{align*}
    &\P\left(
      \sup_{w \in \cW}
      \big| \widehat Z_n^T(w) \big|
      \leq t
      \Bigm\vert \bW_n
    \right) \\
    &\quad\geq
    \P\left(
      \sup_{w \in \cW_\delta}
      \big| \widehat Z_n^T(w) \big|
      \leq t - \varepsilon_n
      \Bigm\vert \bW_n
    \right)
    - \P\left(
      \sup_{|w-w'| \leq \delta_n}
      \left|
      \widehat Z_n^T(w)
      - \widehat Z_n^T(w')
      \right|
      \geq \varepsilon_n
      \Bigm\vert \bW_n
    \right) \\
    &\quad\geq
    \P\left(
      \sup_{w \in \cW_\delta}
      \big| Z_n^T(w) \big|
      \leq t - \varepsilon_n
    \right)
    - O_\P(r_n)
    - \P\left(
      \sup_{|w-w'| \leq \delta_n}
      \left|
      \widehat Z_n^T(w)
      - \widehat Z_n^T(w')
      \right|
      \geq \varepsilon_n
      \Bigm\vert \bW_n
    \right) \\
    &\quad\geq
    \P\left(
      \sup_{w \in \cW}
      \big| Z_n^T(w) \big|
      \leq t - 2\varepsilon_n
    \right)
    - O_\P(r_n) \\
    &\qquad-
    \P\left(
      \sup_{|w-w'| \leq \delta_n}
      \left|
      Z_n^T(w)
      - Z_n^T(w')
      \right|
      \geq \varepsilon_n
    \right)
    - \P\left(
      \sup_{|w-w'| \leq \delta_n}
      \left|
      \widehat Z_n^T(w)
      - \widehat Z_n^T(w')
      \right|
      \geq \varepsilon_n
      \Bigm\vert \bW_n
    \right) \\
    &\quad\geq
    \P\left(
      \sup_{w \in \cW}
      \big| Z_n^T(w) \big|
      \leq t - 2\varepsilon_n
    \right)
    - O_\P(r_n)
    - O_\P(\varepsilon_n^{-1} nh^{-1} \sqrt{\delta_n \log n}) \\
    &\quad\geq
    \P\left(
      \sup_{w \in \cW}
      \big| Z_n^T(w) \big|
      \leq t
    \right)
    - \P\left(
      \left|
      \sup_{w \in \cW}
      \big| Z_n^T(w) \big|
      - t
      \right|
      \leq 2\varepsilon_n
    \right)
    - O_\P(r_n)
    - O_\P(\varepsilon_n^{-1} nh^{-1} \sqrt{\delta_n \log n}).
  \end{align*}
  Combining these uniform upper and lower bounds gives
  \begin{align*}
    &\sup_{t \in \R}
    \left|
    \P\left(
      \sup_{w \in \cW}
      \big| \widehat Z_n^T(w) \big|
      \leq t
      \Bigm\vert \bW_n
    \right)
    -
    \P\left(
      \sup_{w \in \cW}
      \big| Z_n^T(w) \big|
      \leq t
    \right)
    \right| \\
    &\qquad\lesssim_\P
    \sup_{t \in \R}
    \P\left(
      \left|
      \sup_{w \in \cW}
      \big| Z_n^T(w) \big|
      - t
      \right|
      \leq 2\varepsilon_n
    \right)
    + r_n
    + \varepsilon_n^{-1} nh^{-1/2} \delta_n^{1/2} \sqrt{\log n}.
  \end{align*}
  To bound the remaining term,
  we apply the anti-concentration
  result for $Z_n^T$
  from the proof of
  Theorem~\ref{thm:infeasible_ucb}:
  \begin{align*}
    \sup_{t \in \R}
    \P\left(
      \left|
      \sup_{w \in \cW}
      \big| Z_n^T(w) \big|
      - t
      \right|
      \leq \varepsilon
    \right)
    &\lesssim
    \varepsilon
    \sqrt{\log n}.
  \end{align*}
  Therefore
  \begin{align*}
    &\sup_{t \in \R}
    \left|
    \P\left(
      \sup_{w \in \cW}
      \big| \widehat Z_n^T(w) \big|
      \leq t
      \Bigm\vert \bW_n
    \right)
    -
    \P\left(
      \sup_{w \in \cW}
      \big| Z_n^T(w) \big|
      \leq t
    \right)
    \right| \\
    &\qquad\lesssim_\P
    \varepsilon_n \sqrt{\log n}
    + r_n
    + \varepsilon_n^{-1} nh^{-1/2} \delta_n^{1/2} \sqrt{\log n}.
  \end{align*}
  Taking $\varepsilon = r_n / \sqrt{\log n}$
  and then $\delta_n = n^{-2} h r_n^2 \varepsilon_n^2 / \log n$
  yields
  \begin{align*}
    \left|
    \P\left(
      \sup_{w \in \cW}
      \big| \widehat Z_n^T(w) \big|
      \leq t
      \Bigm\vert \bW_n
    \right)
    -
    \P\left(
      \sup_{w \in \cW}
      \big| Z_n^T(w) \big|
      \leq t
    \right)
    \right|
    &\lesssim_\P
    r_n =
    \frac{n^{-1/6}(\log n)^{5/6}}
    {\Dl^{1/3} + (nh)^{-1/6}}.
  \end{align*}
\end{proof}

\begin{proof}[Lemma~\ref{lem:feasible_gaussian_approx}]

  \proofparagraph{Kolmogorov--Smirnov approximation}

  Let $Z_n^T$ and $\widehat Z_n^T$ be defined
  as in the proof of
  Lemma~\ref{lem:distributional_approx_feasible_gaussian}.
  Write
  \begin{align*}
    r_n
    &=
    \frac{n^{-1/6}(\log n)^{5/6}}
    {\Dl^{1/3} + (nh)^{-1/6}}
  \end{align*}
  for the rate of approximation from
  Lemma~\ref{lem:distributional_approx_feasible_gaussian}.
  Then for any $\varepsilon_n > 0$ and
  uniformly in $t \in \R$:
  \begin{align*}
    &\P\left(
      \sup_{w \in \cW}
      \left|
      \widehat Z_n^T(w)
      \right|
      \leq t
      \Bigm\vert \bW_n
    \right)
    \leq
    \P\left(
      \sup_{w \in \cW}
      \left|
      Z_n^T(w)
      \right|
      \leq t
    \right)
    +
    O_\P(r_n) \\
    &\quad\leq
    \P\left(
      \sup_{w \in \cW}
      \left|
      Z_n^T(w)
      \right|
      \leq t - \varepsilon_n
    \right)
    +
    \P\left(
      \left|
      \sup_{w \in \cW}
      \big|
      Z_n^T(w)
      \big|
      -t
      \right|
      \leq \varepsilon_n
    \right)
    +
    O_\P(r_n) \\
    &\quad\leq
    \P\left(
      \sup_{w \in \cW}
      \left| \widehat T_n(w) \right|
      \leq t
    \right)
    +
    \P\left(
      \sup_{w \in \cW}
      \left| \widehat T_n(w) - Z_n^T(w) \right|
      \geq \varepsilon_n
    \right)
    + \P\left(
      \left|
      \sup_{w \in \cW}
      \big|
      Z_n^T(w)
      \big|
      -t
      \right|
      \leq \varepsilon_n
    \right)
    +
    O_\P(r_n) \\
    &\quad\leq
    \P\left(
      \sup_{w \in \cW}
      \left| \widehat T_n(w) \right|
      \leq t
    \right)
    +
    \P\left(
      \sup_{w \in \cW}
      \left| \widehat T_n(w) - Z_n^T(w) \right|
      \geq \varepsilon_n
    \right)
    + \varepsilon_n \sqrt{\log n}
    + O_\P(r_n),
  \end{align*}
  where in the last line we used the anti-concentration result
  from Lemma~\ref{lem:anti-concentration}
  applied to $Z_n^T$,
  as in the proof of
  Lemma~\ref{lem:distributional_approx_feasible_gaussian}.
  The corresponding lower bound is as follows:
  \begin{align*}
    &\P\left(
      \sup_{w \in \cW}
      \left|
      \widehat Z_n^T(w)
      \right|
      \leq t
      \Bigm\vert \bW_n
    \right)
    \geq
    \P\left(
      \sup_{w \in \cW}
      \left|
      Z_n^T(w)
      \right|
      \leq t
    \right)
    -
    O_\P(r_n) \\
    &\quad\geq
    \P\left(
      \sup_{w \in \cW}
      \left|
      Z_n^T(w)
      \right|
      \leq t + \varepsilon_n
    \right)
    -
    \P\left(
      \left|
      \sup_{w \in \cW}
      \big|
      Z_n^T(w)
      \big|
      -t
      \right|
      \leq \varepsilon_n
    \right)
    -
    O_\P(r_n) \\
    &\quad\geq
    \P\left(
      \sup_{w \in \cW}
      \left| \widehat T_n(w) \right|
      \leq t
    \right)
    -
    \P\left(
      \sup_{w \in \cW}
      \left| \widehat T_n(w) - Z_n^T(w) \right|
      \geq \varepsilon_n
    \right)
    - \P\left(
      \left|
      \sup_{w \in \cW}
      \big|
      Z_n^T(w)
      \big|
      -t
      \right|
      \leq \varepsilon_n
    \right)
    -
    O_\P(r_n) \\
    &\quad\geq
    \P\left(
      \sup_{w \in \cW}
      \left| \widehat T_n(w) \right|
      \leq t
    \right)
    -
    \P\left(
      \sup_{w \in \cW}
      \left| \widehat T_n(w) - Z_n^T(w) \right|
      \geq \varepsilon_n
    \right)
    - \varepsilon_n \sqrt{\log n}
    - O_\P(r_n).
  \end{align*}

  \proofparagraph{t-statistic approximation}

  To control the remaining term,
  note that by
  Lemma~\ref{lem:infeasible_t_statistic}
  and Lemma~\ref{lem:studentized_t_statistic},
  \begin{align*}
    \sup_{w \in \cW}
    \left| \widehat T_n(w) - Z_n^T(w) \right|
    &\leq
    \sup_{w \in \cW}
    \left| \widehat T_n(w) - T_n(w) \right|
    + \sup_{w \in \cW}
    \left| T_n(w) - Z_n^T(w) \right| \\
    &\lesssim_\P
    \sqrt{\frac{\log n}{n}}
    \left(
      \sqrt{\log n} + \frac{\sqrt n h^{p \wedge \beta}}
      {\Dl + 1/\sqrt{nh}}
    \right)
    \frac{1}{\Dl + 1/\sqrt{nh}} \\
    &\qquad+
    \frac{
      n^{-1/2} \log n
      + n^{-3/4} h^{-7/8} (\log n)^{3/8} R_n
      + n^{-2/3} h^{-1/2} (\log n)^{2/3}
      + n^{1/2} h^{p \wedge \beta}}
    {\Dl + 1/\sqrt{nh}}
  \end{align*}
  and denote this last quantity by $r_n'$.
  Then for any $\varepsilon_n \gg r_n'$,
  we have
  \begin{align*}
    \sup_{t \in \R}
    \left|
    \P\left(
      \sup_{w \in \cW}
      \left| \widehat T_n(w) \right|
      \leq t
    \right)
    - \P\left(
      \sup_{w \in \cW}
      \left|
      \widehat Z_n^T(w)
      \right|
      \leq t
      \Bigm\vert \bW_n
    \right)
    \right|
    &\lesssim_\P
    \varepsilon_n \sqrt{\log n}
    + r_n
    + o(1).
  \end{align*}

  \proofparagraph{rate analysis}

  This rate can be made $o_\P(1)$
  by some appropriate choice of $\varepsilon_n$
  whenever
  $r_n \to 0$
  and
  $r_n' \sqrt{\log n} \to 0$,
  by Lemma~\ref{lem:slow_convergence}, along with
  a sufficiently slowly diverging sequence $R_n$.
  Explicitly, we require the following.
  \begin{align*}
    \frac{n^{-1/2} (\log n)^{3/2}}{\Dl + 1/\sqrt{nh}}
    &\to 0,
    &\frac{h^{p \wedge \beta} \log n}{\Dl^2 + (nh)^{-1}}
    &\to 0, \\
    \frac{n^{-1/2} (\log n)^{3/2}}
    {\Dl + 1/\sqrt{nh}}
    &\to 0,
    &\frac{n^{-3/4} h^{-7/8} (\log n)^{7/8}}
    {\Dl + 1/\sqrt{nh}}
    &\to 0, \\
    \frac{n^{-2/3} h^{-1/2} (\log n)^{7/6}}
    {\Dl + 1/\sqrt{nh}}
    &\to 0,
    &\frac{n^{1/2} h^{p \wedge \beta} (\log n)^{1/2}}
    {\Dl + 1/\sqrt{nh}}
    &\to 0, \\
    \frac{n^{-1/6}(\log n)^{5/6}}
    {\Dl^{1/3} + (nh)^{-1/6}}
    &\to 0.
  \end{align*}
  Using the fact that $h \lesssim n^{-\varepsilon}$
  for some $\varepsilon > 0$
  and removing trivial statements leaves us with
  \begin{align*}
    \frac{n^{-3/4} h^{-7/8} (\log n)^{7/8}}
    {\Dl + 1/\sqrt{nh}}
    &\to 0,
    &\frac{n^{1/2} h^{p \wedge \beta} (\log n)^{1/2}}
    {\Dl + 1/\sqrt{nh}}
    &\to 0.
  \end{align*}
  Now we analyze these based on the degeneracy type
  and verify that they hold under
  Assumption~\ref{ass:rates}.
  \begin{enumerate}[label=(\roman*)]

    \item No degeneracy:
      if $\Dl > 0$ then we need
      \begin{align*}
        n^{-3/4} h^{-7/8} (\log n)^{7/8}
        &\to 0,
        &n^{1/2} h^{p \wedge \beta} (\log n)^{1/2}
        &\to 0.
      \end{align*}
      These reduce to
      $n^{-6/7} \log n \ll h
      \ll (n \log n)^{-\frac{1}{2(p \wedge \beta)}}$.

    \item Partial or total degeneracy:
      if $\Dl = 0$ then we need
      \begin{align*}
        n^{-1/4} h^{-3/8} (\log n)^{7/8}
        &\to 0,
        &n h^{(p \wedge \beta) + 1/2} (\log n)^{1/2}
        &\to 0.
      \end{align*}
      These reduce to
      $n^{-2/3} (\log n)^{7/3} \ll h
      \ll (n^2 \log n)^{-\frac{1}{2(p \wedge \beta) + 1}}$.

  \end{enumerate}
\end{proof}

\begin{proof}[Theorem~\ref{thm:ucb}]

  \proofparagraph{existence of the conditional quantile}

  We argue as in the proof of
  Lemma~\ref{lem:distributional_approx_feasible_gaussian},
  now also conditioning on the data.
  In particular, using the anti-concentration result from
  Lemma~\ref{lem:anti-concentration},
  the regularity property of $\widehat Z_n^T$
  and the Gaussian process maximal inequality from
  Lemma~\ref{lem:gaussian_process_maximal},
  we see that for any $\varepsilon > 0$,
  \begin{align*}
    \sup_{t \in \R}
    \P\left(
      \left|
      \sup_{w \in \cW}
      \big| \widehat Z_n^T(w) \big|
      - t
      \right|
      \leq 2\varepsilon
      \Bigm\vert \bW_n
    \right)
    &\leq
    8 \varepsilon
    \left(
      1 + \E\left[
        \sup_{w \in \cW}
        \big| \widehat Z_n^T(w) \big|
        \Bigm\vert \bW_n
      \right]
    \right) \\
    &\lesssim \varepsilon \sqrt{\log n}.
  \end{align*}
  Thus letting $\varepsilon \to 0$
  shows that the conditional distribution function of
  $\sup_{w \in \cW} \big|\widehat Z_n^T(w)\big|$
  is continuous,
  and therefore all of its conditional quantiles exist.

  \proofparagraph{validity of the confidence band}

  Define the following (conditional) distribution functions.
  \begin{align*}
    F_Z(t \mid \bW_n)
    &=
    \P\left(
      \sup_{w \in \cW}
      \left| \widehat Z_n^T(w) \right|
      \leq t
      \Bigm\vert \bW_n
    \right),
    &F_T(t)
    &=
    \P\left(
      \sup_{w \in \cW}
      \left| \widehat T_n(w) \right|
      \leq t
    \right),
  \end{align*}
  along with their well-defined right-quantile functions,
  \begin{align*}
    F_Z^{-1}(p \mid \bW_n)
    &=
    \sup
    \big\{
    t \in \R
    \, : \,
    F_Z(t \mid \bW_n)
    = p
    \big\},
    &F_T^{-1}(p)
    &=
    \sup
    \big\{
    t \in \R
    \, : \,
    F_T(t)
    = p
    \big\}.
  \end{align*}
  Note that
  $t \leq F_Z^{-1}(p \mid \bW_n)$
  if and only if
  $F_Z(t \mid \bW_n) \leq p$.
  Take $\alpha \in (0,1)$ and
  define the quantile
  $\widehat q_{1-\alpha} = F_Z^{-1}(1-\alpha \mid \bW_n)$,
  so that
  $F_Z(\widehat q_{1-\alpha} \mid \bW_n) = 1-\alpha$.
  By Lemma~\ref{lem:feasible_gaussian_approx},
  we have that
  \begin{align*}
    \sup_{t \in \R}
    \big|
    F_Z(t \mid \bW_n) - F_T(t)
    \big|
    &=
    o_\P(1)
  \end{align*}
  Thus by Lemma~\ref{lem:slow_convergence},
  this can be replaced by
  \begin{align*}
    \P\left(
      \sup_{t \in \R} \big| F_Z(t \mid \bW_n) - F_T(t) \big|
      > \varepsilon_n
    \right)
    &\leq \varepsilon_n
  \end{align*}
  for some $\varepsilon_n \to 0$.
  Therefore
  \begin{align*}
    \P\left(
      \sup_{w \in \cW}
      \left|
      \widehat T_n(w)
      \right|
      \leq
      \widehat q_{1-\alpha}
    \right)
    &=
    \P\left(
      \sup_{w \in \cW}
      \left|
      \widehat T_n(w)
      \right|
      \leq
      F_Z^{-1}(1-\alpha \mid \bW_n)
    \right) \\
    &=
    \P\left(
      F_Z\left(
        \sup_{w \in \cW}
        \left|
        \widehat T_n(w)
        \right|
        \Bigm\vert \bW_n
      \right)
      \leq
      1 - \alpha
    \right) \\
    &\leq
    \P\left(
      F_T\left(
        \sup_{w \in \cW}
        \left|
        \widehat T_n(w)
        \right|
      \right)
      \leq
      1 - \alpha + \varepsilon_n
    \right)
    + \varepsilon_n
    \leq 1 - \alpha + 3\varepsilon_n,
  \end{align*}
  where in the last line we used the fact that for any
  real-valued random variable $X$ with distribution function $F$,
  we have
  $\big|\P\big(F(X) \leq t\big) - t\big| \leq \Delta$,
  where $\Delta$ is the size of the
  largest jump discontinuity in $F$.
  By taking an expectation and uniform integrability,
  $\sup_{t \in \R} \big| F_Z(t) - F_T(t) \big| = o(\varepsilon_n)$.
  Since $F_Z$ has no jumps,
  we must have $\Delta \leq \varepsilon_n$ for $F_T$.
  Finally a lower bound is constructed in an analogous manner,
  giving
  \begin{align*}
    \P\left(
      \sup_{w \in \cW}
      \left| \widehat T_n(w) \right|
      \leq
      \widehat q_{1-\alpha}
    \right)
    &\geq
    1 - \alpha - 3\varepsilon_n.
  \end{align*}
  Here ends the proof of the theorem.
\end{proof}

\subsubsection{Counterfactual dyadic density estimation}
\begin{proof}[Lemma~\ref{lem:counterfactual_bias}]

  Writing
  $k_{ij} = k_h(W_{ij}^1, w)$,
  $\psi_i = \psi(X_i^1)$,
  $\widehat\psi_i = \widehat\psi(X_i^1)$
  and $\kappa_{ij} = \kappa(X_i^0, X_i^1, X_j^1)$,
  \begin{align*}
    \E\big[\widehat f_W^{\,1 \triangleright 0}(w)\big]
    &=
    \E\left[
      \frac{2}{n(n-1)}
      \sum_{i<j}
      \widehat \psi_i
      \widehat \psi_j
      k_{ij}
    \right] \\
    &=
    \frac{2}{n(n-1)(n-2)}
    \sum_{i < j}
    \sum_{r \notin \{i,j\}}
    \E\left[
      k_{ij}
      \Big(
      \psi_i
      \psi_j
      +\psi_i
      \kappa_{r j}
      +\psi_j
      \kappa_{r i}
      \Big)
    \right]
    + O\left(\frac{1}{n}\right) \\
    &=
    \E\left[
      k_{ij}
      \psi_i
      \psi_j
    \right]
    + O\left(\frac{1}{n}\right)
    =
    \E\big[
      \psi_i
      \psi_j
      \E\left[
        k_h(W_{ij}^1, w)
        \mid X_i^1, X_j^1
      \right]
      \big]
    + O\left(\frac{1}{n}\right) \\
    &=
    \E\big[
      \psi_i
      \psi_j
      f_{W \mid XX}^1(w \mid X_i^1, X_j^1)
      + O_\P\left( h^{p \wedge \beta} \right)
      \big]
    + O\left(\frac{1}{n}\right)
    =
    f_W^{1 \triangleright 0}(w)
    + O\left( h^{p \wedge \beta} + \frac{1}{n}\right)
  \end{align*}
  uniformly in $w \in \cW$, following the proof of
  Lemma~\ref{lem:bias} by H{\"o}lder-continuity
  of $f_{W \mid XX}^1(\,\cdot \mid x_1, x_2)$.
\end{proof}

\begin{proof}[Lemma~\ref{lem:counterfactual_hoeffding}]

  \begin{align*}
    \widehat f_W^{\,1 \triangleright 0}(w)
    &=
    \frac{2}{n(n-1)}
    \sum_{i < j}
    \widehat \psi_i
    \widehat \psi_j
    k_{ij} \\
    &=
    \frac{2}{n(n-1)}
    \sum_{i < j}
    \left(
      \psi_i
      + \frac{1}{n}
      \sum_{r=1}^n \kappa_{r i}
    \right)
    \left(
      \psi_j
      + \frac{1}{n}
      \sum_{r=1}^n \kappa_{r j}
    \right)
    k_{ij}
    + O_\P\left(\frac{1}{n}\right) \\
    &=
    \frac{2}{n(n-1)}
    \sum_{i < j}
    \psi_i
    \psi_j
    k_{ij}
    + \frac{2}{n(n-1)}
    \sum_{i < j}
    \psi_i
    \frac{1}{n}
    \sum_{r \notin \{i,j\}}^n \kappa_{r j}
    k_{ij}
    + \frac{2}{n(n-1)}
    \sum_{i < j}
    \psi_j
    \frac{1}{n}
    \sum_{r \notin \{i,j\}}^n \kappa_{r i}
    k_{ij} \\
    &\quad+
    O_\P\left(\frac{1}{n}\right) \\
    &=
    \frac{2}{n(n-1)(n-2)}
    \sum_{i < j}
    \sum_{r \notin \{i,j\}}
    k_{ij}
    \Big(
    \psi_i
    \psi_j
    +\psi_i
    \kappa_{r j}
    +\psi_j
    \kappa_{r i}
    \Big)
    + O_\P\left(\frac{1}{n}\right) \\
    &=
    \frac{6}{n(n-1)(n-2)}
    \sum_{i < j < r}
    v_{ijr}
    + O_\P\left(\frac{1}{n}\right)
  \end{align*}
  where
  \begin{align*}
    v_{ijr}
    &=
    \frac{1}{3}
    k_{ij} \Big(\psi_i \psi_j +\psi_i \kappa_{r j} +\psi_j \kappa_{r i} \Big)
    + \frac{1}{3}
    k_{ir} \Big(\psi_i \psi_r +\psi_i \kappa_{jr} +\psi_r \kappa_{ji} \Big)
    + \frac{1}{3}
    k_{jr} \Big(\psi_j \psi_r +\psi_j \kappa_{ir} +\psi_r \kappa_{ij} \Big)
  \end{align*}
  So by the Hoeffding decomposition for third-order U-statistics,
  \begin{align*}
    \widehat f_W^{\,1 \triangleright 0}(w)
    &=
    u
    + \frac{3}{n}
    \sum_{i=1}^n
    u_i
    + \frac{6}{n(n-1)}
    \sum_{i=1}^{n-1}
    \sum_{j=i+1}^n
    u_{ij}
    + \frac{6}{n(n-1)(n-2)}
    \sum_{i=1}^{n-2}
    \sum_{j=i+1}^{n-1}
    \sum_{r=j+1}^n
    u_{ijr} \\
    &\quad+
    \frac{6}{n(n-1)(n-2)}
    \sum_{i=1}^{n-2}
    \sum_{j=i+1}^{n-1}
    \sum_{r=j+1}^n
    \big(v_{ijr} - u_{ijr}\big)
    + O_\P\left( \frac{1}{n} \right) \\
    &=
    \E\big[\widehat f_W^{\,1 \triangleright 0}(w) \big]
    + L_n^{1 \triangleright 0}(w)
    + Q_n^{1 \triangleright 0}(w)
    + T_n^{1 \triangleright 0}(w)
    + E_n^{1 \triangleright 0}(w)
    + O_\P\left( \frac{1}{n} \right).
  \end{align*}
  Noting that $\psi_i$, $\kappa_{ij}$
  and $\E[k_{ij} \mid A_i^1, A_j^1]$
  are all bounded and that
  $\E[k_{ij} \mid A_i^1, A_j^1]$
  is Lipschitz in $w$,
  we deduce by
  Lemma~\ref{lem:uprocess_maximal}
  and Proposition~2.3 of
  \citet{arcones1993limit} that
  $\sup_{w \in \cW} |Q_n^{1 \triangleright 0}(w)
  + T_n^{1 \triangleright 0}(w)| \lesssim_\P \frac 1n$.
\end{proof}

\begin{proof}[Lemma~\ref{lem:counterfactual_uniform_consistency}]

  By Lemma~\ref{lem:maximal_vc_inid} we have
  $\sup_{w \in \cW} \big|L_n^{1 \triangleright 0}(w)\big|
  \lesssim_\P \frac{1}{\sqrt n}$.
  Note that in the proof of Lemma~\ref{lem:counterfactual_hoeffding}
  the terms $v_{ijr} - u_{ijr}$ depend only on $V_{ij}, V_{ir}$ and $V_{jr}$
  after conditioning on $\bA_n^1$, $\bX_n^0$ and $\bX_n^1$.
  Thus $E_n^{1 \triangleright 0}(w)$ is a degenerate second-order
  U-statistic and so
  $\sup_{w \in \cW} \big|E_n^{1 \triangleright 0}(w)\big|
  \lesssim_\P \sqrt{\frac{\log n}{n^2h}}$
  by Lemma~\ref{lem:uprocess_maximal}.
\end{proof}

\begin{proof}[Lemma~\ref{lem:counterfactual_sa}]

  Note that
  $L_n^{1 \triangleright 0}(w)
  = \frac 3n \sum_{i=1}^n l_i^{1 \triangleright 0}(w)$
  where $l_i^{1 \triangleright 0}(w)$ depends only on
  $A_i^1$, $X_i^0$ and $X_i^1$.
  Let $\gamma: \cX \times \cX \to \{1, \ldots, |\cX|^2\}$
  be a bijection and
  define $\logistic(x) = \frac{1}{1+e^{-x}}$.
  Let
  $\widetilde A_i = \logistic(A_i^1) + \gamma(X_i^0, X_i^1)$
  so that
  $A_i^1 = \logistic^{-1}\big(\widetilde A_i
  - \lfloor \widetilde A_i \rfloor\big)$
  and
  $(X_i^0, X_i^1) = \gamma^{-1}(\lfloor \widetilde A_i \rfloor)$.
  Thus
  $l_i^{1 \triangleright 0}(w)$ is a bounded-variation function
  of $\widetilde A_i$, uniformly in $w$, and so as in
  Lemma~\ref{lem:strong_approx_Ln} we have that
  on an appropriately enlarged probability space,
  \begin{align*}
    \E\left[
      \sup_{w \in \cW}
      \left|
      \sqrt n L_n^{1 \triangleright 0}(w)
      - Z_n^{L, 1 \triangleright 0}(w)
      \right|
    \right]
    \lesssim
    \frac{\log n}{\sqrt n}
  \end{align*}
  where $Z_n^{L, 1 \triangleright 0}$ is a mean-zero
  Gaussian process with the same covariance structure as
  $\sqrt n L_n^{1 \triangleright 0}$.
  For $E_n^{1 \triangleright 0}(w)$,
  we first construct a strong approximation conditional on
  $\bA_n$ and $\bX_n$ as in Lemma~\ref{lem:conditional_strong_approx_En}
  and deduce an unconditional strong approximation as in
  Lemma~\ref{lem:unconditional_strong_approx_En} to see
  \begin{align*}
    \sup_{w \in \cW}
    \left|
    \sqrt{n^2h} E_n^{1 \triangleright 0}(w)
    - Z_n^{E, 1 \triangleright 0}(w)
    \right|
    \lesssim_\P
    n^{-1/4} h^{-3/8} (\log n)^{3/8} R_n
    + n^{-1/6} (\log n)^{2/3}
  \end{align*}
  where $Z_n^{E, 1 \triangleright 0}$ is a mean-zero
  Gaussian process with the same covariance structure as
  $\sqrt{n^2h} E_n^{1 \triangleright 0}$.
  Arguing as in the proof of Theorem~\ref{thm:strong_approx_fW}
  shows that the Gaussian processes are independent
  and can be summed to yield a single strong approximation.
\end{proof}

\begin{proof}[Lemma~\ref{lem:counterfactual_covariance_structure}]

  Arguing by mean-zero properties and conditional independence,
  \begin{align*}
    \Sigma_n^{1 \triangleright 0}(w,w')
    &=
    \Cov\left[
      \widehat f_W^{\,1 \triangleright 0}(w),
      \widehat f_W^{\,1 \triangleright 0}(w')
    \right] \\
    &=
    \frac{1}{n^2(n-1)^2(n-2)^2}
    \sum_{i \neq j}
    \sum_{r \notin \{i,j\}}
    \sum_{i' \neq j'}
    \sum_{r' \notin \{i',j'\}}
    \E\Big[
      \Big(
      k_{ij} \psi_i \psi_j
      - \E[k_{ij} \psi_i \psi_j]
      + k_{ij}\psi_i \kappa_{r j}
      + k_{ij}\psi_j \kappa_{r i}
      \Big) \\
      &\quad\times
      \Big(
      k_{i' j'}' \psi_{i'} \psi_{j'}
      - \E[k_{ij}' \psi_i \psi_j]
      + k_{i' j'}'\psi_{i'} \kappa_{r' j'}
      + k_{i' j'}'\psi_{j'} \kappa_{r' i'}
      \Big)
      \Big]
    + O\left( \frac{1}{n^{3/2}} + \frac{1}{\sqrt{n^4h}} \right) \\
    &=
    \frac{2}{n^2}
    \E\left[
      k_{ij} \psi_i \psi_j
      k_{ij}' \psi_i \psi_j
    \right]
    + \frac{4}{n}
    \E\left[
      k_{ij} \psi_i \psi_j
      k_{ir}' \psi_i \psi_r
    \right]
    - \frac{4}{n}
    \E\left[
      k_{ij} \psi_i \psi_j
    \right]
    \E\left[
      k_{ij}' \psi_i \psi_j
    \right] \\
    &\quad+
    \frac{4}{n}
    \E\left[
      k_{ij}\psi_i \kappa_{i' j}
      k_{i' j'}' \psi_{i'} \psi_{j'}
    \right]
    + \frac{4}{n}
    \E\left[
      k_{ij} \psi_{i} \psi_{j}
      k_{i' j'}'\psi_{i'} \kappa_{ij'}
    \right]
    + \frac{4}{n}
    \E\left[
      k_{ij} k'_{i' j'}
      \psi_i \psi_{i'}
      \kappa_{r j} \kappa_{r j'}
    \right] \\
    &\quad+
    O\left( \frac{1}{n^{3/2}} + \frac{1}{\sqrt{n^4h}} \right) \\
    &=
    \frac{4}{n}
    \E\left[
      \Big(
      \psi_i
      \E\big[
        k_{ij} \psi_j
        \mid i
        \big]
      + \E\left[
        k_{r j} \psi_r \kappa_{ij}
        \mid i
      \right]
      \Big)
      \Big(
      \psi_i
      \E\big[
        k_{ij}' \psi_j
        \mid i
        \big]
      + \E\left[
        k_{r j}' \psi_r \kappa_{ij}
        \mid i
      \right]
      \Big)
    \right] \\
    &\quad+
    \frac{2}{n^2}
    \E\left[
      k_{ij} k_{ij}'
      \psi_i^2 \psi_j^2
    \right]
    - \frac{4}{n}
    \E\left[
      k_{ij} \psi_i \psi_j
    \right]
    \E\left[
      k_{ij}' \psi_i \psi_j
    \right]
    + O\left( \frac{1}{n^{3/2}} + \frac{1}{\sqrt{n^4h}} \right),
  \end{align*}
  where all indices are distinct.
\end{proof}

\begin{proof}[Lemma~\ref{lem:counterfactual_infeasible_t_statistic}]
  The proof is exactly the same as the proof of
  Lemma~\ref{lem:infeasible_t_statistic}.
\end{proof}

\begin{proof}[Theorem~\ref{thm:counterfactual_infeasible_ucb}]
  This proof proceeds in the same manner as the proof of
  Theorem~\ref{thm:infeasible_ucb}.
\end{proof}

\bibliographystyle{hapalike}
\bibliography{refs}

%% file: preamble_CFU.tex
\usepackage[margin=1in]{geometry}

\usepackage{float}
\usepackage{mathtools}
\usepackage{enumitem}
\usepackage{amsmath}
\usepackage{amssymb}
\usepackage{amsthm}
\usepackage{mathrsfs}
\usepackage{caption}
\usepackage{subcaption}
\usepackage{graphicx}
\usepackage{multirow}
\usepackage{microtype}
\usepackage{titlesec}
\usepackage{xcolor}
\usepackage{mathrsfs}
\usepackage[hang,flushmargin]{footmisc}
\usepackage{stackengine,scalerel}
\usepackage{abstract}

\usepackage[longnamesfirst]{natbib}
\usepackage[onehalfspacing]{setspace}
\usepackage{color,hyperref}
\usepackage[norefs,nocites]{refcheck}
\hypersetup{pdfborder = {0 0 0},colorlinks=true,linkcolor=blue,%
  urlcolor=blue,citecolor=blue}

\usepackage{appendix}

\usepackage[boxruled,
  linesnumbered,
  commentsnumbered,
  procnumbered,
]{algorithm2e}

\DeclareMathOperator*{\argmin}{arg\,min}
\DeclareMathOperator{\IQR}{IQR}

\DeclareMathOperator{\AIMSE}{AIMSE}

\DeclareMathOperator{\Cov}{Cov}

\DeclareMathOperator{\Leb}{Leb}

\DeclareMathOperator{\Var}{Var}

\DeclareMathOperator{\ROT}{ROT}

\newcommand{\bA}{\mathbf{A}}
\newcommand{\bV}{\ensuremath{\mathbf{V}}}
\newcommand{\bX}{\ensuremath{\mathbf{X}}}
\newcommand{\bW}{\ensuremath{\mathbf{W}}}
\newcommand{\bY}{\ensuremath{\mathbf{Y}}}
\renewcommand{\P}{\ensuremath{\mathbb{P}}}
\newcommand{\E}{\ensuremath{\mathbb{E}}}

\newcommand{\Ck}{\ensuremath{C_\mathrm{k}}}

\newcommand{\cW}{\ensuremath{\mathcal{W}}}

\newcommand{\cC}{\ensuremath{\mathcal{C}}}

\newcommand{\CH}{\ensuremath{C_\mathrm{H}}}
\newcommand{\CL}{\ensuremath{C_\mathrm{L}}}

\newcommand{\cH}{\ensuremath{\mathcal{H}}}
\newcommand{\cA}{\ensuremath{\mathcal{A}}}

\newcommand{\I}{\ensuremath{\mathbb{I}}}

\newcommand{\cX}{\ensuremath{\mathcal{X}}}

\newcommand{\R}{\ensuremath{\mathbb{R}}}
\newcommand{\T}{\ensuremath{\mathsf{T}}}

\newcommand{\Dl}{\ensuremath{D_{\textup{lo}}}}
\newcommand{\Du}{\ensuremath{D_{\textup{up}}}}
\newcommand{\cP}{\ensuremath{\mathcal{P}}}
\newcommand{\cPd}{\ensuremath{\mathcal{P}_{\mathrm{d}}}}

\DeclareMathOperator*{\minimize}{minimize}
\DeclareMathOperator*{\subjectto}{subject\ to}

\newcommand{\TV}{\mathrm{TV}}

\newcommand{\diff}[1]{\,\mathrm{d}#1}

\newcommand{\flbeta}{\ThisStyle{%
    \ensurestackMath{%
      \stackengine{-0.5\LMpt}{\SavedStyle \beta}%
      {\SavedStyle {\rule{3.7\LMpt}{0.3\LMpt}}}%
      {U}{c}{F}{F}{S}%
    }%
    \vphantom{\beta}}%
}

\titleformat
{\paragraph} 
[hang] 
{\bfseries\upshape} 
{} 
{1pt} 
{} 
[] 

\titlespacing*{\paragraph}{0pt}{6pt}{0pt}

\newcounter{proofparagraphcounter}

\newlist{inlineroman}{enumerate*}{1}
\setlist[inlineroman]{afterlabel=~,label=(\roman*)}

\newtheoremstyle{break}%
{\topsep}{\topsep}%
{\itshape}{}%
{\bfseries}{}%
{\newline}%
{\thmname{#1}\thmnumber{ #2}%
  \thmnote{ \rm(#3)\addcontentsline{toc}{subsubsection}{{\it#1 #2}: #3}}}

\newtheoremstyle{breakplainproof}
{\topsep}{\topsep}%
{}{}%
{\bfseries}{}%
{\newline}%
{\thmname{#1}\thmnumber{ #2}%
  \thmnote{ \rm(#3)\addcontentsline{toc}{subsubsection}{{\it#1}: #3}}}

\theoremstyle{break}
\newtheorem{theorem}{Theorem}[section]
\newtheorem{lemma}{Lemma}[section]

\newtheorem{assumption}{Assumption}[section]

\theoremstyle{breakplainproof}

\newtheorem*{proof}{Proof}
\AtBeginEnvironment{proof}{\setcounter{proofparagraphcounter}{0}}%
\AtEndEnvironment{proof}{\null\hfill\qedsymbol}%

\theoremstyle{remark}

\DontPrintSemicolon%
\makeatletter%
\renewcommand{\SetKwInOut}[2]{%
  \sbox\algocf@inoutbox{\KwSty{#2}\algocf@typo:}%
  \expandafter\ifx\csname InOutSizeDefined\endcsname\relax%
    \newcommand\InOutSizeDefined{}%
    \setlength{\inoutsize}{\wd\algocf@inoutbox}%
    \sbox\algocf@inoutbox{%
      \parbox[t]{\inoutsize}%
      {\KwSty{#2}\algocf@typo:\hfill}~%
    }%
    \setlength{\inoutindent}{\wd\algocf@inoutbox}%
  \else%
    \ifdim\wd\algocf@inoutbox>\inoutsize%
      \setlength{\inoutsize}{\wd\algocf@inoutbox}%
      \sbox\algocf@inoutbox{%
        \parbox[t]{\inoutsize}%
        {\KwSty{#2}\algocf@typo:\hfill}~%
      }%
      \setlength{\inoutindent}{\wd\algocf@inoutbox}%
    \fi%
  \fi%
  \algocf@newcommand{#1}[1]{%
    \ifthenelse{\boolean{algocf@inoutnumbered}}{\relax}{\everypar={\relax}}{%
      \let\\\algocf@newinout\hangindent=\inoutindent\hangafter=1\parbox[t]%
      {\inoutsize}{\KwSty{#2}%
        \algocf@typo:\hfill}~##1\par%
    }%
    \algocf@linesnumbered%
  }%
}%
\makeatother%
\SetKwInOut{Input}{Input}%
\SetKwInOut{Output}{Output}%

%% file: table.tex
\begin{center}
    \begin{tabular}{|c|c|c|c|c|cc|cc|}
        \hline
        \multirow{2}{*}{$ \pi $}
        & \multirow{2}{*}{Degeneracy type}
        & \multirow{2}{*}{$ \widehat h_{\ROT} $}
        & \multirow{2}{*}{$ p $}
        & \multirow{2}{*}{RIMSE}
        & \multicolumn{2}{|c|}{UCB}
        & \multicolumn{2}{|c|}{PCI} \\
        \cline{6-9}
        & & & &
        & CR & AW
        & CR & AW \\
        \hline
\multirow{2}{*}{$ \left(\frac{1}{2}, 0, \frac{1}{2}\right) $}
& \multirow{2}{*}{Total}
& \multirow{2}{*}{0.161}
 & 2 & 0.00048 & 87.1\% & 0.0028 & 6.5\% & 0.0017 \\
 & &  & 4 & 0.00068 & 95.2\% & 0.0042 & 9.7\% & 0.0025 \\
\hline
\multirow{2}{*}{$ \left(\frac{1}{4}, 0, \frac{3}{4}\right) $}
& \multirow{2}{*}{Partial}
& \multirow{2}{*}{0.158}
 & 2 & 0.00228 & 94.5\% & 0.0112 & 75.6\% & 0.0083 \\
 & &  & 4 & 0.00234 & 94.7\% & 0.0124 & 65.3\% & 0.0087 \\
\hline
\multirow{2}{*}{$ \left(\frac{1}{5}, \frac{1}{5}, \frac{3}{5}\right) $}
& \multirow{2}{*}{None}
& \multirow{2}{*}{0.145}
 & 2 & 0.00201 & 94.2\% & 0.0106 & 73.4\% & 0.0077 \\
 & &  & 4 & 0.00202 & 95.6\% & 0.0117 & 64.3\% & 0.0080 \\
\hline
\end{tabular}
\end{center}